\documentclass[10pt]{amsart}
\usepackage{amsmath}
\usepackage{amsxtra}
\usepackage{amscd}
\usepackage{amsthm}
\usepackage{amsfonts}
\usepackage{amssymb}
\usepackage{eucal}
\usepackage[all]{xy}
\usepackage{tikz-cd}
\usepackage{graphicx}
\usepackage{comment}
\usepackage{epsfig}
\usepackage{psfrag}
\usepackage{mathrsfs}
\usepackage{amscd}
\usepackage{rotating}
\usepackage{lscape}
\usepackage{amsbsy}
\usepackage{verbatim}
\usepackage{moreverb}
\usepackage{url}
\usepackage{mathtools}
\usepackage{hyperref}
\usepackage{csquotes}
\usepackage{mathabx}
\usepackage{booktabs}

\usepackage{xcolor}

\def\acts{\mathrel{\reflectbox{$\righttoleftarrow$}}}

\newcommand{\sssec}{\subsubsection}

\newtheorem{thm}[subsubsection]{Theorem}
\newtheorem{lemma}[subsubsection]{Lemma}
\newtheorem{prop}[subsubsection]{Proposition}
\newtheorem{cor}[subsubsection]{Corollary}
\newtheorem{defn}[subsubsection]{Definition}

\newtheorem{nota}[subsubsection]{Notation}

\newtheorem{claim}[subsubsection]{Claim}

\newtheorem{note}[subsubsection]{Note}

\theoremstyle{remark}
\newtheorem{rmk}[subsubsection]{Remark}
\newtheorem{ex}[subsubsection]{Example}

\newcommand{\mc}[1]{\mathcal{#1}}

\newcommand{\mf}[1]{\mathfrak{#1}}

\newcommand{\te}{\text}
\newcommand{\pa}{\partial}

\newcommand{\ul}{\underline}
\newcommand{\wt}{\widetilde}
\newcommand{\wh}{\widehat}
\newcommand{\C}{\mathbb{C}}
\newcommand{\E}{\mathbb{E}}
\newcommand{\R}{\mathbb{R}}
\newcommand{\Z}{\mathbb{Z}}

\newcommand{\afl}{\mathbb{A}^1}
\newcommand{\pjl}{\mathbb{P}^1}

\usepackage[scr=rsfso]{mathalfa}

\newcommand{\nat}{\natural}
\newcommand{\cP}{\mathcal{P}}
\newcommand{\on}{\operatorname}

\newcommand{\cA}{\mathcal{A}}
\newcommand{\cB}{\mathcal{B}}
\newcommand{\cC}{\mathcal{C}}
\newcommand{\cD}{\mathcal{D}}

\newcommand{\cF}{\mathcal{F}}
\newcommand{\cG}{\mathcal{G}}
\newcommand{\cH}{\mathcal{H}}
\newcommand{\cI}{\mathcal{I}}
\newcommand{\cL}{\mathcal{L}}
\newcommand{\cM}{\mathcal{M}}
\newcommand{\cN}{\mathcal{N}}
\newcommand{\cO}{\mathcal{O}}
\newcommand{\cS}{\mathcal{S}}
\newcommand{\cK}{\mathcal{K}}

\newcommand{\cX}{\mathcal{X}}
\newcommand{\cY}{\mathcal{Y}}

\newcommand{\cJ}{\mathcal{J}}
\newcommand{\cT}{\mathcal{T}}

\newcommand{\Cat}{\on{Cat}}

\newcommand{\CatEc}{\Cat_{E^c}}

\newcommand{\ShvCat}{\on{ShvCat}}
\newcommand{\DGcat}{\on{DGCat}_{\on{cont}}}
\newcommand{\LSCat}{\on{LSCat}}
\newcommand{\LSMon}{\on{LSMon}}
\newcommand{\LSCMon}{\on{LSCMon}}
\newcommand{\LSSymMon}{\on{LSCMon}}
\newcommand{\DGCatMon}{\on{DGMon}}
\newcommand{\DGCatSymMon}{\on{DGCMon}}

\newcommand{\lMod}{\text{-}\on{Mod}}
\newcommand{\lAlg}{\text{-}\on{Alg}}

\newcommand{\Perv}{\on{Perv}}

\newcommand{\Rep}{\on{Rep}}

\newcommand{\Vect}{\on{Vect}}
\newcommand{\Mod}{\on{Mod}}

\newcommand{\QC}{\on{QC}}

\newcommand{\Shv}{\on{Shv}}

\newcommand{\Nilp}{\on{Nilp}}
\newcommand{\Loc}{\on{Loc}}
\newcommand{\LS}{\on{LS}}
\newcommand{\QLisse}{\on{QLisse}}

\newcommand{\Sat}{\on{Sat}}

\newcommand{\nc}{\newcommand}
\nc{\triv}{\on{triv}}
\nc{\ft}{\on{ft}}
\nc{\cl}{\on{cl}}
\nc{\Ab}{\on{Ab}}
\nc{\ULA}{\on{ULA}}
\nc{\TwArr}{\on{TwArr}}
\nc{\fSet}{\on{fSet}}
\nc{\reg}{\on{reg}}
\nc{\laft}{\on{laft}}
\nc{\topl}{\on{top}}
\nc{\sing}{\on{sing}}
\nc{\loc}{\on{loc}}
\nc{\sh}{\on{sh}}
\nc{\Alg}{\mathsf{Alg}}
\nc{\CAlg}{\mathsf{CAlg}}
\nc{\plain}{\on{plain}}
\nc{\constr}{\on{constr}}
\nc{\indconstr}{\on{ind-constr}}
\nc{\indft}{\on{ind-ft}}

\newcommand{\PreStk}{\on{PreStk}}

\newcommand{\PreStklft}{\on{PreStk}_{\on{lft}}}
\nc{\PreStklaft}{\on{PreStk}_{\on{laft}}}

\newcommand{\Stk}{\on{Stk}}
\newcommand{\Stklft}{\on{Stk}_{\on{lft}}}
\newcommand{\NearStklft}{\on{NearStk}_{\on{lft}}}
\newcommand{\Aff}{\on{Aff}}
\newcommand{\Affft}{\on{Aff}_{\on{ft}}}
\newcommand{\Schft}{\on{Sch}_{\on{ft}}}
\newcommand{\Spc}{\on{Spc}}
\nc{\StinftyCat}{\on{DGCat}}
\nc{\inftyCat}{1\on{-Cat}}
\nc{\inftygroup}{\infty\on{-Grpd}}
\nc{\PrL}{\on{Pr}^{\on{L}}}
\nc{\RTens}{\on{RTens}}
\nc{\Dmod}{\on{D-mod}}
\nc{\CMaps}{{\mathcal Maps}}
\nc{\Maps}{\on{Maps}}
\nc{\len}[1]{\prescript{\le n}{}{#1}}
\nc{\cpt}[1]{{#1}^{\on{c}}}
\nc{\et}{\text{\'et}}
\nc{\Bet}[1]{B_{\et} #1}
\nc{\BA}{\Bet{A}}
\nc{\Btet}[1]{B^2_{\et} #1}
\nc{\Biet}[1]{B^i_{\et} #1}
\nc{\Bt}[1]{B^2 #1}
\nc{\Grpd}[1]{\on{Grpd}(#1)}
\nc{\Grp}[1]{\on{Grp}(#1)}
\nc{\Sh}[1]{\on{Sh}(#1)}
\nc{\PSh}{\on{PSh}}
\nc{\clPhi}{\prescript{\cl}{}{\Phi}}
\nc{\Arr}{\on{Arr}}
\nc{\FactGe}{\on{FactGe}_A(X)}
\nc{\chG}{\check{G}}
\nc{\plp}[1]{\prescript{+}{}{#1}}
\nc{\elp}[1]{\prescript{\epsilon}{}{#1}}

\DeclareMathOperator{\Hom}{Hom}
\DeclareMathOperator{\cHom}{\mathcal{H}om}

\DeclareMathOperator{\End}{End}
\DeclareMathOperator{\Aut}{Aut}

\DeclareMathOperator{\Ind}{Ind}

\newcommand{\Fun}{\on{Fun}}
\DeclareMathOperator{\pr}{pr}

\DeclareMathOperator{\op}{op}
\DeclareMathOperator{\colim}{colim}

\nc{\linf}[1]{\prescript{<\infty}{}{#1}}
\newcommand{\id}{\on{id}}

\newcommand{\glob}{\on{glob}}

\makeatletter
\def\lotimes{\@ifnextchar_{\@lotimessub}{\@lotimesnosub}}
\def\@lotimessub_#1{\mathchoice{\mathbin{\mathop{\otimes}^L}_{#1}}%
  {\otimes^L_{#1}}{\otimes^L_{#1}}{\otimes^L_{#1}}}
\def\@lotimesnosub{\mathbin{\mathop{\otimes}^L}}
\makeatother

\newcommand{\bA}{\mathbf{A}}
\newcommand{\bB}{\mathbf{B}}
\nc{\bC}{\mathbf{C}}
\newcommand{\bD}{\mathbf{D}}
\newcommand{\bE}{\mathbf{E}}

\newcommand{\bM}{\mathbf{M}}
\newcommand{\bS}{\mathbf{S}}

\newcommand{\bO}{\mathbf{O}}
\nc{\ba}{\textbf{a}}
\nc{\bbb}{\textbf{b}}
\nc{\bc}{\textbf{c}}
\nc{\bd}{\textbf{d}}
\nc{\be}{\textbf{e}}

\newcommand{\pt}{\on{pt}}
\newcommand{\surj}{\twoheadrightarrow}

\newcommand{\Gr}{\on{Gr}}

\newcommand{\Hecke}{\on{Hecke}}
\newcommand{\Bun}{\on{Bun}}

\newcommand{\Tot}{\on{Tot}}

\DeclareMathOperator{\Spec}{Spec}

\newcommand{\sffH}{\mathsf{H}}
\newcommand{\sfA}{\mathsf{A}}
\newcommand{\sfY}{\mathsf{Y}}

\newcommand{\Mul}{\on{Mul}}
\newcommand{\Conv}{\mathsf{Conv}}
\newcommand{\Convcomm}{\mathsf{Conv}^{\on{comm}}}
\newcommand{\ConvMod}{\mathsf{Conv}^{\on{Mod}}}
\newcommand{\ConvcMod}{\mathsf{Conv}^{\on{commMod}}}

\newcommand{\cHecke}{\mathcal{H}ecke}
\newcommand{\Corr}{\on{Corr}}
\newcommand{\LM}{\mathsf{LM}}

\newcommand{\LMcomm}{\mathsf{LM}^{\on{comm}}}

\newcommand{\Ectens}{\stackrel{\on{f.d.}}{\otimes}}
\newcommand{\Ectensbig}{\stackrel{\on{f.d.}}{\bigotimes}}
\newcommand{\VectAbfd}{\Vect^{\heartsuit,\on{f.d.}}}

\mathchardef\mhyphen="2D

\title{Spectral action in Betti metaplectic Langlands}
\author{JiWoong Park}

\begin{document}

\date{\today}

\begin{abstract}
We construct a spectral action in the context of Betti metaplectic Langlands theory. More precisely, we prove that the category of twisted sheaves on $\Bun_G$ with nilpotent singular support admits a canonical action by the symmetric monoidal category of quasicoherent sheaves on the moduli of twisted $H$-local systems, where $H$ is the metaplectic dual group. This confirms the Betti version of the conjecture of \cite{GL}, originally formulated in the de Rham setting under the name ``metaplectic vanishing conjecture".  

Given a symmetric monoidal category with suitable symmetry, we develop a categorical framework by introducing a notion of twisted factorization homology and studying its categorical action, parallel to \cite{GKRV}. We show that the category of quasicoherent sheaves on the moduli of twisted local systems can be realized as a twisted factorization homology. Using metaplectic factorizable Satake (\cite{Rei}, \cite{GL}, and \cite{ZhaoNew}) we further prove that the metaplectic Hecke action naturally fits into our categorical framework. Combining these results, we obtain the desired spectral action.
\end{abstract}

\maketitle
\tableofcontents

\section{Introduction}
The idea of \textbf{spectral decomposition} plays a key role in various versions of the Langlands program. Notably, in the function field context of the Langlands program, V. Lafforgue’s work \cite{Laf} decomposes the vector space of cusp forms in terms of spectral data, i.e., $L$-parameters. In the geometric context, spectral decomposition is realized through the \textbf{spectral action} \cite{Gai1, NY, GKRV, AGKRRV}, which is the categorical action of the quasi-coherent sheaf category of the moduli of $L$-parameters on the category of automorphic sheaves.
\medskip

The \textbf{metaplectic Langlands program} \cite{GGW, Wei, GG} is a variant of the Langlands program that replaces the role of reductive groups by their topological covering groups. Moreover, this admits a natural geometric analog, called the \textbf{metaplectic geometric Langland program}. One hopes the above theme of spectral decomposition to carry over to the metaplectic context. For example, the metaplectic version of \cite{Laf} was established in \cite{ZhaoNew}. 
\medskip

The goal of this article is to construct the spectral action in the context of Betti metaplectic geometric Langlands: 
\begin{thm}\label{thm: main}
    There is a canonical categorical action of $\QC(\Loc^{\cG_Z}_H)$ on $\cA^q_G$.
\end{thm}
The notation will be explained later throughout the introduction. We comment that $\Loc^{\cG_Z}_H$ is the moduli stack of metaplectic $L$-parameters and $\cA^q_G$ is a geometric incarnation of genuine automorphic forms considered in metaplectic Langlands.
\medskip

The main new ingredients for the proof are:

- the categorical framework relating a factorizable action of a local system of symmetric monoidal categories to an action of (a variant of) factorization homology,

- construction of the factorizable metaplectic Satake functors between derived categories and compatibility data between those, leveraging on the abelian equivalence \cite{Rei, ZhaoNew}.
\medskip

For our purpose, we develop some foundations in the appendix, which may be of independent interest:

- we study some methods of twisting sheaves of categories with algebra and module structures,

- we study twisted sheaves in depth, and establish a six-functor formalism for twisted sheaves.

\subsection{Metaplectic Langlands and passage to geometry}
Topological covering groups of reductive groups appear naturally in various contexts. The classic example of such a covering group - which plays a crucial role in theta correspondence - is the \emph{metaplectic group} \footnote{Metaplectic Langlands borrows its name from it.} $\on{Mp}_{2n}(\R)$, which is the 2-to-1 cover of $\on{Sp}_{2n}(\R)$. More generally, we can replace $\on{Sp}_{2n}$ by a different reductive group and the field $\R$ with either a local field $k$ or the ring of adeles $\mathbb{A}$ of a global field $F$, and hence get a covering group $\wt{G}(k)\rightarrow G(k)$ or $\wt{G}(\mathbb{A})\rightarrow G(\mathbb{A})$. We emphasize that the notation $\wt{G}(\mhyphen)$ is misleading here, because it does \emph{not} mean points of a reductive group. For example, $\on{Mp}_{2n}(\R)$ is not even a linear algebraic group. The groups $\wt{G}(k)$ or $\wt{G}(\mathbb{A})$ only makes sense as a topological group.

\sssec{Metaplectic Langlands}
For the purpose of metaplectic Langlands, we do not consider arbitrary covering groups, but those arising from certain parameters. Classically, one takes \emph{Brylinski-Deligne data} \cite{BD} as the parameter, and this is the approach taken in \cite{GGW, Wei, GG}. More recently, Y. Zhao \cite{Zhao1} introduced another notion of parameter called the \emph{\'etale level}. \footnote{In \emph{op. cit.}, this parameter was called the \emph{\'etale metaplectic cover}. However, since Zhao began referring to it as the \'etale level in the following paper \cite{ZhaoNew}, we adopt that terminology.} There is a natural association of an \'etale level with Brylinski-Deligne data, and hence the latter parameter incorporates the former. The notion of an \'etale level is defined for a group scheme $G\rightarrow S$ of finite type, and we make a choice of a finite abelian group $A$. If $S=\Spec(k)$, an \'etale level gives rise to a central extension \footnote{This needs a mild modification if $k=\R$, but we ignore this point.}
\[1\rightarrow A\rightarrow \wt{G}(k)\rightarrow G(k)\rightarrow 1\]
of topological groups. If $S=\Spec(F)$, an \'etale level gives rise to a central extension \footnote{Likewise, this needs a mild modification which has to do with real places, but we ignore this point.}
\begin{equation}\label{eq: central ext global}
    1\rightarrow A\rightarrow \wt{G}(\mathbb{A})\rightarrow G(\mathbb{A})\rightarrow 1.
\end{equation}
The metaplectic group $\on{Mp}_{2n}(\R)$ is obtained from the \'etale level called `$c_2$ mod 2'.
\medskip

Choose a character $\chi:A\hookrightarrow \C^{\times}$. In the local context, the metaplectic Langlands program studies the \textbf{genuine representations} of the covering group $\wt{G}(k)$. This means a representation $V$ of a topological group $\wt{G}(k)$ with the property
\[a\cdot v=\chi(a)v,\quad a\in A\subset \wt{G}(k),\ v\in V.\]
In the global context, the central extension \eqref{eq: central ext global} is equipped with a canonical splitting along $G(F)$, and hence we can take a quotient space $G(F)\backslash \wt{G}(A)$. In the global context, the metaplectic Langlands program studies the \textbf{genuine automorphic forms} of $G(\mathbb{A})$, which is a suitable class of complex-valued functions on $G(F)\backslash\wt{G}(\mathbb{A})$ satisfying the condition
\begin{equation}\label{eq: genuine condition}
   f(ag)=\chi(a)f(g),\quad a\in A, \ g\in G(F). 
\end{equation}

\sssec{Passage to geometry}
The above situation, more precisely the unramified version of metaplectic Langlands, admits a natural geometric analog. The notion of \'etale level carries over to the geometric context. We choose a complex reductive group $G$ and a smooth projective complex curve $X$, and let the group scheme $G\rightarrow S$ to be the constant group scheme $G\times X\rightarrow X$. Thus, as in the non-metaplectic case, the global field was replaced by the curve $X$. The moduli stack $\Bun_G$ of $G$-bundles on $X$ replaces the role of the double quotient $G(F)\backslash G(\mathbb{A})/G(\cO)$. In the above setup, the notion of \'etale level is equivalent to the notion of factorization gerbes on (the Ran version of) affine Grassmannian of $G$, studied in \cite{Rei, GL}. We switch the terminology and use the term \textbf{metaplectic parameter} in the geometric context. 
\medskip

We review some sheaf-theoretic terminology before proceeding. We work in the Betti setup: sheaves will be Betti sheaves, i.e., complexes of $E$-linear sheaves in the analytic topology, by default.
Given an $A$-\emph{gerbe} \footnote{A gerbe is a 1-categorical higher notion of a torsor. Its definition is reviewed in Sect. \ref{ssec: basic definitions}.} $\cG$ on a topological space (or a complex stack) $\cY$, we have a notion of twisted $E$-linear sheaves \footnote{We refer readers to Sect. \ref{sssec: twisted sheaves brief} for details.} on $\cY$, which gets organized into the (derived) category $\Shv^{\cG}(\cY)$. 
\medskip

Fix an algebraically closed field $E$ of characteristic zero and a finite abelian group $A\hookrightarrow E^{\times}$. The geometric analog of (unramified) genuine automorphic forms should be \emph{twisted} $E$-linear sheaves. To justify this, we observe: although genuine automorphic forms are functions on $G(F)\backslash \wt{G}(\mathbb{A})$, they can also be thought of as \emph{twisted} functions on $G(F)\backslash G(\mathbb{A})$. Indeed, if the $A$-torsor $G(F)\backslash \wt{G}(\mathbb{A})\rightarrow G(F)\backslash G(\mathbb{A})$ was trivial, the functions on the total space satisfying the condition \eqref{eq: genuine condition} would be the same as functions on the base space.
\medskip

An \'etale level gives rise to a covering group of $G(\mathbb{A})$ and hence an $A$-torsor over $G(F)\backslash G(\mathbb{A})$, enabling us to consider genuine automorphic forms. In the geometric context, a metaplectic parameter gives rise to a canonical $A$-\emph{gerbe} over $\Bun_G$, enabling us to consider \textbf{genuine automorphic sheaves}, which are certain twisted sheaves on $\Bun_G$. This analogy can be summarized in the following table:
\[
\begin{tabular}{|l|c|c|}
\hline
               & \textbf{arithmetic}                        & \textbf{geometric}                   \\ \hline
               twisting data  & covering group of $G(\mathbb{A})$ & gerbe on $\Bun_G$           \\ \hline
twisted objects & genuine automorphic forms         & twisted sheaves on $\Bun_G$ \\ \hline

\end{tabular}
\]

\subsection{Spectral decomposition}
A central theme in both arithmetic and geometric Langlands programs is the \emph{spectral decomposition} of the automorphic vector space/category. In both contexts, Hecke operators combined with the (geometric) Satake equivalence provide a large supply of commuting operators - each of which arises from representations of the Langlands dual group $\chG$ - on the vector space of automorphic forms/category of automorphic sheaves. One attempts to understand the vector space/category by means of simultaneous diagonalization under this commuting family of actions. The spectrum, or in other words, simultaneous eigenvalues under such actions amounts to a Galois representation in $\chG$, and hence the idea of spectral decomposition provides us the automorphic-to-Galois direction. This was carried out in the function field setting in \cite{Laf}, leading to the spectral decomposition of cusp forms. This approach can be carried out cleanly in the geometric context, via the \emph{spectral action}.

\sssec{Spectral action}
In the geometric context, the spectrum corresponds to the $\chG$-local systems on the curve $X$, which gets organized into the moduli space $\Loc_{\chG}$. Let $\cA_G$ be the category of automorphic sheaves, or the \textbf{automorphic category} for short, which is a suitable full subcategory of $\Shv(\Bun_G)$.
We have:
\begin{thm}[Spectral action]\label{thm: spectral action non mp}
    There is a canonical action of $\QC(\Loc_{\chG})$ on $\cA_G$.
\end{thm}

Using the language of sheaves of categories developed in \cite{Ga}, the above result spreads the category $\cA_G$ over the stack $\Loc_{\chG}$ as a sheaf of categories. Since $\Loc_{\chG}$ classifies the spectrum of Hecke actions, or in other words $L$-parameters, this categorical action realizes the automorphic-to-Galois direction in the geometric context.

\begin{rmk}
The spectral action is very robust and has avatars in different contexts of geometric Langlands. The first version using D-modules was proved by Gaitsgory \cite{Gai1}, and its restricted variant was proved in \cite{AGKRRV}. The Betti version, which is the most relevant to this paper, was proved in \cite{NY} and later in \cite{GKRV}.
\end{rmk}

\begin{rmk}
    In the context of the Fargues--Scholze program on the geometrization of the local Langlands correspondence \cite{FS}, the spectral action is used to cosntruct $L$-parameters.
\end{rmk}

\begin{rmk}
 The importance of this action in geometric Langlands theory cannot be overstated. Using this action, we get a functor from the spectral side to the automorphic side once we choose an object in $\cA_G$. The right choice is the \emph{vacuum Poincar\'e sheaf}, and this is the starting point of the construction of the geometric Langlands functor.
\end{rmk}

\sssec{Metaplectic situation}
This theme of spectral decomposition carries over to the metaplectic situation. As in \cite{Wei} and \cite{Zhao1}, there is a canonical reductive group $H$ - called the \textbf{metaplectic dual group} - that replaces the role of $\chG$ in the non-metaplectic situation. There are metaplectic Hecke operators arising from representations of $H$ that act on the vector space of genuine automorphic forms/category of genuine automorphic sheaves. This again provides us with a large supply of commuting operators and the spectrum amounts to metaplectic $L$-parameters, or more precisely \emph{twisted} Galois representations in $H$. The twisting here is given by the \textbf{metaplectic dual gerbe} $\cG_Z$, which twists the $L$-group and hence the notion of $L$-parameters. The association of $(H,\cG_Z)$ to an \'etale level is worked out in \cite{Wei} (using Brylinski-Deligne data) and \cite{Zhao1}, and also in \cite{Rei} and \cite{GL} in the geometric context.
In the function field setting, \cite{ZhaoNew} provides a metaplectic analog of the spectral decomposition result in \cite{Laf}. 
\medskip

Consider the geometric setting with a metaplectic parameter $q$. This induces an $A$-gerbe on $\Bun_G$, and we denote the resulting category of twisted sheaves by $\Shv^q(\Bun_G)$. We study a certain full subcategory $\cA^q_G$ of it and consider it as the category of genuine automorphic sheaves, see Sect. \ref{sssec: AqG}. We call $\cA^q_G$ the \textbf{metaplectic automorphic category}.
The spectrum is organized into the moduli space $\Loc^{\cG_Z}_H$ of twisted $H$-local systems on the curve, which we construct in Sect. \ref{sssec: construction of twlocsys}.
\medskip

Theorem \ref{thm: main} admits an interpretation analogous to that following the statement of Theorem \ref{thm: spectral action non mp}. It spreads the category $\cA^q_G$ over the stack $\Loc^{\cG_Z}_H$ as a sheaf of categories. Since $\Loc^{\cG_Z}_H$ classifies the spectrum of metaplectic Hecke actions, or in other words metaplectic $L$-parameters, our main result realizes the automorphic-to-Galois direction in the metaplectic geometric context.

\subsection{Factorization homology}\label{ssec: fact hlgy}
For the rest of the introduction, we discuss some ideas behind the proof of Theorem \ref{thm: main}. To do so, it is helpful to review the proof in the non-twisted case.
The approach of \cite{NY} and \cite{GKRV} to Betti spectral action is crucially based on the categorical framework of \emph{factorization homology} and its categorical action. 

\sssec{Factorizable action and factorization homology}\label{sssec: fact hlgy informal} 
We give an informal account of factorization homology and its relationship to factorizable action. Let $\cC$ be a symmetric monoidal category and $Y$ be a homotopy type. Then we have a symmetric monoidal category $\int_Y \cC$ called the \textbf{factorization homology}, see Sect. \ref{sssec: GKRV fact hlgy} for its definition. For any finite set $I$ and an $I$-tuple of points $\ul{y}\in Y^I$ the factorization homology is equipped with a canonical symmetric monoidal functor \footnote{The colimit expression of factorization homology in Sect. \ref{sssec: GKRV fact hlgy} justifies the name $\on{ins}_{\ul{y}}$ of this functor: this comes from insertion maps.}
\[\on{ins}_{\ul{y}}:\cC^{\otimes I}\rightarrow\int_Y\cC.\]

Let $\cM$ be a category and assume that $\cM$ admits an action of $\cC$ labelled by multiple points of $Y$. More precisely, it means that for any finite set $I$ and $\ul{y}\in Y^I$ we have an action $\on{act}_{\ul{y}}:\cC^{\otimes I}\acts\cM$. We consider the collection of actions $\{\on{act}_{\ul{y}}\}$ where we also vary $I$. 
\medskip

We consider the situation where the collection $\{\on{act}_{\ul{y}}\}$ defines a \textbf{factorizable action} of $\cC$ on $\cM$ along $Y$. Informally, this means the following two structures:

- \textbf{(Local constancy)} For each $I$, the collection $\{\on{act}_{\ul{y}}\}_{\ul{y}\in Y^I}$ varies locally constantly on $\ul{y}$.

- \textbf{(Multiplicity $=\otimes$)} The multiplicity in the tuple $\ul{y}$ corresponds to the symmetric monoidal operation of $\cC$. For example, the actions $\on{act}_{(y,y)}:\cC^{\otimes 2}\acts \cM$ and $\on{act}_{y}:\cC\acts\cM$ are related in the way that we have an equivalence
\[\on{act}_{(y,y)}\simeq \left(\cC\otimes\cC\xrightarrow{\on{mult}}\cC\stackrel{\on{act}_y}{\acts} \cM\right).\]

A factorizable action determines a canonical categorical action of $\int_Y\cC$ on $\cM$ which encodes all the $\on{act}_{\ul{y}}$'s at the same time. To be more precise, there is a canonical equivalence between the actions
\[\on{act}_{\ul{y}}\simeq \left(\cC^{\otimes I}\xrightarrow{\on{ins}_{\ul{y}}}\int_Y\cC\acts \cM\right).\]

\sssec{}
The Hecke action fits into the framework of Sect. \ref{sssec: fact hlgy informal}. Let $X^{\topl,\sing}$ be the homotopy type of $X$. We take $Y=X^{\topl,\sing}$, $\cC=\Rep(\chG)$, and $\cM=\cA_G$. For each $\ul{y}\in X(\C)^I$ the action $\on{act}_{\ul{y}}$ is given by the Hecke action at multiple points $\ul{y}$. A priori, the Hecke action is the action of $\Rep(\chG)^{\otimes I}$ on $\Shv(\Bun_G(X))$. However, in \cite{NY} this action was shown to preserve the full subcategory $\cA_G$, and hence induces an action $\on{act}_{\ul{y}}: \Rep(\chG)^{\otimes I}\acts \cA_G$.
Moreover, this action was shown to vary locally constantly on $\ul{y}$ in \emph{op. cit..} We have the compatibility data pertaining to multiplicities by construction.
\medskip

Therefore, the Hecke gives rise to a canonical action of
$\int_{X^{\topl,\sing}}\Rep(\chG)$ on $\cA_G$. In \cite{BZFN} and later in \cite{GKRV}, this factorization homology is canonically identified with $\QC(\Loc_{\chG})$ as a symmetric monoidal category, and hence the spectral action. From construction, it encodes all the Hecke actions at once.

\subsection{What changes in the metaplectic setting?}
One can attempt to follow the approach of \cite{NY} and \cite{GKRV} using factorization homology to prove Theorem \ref{thm: main}. However, there are some features of metaplectic geometric Langlands theory that make the situation more complicated. The main difference is that what acts on $\cA^q_G$ at each point $x\in X$ has a nontrivial dependence along $X$. This is visible from metaplectic version of factorizable Satake equivalence, which is how the metaplectic dual data $(H,\cG_Z)$ enters the story. 

\sssec{Metaplectic Satake}
At each point $x\in X$, we have the equivariant affine Grassmannian, i.e., the quotient of the affine Grassmannian by the arc group, denoted by $\Hecke^{\loc}_x$. The metaplectic parameter $q$ induces a canonical $A$-gerbe on that space. The corresponding category of twisted perverse sheaves on $\Hecke^{\loc}_x$ is called the \textbf{metaplectic Satake category}, and denoted by $\Sat^q_x$. As in the non-metaplectic case, this category can be endowed with a canonical symmetric monoidal structure. From the metaplectic Satake equivalence, $\Sat^q_x$ is equivalent to $\Rep(H)^{\heartsuit}$ as tensor categories.\footnote{As in the non-metaplectic Satake, there should be a sign modification of commutativity constraint. We ignore this point in the introduction. We refer readers to Sect. \ref{ssec: recollection perv ULA} for the corrected statement. However, this is still valid  as a monoidal equivalence, and it mostly suffices for the application to the spectral action.} \footnote{From now on, we clearly distinguish between abelian categories and DG categories. $\Rep(H)$ is a DG category while $\Rep(H)^{\heartsuit}$ is an abelian category.} However, this identification is not canonical, the reason being that $\Sat^q_x$ does \emph{not} admit a canonical fiber functor. Instead, $\Sat^q_x$ is canonically equivalent to $(\Rep(H))^{\heartsuit}_{(\cG_{Z})_x}$ as a symmetric monoidal category, where the latter is a \emph{twist} of $\Rep(H)^{\heartsuit}$. In fact, this pointwise equivalence fits into the metaplectic Satake equivalence along the curve, i.e., the factorizable metaplectic Satake:

\begin{thm}[\cite{Rei, ZhaoNew}]\label{thm: abelian Satake light}
    There is a canonical symmetric monoidal equivalence   between $\Sat^q_{X}$ and $\Rep(H)^{\Ab,\loc,\ft}_{\cG_Z}(X)$. 
\end{thm}

Since it plays a very important role in our paper, we review this result in detail in Sect. \ref{ssec: recollection perv ULA}. The category $\Sat^q_X$ (resp. $\Rep(H)^{\Ab,\loc,\ft}_{\cG_Z}(X)$) is the relative (to $X$) version of the category $\Sat^q_x$ (resp. $\Rep(H)_{(\cG_Z)_x}$).\footnote{In \cite{ZhaoNew}, the category $\Rep(H)^{\Ab,\loc,\ft}_{\cG_Z}(X)$ was called the category of \emph{twisted $H$-representations}.} The above result has a natural analog relative to $X^I$ for $I\in\fSet$.

\sssec{Factorizable action of a local system}
From the above discussion, it is clear that the paradigm of Sect. \ref{sssec: fact hlgy informal} is not applicable. However, we develop a variant of it replacing the symmetric monoidal category $\cC$ with a \emph{local system} of symmetric monoidal categories. We motivate this below.
\medskip

As a background, we elaborate on what $\Rep(H)_{(\cG_Z)_x}$ is. $\cG_Z$ is an \'etale $Z_H(E)$-gerbe on the curve $X$, and we denote its fiber at $x\in X$ by $(\cG_Z)_x$. The category $\Rep(H)$ has a symmetric monoidal $Z_H(E)$-symmetry, hence a $Z_H(E)$-gerbe yields a symmetric monoidal twist of $\Rep(H)$. The symmetric monoidal category $\Rep(H)_{(\cG_Z)_x}$ is the twist of $\Rep(H)$ by $(\cG_Z)_x$. 
\medskip

Now, the idea is that the family $\{\Rep(H)_{(\cG_Z)_x}\}_{x\in X}$ naturally forms a local system of symmetric monoidal categories. The composition \footnote{Here, we use the language of $\infty$-categories, so the homotopy type $X^{\topl,\sing}$ is an $\infty$-groupoid.}
\begin{equation}\label{eq: loc sys smc RepH}
    X^{\topl,\sing}\xrightarrow{\Phi(\cG_Z)}B^2 Z_H(E)\xrightarrow{\Rep(H)}\DGCatSymMon
\end{equation}
is a local system of symmetric monoidal  DG categories over $X^{\topl,\sing}$. Here $\Phi(\cG_Z)$ is the topological gerbe corresponding to the \'etale gerbe $\cG_Z$, $\DGCatSymMon$ is the category of symmetric monoidal DG categories, and the second arrow is precisely the symmetric monoidal $Z_H(E)$-symmetry of $\Rep(H)$. The fiber of \eqref{eq: loc sys smc RepH} at $x\in X$ is indeed $\Rep(H)_{(\cG_Z)_x}$, and this is the category acting on $\cA^q_G$ at the point $x$.\medskip

Through Sect. \ref{sec: twisted fact hlgy} and \ref{sec: action of twisted fact hlgy} we develop an analog of Sect. \ref{sssec: fact hlgy informal}. Let $\sfA$ be an abelian group, $Y\in\Spc$, and $\cG$ be an $\sfA$-gerbe on $Y$. Given a symmetric monoidal category $\cC$ with an $\sfA$-symmetry, we define a symmetric monoidal category $\int_{(Y,\cG)}\cC$ which we call the \textbf{twisted factorization homology}. \footnote{In the bulk of the paper, we will use the notation $\int^{\CAlg}_{(Y,\cG)}\cC$ instead.} 
We can form a local system of symmetric monoidal categories
\begin{equation}\label{eq: loc sys smc C}
    Y\xrightarrow{\cG}B^2\sfA\xrightarrow{\cC}\DGCatSymMon
\end{equation}
generalizing \eqref{eq: loc sys smc RepH}. We define the notion of a factorizable action of a local system of symmetric monoidal categories, and prove:

\begin{thm}[see Theorem \ref{thm: main categorical thm}]
    A factorizable action of the local system \eqref{eq: loc sys smc C} on a DG category $\cM$ is equivalent to an action of $\int_{(Y,\cG)}\cC$ on $\cM$.
\end{thm}

Moreover, we prove the following result, which implies that the construction of Betti metaplectic spectral action is equivalent to the construction of a factorizable action of \eqref{eq: loc sys smc RepH} on $\cA^q_G$.

\begin{thm}[see Theorem \ref{thm: QC over twlocsys as a tw fact hlgy}]
    The twisted factorization homology of \eqref{eq: loc sys smc RepH} is canonically identified with $\QC(\Loc^{\cG_Z}_H)$.
\end{thm}

\sssec{Factorizable action of \eqref{eq: loc sys smc RepH}}
Through Sect. \ref{sec: metaplectic Satake functor}, \ref{sec: Satake as a nat transf}, and \ref{sec: factorizable action} we construct the factorizable action of \eqref{eq: loc sys smc RepH} on $\cA^q_G$. The idea is to construct the relevant data of the factorizable action \emph{after taking the global section}.
\medskip

Given a homotopy type $Y$, let $\LSMon(Y)$ be the category of local systems of monoidal DG categories, and $\LS(Y)$ be the DG category of $E$-linear local systems.
Denote the local system \eqref{eq: loc sys smc RepH} by $\ul{\Rep(H)}_{X,\cG_Z}$. Suppose we have the desired factorizable action of $\ul{\Rep(H)}_{X,\cG_Z}$ on $\cA^q_G$. One of the requirements for the factorizable action is local constancy, and this is encoded in a map
\begin{equation}\label{eq: map in LSMon}
    \ul{\Rep(H)}_{X,\cG_Z}\rightarrow \ul{\End(\cA^q_G)}_X\quad \te{in }\LSMon(X).
\end{equation}
Indeed, if we evaluate this map at $x\in X$, we get a monoidal map $\Rep(H)_{(\cG_Z)_x}\rightarrow \End(\cA^q_G)$ which is supposed to be the metaplectic Hecke action at $x$. Local constancy is encoded in the enhancement of the collection of those pointwise actions into a map in $\LSMon$. Under the global section functor $\LSMon(Y)\rightarrow \Alg(\LS(Y)\lMod)$, the map \eqref{eq: map in LSMon} turns into a map
\begin{equation}\label{eq: map in Alg LS Mod}
    \Rep(H)^{\loc}_{\cG_Z}(X)\rightarrow \End(\cA^q_G)\otimes\LS(X)\quad \te{in }\Alg(\LS(X)\lMod).
\end{equation}

A factorizable action requires actions parametrized by any $X^I$'s and compatibility between them, and this motivates us to prove:
\begin{thm}[see Theorem \ref{thm: fact Satake nat transf}]\label{thm: fact action after global section}
    There is a compatible family of canonical functors
    \begin{equation*}
        \left(\Rep(H)^{\loc}_{\cG_Z}(X)\right)^{\otimes I}\rightarrow \End(\cA^q_G)\otimes\LS(X^I)\quad \te{in }\Alg(\LS(X^I)\lMod)
    \end{equation*}
    for each $I\in\fSet$.
\end{thm}

This is the `global section' of the desired categorical data of the factorizable action.

\sssec{1-affineness}
The global section functor $\LSMon(Y)\rightarrow \Alg(\LS(Y)\lMod)$ is not an equivalence in general, so the above theorem does not a priori yield the factorizable action of \eqref{eq: loc sys smc RepH}. Nevertheless, we use a special feature of topology of Riemann surfaces to show that we can go backwards - from global sections to local systems - to obtain the desired factorizable action, concluding the proof of Theorem \ref{thm: main}. The crucial input for this backward implication is the 1-affineness result for Betti stacks, see Theorem \ref{thm: 1-affineness}.

\sssec{Factorizable metaplectic Satake}
Let $\Hecke^{\loc}_{X^I}$ be the \emph{local Hecke stack}, which is the relative (to $X^I$) version of $\Hecke^{\loc}_x$. The metaplectic parameter $q$ yields a canonical $A$-gerbe on the local Hecke stack, and hence we have the category $\Shv^q(\Hecke^{\loc}_{X^I})$ of twisted sheaves. As in the non-metaplectic case, this has a natural monoidal structure and acts on $\Shv^q(\Bun_G)$ via the Hecke correspondence. The proof of Theorem \ref{thm: fact action after global section} relies on the following two results:
\begin{thm}[see Theorem \ref{thm: fact Satake nat transf}]\label{thm: fact Satake intro}
    There is a compatible family of canonical functors
    \begin{equation}\label{eq: LS linear met Satake}
        \left(\Rep(H)^{\loc}_{\cG_Z}(X)\right)^{\otimes I}\rightarrow \Shv^q(\Hecke^{\loc}_{X^I})\quad\te{in }\Alg(\LS(X^I)\lMod)
    \end{equation}
    for each $I\in\fSet$.
\end{thm}

\begin{thm}[see Theorem \ref{thm: preservance Nilp metaplectic}]\label{thm: preservance intro}
    The full subcategory $\cA^q_G\subset \Shv^q(\Bun_G)$ is preserved by the Hecke action by the essential image of \eqref{eq: LS linear met Satake}.
\end{thm}

Theorem \ref{thm: preservance intro} is the metaplectic analogue of \cite[Theorem 5.1]{NY} and \cite[Theorem B.5.2]{GKRV}, known as the \emph{preservation of nilpotent singular support}. The proof in the non-metaplectic case carries over to the metaplectic case without much difference.
\medskip

A much more interesting part is Theorem \ref{thm: fact Satake intro}. We call the family of functors \eqref{eq: LS linear met Satake} the \emph{$\LS$-linear} metaplectic Satake functors, to distinguish it from the abelian equivalence (Theorem \ref{thm: abelian Satake light}). The proof of Theorem \ref{thm: fact Satake intro} occupies Sect. \ref{sec: metaplectic Satake functor} and \ref{sec: Satake as a nat transf}, and this is the most technical part of this paper. $\LS$-linear Satake functors look similar to the abelian equivalence, because they have similar origins and these abelian categories sit inside the derived categories:
\[\Rep(H^I)^{\Ab,\loc,\ft}_{\cG^{\boxtimes I}_Z}(X^I)\hookrightarrow\left(\Rep(H)^{\loc}_{\cG_Z}(X)\right)^{\otimes I},\quad \Sat^q_{X^I}\hookrightarrow \Shv^q(\Hecke^{\loc}_{X^I}).\]
In fact, we deduce the theorem using the abelian equivalence. However, the passage from the abelian equivalence to the $\LS$-linear version requires significant efforts. Even for the $I=\{*\}$ case, the passage from $\Rep(H)^{\Ab,\loc,\ft}_{\cG_Z}(X)$ to $\Rep(H)^{\loc}_{\cG_Z}(X)$ involves passing to derived categories and removing the finiteness (of local systems) at once. The presence of twisting by $\cG_Z$ adds extra complexity. We resolve these difficulties using various categorical manipulation, plus the monoidal grading \footnote{The author is grateful to Y. Zhao for teaching this to him.} of $\Shv^q(\Hecke^{\loc}_X)$. For this reason, our construction of functors \eqref{eq: LS linear met Satake} is involved, so constructing the appropriate compatibility data among them poses another challenge. Our approach is to uniquely characterize the $\LS$-linear functors in terms of the abelian equivalence, and transport the compatibility data of abelian Satake equivalences to the desired compatibility data. We refer readers to the introduction of Sect. \ref{sec: metaplectic Satake functor} for a more detailed overview.

\subsection{Acknowledgements}
First and foremost, I would like to thank my advisor, David Ben-Zvi, for his guidance and support over the past few years, and in particular for introducing me to the subject of metaplectic geometric Langlands. I am grateful to Yifei Zhao for illuminating conversations, which greatly shaped my understanding of this topic. I also thank Sam Raskin, Mark Saving, Germ\'an Stefanich, and Philip Thomas for helpful discussions and email correspondence. This work was supported by the Korea Foundation for Advanced Studies Overseas PhD Scholarship and the Provost’s Graduate Excellence Fellowship at the University of Texas at Austin.

\subsection{AI use statement}
The author acknowledges the use of OpenAI's ChatGPT during the development of Sect. \ref{sssec: six functor extension}, specifically to check whether the class of morphisms $\widetilde{E}$ constructed in \cite[Theorem 5.19]{Sch}, in the context of six-functor formalisms, includes a specific type of morphisms. The author independently checked and verified the resulting mathematical claim. In addition, Google Gemini and OpenAI's ChatGPT were used for text editing and sentence polishing.

\section{Preliminaries}
\subsection{Category theory}
We use the language of $\infty$-categories and higher algebra, in the sense of \cite{HTT, HA}. Throughout this text, categories are $(\infty,1)$-categories by default. When we are talking about $(1,1)$-categories, we will use the term \emph{ordinary categories}. 
Let $\Spc$ be the category of homotopy types, i.e., $\infty$-groupoids, $\Cat$ be the category of categories. Fix $E$ to be an algebraically closed field of characteristic zero. Let $\Vect$ be the stable category of chain complexes of $E$-vector spaces.

\sssec{DG categories}
Let $\DGcat$ be the category of $E$-linear presentable stable categories with morphisms being continuous $\Vect$-linear functors. By a \textbf{DG category}, we mean an object in $\DGcat$. We let $\DGCatMon$ be the category of monoidal DG categories and $\DGCatSymMon$ be the category of symmetric monoidal DG categories. We have
\[\DGCatMon:=\Alg(\DGcat)\quad\te{and}\quad \DGCatSymMon:=\CAlg(\DGcat).\]

Let $\Vect^{\on{f.d.}}$ be the category of perfect complexes in $\Vect$.  Let $\on{DGCat}^{\omega}$ be the category of small idempotent complete stable category with an exact $\Vect^{\on{f.d.}}$-action. The $\Ind$-completion is organized into a functor
\[\Ind: \on{DGCat}^{\omega}\rightarrow \DGcat. \]

By a \textbf{non-cocomplete DG category}, we mean an object of $\on{DGCat}^{\omega}$ which is not cocomplete. Note that any category in $\DGcat$ is cocomplete by presentability assumption. $\Vect^{\on{f.d.}}$ is an example of a non-cocomplete DG category.

\sssec{Operads}
We also use the theory of operads \footnote{By our convention, they are $\infty$-operads.} developed in \cite{HA}. The most frequently appearing operads are $\E_1$ and $\E_{\infty}$-operads. Given a symmetric monoidal category $\bC$ we have
\[\Alg(\bC)=\E_1\lAlg(\bC),\quad \CAlg(\bC)=\E_{\infty}\lAlg(\bC).\]
We point out two basic facts about operads that will be used frequently.

\sssec{Adjunction and algebra}\label{sssec: adjunction and algebra}
The following is borrowed from \cite[Part II, Ch.6, Sect. 1.2.3]{GR2}.
Let $\Phi:\bC\rightarrow\bC'$ be a morphism in $\DGCatSymMon$ and $\cP$ be an operad. $\Phi$ induces a functor $\Phi:\cP\lAlg(\bC)\rightarrow \cP\lAlg(\bC')$. Suppose $\Phi$ admits a right adjoint $\Phi^R$ in $\DGcat$. $\Phi^R$ has a natural right-lax symmetric monoidal structure, and hence induces a functor $\Phi^R: \cP\lAlg(\bC')\rightarrow \cP\lAlg(\bC)$. These two functors fit into the adjunction
\[\Phi: \cP\lAlg(\bC)\leftrightarrows \cP\lAlg(\bC'): \Phi^R.\]

\sssec{Tensor product of operads}
In \cite[Sect. 2.2.5]{HA} the tensor product of operads is discussed. Given two operads $\cP$ and $\cP'$, one can associate an operad $\cP\otimes\cP'$ such that there is a functorial equivalence \[\cP\lAlg(\cP'\lAlg(\bC))\simeq (\cP\otimes\cP')\lAlg(\bC)\]
for any symmetric monoidal category $\bC$. This operation endows the category operads with a symmetric monoidal structure. In particular, we have an equivalence
\[\cP\lAlg(\cP'\lAlg(\bC))\simeq \cP'\lAlg(\cP\lAlg(\bC)).\]

\subsection{Stacks and prestacks}
\sssec{Automorphic side}\label{sssec: AG convention}
On the automorphic side, we work with \emph{classical} algebraic geometry over the base field $\C$. Let $\Aff$ be the category of \emph{classical} affine schemes over $\C$, and $\Affft$ is the category of classical affine schemes of finite type over $\C$. Those categories are equipped with \'etale topology, hence we get the categories
\[\PreStk=\PSh(\Aff),\quad \PreStklft=\PSh(\Affft),\quad \Stk=\Shv(\Aff),\quad \NearStklft=\Shv(\Affft).\]
We have an adjunction
\[\on{LKE}:\PreStklft\leftrightarrows \PreStk:\on{Res}\]
where both the left-Kan extension and restriction are considered along $\Affft^{\op}\hookrightarrow \Aff^{\op}$. This left-Kan extension functor realizes $\PreStklft$ as a full subcategory of $\PreStk$.
We have adjunctions 
\[L:\PreStk\leftrightarrows \Stk:\iota\quad \te{and}\quad L_{\ft}:\PreStklft\leftrightarrows \NearStklft: \iota\]
where $\iota$'s are inclusions and $L$ (resp. $L_{\ft}$) is the \'etale sheafification.
\medskip

We follow the convention of \cite{GR} that a scheme is assumed to be separated. (see Remark 3.1.3 of \emph{loc. cit.}) We let $\on{Sch}^{\ft}$ be the category of schemes of finite type over $\C$. We denote $\Spec \C$ by $\pt$. An algebraic stack is assumed to have an affine diagonal.

\sssec{}
From \cite[Sect. 1.3.4]{GL} the stack $\Bet^i(A):=L(B^i A)\in\Stk\subset\PreStk$ actually belongs to $\PreStklft$. Therefore $\Bet^i (A)\in\Stklft$. Its restriction to $\Affft^{\op}$ satisfies descent, hence $\Bet^i (A)|_{\Affft^{\op}}\in\NearStklft$. We can give another description of this object. Consider the prestack $B^i A\in\PreStklft$ and its sheafification $L_{\ft}(B^i A)\in\NearStklft$. Then from \cite[Corollary 2.7.5]{GR} we have a canonical equivalence $\Bet^i (A)|_{\Affft^{\op}} = L_{\ft}(B^i A)$.

\sssec{Spectral side}
The field $E$ is the coefficient field for sheaves on the automorphic side. On the spectral side we work with \emph{derived} prestacks and stacks over the base field $E$. To distinguish algebraic geometry over two different base fields, let us put the superscript $E$ on the categories arising from algebraic geometry over $E$, while we don't put any superscript for complex algebraic geometry. Hence $\Aff^E$ is the category of \emph{derived} affine schemes over $E$, and $\PreStk^E$ is the category of \emph{derived} prestacks over $E$. As opposed to the automorphic side, we would not need to consider their variants with finiteness condition. We denote $\Spec E$ by $\pt$.

\subsection{Sheaf theory}
\sssec{Betti sheaves}
Given a scheme $Y\in \Schft$, we can consider its underlying complex analytic topological space $Y^{\topl}$. Let $\Shv(Y)$ be the DG category of \emph{all} \footnote{We do not impose constructibility, nor any cohomological boundedness.} complexes of sheaves of $E$-vector spaces over $Y^{\topl}$. Given a prestack $\cY\in\PreStklft$, we define
\begin{equation}\label{eq: Shv as a RKE}
    \Shv(\cY):=\lim_{({\Affft}_{/\cY})^{\op}}\Shv(S)\quad \te{in }\DGcat,
\end{equation}
where we use $*$-pullback functors as transition maps. This is an instance of right-Kan extension. More precisely, the definition of $\Shv$ for schemes and the functoriality under $*$-pullbacks gives rise to a functor $\Shv:(\Affft)^{\op}\rightarrow\DGcat$, and we right-Kan extend this functor along  $(\Affft)^{\op}\hookrightarrow(\PreStklft)^{\op}$ to get a functor $\Shv:(\PreStklft)^{\op}\rightarrow \DGcat$. The above limit expression of $\Shv(\cY)$ is exactly the evaluation of the functor $\Shv$ at $\cY$. The $*$-tensor product endows a symmetric monoidal structure on $\Shv(S)$, and the limit \eqref{eq: Shv as a RKE} can be equally understood in $\DGCatSymMon$. This defines the symmetric monoidal structure on $\Shv(\cY)$.
\medskip

Objects in the category $\Shv(\cY)$ will be called \textbf{Betti sheaves}. Our definition of Betti sheaves is different from \cite[Sect. A.1]{GKRV}, in the sense that they use $!$-pullbacks to define the right-Kan extension and work with the \emph{discontinuous} version of $\DGcat$. Nevertheless, our definition and the definition of \emph{op. cit.} agree on algebraic stacks. Indeed, in Sect. A.1.5 of \emph{op. cit.}, it is argued using a smooth descent argument that if $\cY$ is an algebraic stack, then we can restrict the index category $({\Affft}_{/\cY})^{\op}$ of the limit \eqref{eq: Shv as a RKE} to the non-full subcategory $({\Affft}_{/\cY,\on{sm}})^{\op}$ where we admit only smooth $S\rightarrow \cY$'s and smooth maps between those. Therefore, to see those two definitions agree on an algebraic stack $\cY$, we have to show that the two formation of limits
\[\lim_{({\Affft}_{/\cY,\on{sm}})^{\op}}\Shv(S),\]
one using $*$-pullbacks and the other using $!$-pullbacks, agree with each other. This follows from the fact that $*$ and $!$-pullbacks for a smooth map differ by a cohomological shift determined by relative dimension.

\sssec{Local system $\LS$}
On a scheme $Y\in\Schft$ we have a full DG subcategory $\LS(Y)$ of $\Shv(Y)$, consisting of \emph{all} local systems, in the sense that we do not impose any boundedness nor finiteness. $\LS(Y)$ is characterized as Betti sheaves whose cohomology sheaves are local systems. It inherits a symmetric monoidal structure from $\Shv(Y)$. Denote $Y^{\topl,\sing}\in\Spc$ by the underlying homotopy type of the topological space $Y^{\topl}$. Then we have a canonical equivalence
\[\LS(Y)\simeq\lim_{Y^{\topl,\sing}}\Vect \quad\te{in }\DGcat,\]
where we take the diagram $Y^{\topl,\sing}\rightarrow \DGcat$ with constant value $\Vect$. (\cite[Lemma A.4.2]{GKRV})

\sssec{Ind-constructible sheaves}
Another sheaf theory that we use is that of ind-constructible sheaves. Given a scheme $Y\in\Schft$, we consider the non-cocomplete DG category $\Shv^{\on{constr}}(Y)$ of bounded constructible sheaves on $Y$, still using the complex analytic topology. We define 
\[\Shv^{\indconstr}(Y):=\Ind(\Shv^{\on{constr}}(Y)).\]
Similar to the case of Betti sheaves, we extend this definition to prestacks using right-Kan extension and $*$-pullbacks. Hence, for a prestack $\cY\in\PreStklft$, we have
\begin{equation}\label{eq: Shvindconstr via a RKE}
    \Shv^{\indconstr}(\cY):=\lim_{({\Affft}_{/\cY})^{\op}}\Shv^{\indconstr}(S)\quad\te{in }\DGcat,
\end{equation}
using $*$-pullbacks. Using the same logic as before, this agrees with the definition of ind-constructible sheaf category in \cite[Sect. A.1]{GKRV} over algebraic stacks, and we still have a symmetric monoidal structure on $\Shv^{\indconstr}(\cY)$ via $*$-tensor product.
\medskip

Although we are primarily interested in Betti sheaves, ind-constructible sheaves play several important roles, because $\Shv^{\indconstr}$ is `smaller' than $\Shv$, and hence it is easier to control. To be more precise, for an algebraic stack $\cY$, the DG category $\Shv^{\indconstr}(\cY)$ is compactly generated, and its compact objects are constructible, see \cite[Sect. F.1]{AGKRRV}. In contrast, $\Shv(Y)$ is almost never compactly generated.

\sssec{Local systems $\QLisse$ and $\LS^{\indft}$}
Similar to $\LS$ being the category of local systems in $\Shv$, we can talk about local systems in $\Shv^{\indconstr}$. Since the definition of $\Shv^{\indconstr}$ involved (ind)finiteness, so does any notion of a local system therein. There are two versions of it: $\QLisse$ and $\LS^{\indft}$, both defined in \cite[Sect. 1.2]{AGKRRV}. \footnote{In \emph{loc. cit.}, the latter was denoted $\on{IndLisse}$.} The latter is easier to define. Given $Y\in\Schft$, we consider the non-cocomplete DG subcategory $\LS^{\ft}(Y)$ of $\LS(Y)$ consisting of those having finitely many nontrivial cohomology sheaves, each of which being finite type. We define $\LS^{\indft}(Y):=\Ind(\LS^{\ft}(Y))$, and the embedding $\LS^{\ft}(Y)\hookrightarrow \Shv^{\indconstr}(Y)$ induces an embedding $\LS^{\indft}(Y)\hookrightarrow \Shv^{\indconstr}(Y)$. $\LS^{\indft}(Y)$ inherits the symmetric monoidal structure from $\Shv^{\indconstr}(Y)$.
\medskip

$\QLisse(Y)$ is another symmetric monoidal DG subcategory of $\Shv^{\indconstr}(Y)$. Since its definition is a bit more complicated, we refer readers to \emph{loc. cit.} for details but only highlight key features that will be used in this paper. We will use $\QLisse$ only in the case $Y$ is a smooth algebraic curve. In that case, we have symmetric monoidal embeddings
\[\LS^{\indft}(Y)\hookrightarrow \QLisse(Y)\hookrightarrow\LS(Y).\]
If we assume $Y$ is connected, the first arrow is an equivalence unless $Y=\pjl$. If $Y=\pjl$, the second arrow is an equivalence. Moreover, the assignment $Y\mapsto \QLisse(Y)$ satisfies \'etale descent, while it fails for $Y\mapsto\LS^{\indft}(Y)$. Since we will work extensively with sheaf of categories, $\QLisse$ will be much more useful than $\LS^{\indft}$ for our purpose.

\section{Preliminaries on gerbes}\label{sec: prelim on gerbes}

\subsection{Basic definitions}\label{ssec: basic definitions}
We mainly use the notion of gerbes in the following two settings.

\sssec{\'Etale gerbes}
Say that $\sfA$ is a torsion abelian group and $\cY\in \PreStk$. The space of \textbf{\'etale} $\sfA$-\textbf{gerbes} is defined as the mapping space $\on{Ge}_{\sfA}(\cY):=\Maps(\cY, \Btet{\sfA})$ in $\PreStk$. If $\cX\in\Stk$, an \'etale gerbe on $\cY$ is equivalent to a principal $\Bet A$-bundle over $\cY$ in the topos $\Stk$.
In case $\cY\in\NearStklft$,
since $\Btet{\sfA}=L_{\ft}(B^2 \sfA)$ in $\NearStklft$, we conclude that an \'etale gerbe on $\cY$ is the same as a principal $\Bet A$-bundle over $\cX$ in the topos $\NearStklft$.
We have:
\[\pi_i\Maps_{\PreStklft}(\cX,\Bet^j \sfA)=H^{j-i}_{\et}(\cX,\sfA).\]

\sssec{Topological gerbes} Say that $\sfA$ is an abelian group. For each $i\ge 0$ we have $B^i \sfA\in \Spc$. Let $K$ be a homotopy type, i.e., an object in $\Spc$. The space of \textbf{topological} $\sfA$-\textbf{gerbes} is defined as the mapping space $\on{Ge}_{\sfA}(\cY):=\Maps(K, \Bt{\sfA})$ in $\Spc$.
We have
\[\pi_i \Maps_{\Spc}(K,B^j \sfA)=H^{j-i}_{\sing}(K,\sfA).\]

\sssec{Basic operations with gerbes}
We list some basic operations with gerbes. Although we discuss in the context of \'etale gerbes, all the discussions below apply equally to topological gerbes. First of all, the mapping stack $\Maps(\cY,\Bet^2\sfA)$ has a point $\cY\rightarrow\pt\rightarrow\Bet^2 \sfA$, the second map being the neutral point. This gerbe is called the \textbf{trivial gerbe}, and denoted by $\triv_{\cY}$.
\medskip

Suppose we have a group homomorphism $\sfA_1\rightarrow\sfA_2$ between torsion abelian groups. Given an \'etale $\sfA_1$-gerbe $\cG$ on $\cY$ we can inflate it to get an \'etale $\sfA_2$-gerbe. It can be understood as the composition $\cY\xrightarrow{\cG}\Bet^2 \sfA_1\rightarrow \Bet^2 \sfA_2$.
\medskip

Suppose we have \'etale $\sfA$-gerbes $\cG_1$ and $\cG_2$ on $\cY$. We can take the \textbf{sum} $\cG_1+ \cG_2$, which is the \'etale $\sfA$-gerbe on $\cY$ obtained by the composition
\[\cY\xrightarrow{(\cG_1,\cG_2)}\Bet^2\sfA\times\Bet^2\sfA\xrightarrow{\te{add}}\Bet^2\sfA.\]

Suppose we have an \'etale $\sfA_1$-gerbe $\cG_1$ on $\cY_1$ and an \'etale $\sfA_2$-gerbe $\cG_2$ on $\cY_2$. We can form the \textbf{external product} $\cG_1\boxtimes\cG_2$, which is the \'etale $\sfA_1\times\sfA_2$-gerbe on $\cY_1\times\cY_2$ defined as the composition
\[\cY_1\times\cY_2\xrightarrow{\cG_1\times\cG_2}\Bet^2 \sfA_1\times\Bet^2\sfA_2 =\Bet^2(\sfA_1\times\sfA_2).\]

Suppose we have an \'etale $\sfA$-gerbe $\cG_1$ on $\cY_1$ and an \'etale $\sfA$-gerbe $\cG_2$ on $\cY_2$. We can form the \textbf{external sum} $\cG_1\boxplus\cG_2$, which is the \'etale $\sfA$-gerbe on $\cY_1\times\cY_2$ defined as the composition
\[\cY_1\times\cY_2\xrightarrow{\cG_1\times\cG_2}\Bet^2 \sfA\times\Bet^2\sfA \xrightarrow{\on{add}}\Bet^2\sfA.\]

Suppose we have an \'etale $\sfA$-gerbe $\cG$ on $\cY$. We can form the \textbf{total space} $\on{Tot}(\cG)$ of the gerbe, defined as the fiber product along $\cY\xrightarrow{\cG}\Bet^2\sfA\leftarrow\pt$, where the second arrow is the neutral point. There is a canonical projection $\on{Tot}(\cG)\rightarrow\cY$, and there is a canonical $\Bet \sfA$-action on $\on{Tot}(\cG)$ preserving the projection.

\sssec{Category of pairs}
We consider the category
\[\PreStklft+\on{Ge}_{\sfA}:=(\PreStklft)_{/\Bet^2\sfA}.\]
This is the category of pairs $(\cY,\cG)$ where $\cY\in\PreStklft$ and $\cG$ is an \'etale $\sfA$-gerbe on it. A map $(\cY_1,\cG_1)\rightarrow (\cY_2,\cG_2)$ between pairs is the data of $f:\cY_1\rightarrow\cY_2$ together with an isomorphism $f^*\cG_2\simeq \cG_1$. External sum endows this category with a symmetric monoidal structure:
\[(\cY_1,\cG_1)\otimes(\cY_2,\cG_2)= (\cY_1\times\cY_2,\cG_1\boxplus\cG_2).\]
In the topological setup, we can consider the category
\[\Spc+\on{Ge}_{\sfA}:=\Spc_{/B^2 \sfA},\]
also equipped with the symmetric monoidal structure endowed from external sum.

\subsection{Twisting a sheaf of categories}\label{ssec: twisting a sheaf of categories}
The key feature of a gerbe is that it can twist a (sheaf of) categories. We give an informal sketch of this idea, although we provide a rigorous and systematic treatment in Sect. \ref{sec: twisting construction}.

\sssec{An absolute version of twisting}\label{sssec: twisting absolute}
From the definition of gerbe above, a gerbe is a 1-categorical higher notion of a torsor. Given a group $\mathsf{G}$, a $\mathsf{G}$-torsor $\cP$ and a $\mathsf{G}$-space $\mathsf{Y}$, we can twist $\mathsf{Y}$ by the torsor $\cP$ to get a space $\mathsf{Y}_{\cP}$. Similarly, given an abelian group $\mathsf{A}$, a $\mathsf{A}$-gerbe $\cG$ (over a point) and a category $\cC$ with an $\sfA$-symmetry, then we can twist $\cC$ by the gerbe $\cG$ to get a category $\cC_{\cG}$. An $\mathsf{A}$-symmetry on $\cC$ means that we have a compatible family of automorphisms $f_{a,c}\in\Aut_{\cC}(c)$ for $c\in \cC$ and $a\in \sfA$, such that for each $c$ the map $A\rightarrow \Aut_{\cC}(c)$, $a\mapsto f_{a,c}$ is a group homomorphism, and for each $a\in \sfA$ the automorphisms $f_{a,c}$ intertwine every morphisms in $\cC$. As an example, given any $\cY\in\PreStklft$, the categories $\Shv(\cY)$ and $\Shv^{\indconstr}(\cY)$ have natural $E^{\times}$-symmetry, because each sheaf has a scaling automorphism by $E^{\times}$, compatible with morphisms.

\sssec{A relative/sheafy version of twisting}\label{sssec: sheafy twisting}
The above twisting procedure  has a natural sheafy version. Suppose we have a category $\bO$ with finite limits and a site $\cX$, which is either a small \'etale site $Y_{\et}$ of $Y\in\Schft$ or the overcategory $(\Affft)_{/\cY}$ for some $\cY\in\PreStklft$ equipped with \'etale topology. \footnote{We have made specific choices of the site $\cX$ which suffice for our purpose, but the twisting procedure applies for a much more general class of sites, see Sect. \ref{sssec: conditions on cX bO}.} We have categories $\PSh(\cX;\bO)$ and $\Shv(\cX;\bO)$ of presheaves on $\cX$ valued in $\bO$ and sheaves on $\cX$ valued in $\bO$, respectively. We further assume the existence of sheafification functor, i.e., the left adjoint to the inclusion $\Shv(\cX;\bO)\hookrightarrow \PSh(\cX;\bO)$. This is the case when $\bO$ is $\Spc$, $\Cat$ or $\DGcat$ for example. An object in $\bO$ with an $\sfA$-symmetry is an object $\cO\in\bO$ with a monoidal map $B\sfA\rightarrow \Aut_{\bO}(\cO)$ \footnote{When $\bO=\Cat$, we have a notion of a category $\cO$ having an $\sfA$-symmetry. This is equivalent to the definition given in Sect. \ref{sssec: twisting absolute}.}
In this context, a sheaf with an $\sfA$-symmetry is a sheaf $\cF\in \Shv(\cX;\bO)$ such that each section $\cF(x)$ for $x\in\cX$ has an $\sfA$-symmetry compatible with restriction maps. Given such a sheaf $\cF$ with an $\sfA$-symmetry and an \'etale $\sfA$-gerbe $\cG$ over $Y$ (if $\cX=Y_{\et}$) or an \'etale $\sfA$-gerbe $\cG$ over $\cY$ (if $\cX=(\Affft)_{/\cY}$), we can twist $\cF$ by $\cG$ to get $\cF_{\cG}\in \Shv(\cX;\bO)$.

\begin{rmk}
    We use the terms `$\sfA$-symmetry' and `$\sfA$-linearity' interchangeably. The precise meaning of $\sfA$-linearity is reviewed in Sect. \ref{sssec: twisting construction}.
\end{rmk}

\sssec{Twisted sheaves}\label{sssec: twisted sheaves brief}
As an example, take $\cX=(\Affft)_{/\cY}$, $\bO=\DGcat$, and $\mathbf{Shv}\in\Shv(\cX;\bO)$ be
\[\left((\Affft)_{/\cY}\right)^{\op}\rightarrow \DGcat,\quad S\mapsto \Shv(S).\]
This has a natural $E^{\times}$-symmetry, because each of $\Shv(S)$ has a $E^{\times}$-symmetry by scaling, and this structure is compatible with the $*$-pullback functors. Therefore, given a finite subgroup $A\hookrightarrow E^{\times}$ and an \'etale $A$-gerbe on $\cY$ we get a new sheaf of DG categories $\mathbf{Shv}_{\cG}$ on $(\Affft)_{/\cY}$. The \textbf{category of twisted sheaves} is defined to be the global section of this sheaf of DG categories:
\[\Shv^{\cG}(\cY):=\Gamma(\cY; \textbf{Shv}_{\cG}).\]
The above construction equally applies to ind-constructible sheaves, and hence given an \'etale $A$-gerbe over $\cY$ we form  the category $\Shv^{\indconstr,\cG}(\cY)$ by the same procedure.
Twisted sheaves defined as above are main players in this paper. Its more formal definition is given in Sect. \ref{sssec: twisted sheaves}, and we study various properties of twisted sheaves in Sect. \ref{sec: twisted sheaves}.

\sssec{Functoriality of twisted sheaves}\label{sssec: functoriality twshv brief}
Suppose we have a map $f:\cY_1\rightarrow\cY_2$ in $\PreStklft$. We have a $*$-pullback functor $f^*:\Shv(\cY_2)\rightarrow\Shv(\cY_1)$. If $f$ is schematic, we have a $!$-pushforward functor $f_!:\Shv(\cY_1)\rightarrow\Shv(\cY_2)$. We have analogous functoriality for twisted sheaves. However, in the twisted setup, it is necessary to provide a compatibility data of gerbes in order to have a functor between twisted sheaf categories. Let $A\hookrightarrow E^{\times}$ be an abelian group as before and suppose we have \'etale $A$-gerbes $\cG_1$ and $\cG_2$ on $\cY_1$ and $\cY_2$, respectively. Given $f:\cY_1\rightarrow\cY_2$ and an isomorphism $f^*\cG_2\simeq \cG_1$ we have a $*$-pullback functor $f^*:\Shv^{\cG_2}(\cY_2)\rightarrow \Shv^{\cG_1}(\cY_1)$. If $f$ is in addition schematic, we have a $!$-pushforward functor $f_!:\Shv^{\cG_1}(\cY_1)\rightarrow \Shv^{\cG_2}(\cY_2)$. 
In fact, these two functorialities get enhanced into a six-functor formalism for twisted sheaves. To be more precise, we consider the category of correspondences \[\Corr(\PreStklft+\on{Ge}_A,\on{sch})\] where we allow pushforwards along schematic maps. The external sum endows this category with a symmetric monoidal structure, and we have a right-lax symmetric monoidal functor
\[\Shv^{\on{tw}}_{\Corr}:\Corr(\PreStklft+\on{Ge}_A,\on{sch})\rightarrow\DGcat\]
encoding both functorialities.
In particular, the base change formula holds for twisted sheaves given that we have appropriate compatibility data of gerbes, and we have an external tensor product functor
\begin{equation}\label{eq: box tensor}
    \boxtimes:\Shv^{\cG_1}(\cY_1)\otimes\Shv^{\cG_2}(\cY_2)\rightarrow \Shv^{\cG_1\boxplus\cG_2}(\cY_1\times\cY_2).
\end{equation}
We refer the reader to Sect. \ref{ssec: functoriality of twisted sheaves} for details and the construction of $\Shv^{\on{tw}}_{\Corr}$. Finally, all the above discussion equally applies to ind-constructible sheaves.

\sssec{Twisting sheaves of algebras and modules}\label{sssec: twisting example}
For our purpose, it is very important to twist sheaves with additional structures. Consider the site $\cX$ and the category $\bO$ as in Sect. \ref{sssec: sheafy twisting} and suppose $\bO$ has a symmetric monoidal structure whose tensor product commutes with colimits in each variables. In this case, the category $\PSh(\cX;\bO)$ has a pointwise symmetric monoidal structure, and the category $\Shv(\cX;\bO)$ has a symmetric monoidal structure, uniquely determined from the requirement that the sheafification functor $\PSh(\cX;\bO)\rightarrow \Shv(\cX;\bO)$ is symmetric monoidal. Therefore, we can talk about (operadic) algebra objects and module objects in that category. Given such an object/sheaf with an $\sfA$-symmetry compatible with the algebra/module structure, we can twist it by an $\sfA$-gerbe. We refer the reader to Sect. \ref{sec: twisting construction} for details, but we give some examples and non-examples of this.
\medskip

For example, the sheaf $\textbf{Shv}$ considered in Sect. \ref{sssec: twisted sheaves brief} is naturally a commutative algebra object in $\Shv((\Affft)_{/\cY};\DGcat)$. We can consider $\textbf{Shv}$ as a module over itself. The twisting of Sect. \ref{sssec: twisted sheaves brief} preserves the $\textbf{Shv}$-module structure, and hence the twist $\textbf{Shv}^{\cG}$ has a natural $\textbf{Shv}$-module structure. The action map evaluated at each $S\in (\Affft)_{/\cY}$ yields a map $\Shv(S)\otimes \Shv^{\cG}(S)\rightarrow\Shv^{\cG}(S)$. This is in fact a composition 
\[\Shv(S)\otimes \Shv^{\cG}(S)\xrightarrow{\boxtimes}\Shv^{\triv_{S}\boxplus\cG}(S\times S)\rightarrow\Shv^{\cG}(\S).\]
The first map is the box tensor product \eqref{eq: box tensor}, and the second map is a $*$-pullback associated to the map $(S,\cG|_S)\rightarrow (S\times S,\triv_{S}\boxplus \cG|_S)$, where this map is determined by $\Delta:S\rightarrow S\times S$ together with the canonical isomorphism $\Delta^*(\triv_{S}\boxplus \cG|_S)\simeq \cG|_S$.
\medskip

However, the twisting procedure of Sect. \ref{sssec: twisted sheaves brief} does \emph{not} preserve the commutative algebra structure. Indeed, there is no natural way to define a symmetric monoidal structure of $\Shv^{\cG}(\cY)$, because we do not have a canonical isomorphism $\Delta^*(\cG\boxplus\cG)\simeq \cG$.
\medskip

Consider the symmetric monoidal DG category $\Rep(\mathsf{H})$ for an algebraic group $\mathsf{H}$ over $E$. Let $Z_{\mathsf{H}}(E)^{\on{tors}}$ be the torsion subgroup of $Z_{\mathsf{H}}(E)$. Then there is a natural $Z_{\mathsf{H}}(E)^{\on{tors}}$-symmetry on $\Rep(\mathsf{H})$. For each $a\in Z_{\mathsf{H}}(E)$ and $V\in\Rep(\mathsf{H})$ we have an action  $V\xrightarrow{a} V$ of $a\in \mathsf{H}(E)$, and since $a$ lies in the center, this is an automorphism as an $\mathsf{H}$-representation. Therefore, using the notation of Sect. \ref{sssec: twisting absolute}, we let $f_{a,V}= (V\xrightarrow{a} V)$. This is compatible with the symmetric monoidal structure of $\Rep(\mathsf{H})$, for example, we have $f_{a,V\otimes W}=f_{a,V}\otimes f_{a,W}$. Therefore, the twist by a $Z_{\mathsf{H}}(E)^{\on{tors}}$-gerbe preserves the symmetric monoidal structure, and hence a natural symmetric monoidal structure on the twist $\Rep(\mathsf{H})_{\cG}$ for any $Z_{\mathsf{H}}(E)^{\on{tors}}$-gerbe $\cG$. This is symmetric monoidally equivalent to $\QC((\Bet\mathsf{H})_{\cG})$, where the stack $(\Bet\mathsf{H})_{\cG}$ is the twist of $\Bet \mathsf{H}$ using the action $\Bet Z_{\mathsf{H}}(E)\acts \Bet \mathsf{H}$. Here we considered the $Z_{\mathsf{H}}(E)$-gerbe $\cG$ as a $\Bet Z_{\mathsf{H}}(E)$-bundle.

\subsection{Homotopy type of prestacks}
In this subsection, we define the notion of homotopy type for prestacks locally of finite type. For an affine scheme of finite type over $\C$, we can consider its underlying topological space and take its homotopy type. We can formally extend this construction to every prestack of locally finite type via left Kan extension.

\sssec{Construction}
Consider the functor $\Phi^{\on{Aff}}: \Affft\rightarrow \Spc$ given as the composition \[\Affft \xrightarrow{(\mhyphen)^{\topl}} \on{Top} \xrightarrow{(\mhyphen)^{\sing}} \Spc.\] Here the first arrow $(\mhyphen)^{\topl}$ is the passage to the underlying topological space of the scheme and the second arrow $(\mhyphen)^{\sing}$ is the passage to homotopy types.
\medskip

We construct the functor $\Phi: \PreStklft \rightarrow \Spc$ by taking the left Kan extension of the functor $\Phi^{\on{Aff}}: \Affft\rightarrow \Spc$ along the Yoneda embedding $\Affft\hookrightarrow\PreStklft$. From the generality of left Kan extension along the Yoneda embedding, $\Phi$ preserves colimits and $\Phi|_{\Affft}=\Phi^{\Aff}$. (see also Theorem \ref{thm: LKE and colimits})

\sssec{}
Say that $\sfA$ is a finite abelian group. Then $\Phi(\sfA)=\sfA$ where $\sfA$ on the left-hand side is considered as a group scheme, while $\sfA$ on the right-hand side is considered as a discrete topological space. The functor $B(\mhyphen)$ is expressed as a geometric realization, and it commutes with colimits. Therefore, we observe: 

\begin{lemma}\label{lem: Phi for BiA}
    For a torsion abelian group $\sfA$ considered as a group scheme, we have
    $\Phi(B^i \sfA)= B^i \sfA$ for every $i\ge 0$.
\end{lemma}

\sssec{Codescent} The functor $\Phi$ satisfies \'etale codescent, and this enables us to compute $\Phi$ for some stacks of our interest.
\begin{nota}
    Suppose that we are given a category $\cC$ with fiber products and a morphism $U\rightarrow V$ in $\cC$. For each $n\ge 0$ we define
    \[U^n/V:=U\times_V U\times_V\cdots\times_V U\]
    where $U$ appears $(n+1)$-times. We have the \emph{\v{C}ech nerve} $U^{\bullet}/V$ as a simplicial object
    \[\Delta^{\op}\rightarrow \cC,\quad [n]\mapsto U^n/V.\]
\end{nota}
\begin{lemma}
    The functor $\Phi: \PreStklft \rightarrow \Spc$ satisfies codescent, hence it factors through the sheafification functor $L_{\ft}:\PreStklft \rightarrow \NearStklft$. 
\end{lemma}
\begin{proof}
It is enough to show that given an \'etale cover $\wt{S}\rightarrow S$ in $\Affft$, the canonical map
\begin{equation}\label{eq: descent for Phi}
    |\Phi(\wt{S}^{\bullet}/S)|\rightarrow\Phi(S)
\end{equation}
is an equivalence.
If we are given a surjective local homeomorphism $\wt{Z}\rightarrow Z$ of topological spaces in $\on{Top}$, then we have an equivalence \footnote{This follows from \cite[Theorem 1.1]{DI04}.}
\begin{equation}\label{eq: descent for sing}
    |(\wt{Z}^{\bullet}/Z)^{\on{sing}}|\rightarrow Z^{\on{sing}}.
\end{equation}
We emphasize that the \v{C}ech nerve $\wt{Z}^{\bullet}/Z$ in the equation is taken in $\on{Top}$.
Since the functor $\Affft\rightarrow \on{Top}$ commutes with fiber products, that \eqref{eq: descent for Phi} is an equivalence follows from that \eqref{eq: descent for sing} is an equivalence.
\end{proof}

\sssec{}
Consider the following sequence of functors
\[\on{Sch}_{\ft}\xrightarrow{(\mhyphen)^{\topl}}\on{Top}\xrightarrow{(\mhyphen)^{\sing}}\Spc,\]
where the first arrow takes underlying topological space and the second arrow takes homotopy type. We denote the composition $(\mhyphen)^{\topl,\sing}:\on{Sch}_{\ft}\rightarrow \Spc$. Using the above codescent property we prove:

\begin{prop}\label{prop: homotopy type}
    The functor $\Phi: \PreStklft \rightarrow \Spc$ satisfies the following properties.
    \begin{enumerate}
        \item For a finite type scheme $X\in\on{Sch}_{\ft}$, we have a canonical equivalence
        \[\Phi(X)=X^{\topl,\sing}.\]
        \item For a torsion abelian group $\sfA$, we have $\Phi(B^i_{\et}\sfA)=B^i\sfA$.
    \end{enumerate}
\end{prop}

\subsection{Betti realization of an \'etale gerbe}\label{ssec: Betti realization}
Now that we have a notion of the homotopy type of prestacks, we have a bridge between \'etale and Betti gerbes.

\sssec{}
Let $X\in\PreStklft$ and $\sfA$ a torsion abelian group.
The functor $\Phi$ induces a map
\[\Phi: \Maps_{\PreStklft}(X, \Btet{\sfA})\rightarrow \Maps_{\Spc}(\Phi(X), B^2 \sfA),\]
in other words, a functor between the $\infty$-groupoids
\begin{equation}\label{eq: et gerbe to top gerbe}
    \Phi:\{\text{\'etale }\sfA\text{-gerbes on } X\}\rightarrow \{\text{topological }\sfA\text{-gerbes on } X^{\topl,\sing}\}.
\end{equation}

\sssec{}
Given a homotopy type $K\in\Spc$, we have the \textbf{Betti prestack} $K^{\on{pre}}_B:(\Aff^E)^{\op}\rightarrow\Spc$, which is the constant functor with value $K$. 
The \textbf{Betti stack} $K_B\in\Stk^E$ is defined as the \'etale sheafification of the Betti prestack. This construction fits into an adjunction
\begin{equation}\label{eq: Betti functor evaluation adj}
    (\mhyphen)_B: \Spc \leftrightarrows  \Stk^E: \on{ev}_{E},
\end{equation}
where $(\mhyphen)_B$ means taking the Betti stack, and $\on{ev}_{E}$ means the evaluation at $\on{Spec}E$.
\medskip

Consider $\Btet{\sfA}\in \Stk^E$.
Since $\Btet{\sfA}(E)=B^2 \sfA$, the adjunction \eqref{eq: Betti functor evaluation adj} tells us that the Betti functor induces equivalence
\[(\mhyphen)_B:\Maps_{\Spc}(K, B^2 \sfA)\rightarrow \Maps_{\PreStk^E}(K_B, \Btet{\sfA}),\]
in other words, an equivalence between the $\infty$-groupoids
\begin{equation}\label{eq: top gerbe to et gerbe}
    (\mhyphen)_B: \{\text{topological }\sfA\text{-gerbes on }K\}\rightarrow \{\text{etale }\sfA\text{-gerbes on }K_B\}.
\end{equation}

\sssec{}\label{sssec: Betti realization functor}
Composing the functors \eqref{eq: et gerbe to top gerbe} and \eqref{eq: top gerbe to et gerbe} we get the \textbf{Betti realization functor}
\begin{equation}\label{eq: Betti realization}
    \Phi_B: \{\text{etale }\sfA\text{-gerbes on } X\}\rightarrow \{\text{etale }\sfA\text{-gerbes on } (\Phi(X))_B\}
\end{equation}
for $X\in\PreStklft$.
From Proposition \ref{prop: homotopy type}, in case $X$ is a finite-type scheme, then the functor goes as
\begin{equation}\label{eq: Betti realization der sch}
    \Phi_B: \{\text{etale }\sfA\text{-gerbes on } X\}\rightarrow \{\text{etale }\sfA\text{-gerbes on } (X^{\topl, \on{sing}})_B\}.
\end{equation}
We point out that the \'etale $\sfA$-gerbes on $X$ we started with are defined in complex algebraic geometry, while the \'etale $\sfA$-gerbes on $(X^{\topl,\sing})_B$ are defined in algebraic geometry over $E$.

\subsection{Moduli of twisted local systems}\label{ssec: twisted locsys}
In this subsection we construct the moduli stack of twisted local systems.

\sssec{The untwisted case}\label{sssec: untwisted} Let $\mathsf{H}$ be a linear algebraic group over $E$. Given $K\in\Spc$, the moduli space of $\mathsf{H}$-local systems over $K$ is defined to be the mapping stack $\ul{\Maps}(K_B, \Bet \mathsf{H})\in\Stk^E$, and is denoted $\Loc_{\mathsf{H}}(K)$. From its definition, we have
\[\Loc_{\mathsf{H}}(K)(S)=\Maps(S\times K_B, \Bet \mathsf{H}),\quad S\in\Aff^E.\]
If we have a scheme $X\in\on{Sch}_{\ft}$, we define the moduli space $\Loc_{\mathsf{H}}(X)$ of $\mathsf{H}$-local systems over $X$ to be the same as $\Loc_{\mathsf{H}}(X^{\topl,\sing})$.

\sssec{Setup}\label{sssec: BA acts on BH}
To construct the moduli of twisted local systems, we work in the following setup.
Let $\sfA$ be an abelian group, thought of as a group $E$-scheme. Let $\mathsf{H}$ be an algebraic group over $E$, equipped with a central group homomorphism $c:\sfA\rightarrow \mathsf{H}$. Suppose we are given a pair $(Y,\cG)$ where $K\in\Spc$ and $\cG\in \on{Ge}_{\sfA}(K)$.
From the fact that $\sfA\rightarrow \mathsf{H}$ is central, we have a $\Bet{\sfA}$-action on $\Bet{\mathsf{H}}$ in the category $\Stk^E$.

\sssec{Associated bundle construction}\label{sssec: Associated bundle construction}
Suppose that we have a topos $\mc{C}$ and a group object $\mathsf{G}\in\mc{C}$, and an object $V\in\mc{C}$ equipped with an action of $\mathsf{G}$. Let $X\in\mc{C}$ and $\mc{P}$ be a principal $\mathsf{G}$-bundle object over $X$. We define an object $\mc{P}\times^G V\in \mc{C}$ as the geometric realization of the cosimplicial object
\[[n]\mapsto \mc{P}\times G^{\times n}\times V,\]
where the transition maps are given by the action of $G$ on $\mc{P}$ and $V$. In other words, we take the quotient of $\mc{P}\times V$ under the diagonal action of $G$. 
From its construction, the object $\mc{P}\times^G V$ is equipped with a canonical map to $X$.

\sssec{Construction of the moduli space}\label{sssec: construction of twlocsys}
 From the pair $(K,\cG)\in \Spc+\on{Ge}_{\sfA}$, we get an \'etale $\sfA$-gerbe $\mc{G}_B$ over the Betti stack $K_B$.
Now we use the associated bundle construction (Sect. \ref{sssec: Associated bundle construction}) to
\[\mc{C}=\Stk^E,\ G=B_{\et}\sfA,\ \mc{P}=\on{Tot}(\mc{G}_B),\text{ and } V=\Bet \mathsf{H}\]
and get
\begin{equation}\label{eq: BH bundle}
\on{Tot}(\mc{G}_B)\times^{B_{\et}A} \Bet \mathsf{H}\rightarrow K_B.
\end{equation}

We define the moduli space $\Loc^{\mc{G}}_{\mathsf{H}}(K)$ as the relative mapping stack
\[\ul{\Maps}_{/K_B}(K_B, \on{Tot}(\mc{G}_B)\times^{B_{\et}A} (\Bet \mathsf{H})) \in \PreStk^E.\]
This classifies sections of the $\Bet \mathsf{H}$-fibration \eqref{eq: BH bundle}.
Given a test affine scheme $S$, we have
\[\Loc^{\mc{G}}_{\mathsf{H}}(K)(S)=\Maps(K_B\times S,\mc{G}_B\times^{B_{\et}A} (\Bet \mathsf{H}))\times_{\Maps(K_B\times S,K_B)}\pt,\]
where the map
\[\Maps(K_B\times S,\on{Tot}(\mc{G}_B)\times^{B_{\et}A} (\Bet \mathsf{H}))\rightarrow \Maps(K_B\times S,K_B)\]
is induced from the projection \eqref{eq: BH bundle}, and the point
\[\pt\rightarrow \Maps(K_B\times S,K_B)\]  corresponds to the projection to the first factor.

\sssec{}
In case $\mc{G}$ is trivial, $\mc{G}_B$ is trivial, hence the $\Bet \mathsf{H}$-fibration \eqref{eq: BH bundle} gets identified with the projection
\[\Bet \sffH\times K_B\rightarrow K_B,\]
hence, the moduli of its sections recover the mapping space definition of $\Loc_{\mathsf{H}}(K)$.

\sssec{}
If we start with a finite type scheme $X$ over $\C$ and an \'etale $\sfA$-gerbe $\cG$ on it, then we get the topological $\sfA$-gerbe $\Phi(\cG)$ on $X^{\topl,\sing}\in\Spc$ as in \eqref{eq: et gerbe to top gerbe}. Then we can apply the construction in Sect. \ref{sssec: construction of twlocsys} to the pair $(X^{\topl,\sing},\Phi(\cG))$. We denote the resulting moduli space by 
\[\Loc^{\cG}_{\mathsf{H}}(X):=\Loc^{\Phi(\cG)}_{\mathsf{H}}(X^{\topl,\sing}).\]

\sssec{}\label{sssec: moduli space functoriality}
    There is an obvious functoriality of the moduli of twisted local systems. The construction of the moduli space gets upgraded to a functor
    \[(\Spc+\on{Ge}_{\sfA})^{\op}\rightarrow \PreStk^E,\quad (K,\cG)\mapsto \Loc^{\cG}_{\mathsf{H}}(K).\]

\subsection{Metaplectic parameter}
For $G$ and $X$ as before, \cite{Zhao1} defined the space of \textbf{metaplectic parameters} as the mapping space
\[\Maps_{\on{ptd}}(\Bet G\times X, \Bet^4 A(1)).\]
The subscript $\on{ptd}$ means that the map $\Bet G\times X\rightarrow \Bet^4 A(1)$ is equipped with a trivialization after pullback along the neutral map $X\rightarrow \Bet G\times X$. In other words, this is the space of morphisms in the category of pointed $X$-stacks $\Bet G\times X$ to $\Bet^4 A(1)\times X$. We work with a fixed metaplectic parameter $q$ throughout this paper.

\sssec{Metaplectic Hecke action} By means of transgression, a metaplectic parameter canonically induces \'etale $A$-gerbes on various spaces in the automorphic side of geometric Langlands theory. Moreover, those gerbes enjoy canonical compatibility data that enables us to construct the metaplectic Hecke action. We give a brief overview of this, while details can be found in \cite[Sect. 1, Sect. 7.1]{ZhaoNew}.
\medskip

We have the correspondence
\[\begin{tikzcd}
	{X\times\Bun_G} & {\Bet G\times X} \\
	{\Bun_G}
	\arrow["{\on{ev}}", from=1-1, to=1-2]
	\arrow["{\on{pr}}"', from=1-1, to=2-1]
\end{tikzcd}\]
where the horizontal arrow is the evaluation map and the vertical arrow is the projection. Forgetting the pointing of $q$ we have a degree 4 \'etale cochain (valued in $A(1)$) on $\Bet G\times X$. Pulling back along $\on{ev}$ and transgressing along $\on{pr}$, we have a degree 2 \'etale cochain (valued in $A$) on $\Bun_G$. We denote the resulting $A$-gerbe on $\Bun_G$ by $\cG^{\on{glob}}$. We denote the corresponding category of twisted sheaves on $\Bun_G$ by $\Shv^q(\Bun_G)$.
\medskip

The metaplectic parameter $q$ also induces a canonical $A$-gerbe $\cG^I$ on each $I$-legged local Hecke stack $\Hecke^{\loc}_{X^I}$, and we denote the corresponding category of twisted sheaves by $\Shv^q(\Hecke^{\loc}_{X^I})$. As in the non-metaplectic case, the convolution product endows this category with a canonical monoidal structure. A new ingredient in the metaplectic situation is that we need compatibility data between gerbes to make sense of the convolution product. This is because any sheaf operation between twisted sheaves requires compatibility data of gerbes, see Sect. \ref{sssec: functoriality twshv brief}. One can construct the required canonical compatibility data, and hence the monoidal structure on $\Shv^q(\Hecke^{\loc}_{X^I})$. For example, if we consider the convolution diagram
\[\Hecke^{\loc}_{X^I}\times\Hecke^{\loc}_{X^I}\leftarrow \wt{\Hecke}^{\loc,[2]}_{X^I}\rightarrow \Hecke^{\loc}_{X^I},\]
there is a canonical isomorphism between the pullback of $\cG^I\boxplus\cG^I$ from the left and the pullback of $\cG^I$ from the right. 
\medskip

In a similar sense, the Hecke correspondence is equipped with a canonical compatibility data between $\cG^I$ and $\cG^{\on{glob}}$, see Sect. \ref{sssec: preservation metaplectic}. This yields the \textbf{metaplectic Hecke action} of $\Shv^q(\Hecke^{\loc}_{X^I})$ on $\Shv^q(\Bun_G)$. This action is the origin of the Betti metaplectic spectral action. 

\sssec{Metaplectic dual data}\label{sssec: dual data}
From $q$ we extract the data $(H,\cG_Z)$. The \emph{metaplectic dual group} $H$ is a pinned split reductive group over $E$, whose root datum is defined in a combinatorial fashion. From a metaplectic parameter $q$ we extract a discrete invariant: an $A(-1)$-valued $W$-invariant quadratic form on the cocharacter lattice $\Lambda$ of $G$. The root datum of $H$ is determined combinatorially using this quadratic form. The \emph{metaplectic dual gerbe} $\cG_Z$ is an \'etale $Z_H(E)$-gerbe on the curve $X$. We define $\cG_Z$ to be the $Z_H(E)$-gerbe $\prescript{0}{}{\nu}+\theta$ in \cite[Sect. 2.2.6]{ZhaoNew}. Let $\wh{Z}_H$ be the group of characters of $Z_H$. In \emph{op. cit.}, both $\prescript{0}{}{\nu}$ and $\theta$ are defined as $\Z$-linear maps $\wh{Z}_H\rightarrow \on{Ge}_A(X)$, but this gives rise to a $Z_H(E)$-gerbe, because $A\hookrightarrow E^{\times,\on{tors}}$ and we have \footnote{For the equivalence below, we used the assumption that $E$ is algebraically closed.}
\[\Maps_{\Z}(\wh{Z}_H,\on{Ge}_{E^{\times,\on{tors}}}(X))= \on{Ge}_{Z_H(E)}(X).\]

As explained in the introduction, the pointwise metaplectic Satake category $\Sat^q_x$ is non-canonically equivalent to $\Rep(H)^{\heartsuit}$ up to a sign modification, and the dual gerbe $\cG_Z$ measures the dependence of $\Sat^q_x$ on $x\in X$. Now we can hint why we get a different group $H$ in the metaplectic setting, not $\chG$. In the non-metaplectic case, the most immediate evidence of the Satake equivalence is the classification of simple objects in the Satake category: there is a bijection between simple objects and Schubert cells of the affine Grassmannian. This follows from the fact that each Schubert cell supports exactly one equivariant local system. This fails in the metaplectic case: it is \emph{not} true that every Schubert cell supports an equivariant twisted local system. As a consequence, as we incorporate the twisting, there are fewer simple objects in the metaplectic Satake category. The set of Schubert cells that support an equivariant twisted local system is characterized by the quadratic form on $\Lambda$ associated to $q$. This is how the quadratic form enters in the root data of $H$.
\medskip

Finally, we comment more on the terms $\prescript{0}{}{\nu}$ and $\theta$. The construction of $\prescript{0}{}{\nu}$ is purely group-theoretic so that it makes sense for any group scheme $G\rightarrow S$ and it is functorial under the base change. However, the $\theta$ term is where the geometry of the curve $X$ plays the role. This is an inflation of the root gerbe $\omega^{1/2}_X$ of the canonical line bundle. This gets trivialized upon choosing a square root of $\omega_X$. The term $\prescript{0}{}{\nu}$ is in fact obtained as the $\Z$-linear component of the $\E_{\infty}$-map
$\nu:\wh{Z}_H\rightarrow \on{Ge}_{A}(X)$. This map breaks down into two parts \footnote{See \cite[Sect. 4.8.1]{GL} for the precise meaning of this.} : $\nu=\prescript{0}{}{\nu}+\epsilon$ where $\epsilon$ is a group homomorphism $\wh{Z}_H\rightarrow A_{2\te{-tors}}$. Using $A\hookrightarrow E^{\times}$ we get a group homomorphism $\wh{Z}_H\rightarrow \{\pm 1\}$. Hence, $\prescript{0}{}{\nu}$ (resp. $\epsilon$) is the $\Z$-linear (resp. sign) component of $\nu$.

\begin{rmk}\label{rmk: role of eps}
The factor $\epsilon$ plays a role in metaplectic Satake.
To have an equivalence between symmetric monoidal categories, we should alter the commutativity constraint by appropriate sign. On the automorphic side, this is the usual modification according to the parity of the connected component of the affine Grassmannian. On the spectral side, the data $\epsilon$ gives rise to a $\Z/2$-grading on the category of twisted $H$-representations, and we modify the commutativity constraint between odd objects by a sign. We review this in detail in Sect. \ref{ssec: recollection perv ULA}.
\end{rmk}

\section{Twisted factorization homology}\label{sec: twisted fact hlgy}
In this section, we define the notion of \emph{twisted factorization homology} (or \emph{twisted integral} in general), and identify a certain twisted factorization homology with the category of quasi-coherent sheaves on the moduli of twisted local systems. (Theorem \ref{thm: QC over twlocsys as a tw fact hlgy})

\subsection{Construction of twisted factorization homology} We introduce the notion of twisted integral and twisted factorization homology.
\sssec{Integral from \cite{GKRV}}\label{sssec: GKRV fact hlgy}
Suppose that we have a category $\bC$ with colimits. Given a space $Y\in\Spc$ and an object $\bc\in\bC$, we define an \emph{integral}
\[\int^{\bC}_Y \bc := \colim_Y \bc.\]
Here the index category is the $\infty$-groupoid $Y$, and we considered the constant functor $Y\rightarrow \bC$ with value $\bc$.
\medskip

Alternatively, this construction can be understood via the procedure of left Kan extension. We start with the functor $\{*\}\rightarrow\bC$ corresponding to the object $\bc$
and consider its left Kan extension along the functor 
$\{*\}\hookrightarrow \Spc$ corresponding to the object $\pt$.
Then the resulting functor is the \emph{integral}
\[\int^{\bC}_{(\mhyphen)}:\Spc\rightarrow\bC.\]

This operation is also called \emph{tensoring over spaces}, with the commonly used notation $\bc\otimes Y$.
In case $\bC=\DGCatSymMon$ and $\bc\in\bC$, the integral $\int_Y\bc$ deserves the name \emph{factorization homology}, because it coincides with the factorization homology of the $\E_{\infty}$-category $\bc$ over any framed manifold whose homotopy type is $Y$.

\sssec{Input}\label{sssec: fact hlgy setup}
Now we want to define the notion of \textit{twisted integral}. As before, we start with a category $\bC$ with colimits. We fix $B\in\Spc$ and take an object $\sfY\in\Spc_{/B}$ and $\bc\in \Fun(B,\bC)$ as input data. From this input we will define the twisted integral $\int^{\bC}_{\sfY}\bc$.

\sssec{Construction}
Consider the setup in Sect. \ref{sssec: fact hlgy setup}.
The category $\Spc_{/B}$ admits a fully faithful embedding
\begin{equation}\label{eq: pt into Spc}
    B= \pt_{/B}\hookrightarrow \Spc_{/B}.
\end{equation}
On the other hand, we have a functor $\bc:B\rightarrow\bC$.
Now we define the functor
\begin{equation}\label{eq: twisted factorization homology}
    \int^{\bC}_{(\mhyphen)}\bc: \Spc+\on{Ge}_{\sfA}\rightarrow \bC
\end{equation}
as the left Kan extension of the functor $\bc: S\rightarrow\bC$ along the inclusion \eqref{eq: pt into Spc}. We call this \textbf{twisted integral}. As in the untwisted case, we call $\int_{(Y,\cG)} \bc$ a \textbf{twisted factorization homology} if $\bC=\DGCatSymMon$.

\sssec{Some notations for the case $B=B^2\sfA$}\label{sssec: B being B2A}
    The base $B$ relevant to our purpose is $B^2\sfA$ where $\sfA$ is an abelian group. Then an object $\sfY\in \Spc_{/B}$ is the pair $(Y,\cG)$ where $Y\in\Spc$ and $\cG\in\on{Ge}_A(Y)$. For this reason, we may call the category $\Spc_{/B^2\sfA}$ as $\Spc+\on{Ge}_{\sfA}$. The input $\bc:B^2\sfA\rightarrow \bC$ encodes an object of $\bC$ together with a $B\sfA$-action on it; hence we may use the notation \[B\sfA\lMod(\bC):=\Fun(B^2\sfA,\bC).\] Abusing notation, we also denote the underlying object of $\bc\in B\sfA\lMod(\bC)$ as $\bc$, when the $B\sfA$-action is understood. In this situation, the twisted integral will be written as $\int^{\bC}_{(Y,\cG)}\bc$.

\subsection{Left Kan extension and colimits}
We record a result (which follows from the proof in \cite{MORezk}) concerning which colimits does a Kan extension preserve.

\sssec{Setup}
Consider the situation where we have categories $\bA$, $\bB$, and $\bC$ with the functors $f:\bC\rightarrow \bA$ and $g:\bC\rightarrow \bB$. Suppose that $\bA$ has all colimits. Then we have the left Kan extension $\on{LKE}_g f$ of $f$ along $g$, and the following theorem gives us a sufficient condition for a colimit to be preserved by $\on{LKE}_g f$.
\begin{gather*}  
\xy
(-10,0)*+{\bC}="Z";
(10,0)*+{\bB}="W";
(-10,-15)*+{\bA}="T";
{\ar@{->}^{g}"Z";"W"};
{\ar@{->}_{f}"Z";"T"};
{\ar@{->}^{\on{LKE}_g f}"W";"T"};
\endxy
\end{gather*}

\begin{thm}\label{thm: LKE and colimits}
    Suppose that a diagram $F:I\rightarrow B$ has the property that
    \begin{equation}\label{eq: assumption for LKE lemma}
        \colim_{i\in I} \Maps_{\bB}(g(\bc),F(i))=\Maps_{\bB}(g(\bc),\colim_{i\in I}F(i))
    \end{equation}
    for every $\bc\in \bC$. Then we have
    \[\colim_{i\in I}(\on{LKE}_g f)(F(i))=(\on{LKE}_g f)(\colim_{i\in I}F(i)).
    \]
\end{thm}

The rest of this section is devoted to the proof of Theorem \ref{thm: LKE and colimits}.

\sssec{} Define the functor $\wh{f}:\PSh(\bC)\rightarrow \bA$ as the left Kan extension of $f:\bC\rightarrow \bA$ along the Yoneda embedding $\rho: \bC\rightarrow\PSh(\bC)$. Also, define the functor $R_g$
\[R_g:\bB\rightarrow\PSh(\bC),\quad \bbb\mapsto \Maps_{\bB}(g(\mhyphen),\bbb)\]
and the functor $R_f: \bA\rightarrow\PSh(\bC)$ in the same way.

\begin{prop}\label{prop: whf and Rf}
    The functors $\wh{f}$ and $R_f$ defined above fits into the adjunction pair:
    \[\wh{f}:\PSh(\bC)\rightleftarrows \bA:R_f.\]
    In particular, the functor $\wh{f}$ preserves colimits.
\end{prop}
\begin{proof}
    We have to show the equivalence
    \begin{equation}\label{eq: fhat Rf adj}
        \Maps_{\bA}(\wh{f}(X), \ba)=\Maps_{\PSh(\bC)}(X,R_f(\ba))
    \end{equation}
    given $X\in\PSh(\bC)$ and $\ba\in\bA$.
    Since $\bA$ is assumed to have all colimits, the left Kan extension $\wh{f}$ gets identified with the pointwise left Kan extension, hence
    \[\wh{f}(X)=\colim_{\bc\rightarrow X\in \bC_{/X}} f(c)\]
    in $\bA$.
    Also note that 
    \[X=\colim_{\bc\rightarrow X\in \bC_{/X}}\bc\]
    in $\PSh(\bC)$.
    Therefore, it suffices to check \eqref{eq: fhat Rf adj} for representable objects in $\PSh(\bC)$. It follows from the construction of  $R_f$ and the fact that $\wh{f}|_{\bC}=f$.    
\end{proof}

\begin{lemma}\label{lem: Rg as a LKE}
    The functor $R_g: \bB\rightarrow \PSh(\bC)$ is equivalent to the left Kan extension of the Yoneda embedding $\rho:\bC\rightarrow\PSh(\bC)$ along $g:\bC\rightarrow\bB$.
\end{lemma}
\begin{proof}
    Denote the left Kan extension in the lemma as $\on{LKE}_g \rho: \bB\rightarrow\PSh(\bC)$. Since $\PSh(\bC)$ has all colimits, the left Kan extension $\on{LKE}_g \rho$ is computed pointwise, hence
    \[(\on{LKE}_g \rho)(\bbb)=\colim_{\bc, g(\bc)\rightarrow \bbb} \rho(\bc)\]
    for each $\bbb\in\bB$.
    Since colimits in a presheaf category are computed pointwise, we have
    \[(\on{LKE}_g \rho)(\bbb)(\bc')=\colim_{\bc, g(\bc)\rightarrow \bbb}\Maps_{\bC}(\bc',\bc)\]
    for $\bbb\in\bB$ and $\bc\in\bC$. The right-hand side can be rewritten as
    \begin{equation}\label{eq: colim E0 E1}
        \colim_{\bc, g(\bc)\rightarrow \bbb, \bc'\rightarrow \bc}\ast=\colim_{g(\bc')\rightarrow \bbb}\ast = \Maps_{\bB}(g(\bc'),\bbb)=R_g(\bbb)(\bc'),
    \end{equation}
    hence the proof.
    Here, the first equivalence comes from a cofinality argument, which we supply a detail below.
    \medskip

    Given objects $\bbb\in\bB$ and $\bc'\in\bC$, we consider two categories $\bE_0$ and $\bE_1$ where
    \[\bE_1=\Maps_{\bB}(g(\bc),\bbb),\]
    and $\bE_0$ is the category of pairs $(\bc,v,\alpha)$ where
    \[\bc\in\bC,\quad v\in\Maps_{\bB}(g(\bc),\bbb),\quad \alpha\in\Maps_{\bC}(\bc',\bc). \]   
    These are index categories for the colimit appearing in the two sides of the first equivalence of \eqref{eq: colim E0 E1}, hence we are left with constructing a functor $\bE_0\rightarrow \bE_1$ which is cofinal.
    One readily checks the adjunction
    \[\Psi_1:\bE_1\leftrightarrows \bE_0:\Psi_0,\]
    where $\Psi_1$ is defined as 
    \[\Psi_1:\bE_1\rightarrow \bE_0,\quad (g(\bc')\xrightarrow[]{u}\bbb)\mapsto (\bc=\bc', g(\bc')\xrightarrow[]{u} \bbb, \bc'\xrightarrow[]{\id}\bc'), \]
    and $\Psi_0$ is defined as
    \[\Psi_0:\bE_0\rightarrow\bE_1,\quad (\bc,g(\bc)\xrightarrow{v}b,\bc'\xrightarrow{\alpha}\bc)\mapsto (g(\bc')\xrightarrow{vg(\alpha)}\bbb).\]
    A functor admitting a left adjoint is cofinal (\cite[Ch.1 Sect. 2.2.3]{GR}), hence $\Psi_0$ is cofinal.
\end{proof}

\begin{prop}\label{prop: decompose LKE}
    We have an equivalence of functors
    \[\on{LKE}_g f=
    \wh{f}\circ \on{LKE}_g \rho=\wh{f}\circ R_g: \bB\rightarrow \bA.\]
\end{prop}
\begin{proof}
    The second equivalence comes from Lemma \ref{lem: Rg as a LKE}. For the first equivalence, we note that both of the left Kan extensions $\on{LKE}_g f$ and $\on{LKE}_g \rho$ are computed pointwise. Then the equivalence follows from the fact that $\wh{f}$ preserves colimits. (Proposition \ref{prop: whf and Rf}) 
\end{proof}

\sssec{Proof of Theorem \ref{thm: LKE and colimits}}
From the fact that $\on{LKE}_g f$ factors as $\wh{f}\circ R_g$ (Proposition \ref{prop: decompose LKE}) and that $\wh{f}$ preserves colimits (Proposition \ref{prop: whf and Rf}), we know that
\[\colim_{i\in I}(\on{LKE}_g f)(F(i))=(\on{LKE}_g f)(\colim_{i\in I}F(i))\]
holds if
\begin{equation}\label{eq: Rgf colimit}
    \colim_{i\in I}R_g(F(i))=R_g(\colim_{i\in I}F(i)).
\end{equation}

Since colimits in $\PSh(\bC)$ is computed pointwise, the assumption that
\[\colim_{i\in I} \Maps_{\bB}(g(\bc),F(i))=\Maps_{\bB}(g(\bc),\colim_{i\in I}F(i))\]
    for every $\bc\in \bC$ implies \eqref{eq: Rgf colimit}.
\qed

\begin{rmk}
    A dual version of Theorem \ref{thm: LKE and colimits} holds. In the situation where we have categories $\bA$, $\bB$, and $\bC$ with the functors $f:\bC\rightarrow\bA$ and $g:\bC\rightarrow\bB$ and $\bA$ has all limits, we can ask when a limit is preserved by $\on{RKE}_g f$.
    \medskip

    Suppose that a diagram $F:I\rightarrow B$ has the property that
    \begin{equation}\label{eq: assumption for RKE lemma}
        \colim_{i\in I} \Maps_{\bB}
        (F(i),g(\bc))=\Maps_{\bB}(\lim_{i\in I^{\op}}F(i),g(\bc))
    \end{equation}
    for every $\bc\in \bC$. Then we have
    \[\lim_{i\in I}(\on{RKE}_g f)(F(i))=(\on{RKE}_g f)(\lim_{i\in I}F(i)).
    \]
    Its proof is the same, except that instead of considering the $\rho: \bC\hookrightarrow \PSh(\bC)$ we have to consider $\bC\hookrightarrow \PSh(\bC^{\op})^{\op}$.
\end{rmk}
\begin{rmk}
    Theorem \ref{thm: LKE and colimits} and its dual version explains why every left Kan extension along the Yoneda embedding $\bC\hookrightarrow\PSh(\bC)$ preserves colimits, and every right Kan extension along the Yoneda embedding $\bC^{\op}\hookrightarrow \PSh(\bC)^{\op}$ preserves limits.
\end{rmk}
\subsection{Properties of the twisted integral}
In the situation of Sect. \ref{sssec: fact hlgy setup} we constructed the twisted integral as the functor $\Spc_{/B}\rightarrow \bC$.
\begin{prop}\label{prop: twfact colim}
    The functor 
    \begin{equation}\label{eq: twfact functor}
\Spc_{/B}\rightarrow \bC,\quad \sfY\mapsto \int^{\bC}_{\sfY}\bc
    \end{equation}
    preserves colimits.
\end{prop}
\begin{proof}
    We want to apply Theorem \ref{thm: LKE and colimits} to the situation
    \[\bC=B,\quad \bB=\Spc_{/B},\quad \bA=\bC.\]
    Suppose we have a colimit $\colim_{i\in I}\sfY_i =\sfY$ in $\Spc_{/B}$. To apply the theorem, we have to show that for any $b\in B$, the canonical map
    \begin{equation}\label{eq: checking LKE thm condition}
        \colim_{i\in I}\Maps_{\Spc_{/B}}(b,\sfY_i)\rightarrow \Maps_{\Spc_{/B}}(b,\sfY)
    \end{equation}
     is an equivalence. 
     \medskip

     Denote the forgetful functor $\Spc_{/B}\rightarrow \Spc$ by $F_B$.
     Given $\mathsf{K}\in \Spc_{/B}$ we have a canonical equivalence between spaces $\Maps_{\Spc_{/ B}}(b, \mathsf{K})$ and $(\mathsf{K})_b$, where the latter is the notation for the fiber product of
     \[\pr\xrightarrow{b} B\xleftarrow{\mathsf{K}}F_B(\mathsf{K}).\]
     Using the fact that colimits are universal in the topos $\Spc$, the canonical functor
     \[\colim_{i\in I}(\sfY_i)_b\rightarrow (\sfY)_b\]
     is an equivalence, hence \eqref{eq: checking LKE thm condition} is an equivalence. 
\end{proof}

\sssec{Twisted $Y$-diagram}
By the generality on overcategories, the forgetful functor $F_B:\Spc_{/B}\rightarrow \Spc$ is comonadic, hence creates colimits. Given $\sfY\in\Spc_{/B}$ let $Y\in\Spc$ be the underlying space $F_B(\sfY)$. For $y\in \sfY$ we have $(\sfY)_y\in\Spc_{/B}$ as the composition $\pt\xrightarrow{y}Y\xrightarrow{\sfY} B$.
We have a canonical map 
\begin{equation}\label{eq: colimit of fibers}
    \colim_{y\in Y}(\sfY)_y\rightarrow \sfY
\end{equation}
in $\Spc_{/B}$. Using the creation of colimits and that $\colim_{y\in Y}\pt=Y$, we conclude that \eqref{eq: colimit of fibers} is an equivalence.
\medskip

  $(\sfY)_y\in \Spc_{/B}$ actually lies in $\pt_{/B}$ and identified with $\sfY(y)\in B$. In the definition of twisted integral as a left-Kan extension, the functor $B=\pt_{/B}\hookrightarrow \Spc_{/B}$ was fully-faithful, hence \eqref{eq: twfact functor} restricted to $B$ equals $\bc:B\rightarrow \bC$. Therefore,
 The diagram
 \[Y\rightarrow \bC,\quad y\mapsto \int^{\bC}_{(\sfY)_y}\bc\]
 is identified with the diagram
\begin{equation}\label{eq: twisted Y diagram}
    Y\xrightarrow{\sfY}B\xrightarrow{\bc}\bC.
\end{equation}
 Therefore, applying Proposition \ref{prop: twfact colim} to the colimit \eqref{eq: colimit of fibers} we get:
\begin{prop}\label{prop: twisted Y diagram for general B}
    The colimit of the diagram \eqref{eq: twisted Y diagram} is canonically equivalent to the twisted integral $\int^{\bC}_{\sfY}\bc$.
\end{prop}

\sssec{The case $B=B^2\sfA$}
We specialize to the case $B=B^2\sfA$ and use notations in Sect. \ref{sssec: B being B2A}. Suppose that we have a colimit $\colim_{i\in I}Y_i=Y$ in $\Spc$. Using creation of colimits, we get:
\begin{lemma}\label{lem: colim in LKE+GeA}
    Given $(Y,\mc{G})\in\Spc+\on{Ge}_{\sfA}$ and a colimit $\colim_i Y_i =Y$ in $\Spc$, then we have a colimit
    \begin{equation}\label{eq: colim in Spc+GeA}
        \colim_i (Y_i, f^*_i \mc{G})=(Y,\mc{G})
    \end{equation}
in $\Spc+\on{Ge}_{\sfA}$. Here we used the canonical maps $f_i: Y_i\rightarrow Y$ to pullback the gerbe $\mc{G}$.
\end{lemma}

Applying Proposition \ref{prop: twfact colim} to the colimit \eqref{eq: colim in Spc+GeA} we get:
\begin{thm}\label{thm: fact hlgy and colimit}
    Suppose that we have a category $\bC$, an object $\bc\in B\sfA\lMod(\bC)$ and $(Y,\mc{G})\in\Spc+\on{Ge}_{\sfA}$ as in the setup of Sect. \ref{sssec: fact hlgy setup}. Moreover, suppose that we have a colimit $\colim_i Y_i =Y$ in $\Spc$ with canonical maps $f_i: Y_i\rightarrow Y$. Then the canonical map
    \[\colim_i \left(\int^{\bC}_{(Y_i,f^*_i \mc{G})} \bc\right) \rightarrow \int^{\bC}_{(Y,\mc{G})} \bc\]
    in $\bC$ is an equivalence.
\end{thm}

Proposition \ref{prop: twisted Y diagram for general B} in this case reads:
\begin{prop}\label{prop: twisted Y diagram}
    Suppose $(Y,\cG)\in\Spc+\on{Ge}_A$ and $\bc\in B\sfA\lMod(\bC)$ and consider the diagram constructed as the composition
    \begin{equation}\label{eq: twisted Y diagram for B}
        Y\xrightarrow{\mc{G}}B^2 \sfA\xrightarrow{\bc}\bC.
    \end{equation}
    The colimit of this diagram is canonically isomorphic to the twisted integral
    $\int^{\bC}_{(Y,\mc{G})}\bc$.
\end{prop}

\sssec{Comparison to the untwisted case} In case the gerbe $\mc{G}$ is trivial, the twisted $Y$-diagram \eqref{eq: twisted Y diagram} becomes the constant $Y$-diagram with value $\bc$. Therefore we get:

\begin{cor}\label{cor: comparison result}
    Let $Y\in\Spc$ and consider the trivial $\sfA$-gerbe  $(Y,\on{triv}_Y)\in\Spc+\on{Ge}_{\sfA}$. Then the \emph{twisted factorization homology} over $(Y,\on{triv}_Y)$ agrees with the \emph{factorization homology} in the sense of \cite{GKRV} over $Y$ (see Sect. \ref{sssec: GKRV fact hlgy}): we have a canonical equivalence
    \[\int^{\bC}_{(Y,\on{triv}_Y)}\bc=\int^{\bC}_Y \bc.\]
\end{cor}

\subsection{Moduli of twisted local systems and twisted factorization homology}
Let $\sffH$ be a linear algebraic group over $E$ and $\sfA$ be a torsion abelian group equipped with a central group homomorphism $\sfA\rightarrow Z_\sffH$. If we are given a homotopy type $K$ and an $\sfA$-gerbe $\mc{G}$ over it, then as in Subsection \ref{ssec: twisted locsys} we have a prestack $\Loc^{\mc{G}}_\sffH(K)$. In this section, we identify the symmetric monoidal category $\QC(\Loc^{\mc{G}}_\sffH(K))$ with a certain twisted factorization homology.

\sssec{}
The stack $\Bet \sffH\in\Stk^E$ has a canonical lift to $B\sfA\lMod(\Stk^E)$, and hence by taking the $\QC(\mhyphen)$ functor we get $\Rep(\sffH)\in B\sfA\lMod(\DGCatSymMon)$.
We construct a canonical functor
\begin{equation}\label{eq: fact hlgy to QC}
    \int^{\CAlg}_{(K,\mc{G})}\Rep(\sffH) \rightarrow \QC(\Loc^{\mc{G}}_\sffH(K))\quad \te{in }\DGCatSymMon.
\end{equation}

From Proposition \ref{prop: twisted Y diagram} the twisted factorization homology is canonically identified with the colimit of the diagram
\[K\xrightarrow{\cG}B^2\sfA\xrightarrow{\Rep(\sffH)}\DGCatSymMon.\]
To construct \eqref{eq: fact hlgy to QC} it suffices to construct a family of canonical functors
\begin{equation}\label{eq: insertion map at k}
    (\Rep(\sffH))_{\cG_k}\rightarrow \QC(\Loc^{\mc{G}}_\sffH(K))\quad \te{for }k\in K,
\end{equation}
where $\cG_k\in B^2\sfA$ is the pullback of $\cG$ along $\pt\xrightarrow{k}K$, and we are twisting $\Rep(\sffH)$ by $\cG_k$. Note that
\[(\Rep(\sffH))_{\cG_k}= \QC((\Bet \sffH)_{\cG_k}),\]
where on the right-hand side we twist the stack $\Bet \sffH$ by $\cG_k$. Under this identification, the desired functor \eqref{eq: insertion map at k} is the pullback functor along 
\[\Loc^{\mc{G}}_\sffH(K)\rightarrow(\Bet \sffH)_{\cG_k},\]
which is the consequence of functoriality discussed in Sect. \ref{sssec: moduli space functoriality} applied to the canonical map $(\pt,\cG_k)\rightarrow (K,\cG)$.
\medskip

The following is the main result of this subsection:

\begin{thm}\label{thm: QC over twlocsys as a tw fact hlgy}
    Suppose that the homotopy type $K$ has finitely many connected components. The functor \eqref{eq: fact hlgy to QC} is an equivalence in $\DGCatSymMon$. 
\end{thm}

Its special case when $\cG$ is trivial was established in \cite[Theorem 1.5.5]{GKRV}. (see Corollary \ref{cor: comparison result}) As a corollary, we identify $\QC(\Loc^{\cG_Z}_\sffH(X))$ with a certain twisted factorization homology.

\begin{cor}
    We have a canonical equivalence of symmetric monoidal DG categories
    \[\QC(\Loc^{\cG_Z}_\sffH(X))=\int^{\CAlg}_{(X^{\topl,\sing},\Phi(\cG_Z))}\Rep(\sffH).\]
\end{cor}

\sssec{}
We prove Theorem \ref{thm: QC over twlocsys as a tw fact hlgy} in the remainder of this subsection, using the results of \cite{GKRV} combined with the codescent arguments enabled by Theorem \ref{thm: fact hlgy and colimit}. Following \cite[Chapter 4,  Sect. 3.5.1]{GR}, we say that $\cX\in\PreStk^E$ is a quasi-passable prestack if it satisfies \textbf{(1)} the diagonal map $\cX\rightarrow \cX\times \cX$ is affine; \textbf{(2)} the category $\QC(\cX)$ is dualizable.
In the proof of Theorem \ref{thm: QC over twlocsys as a tw fact hlgy}, we actually prove:
\begin{thm}\label{thm: twlocsys quasi passable}
    Assume the situation of Theorem \ref{thm: QC over twlocsys as a tw fact hlgy}. Then $\Loc^{\cG}_\sffH(K)$ is quasi-passable and the functor \eqref{eq: fact hlgy to QC}
    \[\int^{\CAlg}_{(K,\cG)}\Rep(\sffH)\rightarrow \QC(\Loc^{\cG}_\sffH(K))\]
    is an equivalence in $\DGCatSymMon$.
\end{thm}

Using the argument of \cite[Sect. 1.6.6]{GKRV}, we reduce to the case where $K$ is connected.

\sssec{}
We have
\[\Maps_{\Spc+\on{Ge}_{\sfA}}((\pt,\triv),(Y,\cG))=\on{Tot}(\cG),\]
which is the analogue of the following fact: if we are given a principal $G$-bundle $\cP$ over a topological space $X$, then the data of a point $x\in X$ together with trivialization of the fiber $\cP_x$ is equivalent to a choice of a point in the total space of $\cP$. Hence, given $\tau\in \on{Tot}(\cG)$ we get a map
\begin{equation}\label{eq: tau evaluation map}
    \tau^*:\Loc^{\cG}_\sffH(K)\rightarrow \Bet \sffH
\end{equation}
from the functoriality in Sect. \ref{sssec: moduli space functoriality}.
\medskip

If $f: \cX_1\rightarrow \cX_2$ is an affine morphism between prestacks and that $\cX_2$ is quasi-passable, then $\cX_1$ is also quasi-passable. (\cite[Lemma 1.6.2]{GKRV}) Hence, we deduce Theorem \ref{thm: twlocsys quasi passable} once we prove:
\begin{thm}\label{thm: twlocsys affine}
    Let $K$ be a connected homotopy type and $(K,\cG)\in \Spc+\on{Ge}_{\sfA}$, and choose any $\tau\in \on{Tot}(\cG)$.
    Then the map \eqref{eq: tau evaluation map}   
    is affine and the functor \eqref{eq: fact hlgy to QC} is an equivalence.
\end{thm}

\begin{rmk}\label{rmk: twlocsys affine when trivial}
    Its special case where the gerbe $\cG$ is trivial was established in \cite[Theorem 1.6.8]{GKRV}.
\end{rmk}

\sssec{}
For the rest of the subsection, we prove Theorem \ref{thm: twlocsys affine}.  We choose a morphism $f:U\rightarrow K$ in $\Spc$ such that $U$ is also a connected space and the pullback trivializes the gerbe $\cG$. We can choose $f$ to be $\on{Tot}(\cG)\rightarrow K$ for example.
\medskip

The map $U\rightarrow K$ is a surjection on the level of $\pi_0$, hence
\begin{equation}\label{eq: eff epi in Spc}
    K=\colim_{i\in \Delta^{\op}}U^i/K
\end{equation}
in $\Spc$. Therefore by Theorem \ref{thm: fact hlgy and colimit} we have
\begin{equation}\label{eq: colimit of twfact}
    \int^{\CAlg}_{(K,\cG)}\Rep(\sffH)=\colim_{i\in\Delta^{\op}}\int^{\CAlg}_{(U^i/K,f^*_i\cG)}\Rep(\sffH)
\end{equation}
in $\DGCatSymMon$, where we let $f_i:U^i/K\rightarrow K$ the canonical map.

\sssec{}\label{sssec: lifting tau}
From the assumption on $U$ and $K$, both $\on{Tot}(\cG)$ and $\on{Tot}(f^*\cG)$ are connected. We choose a lift $\wt{\tau}\in \on{Tot}(f^*\cG)$ of $\tau\in\on{Tot}(\cG)$, together with the isomorphism $\alpha:\wt{\tau}\cong f^*\tau$. The choice of  $(\wt{\tau},\alpha)$ determines a collection of points
\[\tau_i\in \on{Tot}(f^*\cG)^i/\on{Tot}(\cG),\quad [i]\in \Delta^{\op}\]
equipped with compatibility under transition maps in the \v{C}ech nerve for $\on{Tot}(f^*\cG)\rightarrow \on{Tot}(\cG)$. From the fact that \v{C}ech nerve pulls back, the \v{C}ech nerve is identified with $\{\on{Tot}(f^*_i\cG)\}_{i\in\Delta^{\op}}$, hence $\tau_i\in \on{Tot}(f^*_i\cG)$. Considering the augmented \v{C}ech nerve
\[\on{Tot}(f^*\cG)^{\bullet}/\on{Tot}(\cG)\rightarrow \on{Tot}(\cG),\]
the collection $\{\tau_i\}_{i\in\Delta^{\op}}\cup\{\tau\}$ has a natural compatibility under transition maps.

\sssec{}
Since Theorem \ref{thm: twlocsys affine} holds true for trivial gerbes (Remark \ref{rmk: twlocsys affine when trivial}) and each $f^*_i\cG$ is trivializable, we can apply the result to $\tau_i$ and $(U^i/K,f^*_i\cG)$ and conclude that for each $i\in\Delta$
 the map
    \begin{equation}\label{eq: tau i eval}
        \tau_i:\Loc^{f^*_i\cG}_\sffH (U^i/K)\rightarrow \Bet \sffH
    \end{equation}
    is affine,
    and the canonical map
    \begin{equation}\label{eq: equivalence for pullbacks}
        \int^{\CAlg}_{(U^i/K,f^*_i\cG)}\Rep(\sffH)\rightarrow \QC(\Loc^{f^*_i\cG}_\sffH (U^i/K))
    \end{equation}
    is an equivalence in $\DGCatSymMon$.
\medskip

From the compatibility of $\{\tau^i\}_{i\in\Delta}$, we have a diagram
\begin{equation}\label{eq: diagram twlocsys}
    \Delta\rightarrow \PreStk^E,\quad [i]\mapsto \Loc^{f^*_i\cG}_\sffH (U^i/K),
\end{equation}
equipped with a natural transformation to
\begin{equation*}\label{eq: diagram const BH}
    \Delta\rightarrow \PreStk^E,\quad [i]\mapsto \Bet \sffH
\end{equation*}
given by the collection \eqref{eq: tau i eval}. Passing to limits, we get a morphism
\begin{equation}\label{eq: limit of twlocsys to BH}
    \lim_{i\in\Delta}\Loc^{f^*_i\cG}_\sffH (U^i/K)\rightarrow \Bet \sffH.
\end{equation}

\begin{lemma}\label{lem: limit of twlocsys}
    We have a canonical equivalence
    \begin{equation}\label{eq: limit of twlocsys}
        \lim_{i\in\Delta}\Loc^{f^*_i\cG}_\sffH (U^i/K)=\Loc^{\cG}_\sffH(K)
    \end{equation}
    in $\PreStk^E$. Moreover, the map \eqref{eq: limit of twlocsys to BH} gets identified with
    $\tau^*:\Loc^{\cG}_\sffH(K)\rightarrow \Bet \sffH$
    via \eqref{eq: limit of twlocsys}.
\end{lemma}
\begin{proof}
    The first part is standard, and the second part is a consequence of the compatibility explained in Sect. \ref{sssec: lifting tau}.
\end{proof}

\begin{lemma}[\cite{GKRV}, Lemma 1.6.11]\label{lem: cosifted limit and QC}
    Let $a\mapsto W_a$ be a cosifted diagram of prestacks, affine over some base prestack $Z$. Let $W=\lim_a W_a$ in $\PreStk^E$. Then the map $W\rightarrow Z$ is affine and the induced functor
    \[\colim_a \QC(W_a)\rightarrow \QC(W)\]
    is an equivalence in $\DGCatSymMon$.
\end{lemma}
\begin{rmk}
    In \emph{loc.cit.}, the diagram $a\mapsto W_a$ was assumed to be a sifted diagram, but it was a typo; it actually needs to be a cosifted diagram. It does not affect the validity of \cite[Theorem 1.6.8]{GKRV}.
\end{rmk}

\sssec{}
Using Lemma \ref{lem: cosifted limit and QC} for the cosifted diagram \eqref{eq: diagram twlocsys} and using Lemma \ref{lem: limit of twlocsys}, we conclude that
\[\tau^*:\Loc^{\cG}_\sffH(K)\rightarrow \Bet \sffH\]
is affine and that the canonical functor
\begin{equation}\label{eq: colim of QC twloc}
    \colim_{i\in\Delta^{\op}}\QC(\Loc^{f^*_i\cG}_\sffH (U^i/K))\rightarrow \QC(\Loc^{\cG}_\sffH(K))
\end{equation}
is an equivalence.
\medskip

Now we conclude that the functor \eqref{eq: fact hlgy to QC} is an equivalence for $(K,\cG)$ from combining the fact that the functors \eqref{eq: colimit of twfact}, \eqref{eq: equivalence for pullbacks}, and \eqref{eq: colim of QC twloc} are all equivalences. \qed[Theorem \ref{thm: twlocsys affine}]

\section{Action of twisted factorization homology}\label{sec: action of twisted fact hlgy}
In this section, we study the categorical action of the twisted factorization homology $\int^{\CAlg}_{(Y,\cG)}\cC$, for a symmetric monoidal DG category $\cC$ with an $\mathsf{A}$-symmetry. We prove that this is equivalent to the \emph{factorizable action} of the local system
\[Y\xrightarrow{\cG}B^2\sfA\xrightarrow{\cC}\DGCatSymMon,\]
see Theorem \ref{thm: main categorical thm}.
The notion of a factorizable action of a local system of symmetric monoidal DG categories - which is the main categorical theme of this paper - is also defined in this subsection.
\medskip

In fact, Theorem \ref{thm: main categorical thm} concerns the space of monoidal functors $\int^{\CAlg}_{(Y,\cG)}\cC\rightarrow \cA$ given $\cA\in\DGCatMon$.
 If we have a DG category $\cM$, then we apply the result to $\cA=\End(\cM)$, hence obtaining an alternative description of the categorical action of $\int^{\CAlg}_{(Y,\cG)}\cC$ on $\cM$ in terms of factorization. We will later apply this to the metaplectic Hecke action and construct the spectral action.

\sssec{} 
We give a preview of this section. Theorem \ref{thm: main categorical thm} asserts that given $\cC\in \DGCatMon$, the data of a monoidal functor $\int^{\CAlg}_{(Y,\cG)}\bc\rightarrow \cC$ is equivalent to the data of a natural transformation between certain two functors out of $\fSet$. The idea of the proof is to rewrite both sides using the combinatorics of the twisted arrow category.
\medskip

Given a category $\bD$, we define the category of twisted arrows $\on{TwArr}(\bD)$ as in \cite[Part III, Chapter 9, Sect. 2.4]{GR}. Its objects consist of 1-morphisms $\bd_s\rightarrow \bd_t$ in $\bD$, and a 1-morphism $(\bd^0_s\rightarrow\bd^0_t)\rightarrow (\bd^1_s\rightarrow \bd^1_t)$ is a data of a commutative diagram
\[\begin{tikzcd}
	{\bd^0_s} & {\bd^0_t} \\
	{\bd^1_s} & {\bd^1_t}
	\arrow[from=1-1, to=1-2]
	\arrow[from=1-1, to=2-1]
	\arrow[from=2-1, to=2-2]
	\arrow[from=2-2, to=1-2]
\end{tikzcd}\]

It was shown in \cite[Proposition 5.1]{GHN} that if we have two categories $\bD$ and $\bE$ and functors $F_1,F_2:\bD\rightarrow\bE$, then the space of natural transformations $F_1\rightarrow F_2$ can be expressed as a certain limit over the diagram $\on{TwArr}(\bD)^{\op}$. Applying this to the two functors out of $\fSet$, we express the space of natural transformation as a limit over the diagram $\on{TwArr}(\fSet)^{\op}$. 
\medskip

On the other hand, we can express the underlying monoidal DG category of $\int^{\CAlg}_{(Y,\cG)}\cC$ as a colimit of a certain functor $\on{TwArr}(\fSet)\rightarrow \DGCatMon$ (Theorem \ref{thm: twArr twfact comparison}), hence express the space of monoidal functors $\int^{\CAlg}_{(Y,\cG)}\cC\rightarrow \cA$ as a limit over $\on{TwArr}(\fSet)^{\op}$. After expressing both sides of the equivalence as limits over $\on{TwArr}(\fSet)^{\op}$, the desired equivalence comes from term-by-term equivalences.

\subsection{Twisted factorization homology and twisted arrows}

\sssec{}\label{sssec: cI integrated over YJ}
For a category $\bC$ and an abelian group $\sfA$, we let $B\sfA\lMod(\bC)$ be the functor category $\Fun_{\Cat}(B^2\sfA,\bC)$. This encodes an object $\cC\in\bC$ and a $B\sfA$-action on $\cC$. We have a sequence of forgetful functors
\begin{equation*}
    B\sfA\lMod(\DGCatSymMon)\rightarrow B\sfA\lMod(\DGCatMon)\rightarrow B\sfA\lMod(\DGcat).
\end{equation*}
Suppose we have $\cC\in B\sfA\lMod(\DGCatSymMon)$. Given $\bC$ among $\DGCatSymMon$, $\DGCatMon$, or $\DGcat$, we can consider $\cC$ as an object of $B\sfA\lMod(\bC)$ through the forgetful functor above.

\sssec{}
We construct a functor
\begin{equation}\label{eq: twArr ftr}    \TwArr(\fSet)\rightarrow \bC, \quad (I\rightarrow J)\mapsto \int^{\bC}_{(Y^J,\cG^{\boxtimes J})}\cC^{\otimes I}.
\end{equation}

We first explain what the expression
\begin{equation}\label{eq: twarr termwise}
    \int^{\bC}_{(Y^J,\cG^{\boxtimes J})}\cC^{\otimes I}
\end{equation}
means. We have
\[\cC^{\otimes I}\in BA^{\times I}\lMod(\bC),\]
and then pull back the action of $BA^{\times I}$ along the diagonal map
$\Delta:BA^{\times J}\rightarrow BA^{\times I}$,
which is symmetric monoidal,
to view it as $\cC^{\otimes I}\in BA^{\times J}\lMod(\bC)$.
We have an $A^{\times J}$-gerbe $\cG^{\boxtimes J}$ on $Y^{J}$ defined as
\[Y^J\xrightarrow{\cG^{\times J}} (B^2 \sfA)^{\times J}=B^2 \sfA^{\times J},\]
hence \eqref{eq: twarr termwise} is the twisted integral of $\cC^{\otimes I}$ integrated over $(Y^J, \cG^{\boxtimes J})$, computed in $\bC$.

\sssec{}
We explain what the functor \eqref{eq: twArr ftr} does to morphisms.
Say that we have a morphism $(I_1\xrightarrow{\alpha_1} J_1)\rightarrow (I_2\xrightarrow{\alpha_2} J_2)$ in $\TwArr(\on{fSet})$, i.e., we have a commutative diagram
\[\begin{tikzcd}
	{I_1} & {J_1} \\
	{I_2} & {J_2}
	\arrow["{\alpha_1}", from=1-1, to=1-2]
	\arrow["\beta"', from=1-1, to=2-1]
	\arrow["{\alpha_2}"', from=2-1, to=2-2]
	\arrow["\gamma"', from=2-2, to=1-2]
\end{tikzcd}\]
in $\fSet$. We want to construct a map
\[\int^{\bC}_{(Y^{J_1},\cG^{\boxtimes J_1})}\cC^{\otimes I_1}\rightarrow \int^{\bC}_{(Y^{J_2},\cG^{\boxtimes J_2})}\cC^{\otimes I_2}\]
in $\bC$. To do this, we express both sides as colimit. 
\medskip

Using Proposition \ref{prop: twisted Y diagram} we have
\[\int^{\bC}_{(Y^{J_1},\cG^{\boxtimes J_1})}\cC^{\otimes I_1}=\colim_{\ul{y}\in Y^{J_1}} (\cC^{\otimes I_1})_{(\cG^{\boxtimes J_1})_{\ul{y}}}=\colim_{\ul{y}\in Y^{J_1}} \bigotimes_{i_1\in I_1}\cC_{\cG_{y_{\alpha_1(i_1)}}}\]
and
\[ \int^{\bC}_{(Y^{J_2},\cG^{\boxtimes J_2})}\cC^{\otimes I_2}=
\colim_{\ul{y}'\in Y^{J_2}} (\cC^{\otimes I_2})_{(\cG^{\boxtimes J_2})_{\ul{y}'}}=\colim_{\ul{y}'\in Y^{J_2}} \bigotimes_{i_2\in I_2}\cC_{\cG_{y'_{\alpha_2(i_2)}}}.\]
We have a map $\Delta_{\gamma}:Y^{J_1}\rightarrow Y^{J_2}$ between index categories, and hence it suffices to construct a map
\[\bigotimes_{i_1\in I_1}\cC_{\cG_{y_{\alpha_1(i_1)}}}\rightarrow \bigotimes_{i_2\in I_2}\cC_{\cG_{y'_{\alpha_2(i_2)}}}\]
in $\bC$ for $\ul{y}'=\Delta_{\gamma}(\ul{y})$, or equivalently,
\begin{equation}\label{eq: term by term}
    \bigotimes_{i_1\in I_1}\cC_{\cG_{y_{\alpha_1(i_1)}}}\rightarrow \bigotimes_{i_2\in I_2}\cC_{\cG_{y_{\gamma\alpha_2(i_2)}}}.
\end{equation}

The fact that $\cC\in B\sfA\lMod(\DGCatSymMon)$ means that each $\cC_{\cG_y}$ for $y\in Y$ is a symmetric monoidal category. The desired functor \eqref{eq: term by term} comes from the tensor product of multiplication maps
\begin{equation}\label{eq: for each j1}
    (\cC_{\cG_{y_{j_1}}})^{\otimes \alpha^{-1}_1(j_1)}\rightarrow (\cC_{\cG_{y_{j_1}}})^{\otimes \alpha^{-1}_2\gamma^{-1}(j_1)}
\end{equation}
for each $j_1\in J_1$, induced from
$\alpha^{-1}_1(j_1)\xrightarrow{\beta} \alpha^{-1}_2\gamma^{-1}(j_1)$.
Since $\cC_{\cG_{y_{j_1}}}$ is symmetric monoidal, the multiplication maps \eqref{eq: for each j1} are also symmetric monoidal, and thus defines a map in $\bC$. Taking tensor products, we also get \eqref{eq: term by term} in $\bC$.

\sssec{}
Say that we have $\cC\in B\sfA\lMod(\DGCatSymMon)$. We have the twisted factorization homology
    \[\int^{\on{CAlg}}_{(Y,\cG)}\cC \in\DGCatSymMon\]
    computed in $\DGCatSymMon$. We consider its image under the forgetful functor $\DGCatSymMon\rightarrow \bC$.
        We construct a canonical functor
    \begin{equation}\label{eq: twArr twfact comparison}
        \colim_{(I\rightarrow J)\in\TwArr(\fSet)} \int^{\bC}_{(Y^J,\cG^{\boxtimes J})} \cC^{\otimes I} \rightarrow \int^{\on{CAlg}}_{(Y,\cG)}\cC
    \end{equation}
    in $\bC$, from a collection of functors
    \begin{equation}\label{eq: twArr twfact comparison termwise}
        \int^{\bC}_{(Y^J,\cG^{\boxtimes J})} \cC^{\otimes I} \rightarrow \int^{\on{CAlg}}_{(Y,\cG)}\cC
    \end{equation}
    for $I\xrightarrow{\alpha}J$.
    \medskip

    We construct the functor \eqref{eq: twArr twfact comparison termwise} in $\bC$. The left-hand side can be expressed as a colimit
    \[\colim_{\ul{y}\in Y^J}  \bigotimes_{i\in I}\cC_{\cG_{y_{\alpha(i)}}}. \]
    To construct \eqref{eq: twArr twfact comparison termwise}, we have to  construct a family of functors
    \[\bigotimes_{i\in I}\cC_{\cG_{y_{\alpha(i)}}}\rightarrow \int^{\CAlg}_{(Y,\cG)}\cC\]
    in $\bC$. This is the composition
    \[\bigotimes_{i\in I}\cC_{\cG_{y_{\alpha(i)}}}\rightarrow \left(\int^{\CAlg}_{(Y,\cG)}\cC\right)^{\otimes I} \xrightarrow{\on{mult}} \int^{\CAlg}_{(Y,\cG)}\cC,\]
    where the first arrow is the tensor product of insertion maps
    \[\on{ins}_{y_{\alpha(i)}}:\cC_{\cG_{y_{\alpha(i)}}}\rightarrow\int^{\CAlg}_{(Y,\cG)}\cC\]
    coming from the colimit expression of twisted factorization homology (Proposition \ref{prop: twisted Y diagram}) and the second arrow is the multiplication map in the symmetric monoidal structure of
    $\int^{\CAlg}_{(Y,\cG)}\cC$. Because the twisted factorization homology is computed in $\DGCatSymMon$, both arrows live in $\DGCatSymMon$, hence also in $\bC$.

\begin{thm}\label{thm: twArr twfact comparison}
    The canonical functor \eqref{eq: twArr twfact comparison}
    is an equivalence in $\bC$.
\end{thm}

\begin{proof}
    When $\cG$ is trivial, the equivalence was proved in \cite[Theorem 1.2.4]{GKRV}.
    The idea of the proof is to reduce it to the case where the gerbe is trivial using codescent. 
\medskip

    Denote the total space of gerbe $\cG$ by $P$, which is equipped with a projection map $P\rightarrow Y$. We form a \v{C}ech nerve $\{P^n/Y\}_{n\in\Delta^{\op}}$, and let $\cG^{(n)}$ be the pullback of the gerbe $\cG$ along the map $P^n/Y\rightarrow Y$. From Theorem \ref{thm: fact hlgy and colimit}, we have
    \begin{equation}\label{eq: descent for CAlg twfact}
        \int^{\CAlg}_{(Y,\cG)}\cC = \colim_{n\in \Delta^{\op}} \int^{\CAlg}_{(P^n/Y, \cG^{(n)})}\cC
    \end{equation}
    in $\DGCatSymMon$, but since the forgetful functor $\DGCatSymMon\rightarrow\bC$ preserves sifted colimits, \eqref{eq: descent for CAlg twfact} is still valid in $\bC$.
    \medskip
    
    Since each gerbe $\cG^{(n)}$ is trivializable, we get
    \[\int^{\bC}_{(Y,\cG)}\cC =\colim_{n\in \Delta^{\op}} \int^{\bC}_{(P^n/Y, \cG^{(n)})}\cC= \colim_{n\in \Delta^{\op}} \colim_{(I\rightarrow J)\in \TwArr(\fSet)}\int^{\bC}_{((P^n/Y)^J, (\cG^{(n)})^{\boxtimes J})}\cC^{\otimes I}\]
    in $\bC$.
Switching the order of colimits we get
\begin{equation*}
   \int^{\bC}_{(Y,\cG)}\cC=\colim_{(I\rightarrow J)\in\TwArr(\fSet)}\colim_{n\in \Delta^{\op}} \int^{\bC}_{(P^n/Y)^J, (\cG^{(n)})^{\boxtimes J}}\cC^{\otimes I}, 
\end{equation*}
and from Theorem \ref{thm: fact hlgy and colimit} we get
    \[\colim_{n\in \Delta^{\op}}\int^{\bC}_{((P^n/Y)^J, (\cG^{(n)})^{\boxtimes J})}\cC^{\otimes I} = \int^{\bC}_{(Y^J,\cG^{\boxtimes J})}\cC^{\otimes I},\]
    because the gerbe $(\cG^{(n)})^{\boxtimes J}$ is canonically identified with the pullback of the gerbe $\cG^{\boxtimes J}$ along $(P^n/Y)^J\rightarrow Y^J$.
   Combining these, we get
    \[\int^{\bC}_{(Y,\cG)}\cC = \colim_{(I\rightarrow J)\in \TwArr(\fSet)}\int^{\bC}_{(Y^J,\cG^{\boxtimes J})}\cC^{\otimes I}
    \]
    in $\bC$ as desired.
\end{proof}

\subsection{Local system of categories}\label{ssec: LSCat}
Categorical action of twisted factorization homology can be studied via local system of categories. We review the definition and collect some facts for later.

   \sssec{Definitions}
   Say that we have a homotopy type $K\in\Spc$. We define the \textbf{local system of categories} over $K$ as a functor $K\rightarrow \DGcat$. These form a category 
   \[\LSCat(K):=\Fun(K,\DGcat).\]
   Let
   \[\LSMon(K):=\Fun(K,\DGCatMon)\quad (\te{resp.}\ \on{LSCMon}(K):=\Fun(K,\DGCatSymMon))\]
   to be the category of \textbf{local system of monoidal categories} (resp. \textbf{symmetric monoidal categories}) over $K$.
  
  \sssec{}  
   If we have an operad $\cO$, a category $\bD$ and a symmetric monoidal category $\bE$, then the functor category $\Fun(\bD,\bE)$ attains a pointwise symmetric monoidal structure, and we have a canonical equivalence
   \[\cO\text{-}\Alg(\Fun(\bD,\bE))= \Fun(\bD,\cO\text{-}\Alg(\bE)).\]
Applying this, we get canonical equivalences
\begin{equation}\label{eq: LSMon as alg obj}
    \LSMon(K)=\Alg(\LSCat(K)),\quad \LSCMon(K)=\CAlg(\LSCat(K)).
\end{equation}

\sssec{}\label{sssec: descent of LSMon}
Suppose we are given a map $Z\rightarrow W\in \Spc$ inducing a surjection on $\pi_0$. From the fact that
\[W=\colim_{[n]\in\Delta^{\op}}Z^n/W\]
we get
\[\LSCat(W)=\Fun(W,\DGcat)=\lim_{[n]\in\Delta}\Fun(Z^n/W, \DGcat)=\lim_{[n]\in \Delta}\LSCat(Z^n/W).\] 
The same holds if we replace $\LSCat$ by $\LSMon$ or $\LSCMon$.

\sssec{}\label{sssec: LSCat global section}
We observe that
\begin{equation}\label{eq: LSCat as ShvCat}
    \LSCat(K)=\Fun(K,\DGcat)=\lim_{K^{\op}}\DGcat =\lim_{K^{\op}}\ShvCat(\pt)=\ShvCat(K_B),
\end{equation}
because $K_B =\colim_{K} \pt$ in $\Stk$ and $\ShvCat$ satisfies the \'etale descent. Hence we have an adjunction
\begin{equation}\label{eq: Loc Gamma adjunction}
    \Loc_{K}: \LS(K)\lMod \leftrightarrows \LSCat(K) :\Gamma(K,\mhyphen)
\end{equation}
as in \cite[Sect. 1.3.1]{Ga}, where we invoked $\QC(K_B)=\LS(K)$.
\medskip

The functor $\Loc_K$ is defined as follows: given $\bM\in \LS(K)\lMod$ we form the local system of categories
    \[K\rightarrow \DGcat,\quad k\mapsto \bM\otimes_{\LS(K)}\Vect_k,\]
    where the $\Vect_k$ is a $\LS(K)$-module category given by the symmetric monoidal functor
    $k^*: \LS(K)\rightarrow \Vect$ taking fiber at $k$.
The functor $\Gamma(K,\mhyphen)$ is defined as
     \[\LSCat(K) \rightarrow \LS(K)\lMod, \quad (\cL: K\rightarrow \DGcat) \mapsto \lim_{K}\cL\]
     where $\lim_K\cL$ is understood with the canonical $\LS(K)$-module structure.

\sssec{}
As mentioned in \emph{loc. cit.}, the functor $\Loc_K$ is symmetric monoidal.
Using \eqref{eq: LSMon as alg obj} and Sect. \ref{sssec: adjunction and algebra}, the adjunction \eqref{eq: Loc Gamma adjunction} lifts to adjunctions
\[\Loc_{K}: \Alg(\LS(K)\lMod) \leftrightarrows \LSMon(K) :\Gamma(K,\mhyphen)\]
and
\[\Loc_{K}: \CAlg(\LS(K)\lMod) \leftrightarrows \LSCMon(K) :\Gamma(K,\mhyphen).\]

\sssec{}
The identification \eqref{eq: LSCat as ShvCat} intertwines pullback. Given $f:Y_1\rightarrow Y_2$ in $\Spc$, we have the pullback functor
    \begin{equation}\label{eq: LSCat pullback}
        \LSCat(Y_2)\rightarrow \LSCat(Y_1),\quad (Y_2\xrightarrow{\cL}\DGcat)\mapsto (Y_1\xrightarrow{f}Y_2\xrightarrow{\cL}\DGcat).
    \end{equation}
    On the other hand, $f:Y_1\rightarrow Y_2$ induces a map $f_B: (Y_1)_B\rightarrow (Y_2)_B$, hence we have a pullback of sheaf of categories:
\begin{equation}\label{eq: ShvCat pullback}
    \ShvCat((Y_2)_B)\rightarrow \ShvCat((Y_1)_B).
    \end{equation}
    Then the pullback functor \eqref{eq: LSCat pullback} is identified with \eqref{eq: ShvCat pullback} under the equivalence \eqref{eq: LSCat as ShvCat}.

    \sssec{}\label{sssec: LSCat pull push}
    From the generality of the sheaf of categories, the pullback functor \eqref{eq: LSCat pullback} admits a right adjoint. (see \cite[Sect. 3.1.3]{Ga}) We write this adjunction as
\begin{equation}\label{eq: LSCat pull push}
f^*:\LSCat(Y_2)\leftrightarrows \LSCat(Y_1):f_*.
\end{equation}
    In the notation of \emph{loc. cit.}, $f^*$ corresponds to $\textbf{cores}_{f_B}$ and $f_*$ corresponds to $\textbf{coind}_{f_B}$, respectively. In the special case where $Y_2=\pt$, then this adjunction reads
    \[\ul{(\mhyphen)}_{Y_1}:\DGcat\leftrightarrows \LSCat(Y_1): \Gamma(Y_1,\mhyphen),\]
    where $(\mhyphen)_{Y_1}$ is the functor taking the constant local system on $Y_1$.
    \medskip
    
    The functor $f^*$ is symmetric monoidal, hence using Sect. \ref{sssec: adjunction and algebra} the adjunction \eqref{eq: LSCat pull push} lifts to the adjunctions
    \[f^*:\LSMon(Y_2)\leftrightarrows \LSMon(Y_1):f_*\]
    and
    \[f^*:\LSCMon(Y_2)\leftrightarrows \LSCMon(Y_1):f_*.\]
    
    The following result will play a crucial role:
    \begin{thm}[\cite{PPS}, Corollary 2.13]\label{thm: 1-affineness}
    Let $K$ be a 1-truncated homotopy type. Then the Betti stack $K_B$ is 1-affine.
    \end{thm}

\sssec{Examples arising from gerbes}\label{sssec: examples from gerbes}
Important examples of local system of categories for us arise from a gerbe on a space.
Suppose we have $(Y,\cG)\in\Spc+\on{Ge}_{\sfA}$, and we are given $\cX\in\PreStk^E$ with a $\Bet \sfA$-action. We have the Betti realization
$\cG_B: Y_B\rightarrow \Bet^2 A$, hence a functor
$Y\rightarrow \Btet{\sfA}$, $y\mapsto \cG_{B,y}$,
where $\cG_{B,y}$ is the composition
\[\pt\xrightarrow{y} Y_B\xrightarrow{\cG_B} \Bet^2 A.\]
Therefore, we get a functor
\[Y\rightarrow \PreStk^E,\quad y\mapsto \cX_{\cG_{B,y}},\]
where $\cX_{\cG_{B,y}}$ is the twist of $\cX$ by $\cG_{B,y}$ using the $\Bet A$-action. Taking quasi-coherent sheaves, we get a local system of symmetric monoidal categories
\[\cL_{(Y,\cG,\cX)}:Y\rightarrow \DGCatSymMon,\quad y\mapsto \QC(\cX_{\cG_{B,y}}).\]
Since $\on{Tot}(\cG)\times^{\Bet A}\cX=\colim_{y\in Y} \cX_{\cG_y}$ in $\PreStk^E$,  we have
\[\QC\left(\on{Tot}(\cG)\times^{\Bet A}\cX\right)=\lim_{y\in Y^{\op}} \QC\left(\cX_{\cG_y}\right).\]
This means
\[\Gamma(Y_B, \cL_{(Y,\cG,\cX)})=\QC\left(\on{Tot}(\cG)\times^{\Bet A}\cX\right).\]

More generally, we can start with $(Y,\cG)\in \Spc+\on{Ge}_{\sfA}$ with $\cC\in B\sfA\lMod(\DGcat)$, and form a local system of categories
\begin{equation}\label{eq: LS from gerbe and module}
    \cL_{(Y,\cG,\cC)}: Y\rightarrow \DGcat,\quad y\mapsto \cC_{\cG_y}.
\end{equation}
We have
\[\Gamma(Y_B,  \cL_{(Y,\cG,\cC)})=\lim_{y\in Y} \cC_{\cG_y},\]
and passing to right adjoints, this equals a colimit in $\DGcat$, hence a twisted integral computed in $\DGcat$:
\begin{equation}\label{eq: Gamma of LS form gerbe and module}
    \Gamma(Y_B, \cL_{(Y,\cG,\cC)})=\lim_{y\in Y} \cC_{\cG_y}=
\colim_{y\in Y} \cC_{\cG_y} =\int^{\plain}_{(Y,\cG)}\cC.
\end{equation}

If $\cC$ admits a lift to $B \sfA\lMod(\DGCatMon)$ (resp. $B\sfA\lMod(\DGCatSymMon)$), then the local system of categories \eqref{eq: LS from gerbe and module} admits a lift to $\LSMon(Y)$ (resp. $\LSCMon(Y)$). $\cL_{(Y,\cG,\cX)}$ introduced earlier is a special case of $\cL_{(Y,\cG,\cC)}$ for $\cC=\QC(\cX)\in B\sfA\lMod(\DGCatSymMon)$.

\sssec{Some notations concerning local systems}\label{sssec: some notations}
    We find it convenient to abuse the notation and denote a local system of categories by its category of global sections. For example:
    \begin{enumerate}
        \item consider the constant local system of categories over $Y$ with value $\cC\in\DGcat$. We use the notation $\cC\otimes \LS(Y)$ to denote this constant local system; 
        \item in view of \eqref{eq: Gamma of LS form gerbe and module}, we understand $\int^{\plain}_{(Y,\cG)}\cC$ as the local system $\cL_{(Y,\cG,\cC)}$.
    \end{enumerate}
 Given $I\in\fSet$ we have 
    \begin{equation}\label{eq: box tensor product and Gamma}
        (\cL_{(Y,\cG,\cC)})^{\boxtimes I}=  \cL_{(Y^I,\cG^{\boxtimes I},\cC^{\otimes I})}
    \end{equation}
    in $\LSCat(Y^I)$ and
    \[\Gamma\left((Y^I)_B,\cL_{(Y^I,\cG^{\boxtimes I},\cC^{\otimes I})}\right)=\left(\int^{\plain}_{(Y,\cG)}\cC\right)^{\otimes I}.\]
    Hence, similarly to the previous case, we understand $(\int^{\plain}_{(Y,\cG)}\cC)^{\otimes I}$ as the local system \eqref{eq: box tensor product and Gamma}.

\subsection{Functors out of twisted factorization homology}
The goal of this subsection is to understand what it means to have a monoidal functor out of twisted factorization homology. The main result Theorem \ref{thm: main categorical thm} means that given $\cA\in B\sfA\lMod(\DGCatSymMon)$ and a monoidal category $\cC$, a monoidal functor 
\[\int^{\CAlg}_{(Y,\cG)}\cC\rightarrow \cA\]
is equivalent to a compatible family of $\cG$-twisted local system of monoidal functors $\cC^{\otimes I}\rightarrow\cA$ on $Y^I$ for each finite set $I$.

\begin{rmk}
    Although Theorem \ref{thm: main categorical thm} concerns the case $\cA\in \DGCatMon$ and monoidal functors $\int^{\CAlg}_{(Y,\cG)}\cC\rightarrow\cA$, the same conclusion holds if we replace `monoidal' by `plain' or `symmetric monoidal'. We state it for the monoidal case because this is what is relevant to the construction of spectral action. 
\end{rmk}

\sssec{}
Suppose $\cC\in B\sfA\lMod(\DGCatSymMon)$ and $\cA\in \DGCatMon$. As observed in Sect. \ref{sssec: examples from gerbes}, the twisted integral $\int^{\on{plain}}_{(Y,\cG)}\cC$
is an object of $\CAlg(\LS(Y)\lMod)$. From the limit expression in \eqref{eq: Gamma of LS form gerbe and module} of twisted integral, for each $y\in Y$ we have a functor
\[\int^{\plain}_{(Y,\cG)}\cC\xrightarrow{\on{ev}_y}\cC_{\cG_y}.\]
This is symmetric monoidal, because the limit can be taken equally in $\DGCatSymMon$.

\sssec{}\label{sssec: fSet LSMon}
Consider the category $\fSet + \LSMon$ of the tuple $[I, \bA]$ of $I\in \fSet$ and $\bA\in \LSMon(Y^I)$. A morphism $[I_1, \bA_1]\rightarrow [I_2,\bA_2]$ is the data of $I_1\rightarrow I_2$ in $\fSet$ together with $\bA_1\rightarrow \bA_2$ in $\LSMon(Y^{I_2})$, where $\bA_1$ is thought of as an object of $\LSMon(Y^{I_2})$ under the pullback 
\[\LSMon(Y^{I_1})\rightarrow \LSMon(Y^{I_2})\]
along the diagonal $Y^{I_2}\rightarrow Y^{I_1}$. More precisely, this category, together with the forgetful functor $\fSet + \LSMon\rightarrow \fSet$ is the coCartesian fibration corresponding to the functor
\[\fSet\rightarrow \Cat,\quad I\mapsto \LSMon(Y^I).\]

\sssec{}\label{sssec: main cat thm prep}
Given $\cB\in \LSSymMon(Y)$ we construct a functor
\[\fSet\rightarrow \fSet + \LSMon,\quad I\mapsto [I, \cB^{\boxtimes I}].\]
Taking $\cB=\cL_{(Y,\cG,\cC)}$, we get the functor
\begin{equation}\label{eq: funct CX}
    \fSet\rightarrow \fSet + \LSMon, \quad I\mapsto \left[I, \left(\int^{\on{plain}}_{(Y,\cG)}\cC\right)^{\otimes I}\right].
\end{equation}

On the other hand, we have the functor
\begin{equation}\label{eq: funct AG}
    \fSet\rightarrow \fSet + \LSMon, \quad I\mapsto \left[ I, \cA\otimes\LS(Y^I) \right].
\end{equation}
We have the following assertion, which will be proved in the next subsection:
\begin{thm}\label{thm: main categorical thm}
    Let $(Y,\cG)\in \Spc+\on{Ge}_{\sfA}$, $\cC\in B\sfA\lMod(\DGCatSymMon)$ and $\cA\in\DGCatMon$. Consider the space of natural transformations from the functor \eqref{eq: funct CX} to \eqref{eq: funct AG} that is identified with $\fSet\xrightarrow{\id}\fSet$ after composing $\fSet+\LSMon\rightarrow \fSet$. This is equivalent to the space of maps $\int^{\CAlg}_{(Y,\cG)}\cC\rightarrow \cA$ in $\DGCatMon$. 
\end{thm}

\sssec{}\label{sssec: respects the projection}
For brevity, we say that a natural transformation $\eqref{eq: funct CX} \rightarrow \eqref{eq: funct AG}$ \textbf{respects the projection} if it satisfies the condition in the statement of the above theorem.

\sssec{Factorizable action}
Consider the special case where $\cM\in\DGcat$ and $\cA=\End(\cM)$. By the \textbf{factorizable action} of the local system
\[Y\xrightarrow{\cG}B^2\sfA\xrightarrow{\cC}\DGCatSymMon\]
of symmetric monoidal DG categories on the DG category $\cM$, we mean a natural transformation $\eqref{eq: funct CX}\rightarrow\eqref{eq: funct AG}$ respecting the projection. To justify the terminology, we observe that such a data includes a map
\begin{equation}\label{eq: LSMon map each I}
    \left(\int^{\on{plain}}_{(Y,\cG)}\cC\right)^{\otimes I}\rightarrow \End(\cM)\otimes\LS(Y^I)\quad\te{in }\LSMon(Y^I).
\end{equation}
Given a tuple $\ul{y}=(y_i)_{i\in I}\in Y^I$, we can evaluate this map at $\ul{y}$ to get
\[\bigotimes_{i\in I}\cC_{\cG_{y_i}}\rightarrow \End(\cM)\quad\te{in }\DGCatMon,\]
in other words, an action $\on{act}_{\ul{y}}:\bigotimes_{i\in I}\cC_{\cG_{y_i}}\acts\cM$. This leads us to the twisted analog of the situation of Sect. \ref{sssec: fact hlgy informal}. Local constancy is encoded in the fact that the collection $\{\on{act}_{\ul{y}}\}_{\ul{y}\in Y^I}$ in fact came from a map \eqref{eq: LSMon map each I}. The compatibility between the multiplicity of $\ul{y}$ and the symmetric monoidal operation of $\cC$ is encoded in the fact that the collection of maps \eqref{eq: LSMon map each I} should intertwine the functoriality in $I\in\fSet$.

\subsection{Proof of Theorem \ref{thm: main categorical thm}}
\sssec{}
We construct a map
\[\Maps_{\DGCatMon}\left(\int^{\CAlg}_{(Y,\cG)}\cC, \cA\right)\rightarrow \{\text{natural transformations } \eqref{eq: funct CX}\rightarrow \eqref{eq: funct AG} \te{ that respects the projection}\}.\]
For this purpose, it suffices to cover the case where \[\int^{\CAlg}_{(Y,\cG)}\cC\rightarrow \cA\]
is the identity. More precisely, it suffices to construct a natural transformation from \eqref{eq: funct CX} to the functor
\[\fSet\rightarrow\fSet + \LSMon,\quad I\mapsto \left[I, \int^{\CAlg}_{(Y,\cG)}\cC\otimes \LS(Y^I)\right].\]
This amounts to constructing a compatible family of functors
\begin{equation}\label{eq: plain to CAlg}
    \left(\int^{\on{plain}}_{(Y,\cG)}\cC\right)^{\otimes I}\rightarrow \int^{\CAlg}_{(Y,\cG)}\cC \otimes \LS(Y^I)
\end{equation}
in $\LSMon(Y^I)$.
\medskip

The left-hand side of \eqref{eq: plain to CAlg} is the local system of monoidal categories
\[Y^I\rightarrow \DGCatMon,\quad \ul{y}=(y_i)_{i\in I}\mapsto \bigotimes_{i\in I}\cC_{\cG_{y_i}},\]
whereas the right-hand side of \eqref{eq: plain to CAlg} is the constant local system of monoidal categories over $Y^I$ with value $\int^{\CAlg}_{(Y,\cG)}\cC$.
The functor \eqref{eq: plain to CAlg} comes from the family of monoidal functors
\[\bigotimes_{i\in I}\cC_{\cG_{y_i}}\xrightarrow{\otimes_{i\in I}\on{ins}_{y_i}} \left(\int^{\CAlg}_{(Y,\cG)}\cC\right)^{\otimes I}\xrightarrow{\on{mult}}\int^{\CAlg}_{(Y,\cG)}\cC\]
parametrized by $\ul{y}\in Y^I$.
Note that all the arrows are symmetric monoidal, hence monoidal.

\sssec{}
    From \cite[Proposition 5.1]{GHN} (see also \cite[Lemma 1.3.12]{GKRV}), given two categories $\bD$, $\bC$ and functors
    \[\Phi_1, \Phi_2: \bD\rightarrow \bC,\]
    the space of natural transformations from $\Phi_1\rightarrow\Phi_2$ gets identified with 
    \[\lim_{(\bd_s\rightarrow \bd_t)\in\TwArr(\bD)^{\op}} \Maps_{\bC}(\Phi_1(\bd_s),\Phi_2(\bd_t)).\]

    Hence, the space of natural transforms from \eqref{eq: funct CX} to \eqref{eq: funct AG} is the limit
    \[\lim_{(I\rightarrow J)\in \TwArr(\fSet)^{\op}} \Maps_{\fSet + \LSMon}\left(\left[I,\left(\int^{\plain}_{(Y,\cG)}\cC\right)^{\otimes I}\right], \left[J, \cC\otimes\LS(Y^J)\right]\right).\]
    Under the diagonal map $Y^J\rightarrow Y^I$, the local system
    \[\left(\int^{\plain}_{(Y,\cG)}\cC\right)^{\otimes I}=\cL_{(Y^I,\cG^{\boxtimes I},\cC^{\otimes I})}\in\LSMon(Y^I)\]
    restricts to the local system
    \begin{equation}\label{eq: local system restricted}
        \cL_{(Y^J,\cG^{\boxtimes J},\cC^{\otimes I})}\in \LSMon(Y^J),
    \end{equation}
    where the $(B\sfA)^{\times J}$-action on $\cC^{\otimes I}$ is as in Sect. \ref{sssec: cI integrated over YJ}. To be consistent with the notational convention in Sect. \ref{sssec: some notations}, we use the notation
    \[\int^{\plain}_{(Y^J,\cG^{\boxtimes J})}\cC^{\otimes I}\]
    in place of \eqref{eq: local system restricted}. Therefore, the space of natural transformations $\eqref{eq: funct CX}\rightarrow \eqref{eq: funct AG}$ that respects the projection is the limit
    \[\lim_{(I\rightarrow J)\in \TwArr(\fSet)^{\op}} \Maps_{\LSMon(Y^J)}\left(\cL_{(Y^J,\Delta^*\cG^{\boxtimes I},\cC^{\otimes I})},  \cC\otimes\LS(Y^J)\right).\]
    
    \medskip
    
    On the other hand, from Theorem \ref{thm: twArr twfact comparison} we get
    \[\Maps_{\DGCatMon}\left(\int^{\CAlg}_{(Y,\cG)} \cC,\cA\right) = \lim_{(I\rightarrow J)\in\TwArr(\fSet)^{\op}}\Maps_{\DGCatMon}\left(\int^{\Alg}_{(Y^J,\cG^{\boxtimes J})}\cC^{\otimes I}, \cA\right).\]
    Therefore, to prove Theorem \ref{thm: main categorical thm}, it suffices to show that for any $I\rightarrow J$ in $\fSet$, there is a canonical equivalence between
    \begin{equation}\label{eq: plain side}
        \Maps_{\LSMon(Y^J)}\left(\int^{\plain}_{(Y^J,\cG^{\boxtimes J})}\cC^{\otimes I},  \cA\otimes\LS(Y^J)\right)
    \end{equation}
    and
    \begin{equation}\label{eq: Mon side}
        \Maps_{\DGCatMon}\left(\int^{\Alg}_{(Y^J,\cG^{\boxtimes J})}\cC^{\otimes I}, \cA\right).
    \end{equation}

    \sssec{}\label{sssec: functor from Mon to plain}
We construct a canonical functor from \eqref{eq: Mon side} to \eqref{eq: plain side}. It suffices to construct a functor
\[\int^{\plain}_{(Y^J,\cG^{\boxtimes J})}\cC^{\otimes I}\rightarrow \left(\int^{\Alg}_{(Y^J,\cG^{\boxtimes J})}\cC^{\otimes I}\right)\otimes \LS(Y^J)\]
in $\LSMon(Y^J)$. 
The left-hand side is the local system
\[Y^J\rightarrow \DGCatMon,\quad \ul{y}\in Y^J\mapsto (\cC^{\otimes I})_{(\cG^{\boxtimes J})_{\ul{y}}}\]
while the right-hand side is the constant local system with value 
\[ \int^{\Alg}_{(Y^J,\cG^{\boxtimes J})}\cC^{\otimes I}.
\]

The desired functor comes from the family of monoidal functors
\begin{equation}\label{eq: comparison as Alg}
    (\cC^{\otimes I})_{(\cG^{\boxtimes J})_{\ul{y}}}\xrightarrow{\on{ins}_{\ul{y}}}\int^{\Alg}_{(Y^J,\cG^{\boxtimes J})}\cC^{\otimes I}
\end{equation}
parametrized by $\ul{y}\in Y^J$. 
In the remainder of this subsection, we prove:

\begin{prop}\label{prop: Mon side and plain side}
    The canonical functor from \eqref{eq: Mon side} to \eqref{eq: plain side} 
    is an equivalence.
    \end{prop}

\sssec{}\label{sssec: trivializable case}
We first confirm the equivalence in the case $\cG=\triv_Y$. In this case \[\int^{\plain}_{(Y^J,\triv_{Y^J})}\cC^{\otimes I}=\ul{(\cC^{\otimes I})}_{Y^J}\in\LSMon(Y^J)\]
and using Sect. \ref{sssec: LSCat pull push} the space \eqref{eq: plain side} is equivalent to 
\[\Maps_{\DGCatMon}\left(\cC^{\otimes I}, \cC\otimes \LS(Y^J)\right).\]
We conclude as
\begin{align*}
    \Maps_{\DGCatMon}\left(\cC^{\otimes I}, \cC\otimes \LS(Y^J)\right)&=\Maps_{\DGCatMon}\left(\cC^{\otimes I}, \lim_{Y^J}\cC\right)\\
    &=\Maps_{\DGCatMon}\left(\colim_{Y^J}\cC^{\otimes I},\cC\right)=\Maps_{\DGCatMon}\left(\int^{\Alg}_{(Y^J,\triv_{Y^J})}\cC^{\otimes I}, \cC\right).
\end{align*}
The fact that the canonical functor is an equivalence for $\cG=\triv_Y$ implies the same is true for the case $\cG$ is trivializable.

\sssec{}
    Now we cover the general case of $(Y,\cG)$. Choose $P\rightarrow Y$ in $\Spc$ which induces a surjection on $\pi_0$ and trivializes $\cG$, e.g., $P=\on{Tot}(\cG)$.
    We have
    \begin{equation}\label{eq: descent LSMon}
        \LSMon(Y^J)=\lim_{[n]\in\Delta}\LSMon((P^n/Y)^J)
    \end{equation}
    as in Sect. \ref{sssec: descent of LSMon}.
    Under this equivalence, $\cC\otimes\LS(Y^J) \in \LSMon(Y^J)$ corresponds to the collection
    \[\cC\otimes\LS((P^n/Y)^J)\in\LSMon((P^n/Y)^J))\]
    On the other hand,
    \[\int^{\plain}_{(Y^J, \cG^{\boxtimes J})}\cC^{\otimes I}\in \LSMon(Y^J)\]
    corresponds to the collection
    \[\int^{\plain}_{((P^n/Y)^J,(\cG^{(n)})^{\boxtimes J})}\cC^{\otimes I}\in \LSMon((P^n/Y)^J),\]
    where $\cG^{(n)}$ is the pullback of $\cG$ along $P^n/Y\rightarrow Y$.
    \medskip
    
Using the equivalence \eqref{eq: descent LSMon}, we have an equivalent between \eqref{eq: plain side} and
\begin{equation}\label{eq: lim mapping space}
\lim_{[n]\in\Delta} \Maps_{\LSMon((P^n/Y)^J)}\left(\int^{\plain}_{((P^n/Y)^J,(\cG^{(n)})^{\boxtimes J})}\cC^{\otimes I}, \cC\otimes \LS((P^n/Y)^J)\right).
\end{equation}
Each $\cG^{(n)}$ is trivializable, hence applying the result of Sect. \ref{sssec: trivializable case} to $(P^n/Y,\cG^{(n)})$ we have an equivalence
\begin{align*}
    &\Maps_{\LSMon((P^n/Y)^J)}\left(\int^{\plain}_{((P^n/Y)^J,(\cG^{(n)})^{\boxtimes J})}\cC^{\otimes I}, \cC\otimes \LS((P^n/Y)^J)\right)\\
    &= \Maps_{\DGCatMon}\left(\int^{\Alg}_{( (P^n/Y)^J, (\cG^{(n)})^{\boxtimes J})}\cC^{\otimes I}, \cC\right).
\end{align*}
Therefore, \eqref{eq: lim mapping space} is equivalent to
\begin{equation}\label{eq: lim mapping Alg side}
    \lim_{[n]\in\Delta}\Maps_{\DGCatMon}\left(\int^{\Alg}_{( (P^n/Y)^J, (\cG^{(n)})^{\boxtimes J})}\cC^{\otimes I}, \cC\right),
\end{equation}
which is equivalent to \eqref{eq: Mon side} via Theorem \ref{thm: fact hlgy and colimit}.
\medskip

In conclusion, we proved a chain of canonical equivalences
\[\eqref{eq: plain side} \rightarrow \eqref{eq: lim mapping space}\rightarrow \eqref{eq: lim mapping Alg side} \rightarrow \eqref{eq: Mon side},\]
hence the proof of Proposition \ref{prop: Mon side and plain side} and thus Theorem \ref{thm: main categorical thm}.

\section{Factorizable metaplectic Satake functor}\label{sec: metaplectic Satake functor}
In previous sections, we have identified $\QC(\Loc^{\cG_Z}_H(X))$ with twisted factorization homology \[\int^{\CAlg}_{(X^{\topl,\sing},\Phi(\cG_Z))}\Rep(H),\] and developed a categorical framework for the action of twisted factorization homology. From now on, we prove that the multiple point metaplectic Hecke action fits into our categorical framework, i.e., that we have a canonical natural transformation $\eqref{eq: funct CX}\rightarrow \eqref{eq: funct AG}$ arising from the metaplectic Hecke action where 
\[\sfA=Z_H(E),\quad \cC=\Rep(H),\quad  (Y,\cG)=(X^{\topl,\sing},\Phi(\cG_Z)), \te{ and } \cA=\End(\cA^q_G).\]
The desired natural transform, when evaluated at each $I\in\fSet$, should take a form
\[\left(\int^{\plain}_{(X^{\topl,\sing},\Phi(\cG_Z))}\Rep(H)\right)^{\otimes I}\rightarrow \End(\cA^q_G)\otimes \LS(X^I) \quad \te{in }\LSMon(X^I),\]
and after taking global sections it has a form\footnote{The reason why the global section of the source of the above functor is the source of the below functor will be explained later in the proof of Proposition \ref{prop: equiv LS global}.}
\[\left(\Rep(H)^{\loc}_{\cG_Z}(X)\right)^{\otimes I}\rightarrow\End(\cA^q_G)\otimes \LS(X^I) \quad\te{in }\Alg(\LS(X^I)\lMod).\]

We can break this into two parts: the `hard' part \footnote{Here, $\Hecke^{\loc}_{X^I}$ stands for the $I$-legged local Hecke stack. Its definition will be recalled in Sect. \ref{sssec: stacks on the A side}.}
\[\left(\Rep(H)^{\loc}_{\cG_Z}(X)\right)^{\otimes I}\rightarrow \Shv^q(\Hecke^{\loc}_{X^I})\quad \te{in }\Alg(\LS(X^I)\lMod)\]
and the `easy' part \footnote{Strictly speaking, we do not have a natural transform of this form. To make sense of this functor, we should restrict the source category to a certain full subcategory. We ignore this point here for simplicity.}
\[\Shv^q(\Hecke^{\loc}_{X^I})\rightarrow \End(\cA^q_G)\otimes \LS(X^I) \quad \te{in }\Alg(\LS(X^I)\lMod).\]
Obviously, the first part is harder because it involves geometric Satake, while the second part is entirely on the automorphic side. While this was a description after evaluation on each $I$, we also have to construct functoriality in $\fSet$ in order to upgrade it into a natural transformation.
\medskip

In this section and the following Section \ref{sec: Satake as a nat transf} we do the hard part of the construction. We construct the functor on each $I$ as $\LS(X^I)$-linear functors \footnote{We do not construct it as an $\E_1$-functor. The $\E_1$-structure and functoriality will be established simultaneously in Sect. \ref{sec: Satake as a nat transf}.}, and the next Section upgrades the collection of functors for each $I$ into a desired natural transformation.

\sssec{}\label{sssec: explaining main results}
Now we explain the main result of this section.
We construct a collection of $\LS(X^I)$-linear functors
\[\cS_{X,I}:\left(\Rep(H)^{\loc}_{\cG_Z}(X)\right)^{\otimes I}\rightarrow \Shv^q(\Hecke^{\loc}_{X^I})\]
for each $I\in\fSet$ (Sect. \ref{sssec: construction of S QL I}). The construction of $\cS_{X,I}$ relies on the pre-existing result of perverse ULA version established in \cite{Rei} and \cite{ZhaoNew}. What we do here is to reformulate their results through various categorical manipulations in a form that is best suited to our purpose. Let us call the functors $\cS_{X,I}$ above the $\LS$-linear versions of factorizable Satake, to distinguish them from perverse ULA versions. 
\medskip

Our construction of $\LS$-linear Satake functors will be a bit involved and also indirect in the sense that we have access to $\cS_{X,I}$ for general $I$'s only through the singleton case. Therefore, not only do we not have a clear idea of what this functor looks like, but it would also be difficult to establish various compatibilities needed to upgrade it to a natural transformation. To overcome these difficulties altogether, we prove a certain compatibility of $\cS_{X,I}$ with the perverse ULA Satake equivalence (Theorem \ref{thm: characterization for I}) and moreover prove that such a compatibility uniquely characterizes the functor (Lemma \ref{lem: the characterization lemma}).

\sssec{}\label{sssec: technical aspects LS linear}
We discuss some technical aspects of the construction of $\cS_{X,I}$'s, which should clarify the need for the more involved manipulations that appear later. The (metaplectic) perverse ULA Satake equivalence yields us a \emph{monoidal} equivalence
\[\cS^{\Ab}_{X,I}:\Rep(H^I)^{\Ab,\loc,\ft}_{\cG^{\boxtimes I}_Z}(X^I)\xrightarrow{\sim}\Sat^q_{X^I},\]
where the left-hand side is the category of twisted $H^I$-representation on (abelian) local systems of finite type on $X^I$ and the right-hand side is a twisted Satake category over the affine Grassmannian over $X^I$. (we review their definition in Sect. \ref{ssec: recollection perv ULA}) The source (resp. the target) is a full subcategory of the source (resp. the target) of $\cS_{X,I}$:
\[\Rep(H^I)^{\Ab,\loc,\ft}_{\cG^{\boxtimes I}_Z}(X^I)\subseteq \left(\Rep(H)^{\loc}_{\cG_Z}(X)\right)^{\otimes I}\quad\te{and}\quad \Sat^q_{X^I}\subseteq \Shv^q(\Hecke^{\loc}_{X^I}).\]

Since $\left(\Rep(H)^{\loc}_{\cG_Z}(X)\right)^{\otimes I}$ can be understood as the category of twisted $H^I$-representations on (derived) local systems on $X^I$, the difference between $\Rep(H^I)^{\Ab,\loc,\ft}_{\cG^{\boxtimes I}_Z}(X^I)$ and $\left(\Rep(H)^{\loc}_{\cG_Z}(X)\right)^{\otimes I}$ does not seem to be that significant: we are changing the notion of local systems from $\LS^{\Ab,\ft}$ to $\LS$, where the former means abelian local systems of finite type. So, at the first glance, one might guess that the construction  of $\cS_{X,I}$ would follow immediately from the perverse Satake equivalence.

\sssec{}\label{sssec: subtleties LS linear}
However, the passage turns out to be more subtle. Comparing the source categories in two setups, we note the following two differences:
\begin{enumerate}
    \item perverse ULA version concerns abelian $H^I$-representations while the desired functor $\cS_{X,I}$ concerns derived $H^I$-representations;
    \item a finiteness condition was imposed on the spectral side of the perverse ULA version while it is not the case for the spectral side of \eqref{eq: fact Satake I singleton}. 
\end{enumerate}

For the first point, we wish to apply the universal property of the derived functors. However, we run into a problem here, because the perverse ULA functor $\cS^{\Ab}_{X,I}$ lands in the category of perverse sheaves while $\Shv^q(\Hecke^{\loc}_{X^I})$ does \emph{not} have a perverse $t$-structure. To fix this, we are led to consider the $\Ind$-constructible version $\Shv^{\on{ind-constr},q}(\Hecke^{\loc}_X)$ as an intermediate. Note that the convolution operation can be defined for ind-constructible sheaves in the same way so that we have a convolution monoidal structure on $\Shv^{\on{ind-constr},q}(\Hecke^{\loc}_{X^I})$ \footnote{This still uses $*$-pullback and $!$-pushforward.}. Moreover, we have a canonical monoidal functor
\begin{equation}\label{eq: indconstr to Betti}
    \Shv^{\on{ind-constr},q}(\Hecke^{\loc}_{X^I})\rightarrow \Shv^{q}(\Hecke^{\loc}_{X^I})
\end{equation}
using the fact that the source category is compactly generated, and its compact objects belong to the full subcategory of constructible sheaves in $\Shv^q(\Hecke^{\loc}_{X^I})$ (\cite[Sect. F.1.1, Sect. F.2.2]{AGKRRV}). Roughly speaking, the plan of our construction is to factor \eqref{eq: fact Satake I singleton} through \eqref{eq: indconstr to Betti}. Since $\Shv^{\indconstr,q}(\Hecke^{\loc}_{X^I})$ has a perverse $t$-structure and whose heart contains $\Sat^{q}_{X^I}$, we can construct a derived functor into $\Shv^{\indconstr,q}(\Hecke^{\loc}_{X^I})$, and then compose with \eqref{eq: indconstr to Betti} to land in the right target.
\medskip

The second point is more subtle. In many cases, removing the finiteness condition amounts to passing to Ind-completion. However, this does not work here, because the Ind-completion of $\LS^{\Ab,\ft}(X^I)$ is \emph{not} in general equivalent to $\LS(X^I)^{\heartsuit}$. Therefore, rather than understanding the passage from $\LS^{\Ab,\ft}(X^I)$ to $\LS(X^I)$ as an Ind-completion, it is better to understand it as \emph{tensoring up} with the category $\LS(X^I)$. To make this idea precise, we need to make an \emph{opposite twist} on both sides to $\cS_{X,I}$, which means that we twist both sides of $\cS_{X,I}$ by the negative of $\cG_Z$ so that the source category gets turned into 
\[(\Rep(H)^{\loc}(X))^{\otimes I}=\Rep(H^I)\otimes \LS(X^I).\]
In fact, this manipulation plays a crucial role in the proof of perverse ULA Satake (see Sect. \ref{sssec: opposite twist}). While the opposite twist can be done for every $I$'s on the level of abelian categories, we know how to do it for derived categories only when $I$ is a singleton.\footnote{We do not know how to twist $\Shv^{\indconstr,q}(\Hecke^{\loc}_{X^I})$ in general.} For this reason, we construct $\cS_{X,I}$ for $I=\{*\}$ following the aforementioned strategy and then construct $\cS_{X,I}$ for general $I$'s using the external convolution product.\footnote{This is one of the reasons why our construction of $\cS_{X,I}$ is indirect.}

\sssec{The role of $\QLisse$ sheaves}
As discussed above, our construction of $\cS_{X,I}$, say for $I=\{*\}$, should factor through the category of $\Ind$-constructible sheaves on $\Hecke^{\loc}_{X}$. Since this category does not have a $\LS(X)$-linear structure, it does not make sense to construct a $\LS(X)$-linear functor into it. What we actually construct is a $\QLisse(X)$-linear version:
\[\cS^{\QLisse}_{X,\{*\}}:\Rep(H)^{\QLisse}_{\cG_Z}(X)\rightarrow \Shv^{\indconstr,q}(\Hecke^{\loc}_{X}),\]
where the source is the category of twisted $H$-representations on $\QLisse$ sheaves. (see Sect. \ref{sssec: list of sheaves}) The desired $\LS(X)$-linear functor $\cS_{X,\{*\}}$ is obtained by composing $\cS^{\QLisse}_{X,\{*\}}$ with \eqref{eq: indconstr to Betti} and using the $\LS(X)$-linearity on $\Shv^q(\Hecke^{\loc}_X)$. Similarly, the construction of $\LS(X^I)$-linear functor $\cS_{X,I}$ follows from the $\QLisse(X)^{\otimes I}$-linear version
\[\cS^{\QLisse}_{X,I}: \left(\Rep(H)^{\QLisse}_{\cG_Z}(X)\right)^{\otimes I}\rightarrow \Shv^{\indconstr,q}(\Hecke^{\loc}_{X^I}).\]

Although what we ultimately need to use the twisted factorization homology is the $\LS$-linear functors $\cS_{X,I}$, what actually plays a prominent role in the construction of those input is the $\QLisse$-linear versions $\cS^{\QLisse}_{X,I}$. The above discussion provides one of the reasons for that. Furthermore, the $\QLisse$-linear version interacts very nicely with the perverse ULA version, in the sense that $\cS^{\QLisse}_{X,I}$ is uniquely determined by the perverse ULA Satake $\cS^{\Ab}_{X,I}$.\footnote{We discussed the importance of this property in Sect. \ref{sssec: explaining main results}.} We refer the readers to Lemma \ref{lem: the characterization lemma} for the precise statement. This has to do with the fact that $\QLisse(X)$ is compactly generated by $\LS^{\Ab,\ft}(X)$ if $X\ne \pjl$. Note that the analogous statement for $\LS(X)$ is false.

\sssec{}
Let us describe briefly the structure of this section.
In Sect. \ref{ssec: recollection perv ULA} we review the perverse ULA Satake equivalence. In Sect. \ref{ssec: LS linear fact Satake} we construct the functor $\cS_{X,I}$ for the case $I$ is a singleton, assuming an existence of certain graded lift of the target (see Theorem \ref{thm: A side grading} for a precise statement). This graded lift will be constructed in Sect. \ref{ssec: A side grading}. We then proceed to the construction of $\cS_{X,I}$ for general $I$'s together with their compatibility with perverse ULA Satake equivalence. Sect. \ref{ssec: various versions of local systems} and \ref{ssec: full compatibility with extensions} both serve as a preparation for Sect. \ref{ssec: characterization lemmas} and \ref{ssec: LS linear general}. Since various versions of local systems will be used in our argument, we record some facts about them in Sect. \ref{ssec: various versions of local systems}. For the purpose of studying the compatibility between the perverse ULA and $\LS$-linear versions, it is useful to study how extensions in an abelian category and a DG category interact with each other, and this is done in Sect. \ref{ssec: full compatibility with extensions}. After developing these tools, in Sect. \ref{ssec: characterization lemmas} we provide a characterization of functors in a form\footnote{This means that the functors having the same source and linearity.} $\cS^{\QLisse}_{X,I}$. Finally, in Sect. \ref{ssec: LS linear general} we proceed to the construction of $\cS_{X,I}$ for general $I$'s and establish compatibility with perverse ULA Satake (Theorem \ref{thm: characterization for I}). From the characterization established before, this compatibility uniquely characterizes the $\QLisse$-linear version of Satake functors, and therefore gives us a good control of these functors. This fact will play a crucial role in Sect. \ref{sec: Satake as a nat transf}.

\begin{nota}\label{nota: introducing oneCatEc}
    Let $\VectAbfd$ be the abelian category of finite-dimensional $E$-vector spaces and let
\[\CatEc:=\Vect^{\heartsuit,\on{f.d.}}\lMod(\Cat).\]
Using the Cartesian monoidal structure of $\Cat$ we can equip $\CatEc$ with the natural symmetric monoidal structure $\Ectensbig$ by means of relative tensor product.\footnote{By adding $\on{f.d.}$ on the top, we distinguish it from the Lurie tensor product.} We will sometimes write the maps in a form
$C_1\Ectens \cdots\Ectens C_r\rightarrow C$
simply as $C_1\times\cdots\times C_r\rightarrow C$, where the $\Vect^{\heartsuit,\on{f.d.}}$-multilinear structure is implicitly understood.
\end{nota}

\sssec{Subtlety with $\Shv^q(\Hecke^{\loc}_{X^I})$}\label{sssec: subtlety Hecke}
Since the category $\Shv^q(\Hecke^{\loc}_{X^I})$ plays a major role in the remainder of this paper, let us clarify on the definition of it. The expression $\Shv^q(\Hecke^{\loc}_{X^I})$ does not make a literal sense because $\Hecke^{\loc}_{X^I}$ does \emph{not} lie in $\PreStklft$. We define $\Shv^q(\Hecke^{\loc}_{X^I})$ as a colimit of twisted sheaf categories on \emph{finite} pieces of the local Hecke stack, analogous to the approach of
\cite[Sect. 9.2]{Ach} defining the Satake category. More precisely, we define
\[\Shv^q(\Hecke^{\loc}_{X^I})=\colim_{\ul{\lambda}\in (\Lambda^+)^I} \Shv^q(\Hecke^{\loc,\le\ul{\lambda}}_{X^I})\quad\te{in }\DGcat.\]
We explain some notations:
\begin{enumerate}
    \item Let $(L^+G)_I\rightarrow X^I$ be the $I$-legged arc group, and $\Gr^{\le \ul{\lambda}}_{X^I}\rightarrow X^I$ is the Schubert variety of the $I$-legged affine Grassmannian, associated to $\ul{\lambda}$. 
    The action of $(L^+G)_I$ on $\Gr^{\le \ul{\lambda}}_{X^I}$ factors through the quotient $(L^+G)_I\rightarrow (L^{n(\ul{\lambda})}G)_I$
where the target means the ($I$-legged version of) $n(\ul{\lambda})$-th jet group scheme of $G$. We let
\[\Shv^q(\Hecke^{\loc,\le\ul{\lambda}}_{X^I}):= \Shv^q((L^{n} G)_I\backslash \Gr^{\le \ul{\lambda}}_{X^I})\]
for any $n\ge n(\ul{\lambda})$.

    \item $\Lambda^+\subseteq \Lambda$ is the set of dominant cocharacters of $G$, and this is considered as a poset in a natural way. Hence, we view $(\Lambda^+)^I$ as a category. In the colimit expression, each transition map is induced by the $!$-pushforward along the closed embedding between $\Gr^{\le\ul{\lambda}}_{X^I}$'s.
\end{enumerate}
For the category of twisted sheaves on the variant of the local Hecke stack, such as the iterated Hecke stack $\wt{\Hecke}^{\loc,[k]}_{X^I}$ we use a similar approach of taking colimit over the categories arising from the finite pieces. 
\medskip

We will mostly ignore this subtlety and work as if $\Shv^q(\Hecke^{\loc}_{X^I})$ and other variants were a part of the six-functor formalism $\Shv^{\on{tw}}_{\Corr}$
This is justified by the fact that various operations, e.g., convolution product, we perform on those spaces can be approximated by the operations on finite pieces. For example, the convolution product
\[\Shv^q(\Hecke^{\loc}_{X^I})\otimes \Shv^q(\Hecke^{\loc}_{X^I})\rightarrow \Shv^q(\Hecke^{\loc}_{X^I})\]
is approximated by
\[\Shv^q(\Hecke^{\loc,\le\ul{\lambda}}_{X^I})\otimes \Shv^q(\Hecke^{\loc,\le\ul{\mu}}_{X^I})\rightarrow \Shv^q(\Hecke^{\loc,\le\ul{\lambda}+\ul{\mu}}_{X^I}).\]

We take the same approach for $\Ind$-constructible sheaves.

\subsection{Recollection on the perverse ULA Satake equivalence}\label{ssec: recollection perv ULA}
The construction of the factorizable Satake functor in the metaplectic situation was studied in \cite{Rei, GL, ZhaoNew}. The version most relevant to us is that of \cite[Theorem 2.3.2]{ZhaoNew}, which we review in this subsection. To fit our purpose, we state the result using sheaves in analytic topology.
Although the result of \emph{loc. cit.} was stated in the context of \'etale sheaves, the same proof can be applied to the analytic context. (see Remark 2.3.4 of \emph{op. cit.})
\medskip

The perverse ULA version asserts a symmetric monoidal equivalence between the metaplectic version of Satake category over $X^I$ and the category of twisted $H^I$-representations on local systems on $X^I$.

\sssec{Stacks on the automorphic side}\label{sssec: stacks on the A side}
For each $I\in\fSet$ we consider the following three stacks:
\[\Hecke^{\loc}_{X^I}, \quad \Hecke_{X^I}, \ \te{and} \ \Gr_{X^I},\]
called the \textbf{local Hecke stack}, \textbf{Hecke stack}, and Beilinson-Drinfeld \textbf{affine Grassmannian}. We call them e.g., $I$-legged local Hecke stack, if we want to make $I$ explicit. We review the functor of points definition of those stacks. 
\medskip

Each of them is a stack over $X^I$. Given an $S$-point $\ul{x}$ of $X^I$, let $\Gamma_{\ul{x}}\subseteq X\times S$ be the graph of $\ul{x}$ and $D_{\ul{x}}$ be the formal completion of $X\times S$ along $\Gamma_{\ul{x}}$, and $D^{\circ}_{\ul{x}}= D_{\ul{x}}\setminus \Gamma_{\ul{x}}$.
\begin{enumerate}
    \item An $S$-point of $\Hecke^{\loc}_{X^I}$ is a tuple $(\ul{x}, \cP_1,\cP_2, \alpha)$ where $\ul{x}$ is an $S$-point of $X^I$, each $\cP_i$ is a $G$-bundle over $D_{\ul{x}}$, and $\alpha$ is an isomorphism of $\cP_1$ and $\cP_2$ over $D^{\circ}_{\ul{x}}$.
    
    \item An $S$-point of $\Hecke_{X^I}$ is a tuple $(\ul{x}, \cP_1,\cP_2, \alpha)$ where $\ul{x}$ is an $S$-point of $X^I$, each $\cP_i$ is a $G$-bundle over $X\times S$, and $\alpha$ is an isomorphism of $\cP_1$ and $\cP_2$ over $X\times S\setminus \Gamma_{\ul{x}}$.

    \item An $S$-point of $\Gr_{X^I}$ is a tuple $(\ul{x},\cP,\alpha)$ where $\ul{x}$ is an $S$-point of $X^I$, each $\cP$ is a $G$-bundle over $X\times S$, and $\alpha$ is a trivialization of $\cP$ over $X\times S\setminus \Gamma_{\ul{x}}$.
\end{enumerate}
By construction, for each $I\in\fSet$ we have canonical maps \[\Gr_{X^I}\rightarrow\Hecke_{X^I}\rightarrow\Hecke^{\loc}_{X^I}\rightarrow X^I.\]

\sssec{automorphic side}\label{sssec: A side of Satake}
 The metaplectic parameter $q$ canonically induces an $A$-gerbe $\cG^I$ on $\Hecke^{\loc}_{X^I}$. We denote the resulting category of twisted sheaves on $\Hecke^{\loc}_{X^I}$ as $\Shv^q(\Hecke^{\loc}_{X^I})$.
We define $\Sat^q_{X^I}$ as the full subcategory of $\Shv^q(\Hecke^{\loc}_{X^I})$ spanned by the objects whose pullback to $\Gr_{X^I}$ lies in the category
\begin{equation}\label{eq: Perv ULA subcat}
    \Perv^{\ULA,q}(\Gr_{X^I})[-|I|]\subset \Shv^q(\Gr_{X^I}). 
\end{equation}
Here the notation $\Perv^{\ULA}$ means the category of (twisted) perverse sheaves on $\Gr_{X^I}$ that are universally locally acyclic (ULA) with respect to the map $\Gr_{X^I}\rightarrow X^I$. The cohomological shift in \eqref{eq: Perv ULA subcat} makes it convenient to state the compatibility pertaining maps $I\rightarrow J$. 
\medskip

This construction of perverse ULA satake category upgrades to a sheaf. The assignment 
\begin{equation}\label{eq: sheaf of Satake categories}
  \textbf{Sat}^q_I: (X_{\et})^{\op}\rightarrow \CatEc,\quad X'\mapsto \Sat^q_{(X')^I} 
\end{equation}
is a sheaf of $\VectAbfd$-linear categories on the small \'etale site of $X$. The category $\Sat^q$ has a \emph{convolution} monoidal structure (\cite[Sect. 1.3.5]{ZhaoNew}) lifting \eqref{eq: sheaf of Satake categories} as a sheaf of $\VectAbfd$-linear monoidal categories, where we use $*$-pullbacks and $!$-pushforward in the formation of convolution product.
In fact, this monoidal structure canonically lifts to a symmetric monoidal structure, by identifying the convolution product with a \emph{fusion} product.
(\cite[Sect. 1.3.12]{ZhaoNew})

\sssec{Convolution product}\label{sssec: convolution product}
For future reference, we briefly review the convolution product of the Satake category. We refer the reader to \cite[Sect. 1.3.5]{ZhaoNew} for details that we omit. However, this construction is completely parallel to the non-metaplectic case, apart from the necessary compatibilities among gerbes required for the pull–push constructions. Thus, readers already familiar with the non-metaplectic setting, as in \cite[Sect. 5]{Zhu}, may safely take these constructions on faith.
\medskip

The convolution monoidal structure on $\Sat^q_{X^I}$ is in fact a restriction of the convolution monoidal structure on $\Shv^q(\Hecke^{\loc}_{X^I})$, i.e., the monoidal operation of the latter category preserves the full subcategory $\Sat^q_{X^I}$.
To construct the convolution monoidal structure on $\Shv^q(\Hecke^{\loc}_{X^I})$ we consider the \textbf{iterated local Hecke stack} $\wt{\Hecke}^{\loc, [k]}_{X^I}$ for each $k\ge 0$, whose $S$-point is a tuple 
\[(\ul{x}, \cP^0,\cdots, \cP^k, \alpha_1,\cdots,\alpha_k),\] 
where $\ul{x}$ is an $S$-point of $X^I$, and each $\cP_i$ is a $G$-bundle over $D_{\ul{x}}$, and each $\alpha_a$ is an isomorphism between $\cP_{a-1}$ and $\cP_a$ over $D^{\circ}_{\ul{x}}$. Then we have a correspondence
\begin{equation}\label{eq: convolution diagram}
    \begin{tikzcd}
	{\wt{\Hecke}^{\loc,[k]}_{X^I}} & {\prod^k_{a=1} \Hecke^{\loc}_{X^I}} \\
	{\Hecke^{\loc}_{X^I}.}
	\arrow["{\prod p_a}", from=1-1, to=1-2]
	\arrow["m"', from=1-1, to=2-1]
\end{tikzcd}
\end{equation}
Moreover, it can be shown that the gerbe on $\Hecke^{\loc}_{X^I}$ (denoted $\cG^I$ in Sect. \ref{sssec: A side of Satake}) is compatible with the correspondence, so that we have a functor
\begin{equation}\label{eq: product kth level}
    \bigotimes^k_{a=1}\Shv^q(\Hecke^{\loc}_{X^I})\rightarrow \Shv^q(\Hecke^{\loc}_{X^I})
\end{equation}
using $*$-pullback and $!$-pushforward. The collection of \eqref{eq: product kth level} defines a monoidal structure on $\Shv^q(\Hecke^{\loc}_{X^I})$, called the \textbf{convolution monoidal structure}.
\medskip

We can also consider an external convolution product using a similar idea. Suppose we have an \emph{ordered} partition $I=I_1\coprod I_2\coprod\cdots \coprod I_k$. 
We consider a stack $\wt{\Hecke}^{\loc,I_1,\cdots,I_k}_{X^I}$ whose $S$-point is a tuple
\[(\ul{x}, \cP^0,\cdots,\cP^k,\alpha_1,\cdots,\alpha_k)\]
where each $\cP^i$ is a $G$-bundle over $D_{\ul{x}}$, where $\alpha_a$ is an isomorphism between $\cP^{a-1}$ and $\cP^a$ over $D_{\ul{x}}\setminus\Gamma_{\ul{x}_a}$, where $\ul{x}_a$ is an $S$-point of $X^{I_a}$ obtained by a restriction of $\ul{x}$. Then we have a correspondence
\[\begin{tikzcd}
	{\wt{\Hecke}^{\loc,I_1,\cdots,I_k}_{X^I}} & {\prod^k_{a=1} \Hecke^{\loc}_{X^{I_a}}} \\
	{\Hecke^{\loc}_{X^I}.}
	\arrow["{\prod p_a}", from=1-1, to=1-2]
	\arrow["m"', from=1-1, to=2-1]
\end{tikzcd}\]
Moreover, it can be shown that the gerbes $\cG^I$ and $\cG^{I_a}$ on $\Hecke^{\loc}_{X^I}$ and $\Hecke^{\loc}_{X^{I_a}}$ are compatible with the correspondence (\cite[Sect. 1.3.5]{ZhaoNew}), so that we have a functor 
\begin{equation}\label{eq: external convolution product}
    \Shv^q(\Hecke^{\loc}_{X^{I_1}})\otimes\cdots\otimes \Shv^q(\Hecke^{\loc}_{X^{I_k}})\rightarrow\Shv^q(\Hecke^{\loc}_{X^{I}})
\end{equation}
using $*$-pullback and $!$-pushforward. This functor is called the \textbf{external convolution product}, and this restricts to Satake categories to yield a functor
\[\Sat^q_{X^{I_1}}\times\cdots\times \Sat^q_{X^{I_k}}\rightarrow \Sat^q_{X^I}.\]

\sssec{Symmetric monoidal structure}
The convolution monoidal structure on $\Sat^q_{X^I}$ admits a natural lift to a symmetric monoidal structure by means of \textbf{fusion product} (\cite[Sect. 1.3.12]{ZhaoNew}), and therefore the sheaf \eqref{eq: sheaf of Satake categories} of Satake categories lifts to a sheaf of $\VectAbfd$-linear symmetric monoidal categories:
\begin{equation*}
  \mathbf{Sat}^q_I: (X_{\et})^{\op}\rightarrow \CAlg(\CatEc),\quad X'\mapsto \Sat^q_{(X')^I}.
\end{equation*}
However, as in the case of the non-metaplectic geometric Satake, to get a symmetric monoidal equivalence we have to alter the commutativity constraint by a sign according to a parity of connected components of $\Gr_{X^I}$ (\cite[Proposition 5.2.6]{Zhu}). Following the notation of \cite[Sect. 1.3.16]{ZhaoNew}, we denote the resulting sheaf of $\VectAbfd$-linear symmetric monoidal categories
\begin{equation}\label{eq: sheaf of plp Sat}
    \plp{\mathbf{Sat}^q_I}: (X_{\et})^{\op}\rightarrow \CAlg(\CatEc),\quad X'\mapsto \plp{\Sat^q_{(X')^I}}.   
\end{equation}

\sssec{Twisted $H$-representations}\label{sssec: twisted H reps}
The other side of the equivalence is concerned with \emph{twisted} $H$\emph{-representations} on local systems. Given a scheme $Y$ of finite type over $\C$ we have the $E$-linear abelian category $\LS^{\Ab,\ft}(Y)$ of $E$-linear finite type local systems concentrated in degree $0$. Given a reductive group $\mathsf{H}$ over $E$, we can consider the constant local system with value $\cO_{\mathsf{H}}$, and denote it by the same notation for convenience. We have $\cO_{\mathsf{H}}\in \Ind(\LS^{\Ab,\ft}(Y))$ as a coalgebra object, hence we can consider its comodules. We consider the full subcategory 
\[\Rep(\mathsf{H})^{\Ab,\loc,\ft}(Y)\subset \cO_{\mathsf{H}}\te{-CoMod}(\Ind(\LS^{\Ab,\ft}(Y)))\]
spanned by the objects whose image under the forgetful functor to $\Ind(\LS^{\Ab,\ft}(Y))$ lies in $\LS^{\Ab,\ft}(Y)$. In other words, an object in $\Rep(\mathsf{H})^{\Ab,\loc,\ft}(Y)$ is the data of
\begin{equation}\label{eq: objects of lisse repn}
    \cF\in \LS^{\Ab,\ft}(Y) \ \te{ and } \ \cF\xrightarrow{\on{coact}_{\cF}} \cF\otimes\cO_{\mathsf{H}}
\end{equation}
where the latter map defines a coaction. This can be thought of as the category of $\mathsf{H}$-representations on finite type local systems. Note that $\Rep(\mathsf{H})^{\Ab,\loc,\ft}(Y)$ has a natural symmetric monoidal structure such that the forgetful map
\[\Rep(\mathsf{H})^{\Ab,\loc,\ft}(Y)\rightarrow \LS^{\Ab,\ft}(Y)\]
is symmetric monoidal. 
\medskip

For any $I\in\fSet$, we have a sheaf 
\begin{equation}\label{eq: sheaf of H lisse repns}
 \mathbf{Rep}(\mathsf{H})^{\Ab,\loc,\ft}_I:
    X_{\et}^{\op}\rightarrow \CAlg(\CatEc),\quad X'\mapsto \Rep(\mathsf{H})^{\Ab,\loc,\ft}(X'^I).
\end{equation}
of $\VectAbfd$-linear symmetric monoidal categories on $X_{\et}$: the small \'etale site of $X$. We observe that this sheaf of symmetric monoidal categories has a natural $Z_{\mathsf{H}}(E)^{\on{tors}}$-linear structure in the sense of Sect. \ref{sssec: twisting construction} and \ref{sssec: twisting algebra construction}. Indeed, \eqref{eq: sheaf of H lisse repns} is a symmetric monoidal subsheaf of
\begin{equation}\label{eq: Rep H I}
    \mathbf{Rep}(\mathsf{H})^{\loc}_I: X^{\op}_{\et}\rightarrow \CAlg(\DGcat),\quad X'\mapsto \Rep(\mathsf{H})\otimes \LS(X'^I),
\end{equation}
because we have an identification
\[\Rep(\mathsf{H})\otimes \LS(X'^I)=\cO_{\mathsf{H}}\te{-CoMod}(\LS(X'^I))\]
and we have $\Ind(\LS^{\Ab,\ft}(X'^I))$ as a symmetric monoidal subcategory of $\LS(X'^I)$, and moreover the sheaf \eqref{eq: Rep H I} admits a symmetric monoidal twist by a $Z_{\mathsf{H}}(E)^{\on{tors}}$-gerbe. (see Sect. \ref{sssec: twisted representations}) Therefore, using the twisting procedure of Sect. \ref{sssec: twisting algebra construction}, given a $Z_{\mathsf{H}}(E)^{\on{tors}}$-gerbe $\cG$ on $X$ we get the twist of $\mathbf{Rep}(\mathsf{H})^{\Ab,\loc,\ft}$ as a sheaf of $\VectAbfd$-linear symmetric monoidal categories
\begin{equation}\label{eq: sheaf of twisted H lisse repns}
 \mathbf{Rep}(\mathsf{H})^{\Ab,\loc,\ft}_{I,\cG}:
    X_{\et}^{\op}\rightarrow \CAlg(\CatEc),\quad X'\mapsto \Rep(\mathsf{H})^{\Ab,\loc,\ft}_{\cG}(X'^I).
\end{equation}

\begin{rmk}
    The twist of \eqref{eq: Rep H I} by a $Z_{\mathsf{H}}(E)^{\on{tors}}$-gerbe can be described in a more elementary way. Given any $a\in Z_{\mathsf{H}}(E)^{\on{tors}}$ we have $\cO_H\xrightarrow{\on{ev}_a}\ul{E}_Y$ in $\Ind(\LS^{\Ab,\ft}(Y^I))$ given by evaluation of functions at $a$, and the composition
\[\cF\xrightarrow{\on{coact}_{\cF}}\cF\otimes \cO_{\mathsf{H}}\xrightarrow{\cF\otimes\on{ev}_a}\cF\]
defines an automorphism for each object \eqref{eq: objects of lisse repn}. It is easy to check that the collection of these automorphisms for every object $(\cF,\on{coact}_{\cF})$ defines an automorphism of the identity functor of $\Rep(\mathsf{H})^{\Ab,\loc,\ft}(Y)$ as a symmetric monoidal functor.
\end{rmk}

\begin{rmk}\label{rmk: twisted Rep H containing twisted Rep H lisse}
    Let 
    \begin{equation*}
       \mathbf{Rep}(\mathsf{H})^{\loc}_{I,\cG}: X^{\op}_{\et}\rightarrow \CAlg(\DGcat),\quad X'\mapsto \Rep(\mathsf{H})^{\loc}_{\cG}(X'^I) 
    \end{equation*}
    be the sheaf obtained as the symmetric monoidal twist of \eqref{eq: Rep H I} by $\cG$.
    The procedure of twisting preserves the subsheaf relation, and therefore the sheaf above has \eqref{eq: sheaf of twisted H lisse repns} as a symmetric monoidal subsheaf. When $I=\{*\}$, this sheaf is simply written as $\mathbf{Rep}(\mathsf{H})^{\loc}_{\cG}$, and this is the category of twisted $\mathsf{H}$-representations, which will play a key role in our paper (see also Sect. \ref{sssec: twisted representations}).
\end{rmk}

\sssec{Spectral side}\label{sssec: B side of Satake}
For each $I\in\fSet$ we apply the above twisting construction \eqref{eq: sheaf of twisted H lisse repns} to
\[\mathsf{H}=H^I,\quad Z=X,\ \te{ and } \ \cG=(\Delta_I)_*\cG_Z,\]
 where $(\Delta_I)_*\cG_Z$ means the inflation of the $Z_H(E)^{\on{tors}}$-gerbe $\cG_Z$ to a $Z_{H^I}(E)^{\on{tors}}$-gerbe using the diagonal map $\Delta_I:Z_H(E)^{\on{tors}}\rightarrow Z_{H^I}(E)^{\on{tors}}$. We denote the resulting sheaf of categories as
\begin{equation}\label{eq: sheaf of twisted lisse repns}
\mathbf{Rep}(H^I)^{\Ab,\loc,\ft}_{\cG^{\boxtimes I}_Z}: (X_{\et})^{\op}\rightarrow \CAlg(\CatEc),\quad X'\mapsto \Rep(H^I)^{\Ab,\loc,\ft}_{\cG_Z^{\boxtimes I}}((X')^I)
\end{equation}
and its global section by 
\begin{equation}\label{eq: Ith product twisted representation}
    \Rep(H^I)^{\Ab,\loc,\ft}_{\cG^{\boxtimes I}_Z}(X^I).
\end{equation}

The use of $\cG^{\boxtimes I}_Z$ in the subscript is justified by the following alternative description of \eqref{eq: Ith product twisted representation}. We can start with a sheaf
\[\mathbf{Rep}(H^I)^{\Ab,\loc,\ft}:(X^I)^{\op}_{\et}\rightarrow \CAlg(\CatEc),\quad Y\mapsto \Rep(H^I)^{\Ab,\loc,\ft}(Y)\]
of $E$-linear symmetric monoidal categories small \'etale site of $X^I$, and twist it by the $Z_{H^I}(E)^{\on{tors}}$-gerbe $\cG^{\boxtimes I}_Z$ on $X^I$ obtained as the $I$-fold external product of $\cG_Z$. Using the map
\[X_{\et}\rightarrow (X^I)_{\et},\quad X'\mapsto (X')^I\]
of sites we can restrict this sheaf to the small \'etale site of $X$, and this agrees with \eqref{eq: sheaf of twisted lisse repns}. That those two agree with each other can be verified easily using Sect. \ref{ssec: restriction to triv}.

\sssec{Sign modification on the spectral side}\label{sssec: spectral sign}
As mentioned in Remark \ref{rmk: twist by sign}, the spectral side also has to go under the sign modification. In general, if we have a symmetric monoidal category $\cC$ and this admits a symmetric monoidal $\Z/2$-graded structure \footnote{See Sect. \ref{sssec: graded monoidal structure} for the precise meaning of this.}, then we can get a new symmetric monoidal category by modifying the commutativity constraint by a sign. Indeed, we have a decomposition $\cC=\cC_{\on{even}}\oplus \cC_{\on{odd}}$, and we alter the commutativity constraint between objects in $\cC_{\on{odd}}$ by a sign. Note that the passage from $\Sat^q_{X^I}$ to $\plp{\Sat}^q_{X^I}$ is an example of this.
\medskip

Consider the general setup of Sect. \ref{sssec: twisted H reps} where we consider the sheaf $\mathbf{Rep}(\sffH)^{\Ab,\loc,\ft}_I$ in \eqref{eq: sheaf of H lisse repns}. Given a 2-torsion element $a\in Z_{\sffH}(E)$, its action on $\sffH$-representations defines an order 2 automorphism of every object in $\Rep(\sffH)^{\Ab,\loc,\ft}(X'^I)$, compatible to every morphism. This gives rise to the decomposition of $\Rep(\sffH)^{\Ab,\loc,\ft}(X'^I)$ into even and odd parts. This is in fact a symmetric monoidal $\Z/2$-grading.
Moreover, the order 2 automorphism is compatible with the symmetric monoidal $Z_{\sffH}(E)^{\on{tors}}$-linear structure on $\Rep(\sffH)^{\Ab,\loc,\ft}(X'^I)$, and therefore the symmetric monoidal $\Z/2$-grading is preserved after twisting by a $Z_{\sffH}(E)^{\on{tors}}$-gerbe. In conclusion, the action of $a$ gives rise to a symmetric monoidal $\Z/2$-grading on the sheaf \eqref{eq: sheaf of twisted H lisse repns}. Now we can modify the commutativity constraint by a sign. We denote it by
\begin{equation}\label{eq: sheaf of twisted H lisse repns sgn}
 \prescript{a}{}{\mathbf{Rep}}(\mathsf{H})^{\Ab,\loc,\ft}_{I,\cG}:
    X_{\et}^{\op}\rightarrow \CAlg(\CatEc),\quad X'\mapsto \prescript{a}{}{\Rep}(\mathsf{H})^{\Ab,\loc,\ft}_{\cG}(X'^I).
\end{equation}

From $q$ we get a group homomorphism $\wh{Z}_H\rightarrow A_{2\on{-tors}}$ and hence $\wh{Z}_H\rightarrow \{\pm 1\}$, see Sect. \ref{sssec: dual data}. This gives rise to a 2-torsion element of $Z_H(E)$, still denoted by $\epsilon$.
Specializing to the situation of metaplectic Satake, for each \eqref{eq: sheaf of twisted lisse repns} we apply the above construction to the 2-torsion element $(\epsilon)_I\in Z_H(E)^I=Z_{H^I}(E)$, and denote the resulting sheaf by
\begin{equation}\label{eq: sheaf of twisted lisse repns sgn}
\elp{\mathbf{Rep}}(H^I)^{\Ab,\loc,\ft}_{\cG^{\boxtimes I}_Z}: (X_{\et})^{\op}\rightarrow \CAlg(\CatEc),\quad X'\mapsto \elp{\Rep}(H^I)^{\Ab,\loc,\ft}_{\cG_Z^{\boxtimes I}}((X')^I).
\end{equation}

\sssec{}
The perverse ULA version of metaplectic factorizable Satake can be stated as follows:

\begin{thm}[\cite{ZhaoNew}, Theorem 2.3.2]\label{thm: factorizable Satake}
    For each $I\in\fSet$ there is a canonical equivalence between sheaves of $\VectAbfd$-linear symmetric monoidal categories \eqref{eq: sheaf of plp Sat} and \eqref{eq: sheaf of twisted lisse repns sgn}.
\end{thm}

\begin{rmk}\label{rmk: action vs linearity}
    In \emph{op. cit.}, both sheaves are thought of as sheaf of symmetric monoidal $E$-linear abelian categories. We chose a slightly different way of presentation because $\VectAbfd$-linearity is more useful than $E$-linearity when it comes to the relationship to derived categories, i.e., the construction of $\LS$-linear factorizable Satake functors. However, it is easy to recover the $E$-linearity from the $\VectAbfd$-action. It is easy to  see that both \eqref{eq: sheaf of plp Sat} and \eqref{eq: sheaf of twisted lisse repns} are valued in abelian categories, and that the $\VectAbfd$-action is exact. In general, if we are given an abelian category $\cA$ with an exact $\VectAbfd$-action, then the passage to $\Ind$-completion gives rise to a presentable action of $\Vect^{\heartsuit}$ on $\Ind(\cA)$. Hence, $\Ind(\cA)$ naturally gets enriched in $\Vect^{\heartsuit}$, and therefore the $E$-linearity on $\cA$.
\end{rmk}

\sssec{}\label{sssec: monoidal formulation of PULA Satake}
For the purpose of this paper, most of the times it will be sufficient to view this equivalence as a monoidal equivalence. Since $\mathbf{Sat}^q$ and $\plp{\mathbf{Sat}^q}$ have the same monoidal structure, this yields us an equivalence
\begin{equation}\label{eq: perverse ULA Satake sheafy}
    \textbf{S}^{\Ab}_{I}: \mathbf{Rep}(H)^{\Ab,\loc,\ft}_{I,\cG_Z^{\boxtimes I}} \rightarrow \textbf{Sat}^q_{I}\quad\te{in }\Shv(X_{\et};\Alg(\CatEc)).
\end{equation}
We denote its section at each $X'\in X_{\et}$ by
\[\cS^{\Ab}_{X',I}: \Rep(H)^{\Ab,\loc,\ft}_{\cG_Z^{\boxtimes I}}((X')^I) \rightarrow \Sat^q_{(X')^I} \quad \te{in }\Alg(\CatEc).\]

This equivalence enjoys a compatibility pertaining to morphisms in $\fSet$. To state this, we construct a certain monoidal functor on each side, given a map $I\rightarrow J$ in $\fSet$.

\sssec{Automorphic side}\label{sssec: A side functoriality} Given $I\rightarrow J$ in $\fSet$, a map
\begin{equation}\label{eq: fSet map A side}
    \mathbf{Sat}^q_{I}\rightarrow\mathbf{Sat}^q_{J}\quad \te{in }\Shv(X_{\et};\Alg(\CatEc))
\end{equation}
 was constructed in \cite[Remark 1.3.11]{ZhaoNew}. In fact, \emph{loc. cit.} constructs a monoidal functor
 \begin{equation}\label{eq: fSet map A side on Shv}
     \Shv^q(\Hecke^{\loc}_{X^I})\rightarrow \Shv^q(\Hecke^{\loc}_{X^J}),
 \end{equation}
and then shows that it restricts to the perverse ULA categories and hence we get \eqref{eq: fSet map A side}. The monoidal functor \eqref{eq: fSet map A side on Shv} is obtained by the $*$-pullback and $!$-pushforward along the correspondence
\[\begin{tikzcd}
	{\Hecke^{\loc}_{X^{I\rightarrow J}}} & {\Hecke^{\loc}_{X^I}} \\
	{\Hecke^{\loc}_{X^J},}
	\arrow[from=1-1, to=1-2]
	\arrow[from=1-1, to=2-1]
\end{tikzcd}\]
where an $S$-point of $\Hecke^{\loc}_{X^{I\rightarrow J}}$ is a tuple $(\ul{x}, \cP^0,\cP^1,\alpha)$ where $\ul{x}$ is an $S$-point of $X^J$ whose restriction to $X^I$ is denoted $\ul{x}_I$, and $\cP^i$ are $G$-bundles over $D_{\ul{x}}$, and $\alpha$ is an isomorphism between $\cP^0$ and $\cP^1$ over $D_{\ul{x}}\setminus \Gamma_{\ul{x}_I}$.
\medskip

An alternative way of describing \eqref{eq: fSet map A side on Shv} makes use of a unique factorization of $I\rightarrow J$ into a surjection $I\surj J'$ and injection $J'\subseteq J$. In the first case, \eqref{eq: fSet map A side on Shv} is a restriction along the diagonal map $\Hecke^{\loc}_{X^{J'}}\rightarrow \Hecke^{\loc}_{X^I}$, and in the second case, \eqref{eq: fSet map A side on Shv} is obtained by inserting the monoidal unit of $\Shv^q(\Hecke^{\loc}_{X^{J\setminus J'}})$ into the external convolution product functor
\[\Shv^q(\Hecke^{\loc}_{X^{J'}})\otimes \Shv^q(\Hecke^{\loc}_{X^{J\setminus J'}})\rightarrow \Shv^q(\Hecke^{\loc}_{X^{J}}).\]

\begin{rmk}
    In fact, the monoidal functor \eqref{eq: fSet map A side} has a natural symmetric monoidal lift, and by altering commutativity constraints on both sides we get a symmetric monoidal functor
    \begin{equation}\label{eq: functoriality plpSat}
        \plp{\mathbf{Sat}^q_I}\rightarrow \plp{\mathbf{Sat}^q_J} \quad\te{in }\Shv(X_{\et};\CAlg(\CatEc)).
    \end{equation}
\end{rmk}

\sssec{Spectral side}
On the spectral side, we can also construct 
\begin{equation}\label{eq: fSet map B side}
    \mathbf{Rep}(H^I)^{\Ab,\loc,\ft}_{\cG^{\boxtimes I}_Z}\rightarrow \mathbf{Rep}(H^J)^{\Ab,\loc,\ft}_{\cG^{\boxtimes J}_Z}\quad\te{in }\Shv(X_{et};\CAlg(\CatEc)).
\end{equation}
This is the composition
\begin{equation}\label{eq: sequence sheaf of twisted repns}
    \mathbf{Rep}(H^I)^{\Ab,\loc,\ft}_{I,(\Delta_I)_*\cG_Z}\rightarrow \mathbf{Rep}(H^I)^{\Ab,\loc,\ft}_{J,(\Delta_I)_*\cG_Z}\rightarrow \mathbf{Rep}(H^J)^{\Ab,\loc,\ft}_{J,(\Delta_J)_*\cG_Z}
\end{equation}
of maps between sheaves of symmetric monoidal categories. The sequence \eqref{eq: sequence sheaf of twisted repns} arises as a twist of the sequence
\begin{equation}\label{eq: sequence sheaf of repns}
    \mathbf{Rep}(H^I)^{\Ab,\loc,\ft}_{I}\rightarrow \mathbf{Rep}(H^I)^{\Ab,\loc,\ft}_{J}\rightarrow \mathbf{Rep}(H^J)^{\Ab,\loc,\ft}_{J},
\end{equation}
where the first arrow is the restriction of $H^I$-representations along the diagonal $X^J\rightarrow X^I$, and the second arrow restricts $H^I$-representations to $H^J$-representations along the diagonal $H^J\rightarrow H^I$. It is easy to check that the maps of \eqref{eq: sequence sheaf of repns} intertwine (symmetric monoidal) twisting data applied to each term of the sequence, and therefore by the functoriality of twisting construction we get \eqref{eq: sequence sheaf of twisted repns} after twist. This completes the construction of \eqref{eq: fSet map B side}.

\begin{rmk}
    Using Remark \ref{rmk: twisted Rep H containing twisted Rep H lisse}, the functor \eqref{eq: fSet map B side} admits an alternative description. The source and the target can be realized as a symmetric monoidal subsheaf of
    \[\mathbf{Rep}(H^I)^{\loc}_{I,(\Delta_I)_*\cG_Z}\quad\te{and}\quad \mathbf{Rep}(H^J)^{\loc}_{J,(\Delta_J)_*\cG_Z},\]
    which are in turn canonically identified with
    \[\left(\mathbf{Rep}(H)^{\loc}_{\cG_Z}\right)^{\otimes I}\quad \te{and}\quad \left(\mathbf{Rep}(H)^{\loc}_{\cG_Z}\right)^{\otimes J}.\]
    The $\E_{\infty}$-structure on $\mathbf{Rep}(H)^{\loc}_{\cG_Z}$ defines an $\E_{\infty}$-map
    \[\left(\mathbf{Rep}(H)^{\loc}_{\cG_Z}\right)^{\otimes I}\rightarrow \left(\mathbf{Rep}(H)^{\loc}_{\cG_Z}\right)^{\otimes J},\]
    and this restricts to the composition \eqref{eq: sequence sheaf of twisted repns}.
\end{rmk}

\begin{rmk}
    The functor \eqref{eq: fSet map B side} has an obvious analog after sign modification:
    \begin{equation}\label{eq: fSet map B side sgn}
    \elp{\mathbf{Rep}}(H^I)^{\Ab,\loc,\ft}_{\cG^{\boxtimes I}_Z}\rightarrow \elp{\mathbf{Rep}}(H^J)^{\Ab,\loc,\ft}_{\cG^{\boxtimes J}_Z}\quad\te{in }\Shv(X_{et};\CAlg(\CatEc)).
\end{equation}
\end{rmk}

\sssec{}
Having constructed monoidal functors on both sides, we can make precise what compatibility the equivalence of Theorem \ref{thm: factorizable Satake} has under morphisms in $\fSet$. Given $I\rightarrow J$ in $\fSet$, we have a canonical commutative diagram of monoidal functors
\begin{equation}\label{eq: fSet compatibility diagram}
    \begin{tikzcd}
	{\textbf{Rep}(H^I)^{\Ab,\loc,\ft}_{\cG^{\boxtimes I}_Z}}& {\mathbf{Sat}^q_{I}} \\
	  {\mathbf{Rep}(H^J)^{\Ab,\loc,\ft}_{\cG^{\boxtimes J}_Z}} & {\mathbf{Sat}^q_{J}}
	\arrow["\mathbf{S}^{\Ab}_{I}",from=1-1, to=1-2]
	\arrow["{\eqref{eq: fSet map B side}}"', from=1-1, to=2-1]
	\arrow["{\eqref{eq: fSet map A side}}", from=1-2, to=2-2]
	\arrow["\mathbf{S}^{\Ab}_{J}"',from=2-1, to=2-2]
\end{tikzcd}
\end{equation}
and the commutativity data of this diagram is compatible with the composition $I\rightarrow J\rightarrow K$.

\begin{rmk}\label{rmk: sym mon version functoriality}
    We have an obvious symmetric monoidal version of this commutative diagram replacing $\mathbf{Sat}^q_{?}$ by $\plp{\mathbf{Sat}^q_?}$ and replacing \eqref{eq: fSet map A side} by \eqref{eq: functoriality plpSat} and \eqref{eq: fSet map B side} by \eqref{eq: fSet map B side sgn}.
\end{rmk}

\sssec{Grading}\label{sssec: grading overview}
A crucial feature of the equivalence of Theorem \ref{thm: factorizable Satake} - which plays a key role in its proof - is the canonical lift to an equivalence between graded monoidal sheaves \footnote{We refer readers to Sect. \ref{ssec: twisting from grading} for the precise meaning of this.}. Let $\wh{Z}_H$ be the group of characters of $Z_H$. Each of the sheaves \eqref{eq: sheaf of Satake categories} and \eqref{eq: sheaf of twisted lisse repns} of monoidal categories has a canonical lift as a sheaf of $(\wh{Z}_H)^I$-graded monoidal categories.
On the spectral side \eqref{eq: sheaf of twisted lisse repns} this graded monoidal structure comes from the $(\wh{Z}_H)^I$-linearity of $\mathbf{Rep}(H^I)^{\Ab,\loc,\ft}$ discussed previously. The construction of the monoidal $(\wh{Z}_H)^I$-grading on the automorphic side \eqref{eq: sheaf of Satake categories} is more subtle, and involves the study of $(L^+G_{\on{ad}})_X$-equivariance on $\Gr_{X}$. We refer to \cite[Sect. 5]{ZhaoNew} for details, although we review it in Sect. \ref{ssec: A side grading}.
\medskip

This monoidal grading is compatible with maps in $\fSet$. A map $I\rightarrow J$ in $\fSet$ induces a map $(\wh{Z}_H)^I\rightarrow (\wh{Z}_H)^J$. Call this map $\phi$ for convenience. Both of the monoidal functors \eqref{eq: fSet map A side} and \eqref{eq: fSet map B side} have a monoidal $(\wh{Z}_H)^I$-graded structure on the source and a monoidal $(\wh{Z}_H)^J$-graded structure on the target. Given $a\in (\wh{Z}_H)^I$, both monoidal functors send an 
$a$-graded piece of the source to the $\phi(a)$-graded piece of the target. 

\begin{rmk}
    The graded monoidal lifts and compatibilities mentioned above naturally enhance to that of graded symmetric monoidality. 
\end{rmk}

\begin{rmk}\label{rmk: induced grading}
    By a chain of group homomorphisms \footnote{$\epsilon$ is the sign component of $\nu$ introduced in Sect. \ref{sssec: dual data}.}
    \[(\wh{Z}_H)^I\xrightarrow{\on{sum}}\wh{Z}_H\xrightarrow{\epsilon}\{\pm 1\},\]
    we have a symmetric monoidal $\wh{Z}_H$-graded (and $\Z/2$-graded) structure on $\Sat^q_{X^I}$.
\end{rmk}

\sssec{}\label{sssec: opposite twist}
The importance of this symmetric monoidal grading in the proof of Theorem \ref{thm: factorizable Satake} is that we can prove the equivalence after untwisting both sides. Tannakian reconstruction theorem plays a key role in the proof of geometric Satake and serves as the definition of dual group, but the category $\plp{\Sat^q_{X^I}}$ does \emph{not} admit a fiber functor to $\LS^{\Ab,\ft}(X^I)$. Instead, the \emph{opposite} twist of it, or more precisely $\left(\plp{\Sat^q_{X^I}}\right)_{-(\cG_Z+\epsilon)}$, admits a fiber functor to $\LS^{\Ab,\ft}(X^I)$ whose Tannakian group is $H^I$. 
\medskip

Let us elaborate on what the twist of $\plp{\Sat}^q_{X^I}$ by $-(\cG_Z+\epsilon)$ means. Here $-(\cG_Z+\epsilon)$ is a $\E_{\infty}$-map $\wh{Z}_H\rightarrow \on{Ge}_A(X)$, where we start from the $\Z$-linear map $-\cG_Z$ and tweak the commutativity constraint according to $\epsilon$. Using the method of twisting by grading (Sect. \ref{sssec: twist by grading}), the data $-(\cG_Z+\epsilon)$ twists a symmetric monoidal $\wh{Z}_H$-graded sheaf $\Sat^q_{X^I}$.
\medskip

We have a canonical \emph{monoidal} equivalence
\begin{equation}\label{eq: twist and sign}
    (\plp{\Sat}^q_{X^I})_{-(\cG_Z+\epsilon)}\simeq (\plp{\Sat}^q_{X^I})_{-\cG_Z},
\end{equation}
and the difference between those is precisely the sign modification. The data of $\epsilon$ determines a symmetric monoidal $\Z/2$-grading on $\plp{\Sat}^q_{X^I}$ (see Remark \ref{rmk: induced grading}), and this structure is preserved under the twisting by $-\cG_Z$, inducing the $\Z/2$-grading on the right-hand side of \eqref{eq: twist and sign}. Now, the left-hand side of \eqref{eq: twist and sign} is obtained from the right-hand side by the sign modification.

\begin{rmk}
    There were two sign modifications made in the formation of the left-hand side of \eqref{eq: twist and sign}: the passage from $\Sat^q$ to $\plp{\Sat}^q$ and the twisting by $\epsilon$. Each of them go through a symmetric monoidal $\Z/2$-graded structure, but the graded structures used in those two twisting are \emph{different} from each other.
\end{rmk}

\begin{rmk}
    We earlier described the construction of $\elp{\textbf{Rep}}(H^I)^{\Ab,\loc,\ft}_{(\cG_Z)^{\boxtimes I}}$ as a two step procedure: start from $\textbf{Rep}(H^I)^{\Ab,\loc,\ft}$, twist by $\cG_Z$, and modify commutativity constraint using $\epsilon$. This in fact can be done in one shot, using the symmetric monoidal $\wh{Z}_H$-grading on $\textbf{Rep}(H^I)^{\Ab,\loc,\ft}$ and twisting by $\cG_Z+\epsilon$.
\end{rmk}

\sssec{}
The commutative diagram \eqref{eq: fSet compatibility  diagram} is obtained as a $\cG_Z$-twist of the commutative diagram of monoidal functors
\begin{equation}\label{eq: fSet compatibility diagram untwist}
    \begin{tikzcd}
	{\left(\Sat^q_{X^I}\right)_{-\cG_Z}} & {\Rep(H^I)^{\Ab,\loc,\ft}(X^I)} \\
	{\left(\Sat^q_{X^J}\right)_{-\cG_Z}} & {\Rep(H^J)^{\Ab,\loc,\ft}(X^J).}
	\arrow["\sim",from=1-1, to=1-2]
	\arrow[from=1-1, to=2-1]
	\arrow[from=1-2, to=2-2]
	\arrow["\sim",from=2-1, to=2-2]
\end{tikzcd}
\end{equation}
The symmetric monoidal version of \eqref{eq: fSet compatibility  diagram} - mentioned in Remark \ref{rmk: sym mon version functoriality} - is obtained as a $\cG_Z+\epsilon$-twist of the commutative diagram of symmetric monoidal functors
\begin{equation}\label{eq: fSet compatibility diagram untwist symmon}
    \begin{tikzcd}
	{\left(\plp{\Sat}^q_{X^I}\right)_{-(\cG_Z+\epsilon)}} & {\Rep(H^I)^{\Ab,\loc,\ft}(X^I)} \\
	{\left(\plp{\Sat}^q_{X^J}\right)_{-(\cG_Z+\epsilon)}} & {\Rep(H^J)^{\Ab,\loc,\ft}(X^J).}
	\arrow["\sim",from=1-1, to=1-2]
	\arrow[from=1-1, to=2-1]
	\arrow[from=1-2, to=2-2]
	\arrow["\sim",from=2-1, to=2-2]
\end{tikzcd}
\end{equation}
     The diagram \eqref{eq: fSet compatibility diagram untwist symmon} considered as the diagram of monoidal functors equals the diagram \eqref{eq: fSet compatibility diagram untwist}.

\subsection{A $\mathbf{LS}$-linear version of factorizable Satake: the singleton case}\label{ssec: LS linear fact Satake}
Now we start the construction of a collection of canonical functors 
\begin{equation*}
    \cS_{X,I}:\left(\Rep(H)^{\loc}_{\cG_Z}(X)\right)^{\otimes I}\rightarrow \Shv^q(\Hecke^{\loc}_{X^I})\quad \te{in } \LS(X^I)\lMod
\end{equation*}
for $I\in\fSet$.
The case $I=\{\ast\}$, i.e, a $\LS(X)$-linear functor
\begin{equation}\label{eq: fact Satake I singleton}
    \cS_{X,\{*\}}:\Rep(H)^{\loc}_{\cG_Z}(X)\rightarrow \Shv^q(\Hecke^{\loc}_{X})
\end{equation}
is the most interesting part, and the rest is done by external convolution product (see Sect. \ref{sssec: convolution product}). In the case of a singleton, we will do a little bit more: construct \eqref{eq: fact Satake I singleton} as a $\LS(X)$-linear monoidal functor \footnote{The monoidality we establish here is not strictly needed for future use, because it would follow from the $\Conv$-algebra formalism to be introduced in Sect. \ref{sec: Satake as a nat transf}.}.

\begin{rmk}\label{rmk: characterization of LS linear Satake}
    The construction of the functors $\{\cS_{X,I}\}$ involves a bit of technicalities. However, we emphasize that those functors are completely controlled by the perverse ULA version of factorizable Satake. This point will be elaborated on later in Sect. \ref{ssec: LS linear general}.
\end{rmk}

\sssec{}
Having explained the ideas and subtlety around the construction of $\cS_{X,\{*\}}$ in the introduction to this section (see Sect. \ref{sssec: technical aspects LS linear} and \ref{sssec: subtleties LS linear}), we proceed to the construction. We need:

\begin{thm}\label{thm: A side grading}
    We have a canonical $\wh{Z}_H$-graded $\Shv^{\indconstr}(X)$-linear monoidal lift of the category $\Shv^{\indconstr,q}(\Hecke^{\loc}_X)$.
\end{thm}

The construction of such a lift for constructible sheaves was first considered in \cite[Remark 5.3.7]{ZhaoNew}. We provide the proof of this assertion in Sect. \ref{ssec: A side grading} using the ideas of \emph{loc. cit.}.

\begin{note}
We will first give a slightly imprecise construction of $\cS_{X,\{*\}}$ ignoring the subtlety with the definition of $\Shv^q(\Hecke^{\loc}_X)$ mentioned in Sect. \ref{sssec: subtlety Hecke}, and fix it later. This presentation is intended to keep readers from becoming bogged down in minor details. It suffices for readers to accept the initial construction. The logical gap in our first construction lies in Sect. \ref{sssec: problematic} and we provide a fix in Sect. \ref{sssec: ind nature precise}.    
\end{note}

\sssec{}
We consider the inclusion functor
\[\Sat^q_X\hookrightarrow \Shv^{\indconstr,q}(\Hecke^{\loc}_X).\]
By construction, the $\wh{Z}_H$-graded monoidal structure in Theorem \ref{thm: A side grading} restricts to the $\wh{Z}_H$-graded monoidal structure on $\Sat^q_X$ discussed in Sect. \ref{sssec: grading overview}. Therefore, we can untwist both sides to get a monoidal functor
\[\left(\Sat^q_X\right)_{-\cG_Z}\rightarrow \left(\Shv^{\indconstr,q}(\Hecke^{\loc}_X)\right)_{-\cG_Z}.\]
Thus, we have a composition
\[\Rep^{\Ab,\loc,\ft}(H)(X)\xrightarrow{\simeq} \left(\Sat^q_X\right)_{-\cG_Z}\rightarrow \left(\Shv^{\indconstr,q}(\Hecke^{\loc}_X)\right)_{-\cG_Z}\]
of monoidal functors. Let $\Rep(H)^{\heartsuit,\on{f.d.}}$ be the abelian category of finite dimensional $H$-representations. We have a monoidal functor
\[\Rep(H)^{\heartsuit,\on{f.d.}}\rightarrow \Rep^{\Ab,\loc,\ft}(H)(X)\]
by pulling back local systems along $X\rightarrow \pt$, and hence we get a monoidal functor
\begin{equation}\label{eq: functor from Ab ft Rep H}
    \Rep(H)^{\heartsuit,\on{f.d.}}\rightarrow \left(\Shv^{\indconstr,q}(\Hecke^{\loc}_X)\right)_{-\cG_Z}.
\end{equation}

\sssec{}\label{sssec: problematic}
We have a perverse $t$-structure on the target, and \eqref{eq: functor from Ab ft Rep H} lands in the perverse heart of the target, i.e., we have
\begin{equation}\label{eq: functor from Ab ft btw Abs}
    \Rep(H)^{\heartsuit,\on{f.d.}}\rightarrow \left(\Shv^{\on{ind-constr},q}(\Hecke^{\loc}_{X})\right)_{-\cG_Z}^{\heartsuit}.
\end{equation}

    The perverse $t$-structure on $\left(\Shv^{\on{ind-constr},q}(\Hecke^{\loc}_{X'})\right)_{-\cG_Z}$ is left-complete, compatible with filtered colimits, and its non-positive part has filtered colimits. Its abelian heart is $\Ind$-complete. Hence, we can $\Ind$-complete \eqref{eq: functor from Ab ft btw Abs} to
\[\Rep(H)^{\heartsuit}\rightarrow \left(\Shv^{\on{ind-constr},q}(\Hecke^{\loc}_{X})\right)_{-\cG_Z}^{\heartsuit}.\]
Using the left-completeness of the $t$-structure we can use \cite[Theorem 1.3.3.2]{HA} to get
\[\Rep(H)^b\rightarrow \left(\Shv^{\on{ind-constr},q}(\Hecke^{\loc}_{X})\right)_{-\cG_Z}\]
between presheaves, where $\Rep(H)^b$ is the bounded derived category. Pulling it back along the inclusion $\Rep(H)^c\hookrightarrow \Rep(H)^b$ and taking $\Ind$-completion we get
\begin{equation}\label{eq: desired map from Rep H}
    \Rep(H)\rightarrow \left(\Shv^{\on{ind-constr},q}(\Hecke^{\loc}_{X})\right)_{-\cG_Z}\quad\te{in }\DGCatMon.
\end{equation}

\sssec{}
Since the $\wh{Z}_H$-grading on $\Shv^{\on{ind-constr},q}(\Hecke^{\loc}_{X})$ is a grading as a $\QLisse(X)$-linear monoidal category, its twist by $-\cG_Z$ also has a canonical $\QLisse(X)$-linear monoidal structure. Hence, from \eqref{eq: desired map from Rep H} we get 
\begin{equation*}
    \Rep(H)\otimes \QLisse(X)\rightarrow \left(\Shv^{\on{ind-constr},q}(\Hecke^{\loc}_{X})\right)_{-\cG_Z}
\end{equation*}
in $\Alg(\QLisse(X)\lMod)$.

\sssec{}\label{sssec: list of sheaves}
In fact, we can do the same construction functorially for $X'\in X_{\et}$ to get a map between sheaves. Let us expand on this. Consider the sheaves
\begin{gather*}
    \mathbf{QLisse}: X^{\op}_{\et}\rightarrow \DGCatSymMon,\quad X'\mapsto \QLisse(X'), \\ \mathbf{Rep}(H)^{\QLisse}: X^{\op}_{\et}\rightarrow \DGCatSymMon,\quad X'\mapsto \Rep(H)\otimes\QLisse(X'), \\
\mathbf{Shv}^{\indconstr,q}(\Hecke^{\loc}):X^{\op}_{\et}\rightarrow \DGCatMon,\quad X'\mapsto \Shv^{\indconstr,q}(\Hecke^{\loc}_{X'}).
\end{gather*}
The latter two naturally live in the category of $\mathbf{QLisse}$-linear algebra objects in $\Shv(X_{\et}, \DGcat)$, and we have obtained a map
\[\mathbf{Rep}(H)^{\QLisse}\rightarrow \mathbf{Shv}^{\indconstr,q}(\Hecke^{\loc})_{-\cG_Z}\]
in that category. For the rest of the argument, let us omit $\DGcat$ to unburden the notations regarding sheaf of categories. We can twist $\cG_Z$ back to get
\begin{equation}\label{eq: functor 1}
    \mathbf{Rep}(H)^{\QLisse}_{\cG_Z}\rightarrow \mathbf{Shv}^{\indconstr,q}(\Hecke^{\loc})\quad\te{in }\Alg(\mathbf{QLisse}\lMod(\Shv(X_{\et}))).
\end{equation}

On the other hand, we have
\[\mathbf{Shv}^q(\Hecke^{\loc}): X^{\op}_{\et}\rightarrow \DGCatMon,\quad X'\mapsto \Shv^q(\Hecke^{\loc}_{X'})\]
as an object of $\Alg(\mathbf{LS}\lMod(\Shv(X_{\et})))$, and using the canonical symmetric monoidal map $\mathbf{QLisse}\rightarrow \mathbf{LS}$,
we can also consider it as an object of $\Alg(\mathbf{QLisse}\lMod(\Shv(X_{\et})))$. Then \eqref{eq: indconstr to Betti} gets upgraded to a map
\begin{equation}\label{eq: functor 2}
    \mathbf{Shv}^{\indconstr,q}(\Hecke^{\loc})\rightarrow \mathbf{Shv}^q(\Hecke^{\loc})\quad \te{in }\Alg(\mathbf{QLisse}\lMod(\Shv(X_{\et}))).
\end{equation}

Combining the two functors \eqref{eq: functor 1} and \eqref{eq: functor 2} above we get
\begin{equation}\label{eq: Rep QLisse to Shv all}
\mathbf{S}^{\QLisse}_{\{*\}}:\mathbf{Rep}(H)^{\QLisse}_{\cG_Z}\rightarrow \mathbf{Shv}^{q}(\Hecke^{\loc})\quad \te{in } \Alg(\mathbf{QLisse}\lMod(\Shv(X_{\et}))).
\end{equation}
 Since
\[\mathbf{Rep}(H)^{\loc}_{\cG_Z}=\mathbf{Rep}(H)^{\QLisse}_{\cG_Z}\otimes_{\mathbf{QLisse}}\mathbf{LS} \quad\te{in } \CAlg(\mathbf{QLisse}\lMod(\Shv(X_{\et}))),\]
 using the $\mathbf{LS}$-linear monoidality on the target of \eqref{eq: Rep QLisse to Shv all} we get
\begin{equation}\label{eq: fact Satake singleton sheafy}
\mathbf{S}_{\{*\}}:\mathbf{Rep}(H)^{\loc}_{\cG_Z}\rightarrow \mathbf{Shv}^{q}(\Hecke^{\loc})\quad\te{in } \Alg(\mathbf{LS}\lMod(\Shv(X_{\et}))),
\end{equation}
 hence we finish the construction of \eqref{eq: fact Satake I singleton} by taking global sections.

\sssec{Filling the gap in Sect. \ref{sssec: problematic}}\label{sssec: ind nature precise}
In fact, we have \emph{not} verified the assertion (made in Sect. \ref{sssec: problematic}) pertaining to the perverse $t$-structure on the target of \eqref{eq: functor from Ab ft Rep H}. By Proposition \ref{prop: t str on untwisted indconstr Hecke loc}, we have a perverse $t$-structure on each finite piece 
\[\left(\Shv^{\indconstr,q}(\Hecke^{\loc,\le\lambda}_X)\right)_{-\cG_Z}\]
for $\lambda\in\Lambda^+$ and they satisfy the properties listed in Sect. \ref{sssec: problematic}. However, we found it difficult to analyze the $t$-structure of the colimit of such. \footnote{In fact, we do not even know if there is a canonical $t$-structure on the colimit.}
\medskip

To circumvent this difficulty, we perform the construction of Sect. \ref{sssec: problematic} on each finite piece. Let $\Rep(H)^{\heartsuit,\le\lambda}$ be the subcategory of $\Rep(H)^{\heartsuit}$ generated by irreducible representations of highest weights $\le \lambda$ under arbitrary colimits. This is the $\Ind$-completion of the subcategory $\Rep(H)^{\heartsuit,\on{f.d.},\le\lambda}$ of finite-dimensional representations. Then the functor \eqref{eq: functor from Ab ft Rep H} factors through
\[\Rep(H)^{\heartsuit,\on{f.d}.\le\lambda}\rightarrow \left(\Shv^{\on{ind-constr},q}(\Hecke^{\loc,\le\lambda}_{X'})\right)^{\heartsuit}_{-\cG_Z},\]
and by the procedure of Sect. \ref{sssec: problematic} we obtain
\begin{equation}\label{eq: derived functor on the finite piece}
    D(\Rep(H)^{\heartsuit,\le\lambda})\rightarrow \left(\Shv^{\on{ind-constr},q}(\Hecke^{\loc,\le\lambda}_{X'})\right)_{-\cG_Z} \quad\te{in }\DGcat.
\end{equation}

We have
\[\Rep(H)=\colim_{\lambda\in\Lambda^+} D(\Rep(H)^{\heartsuit,\le\lambda})\quad\te{in }\DGcat,\]
and hence the desired functor \eqref{eq: desired map from Rep H} is obtained by taking colimits on both sides of \eqref{eq: derived functor on the finite piece}.

\subsection{Grading on the automorphic side}\label{ssec: A side grading}
In this subsection we prove Theorem \ref{thm: A side grading}. We follow the approach of \cite{ZhaoNew} using the $(L^+G_{\on{ad}})_X$-equivariance of the pullback of the $A$-gerbe $\cG^{\{1\}}$ on $\Hecke^{\loc}_X$ to $\Gr_X$. Although they only needed the monoidal $\wh{Z}_H$-grading for the abelian category $\Sat^q_X$, it was pointed out in Remark 5.3.7 of \emph{loc. cit.} that the construction can be extended to the derived category. What we prove here is slightly different in the sense that we use ind-constructible sheaves while \emph{loc. cit.} uses constructible sheaves.\footnote{Over a stack, the category of ind-constructible sheaves is in general different from the $\Ind$-completion of the category of constructible sheaves.}  However, the key arguments were already provided there.

\sssec{$(L^+G_{\on{ad}})_X$-equivariance} 
We review the $(L^+G_{\on{ad}})_X$-equivariance of $\cG^{\{1\}}$ studied in Sect. 5.1 of \emph{op. cit.}.
\medskip

For $\theta\in\pi_1(G)$ let $\Hecke^{\loc,\theta}_X$ (resp. $\Gr_{X}$) be the corresponding connected component of $\Hecke^{\loc}_X$ (resp. $\Gr_X$), and consider the restriction $\cG^{\theta}$ of $\cG^{\{1\}}$ to $\Hecke^{\loc,\theta}_X$. We have a canonical map
$\wh{Z}_H\rightarrow \pi_1(G)$, and a lift $\theta^{\sharp}\in \wh{Z}_H$ of $\theta\in\pi_1(G)$ along this map determines a canonical $(L^+G_{\on{ad}})_X$-equivariant structure on the pullback of $\cG^{\{1\}}$ to $\Gr^{\theta}_X$. We denote the resulting $A$-gerbe on $(L^+G)_{\on{ad}}\backslash \Gr^{\theta}_X$ by $\cG^{\theta^{\sharp}}$.

\sssec{Construction of graded pieces}
For each $\theta^{\sharp}\in\wh{Z}_H$ we consider the functor
\[\Sigma'_{\theta^{\sharp}}: \Shv^{\indconstr,\cG^{\theta^{\sharp}}}((L^+G_{\on{ad}})_X\backslash\Gr^{\theta}_X)\rightarrow \Shv^{\indconstr,q}(\Hecke^{\loc}_X)\]
given as the composition of the $*$-pullback along
\[(L^+G)_X\backslash \Gr^{\theta}_X\rightarrow (L^+G_{\on{ad}})_X\backslash \Gr^{\theta}_X\]
and the $!$-pushforward along
\[(L^+G)_X\backslash \Gr^{\theta}_X = \Hecke^{\loc,\theta}_X\hookrightarrow \Hecke^{\loc}_X.\]

Consider its essential image $\on{Im}(\Sigma'_{\theta^{\sharp}})$. It consists of $\cG^{\theta^{\sharp}}$-twisted $(L^+G_{\on{ad}})_X$-equivariant sheaves supported on $\Gr^{\theta}_X$. Consider the full subcategory
\[\on{Im}(\Sigma'_{\theta^{\sharp}})^{\on{fin.ext}}\subseteq \Shv^{\indconstr,q}(\Hecke^{\loc}_X)\]
spanned by objects that can be obtained as a finite extension of objects in $\on{Im}(\Sigma'_{\theta^{\sharp}})$. Consider the full subcategory 
\[\on{Im}(\Sigma'_{\theta^{\sharp}})^{\on{fin.ext},c}\subseteq \on{Im}(\Sigma'_{\theta^{\sharp}})^{\on{fin.ext}}\]
spanned by objects that are compact in $\Shv^{\indconstr,q}(\Hecke^{\loc}_X)$.
Taking $\Ind$-completion we get a fully-faithful functor
\[\Sigma_{\theta^{\sharp}}:\Ind\left(\on{Im}(\Sigma'_{\theta^{\sharp}})^{\on{fin.ext},c}\right) \rightarrow \Shv^{\indconstr,q}(\Hecke^{\loc}_X). \]
Combining these functors for all $\theta^{\sharp}\in\wh{Z}_H$ we obtain a canonical functor
\[\Sigma: \bigoplus_{\theta^{\sharp}\in\wh{Z}_H}\Ind\left(\on{Im}(\Sigma'_{\theta^{\sharp}})^{\on{fin.ext},c}\right)\rightarrow \Shv^{\indconstr,q}(\Hecke^{\loc}_X).\]

In the rest of this subsection, we prove the following sharper version of Theorem \ref{thm: A side grading}:
\begin{thm}\label{thm: A side grading sharper}
    The functor $\Sigma$ is an equivalence in $\DGcat$, and this realizes a $\Shv^{\indconstr}(X)$-linear monoidal $\wh{Z}_H$-graded structure of $\Shv^{\indconstr,q}(\Hecke^{\loc}_X)$.
\end{thm}

\sssec{Orthogonality}
The proof relies on the following key observation:
\begin{prop}[\cite{ZhaoNew}, Proposition 5.3.3]\label{prop: orthogonality}
    The categories $\{\on{Im}(\Sigma'_{\theta^{\sharp}})\}$ for different $\theta^{\sharp}$'s are mutually orthogonal in the category $\Shv^{\indconstr,q}(\Hecke^{\loc}_X)$.
\end{prop}

\sssec{}
We recall the following generalities about compactly generated DG categories. 

\begin{lemma}[\cite{GR}, Corollary 7.2.7]\label{lem: colim cpt gen}
    Suppose we are given a diagram
    \[I\rightarrow\DGcat,\quad i\mapsto \bC_i\]
    where each $\bC_i$ is compactly generated and each of the functors $\bC_i\rightarrow\bC_j$ preserves compactness. Then each insertion map $\bC_i\xrightarrow{\on{ins}_i}\colim_i \bC_i$
    preserves compactness, and the category $\colim_i \bC_i$ is compactly generated by the essential images of the functors
    \[\bC^c_i\hookrightarrow \bC_i\xrightarrow{\on{ins}_i}\colim_i \bC_i.\]
\end{lemma}

\begin{lemma}\label{lem: cpt gen full faithful}
    Suppose we have a map $F:\bC\rightarrow\bD$ in $\DGcat$ where $\bC$ is compactly generated and $F$ preserves compacts. If the collection $\{\bc_{\alpha}\}$ compactly generates $\bC$, then $F$ is fully faithful if and only if its restriction along $\bC^{\circ}\hookrightarrow\bC$ is fully faithful, where $\bC^{\circ}$ is the full subcategory of $\bC$ spanned by $\{\bc_{\alpha}\}$.
\end{lemma}

\sssec{}\label{sssec: setup orthogonal decomposition}
Using the lemmas above, we can prove the following general criterion for the orthogonal decomposition of a compactly generated DG category. 
Suppose that we have $\bC\in\DGcat$ which is compactly generated, and that $\bC^c$ has a collection $\{\bC^c_{\alpha}\}_{\alpha}$ of full subcategories satisfying
\begin{enumerate}
    \item each $\bC^c_{\alpha}$ is stable under the $\Vect^{\on{f.d.}}$-action restricted from the $\Vect$-action on $\bC$;
    \item each $\bC^c_{\alpha}$ is stable under translations and finite extensions;
    \item $\bC^c_{\alpha}$'s are mutually orthogonal to each other in $\bC$;
    \item every object in $\bC^c$ can be expressed as a finite extension of objects from $\{\bC^c_{\alpha}\}_{\alpha}$.
\end{enumerate}

\begin{prop}\label{prop: orth decomp of cpt gen cat}
Under the hypothesis of Sect. \ref{sssec: setup orthogonal decomposition} we have:
    \begin{enumerate}
    \item each $\{\bC^c_{\alpha}\}_{\alpha}$ belongs to $\on{DGCat}^{\omega}$; 
    \item and the natural functor
    \[\Sigma:\bigoplus_{\alpha}\Ind(\bC^c_{\alpha})\rightarrow \bC\quad\te{in }\DGcat\]
    is an equivalence.
\end{enumerate}
\end{prop}

\sssec{}
Using Lemma \ref{lem: colim cpt gen} we prove that $\Shv^{\indconstr,q}(\Hecke^{\loc}_X)$ is compactly generated. 
Moreover, the compact objects therein are constructible, using the same argument as in \cite[Sect. F.2.2]{AGKRRV}.
\medskip

We apply Proposition \ref{prop: orth decomp of cpt gen cat} to the situation where the compactly generated DG category $\bC$ is $\Shv^{\indconstr,q}(\Hecke^{\loc}_X)$ and we use the collection $\{\bC^c_{\theta^{\sharp}}\}_{\theta^{\sharp}\in\wh{Z}_H}$ where
\[\bC^c_{\theta^{\sharp}}=\on{Im}(\Sigma'_{\theta^{\sharp}})^{\on{fin.ext},c}.\]
If we verify all the hypotheses of Proposition \ref{prop: orth decomp of cpt gen cat}, then we are done proving that the functor $\Sigma$ is an equivalence in $\DGcat$. The first two hypotheses obviously hold true, and the third hypothesis is verified by the orthogonality established in Proposition \ref{prop: orthogonality}. Now we are left with verifying the last hypothesis. 

\sssec{Verification of the last hypothesis}
To verify the last hypothesis, it suffices to show that every compact object in $\Shv^{\indconstr,q}(\Hecke^{\loc}_X)$ can be expressed as a finite extension of objects in $\on{Im}(\Sigma'_{\theta^{\sharp}})$ for varying $\theta^{\sharp}$'s. Indeed, once we prove this, from the orthogonality (Proposition \ref{prop: orthogonality}) a compact object in $\Shv^{\on{constr},q}(\Hecke^{\loc}_X)$ is a direct sum of objects in $\on{Im}(\Sigma'_{\theta^{\sharp}})^{\on{fin.ext}}$ for $\theta^{\sharp}$'s, and since the category of compact objects is closed under retracts, each of the summands is compact, hence lies in $\on{Im}(\Sigma'_{\theta^{\sharp}})^{\on{fin.ext},c}$. Since compact objects are constructible, it suffices to prove the same property for constructible sheaves.
\medskip

 Every object of $\Shv^{\constr,q}(\Hecke^{\loc}_X)$ can be obtained as a finitely iterated extension of objects in the form $h_!\cL$, where we have $\lambda\in\Lambda^+$ and $h:Y\hookrightarrow \Gr^{\lambda}_X$ is a $L^+G$-equivariant locally closed embedding of a smooth variety and $\cL$ is a $\cG$-twisted $L^+G$-equivariant abelian local system of finite rank over $Y$. Hence, it suffices to show that such $h_!\cL$ belongs to $\on{Im}(\Sigma'_{?})$ for some $?\in \wh{Z}_H$.
 From \cite[Cor 3.1.10]{ZhaoNew} the existence of a nonzero $\cL$ implies $\lambda\in\Lambda^{\sharp,+}$. Using the identification $\wh{Z}_H= \Lambda^{\sharp}/\Lambda^{\sharp,r}$ we have $\theta^{\sharp}\in\wh{Z}_H$ corresponding to $\lambda$. We claim that $h_!\cL\in \on{Im}(\Sigma'_{\theta^{\sharp}})$. It is sufficient to equip $\cL$ with a $\cG^{\theta^{\sharp}}$-twisted $(L^+G_{\on{ad}})_X$-equivariance. We have a commutative diagram
\[\begin{tikzcd}
	{(L^+G)_X\backslash \Gr^{\lambda}_X} & X \\
	{(L^+G_{\on{ad}})_X\backslash \Gr^{\lambda}_X.}
	\arrow["p", from=1-1, to=1-2]
	\arrow["q"', from=1-1, to=2-1]
	\arrow["{p'}"', from=2-1, to=1-2]
\end{tikzcd}\]
For such $\lambda\in\Lambda^{\sharp,+}$ the restriction $\cG$ to the Schubert cell labelled by $\lambda$ canonically descends to $X$: there is a gerbe $\cG_{\varpi^{\lambda}}$ on $X$ and a canonical isomorphism
\begin{equation}\label{eq: descending gerbe on Schubert cell}
    \cG\cong p^*\cG_{\varpi^{\lambda}}\quad \te{over } (L^+G)_X\backslash \Gr^{\lambda}_X
\end{equation}
(\cite[Corollary 3.1.12]{ZhaoNew}), and by Lemma 5.2.7 of \emph{op. cit.} we have a canonical isomorphism
\begin{equation}\label{eq: descending Schubert cell plus BZ}
    (p')^*\cG_{\varpi^{\lambda}}\simeq \cG^{\theta^{\sharp}}\quad \te{over } (L^+G_{\on{ad}})_X\backslash \Gr^{\lambda}_X.
\end{equation}
The equivalence \eqref{eq: descending Schubert cell plus BZ} lifts the equivalence \eqref{eq: descending gerbe on Schubert cell} in the sense that those two are compatible with the data $q^*\cG^{\theta^{\sharp}}\simeq \cG$ in an obvious way.
\medskip

From the $(L^+G)_X$-invariance of $h$, the subvariety $Y$ should be in a form $\pi^{-1}(Y_0)$ under the projection $\pi:\Gr^{\lambda}_X\rightarrow X$ and a locally closed smooth subscheme $Y_0\subseteq X$. By pullback, we get a commutative diagram
\[\begin{tikzcd}
	{\LS^{\Ab,\ft,\cG}((L^+G)_{Y_0}\backslash Y)} & {\LS^{\Ab,\ft,\cG_{\varpi_{\lambda}}}(Y_0)} \\
	{\LS^{\Ab,\ft,\cG^{\theta^{\sharp}}}((L^+G_{\on{ad}})_{Y_0}\backslash Y).}
	\arrow[from=1-2, to=1-1]
	\arrow[from=1-2, to=2-1]
	\arrow[from=2-1, to=1-1]
\end{tikzcd}\]
The horizontal arrow is an equivalence because $L^+G\backslash Y$ is a quotient of $Y_0$ by a connected group scheme over $Y_0$.
Therefore we conclude the essential surjectivity of the vertical arrow. This completes the verification of the last hypothesis.

\sssec{}\label{sssec: lift of Shv Hecke}
From the arguments above, we have proved that the functor $\Sigma$ is an equivalence in $\DGcat$. Now we discuss how the equivalence $\Sigma$ equips $\Shv^{\indconstr,q}(\Hecke^{\loc}_X)$ with a $\Shv^{\indconstr}(X)$-linear $\wh{Z}_H$-graded monoidal structure.
\medskip

The category $\Shv^{\indconstr,q}(\Hecke^{\loc}_X)$ has a natural $\Shv^{\indconstr}(X)$-linear monoidal structure. Hence the constant functor
\begin{equation}\label{eq: constant graded monoidal functor}
    \wh{Z}_H\rightarrow \Shv^{\indconstr}(X)\lMod
\end{equation}
with value $\Shv^{\indconstr,q}(\Hecke^{\loc}_X)$ has a right-lax monoidal structure. We will equip $\Shv^q(\Hecke^{\loc}_X)$ with a desired graded linear monoidal structure by realizing $\Sigma$ as a `subfunctor' of \eqref{eq: constant graded monoidal functor}. To do so, it suffices to prove the following three claims:
\begin{enumerate}
    \item each graded piece $\Ind(\on{Im}(\Sigma'_{\theta^{\sharp}})^{\on{fin.ext,c}})$ is stable under the $\Shv^{\indconstr}(X)$-action;
    \item the monoidal unit of $\Shv^{\indconstr,q}(\Hecke^{\loc}_X)$ lies in the $0$-th graded piece;
    \item given $\theta^{\sharp}_1,\theta^{\sharp}_2\in\wh{Z}_H$, the convolution product of $\Shv^{\indconstr,q}(\Hecke^{\loc}_X)$ restricts to
    \[\Ind(\on{Im}(\Sigma'_{\theta^{\sharp}_1})^{\on{fin.ext,c}})\times \Ind(\on{Im}(\Sigma'_{\theta^{\sharp}_2})^{\on{fin.ext,c}})\rightarrow \Ind(\on{Im}(\Sigma'_{\theta^{\sharp}_1+\theta^{\sharp}_2})^{\on{fin.ext,c}}).\]
\end{enumerate}

\begin{lemma}\label{lem: Im containment}
    We have a containment
    \[\on{Im}(\Sigma'_{\theta^{\sharp}})\subseteq\Ind(\on{Im}(\Sigma'_{\theta^{\sharp}})^{\on{fin.ext,c}}).\]
\end{lemma}
\begin{proof}
    Since we have established the orthogonal decomposition of $\Shv^{\indconstr,q}(\Hecke^{\loc}_X)$, it suffices to show that $\on{Im}(\Sigma'_{\theta^{\sharp}})$ is \emph{right} orthogonal to all the other graded pieces $\Ind(\on{Im}(\Sigma'_{\theta^{\sharp }_{
    \circ
    }})^{\on{fin.ext,c}})$ for $\theta^{\sharp}_{\circ}\ne \theta^{\sharp}$. However, this is a direct consequence of Proposition \ref{prop: orthogonality}.
\end{proof}
\begin{prop}\label{prop: graded pieces alternative}
 The subcategory
    $\Ind(\on{Im}(\Sigma'_{\theta^{\sharp}})^{\on{fin.ext,c}})$ is identical to the subcategory of $\Shv^{\indconstr,q}(\Hecke^{\loc}_X)$ generated by $\on{Im}(\Sigma'_{\theta^{\sharp}})$ under colimits.
\end{prop}
\begin{proof}
    For convenience, denote the latter category by $\bD_{\theta^{\sharp}}$. There is an obvious inclusion 
    \[\Ind(\on{Im}(\Sigma'_{\theta^{\sharp}})^{\on{fin.ext,c}})\subseteq \bD_{\theta^{\sharp}}.\]
    On the other hand, by Lemma \ref{lem: Im containment} and the fact that 
    \[\Ind(\on{Im}(\Sigma'_{\theta^{\sharp}})^{\on{fin.ext,c}})\hookrightarrow \Shv^{\indconstr,q}(\Hecke^{\loc}_X)\]
    lies in $\DGcat$ we also get the reverse inclusion.
\end{proof}

\sssec{}
Now we prove the claims made in Sect. \ref{sssec: lift of Shv Hecke} using the alternative description of graded pieces provided by Proposition \ref{prop: graded pieces alternative}.
\medskip

The second claim follows immediately from the fact that the monoidal unit belongs to $\on{Im}(\Sigma'_{0})$.
The first claim follows from the fact that the functor $\Sigma'_{\theta^{\sharp}}$ has a natural $\Shv^{\indconstr}(X)$-linear structure, and hence $\on{Im}(\Sigma'_{\theta^{\sharp}})$ is stable under the action of
$\Shv^{\indconstr}(X)$.
For the third claim, it suffices to show that the convolution product restricts to
\[\on{Im}(\Sigma'_{\theta^{\sharp}_1})\times \on{Im}(\Sigma'_{\theta^{\sharp}_2})\rightarrow \on{Im}(\Sigma'_{\theta^{\sharp}_1+\theta^{\sharp}_2}).\]
However, this follows from the compatibility of the gerbes $\{\cG_{\theta^{\sharp}}\}$ with convolution product established in \cite[Lemma 5.2.3]{ZhaoNew}.
\medskip

Having verified all the claims made in Sect. \ref{sssec: lift of Shv Hecke}, we complete the proof of Theorem \ref{thm: A side grading sharper}, and hence Theorem \ref{thm: A side grading}. \qed

\subsection{Various versions of local systems}\label{ssec: various versions of local systems}
This Subsection serves as a preparation for Sect. \ref{ssec: characterization lemmas} and \ref{ssec: LS linear general}. We introduce other notions of local system and study their relationships with $\QLisse$ and $\LS$. 

\sssec{}
For $Y\in\Schft$, we can consider the category $\LS^{\ft}(Y)$, which is the full subcategory of $\LS(Y)$ consisting of objects having finitely many non-trivial cohomology sheaves, each of which having finite rank. Equivalently, $\LS^{\ft}(Y)$ is the category of dualizable objects in $\LS(Y)$. We take its $\Ind$-completion to obtain a DG category \footnote{In \cite[Sect. 1.3]{AGKRRV}, the same category was denoted as $\on{IndLisse}(Y)$.}
\[\LS^{\indft}(Y):=\Ind(\LS^{\ft}(Y)).\]
Note that it was not necessary to let $Y\in\Schft$; we can define both categories given a homotopy type.
\medskip

We have a sequence of canonical functors\footnote{The first functor was discussed in \cite[Sect. 1.3.2]{AGKRRV}, and the second functor was discussed in Remark G.4.4 of \emph{op. cit.}.}
\begin{equation}\label{eq: comparison sequence}
    \LS^{\indft}(Y)\rightarrow \QLisse(Y)\rightarrow \LS(Y).
\end{equation}

\begin{lemma}[\cite{AGKRRV}, Theorem E.2.8]\label{lem: indft vs QLisse}
    If $Y'$ is a connected algebraic curve other than $\pjl$, then the first arrow of \eqref{eq: comparison sequence} is an equivalence.
\end{lemma}
\begin{lemma}[\cite{AGKRRV}, Sect. E.2.6]\label{lem: QLisse vs LS}
    If $Y'=\pjl$, then the second arrow of \eqref{eq: comparison sequence} is an equivalence.
\end{lemma}

\begin{lemma}\label{lem: dualizability}
    For $Y\in\Spc$, the categories $\LS^{\indft}(Y)$ and $\LS(Y)$ are compactly generated, and hence dualizable.
    For a smooth curve $X'$ in $\Schft$, the category $\QLisse(X')$ is compactly generated, and thus dualizable.
\end{lemma}
\begin{proof}
    The only non-trivial part is the compact generation of $\QLisse(X')$.
    It suffices to prove the claim for the case $X'$ is connected. If $X'\ne\pjl$, this follows from Lemma \ref{lem: indft vs QLisse}. If $X'=\pjl$, this follows from Lemma \ref{lem: QLisse vs LS}.
\end{proof}

\begin{lemma}\label{lem: QL to LS containment}
    For a smooth curve $X'$ in $\Schft$, the canonical functor $\QLisse(X')\rightarrow \LS(X')$ is fully faithful.
\end{lemma}
\begin{proof}
    It suffices to prove this for a connected $X'$. If $X'=\pjl$, then this functor is an equivalence. If $X'\ne\pjl$, since $\LS^{\indft}(Y)\rightarrow \QLisse(Y)$ is an equivalence, it suffices to show the containment
    $\LS^{\ft}(X')\subseteq \LS(X')^c$.
    This holds true in general when we replace $X'$ by any topological space $K$ homotopy equivalent to a finite cell complex. In that case, the constant sheaf is compact in $\LS(K)$, and hence dualizable objects are compact in that category. \end{proof}

    \begin{lemma}\label{lem: Ab ft on product into tensor of ind ft}
        Suppose we have topological spaces $Y_1,\cdots,Y_r$, each homotopy equivalent to a finite cell complex. There is a canonical embedding
        \[\LS^{\Ab,\ft}(Y_1\times\cdots\times Y_r)\hookrightarrow \LS^{\indft}(Y_1)\otimes\cdots\otimes \LS^{\indft}(Y_r).\]
    \end{lemma}
    \begin{proof}
        Let us prove the case $r=2$ for notational convenience. The other cases are proved with the same logic.
        We have an equivalence $\LS(Y_1\times Y_2)=\LS(Y_1)\otimes \LS(Y_2)$. From the assumptions on $Y_i$ we get a fully-faithful functor
        \[\LS^{\indft}(Y_i)\hookrightarrow \LS(Y_i)\quad \te{in }\DGcat\]
        as in the proof of Lemma \ref{lem: QL to LS containment}.
        Since the categories $\LS^{\indft}(Y_i)$ and $\LS(Y_i)$ are all dualizable, (Lemma \ref{lem: dualizability}) we have a fully-faithful \footnote{The fully faithfulness follows from the following general fact: if $C_1\rightarrow C_2$ in $\DGcat$ is fully faithful and $D$ is a dualizable DG category, then $C_1\otimes D \rightarrow C_2\otimes D$ is fully faithful.} embedding
        \begin{equation}\label{eq: tensor prod indft into tensor prod LS}
            \LS^{\indft}(Y_1)\otimes \LS^{\indft}(Y_2)\hookrightarrow \LS(Y_1)\otimes \LS(Y_2).
        \end{equation}
        On the other hand, we have a fully-faithful embedding
        \[\LS^{\Ab,\ft}(Y_1\times Y_2)\hookrightarrow \LS(Y_1)\otimes \LS(Y_2),\]
 hence, it suffices to show that all objects of $\LS^{\Ab,\ft}(Y_1\times Y_2)$ belong to $\LS^{\indft}(Y_1)\otimes \LS^{\indft}(Y_2)$.
 \medskip
        
        The DG subcategory $\LS^{\indft}(Y_1)\otimes \LS^{\indft}(Y_2)$ is closed under extensions in $\LS(Y_1)\otimes \LS(Y_2)$, hence it suffices to show that every irreducible local system of finite type over $Y_1\times Y_2$ belong to $\LS^{\indft}(Y_1)\otimes \LS^{\indft}(Y_2)$. However, such a local system takes a form $\cL_1\boxtimes \cL_2$ where each $\cL_i$ is an irreducible local system of finite type over $Y_i$\footnote{Here we used the assumption that $E$ is algebraically closed.},  and therefore it lies in $\LS^{\indft}(Y_1)\otimes \LS^{\indft}(Y_2)$.
    \end{proof}

    \begin{prop}\label{prop: Ab ft I into QLisse tensor}
        For a smooth curve $X'$ in $\Schft$, the full subcategory
        \[\LS^{\Ab,\ft}(X'^I)\hookrightarrow \LS(X'^I)\]
        lies in the full subcategory
        \[\QLisse(X')^{\otimes I}\hookrightarrow \LS(X')^{\otimes I}=\LS(X'^I).\]
    \end{prop}
    \begin{proof}
        If $X'\ne\pjl$, then $\QLisse(X')=\LS^{\indft}(X')$, and hence this is a consequence of Lemma \ref{lem: Ab ft on product into tensor of ind ft} and its proof. If $X'=\pjl$, then this is trivial since $\QLisse(X')=\LS(X')$.
    \end{proof}

\sssec{}
From Proposition \ref{prop: Ab ft I into QLisse tensor} we get a chain of fully-faithful functors 
\begin{equation}\label{eq: chain inclusion}
    \LS^{\Ab,\ft}(X'^I)\hookrightarrow \QLisse(X')^{\otimes I}\hookrightarrow \LS(X')^{\otimes I}.
\end{equation}
All the inclusion maps are symmetric monoidal because both $\LS^{\Ab,\ft}(X'^I)$ and $\QLisse(X')^{\otimes I}$ are closed in $\LS(X'^I)$ under the monoidal operation.

\sssec{}
For each $I\in\fSet$ we have a sheaf
\[\mathbf{LS}^{\Ab,\ft}_I:(X_{\et})^{\op}\rightarrow\CatEc,\quad X'\mapsto \LS^{\Ab,\ft}(X'^I).\]
The collection of \eqref{eq: chain inclusion} defines a chain of maps
\begin{equation}\label{eq: chain of sheaves}
    \textbf{LS}^{\Ab,\ft}_I\rightarrow\textbf{QLisse}^{\otimes I}\rightarrow\textbf{LS}^{\otimes I}\quad \te{in }\CAlg(\Shv(X_{\et};E^c)).
\end{equation}
Note that the $I$-fold tensor product of sheaves $\textbf{QLisse}$ and $\textbf{LS}$ are calculated pointwise, i.e., we have
\[\Gamma(X', \textbf{QLisse}^{\otimes I})=\QLisse(X')^{\otimes I},\quad \Gamma(X', \textbf{LS}^{\otimes I})=\LS(X')^{\otimes I}.\]
This is because of the dualizability in Lemma \ref{lem: dualizability}. 

\subsection{Full compatibility with extensions}\label{ssec: full compatibility with extensions}
The interplay of abelian categories and DG categories will be important for us (see the introduction to Sect. \ref{ssec: characterization lemmas}). In particular, we will often encounter the situation where we have an $E$-linear abelian category $\cA$ and a (\emph{not} necessarily cocomplete) DG category $\cC$ together with a fully faithful functor $F: \cA\hookrightarrow \cC$. Each category has their own notion of extension: for abelian category it is a short exact sequence, and for a stable category it is a cofiber sequence (or equivalently a fiber sequence). In this subsection, we will record some facts about the interaction of $F$ with those two notions of extensions.

\sssec{}
Consider $F:\cA\hookrightarrow \cC$ as before. We say that $F$ is \textbf{fully compatible with extensions} if it satisfies the following three properties:
\begin{enumerate}
    \item (\emph{closure}) $\cA$ is closed under the extension in $\cC$;
    \item (\emph{preservation}) an exact sequence in $\cA$ is sent to a cofiber sequence in $\cC$ under $F$;
    \item (\emph{reflection}) if a cofiber sequence in $\cC$ satisfies the condition that every object lies in $\cA$, then it is an exact sequence in $\cA$.
\end{enumerate}

\begin{ex}\label{ex: t str full compatibility}
    If $\cC$ is equipped with a $t$-structure and $F$ realizes $\cA$ as a Serre subcategory of $\cC^{\heartsuit}$, then $F$ is fully compatible with extensions.
\end{ex}

\begin{lemma}\label{lem: full compatibility propagate}
    Suppose we have an $E$-linear abelian category $\cA$, a (not necessarily cocomplete) DG category $\cC$, and its full DG subcategory $\cD$. 
    Suppose we have a fully faithful functor $F:\cA\rightarrow\cD$ and hence $G:\cA\rightarrow \cC$ as the composition. In this situation, $F$ is fully compatible with extensions if and only if $G$ is fully compatible with extensions.
\end{lemma}

\begin{cor}\label{cor: Abft QL full compatibility}
    For a smooth curve $X'$ in $\Schft$, the fully faithful functor $\LS^{\Ab,\ft}(X'^I)\hookrightarrow \QLisse(X')^{\otimes I}$ constructed in Proposition \ref{prop: Ab ft I into QLisse tensor} is fully compatible with extensions.
\end{cor}
\begin{proof}
This is an immediate consequence of Lemma \ref{lem: full compatibility propagate} and Proposition \ref{prop: Ab ft I into QLisse tensor}.
We have a commutative diagram
\begin{equation}\label{eq: Abft QL LS diagram}
    \begin{tikzcd}
	{\LS^{\Ab,\ft}(X'^I)} & \\
	{\QLisse(X')^{\otimes I}} & {\LS(X'^I).}
	\arrow[hook, from=1-1, to=2-1]
	\arrow[hook', from=1-1, to=2-2]
	\arrow[hook, from=2-1, to=2-2]
\end{tikzcd}
\end{equation}
The diagonal arrow is fully compatible with extensions because it is an instance of Example \ref{ex: t str full compatibility}, where the natural $t$-structure on $\LS(X'^I)$ is considered. Moreover, the horizontal arrow lies in $\DGcat$, hence is an inclusion of a DG subcategory.
\end{proof}

\begin{cor}\label{cor: Sat full compatibility}
    For a smooth curve $X'$ in $\Schft$, the inclusion
    \[\Sat^q_{X'^I}\hookrightarrow \Shv^q(\Hecke^{\loc}_{X'^I})\]
    is fully compatible with extensions.
\end{cor}
\begin{proof}
    We consider the commutative diagram
    \[\begin{tikzcd}
	{\Sat^q_{X'^I}} & \\
	{\Shv^{\on{constr},q}(\Hecke^{\loc}_{X'^I})} & {\Shv^{q}(\Hecke^{\loc}_{X'^I}),}
	\arrow[hook, from=1-1, to=2-1]
	\arrow[hook', from=1-1, to=2-2]
	\arrow[hook, from=2-1, to=2-2]
\end{tikzcd}\]
where the superscript $\on{constr}$ means the (non-cocomplete) DG category of constructible sheaves.
The vertical arrow is fully compatible with extensions, because it is an instance of Example \ref{ex: t str full compatibility}. Indeed, we can consider the perverse $t$-structure on the category of constructible sheaves, and the category of perverse ULA sheaves is a Serre subcategory of the category of perverse sheaves. Since the horizontal arrow is an inclusion of DG subcategory, we conclude from Lemma \ref{lem: full compatibility propagate} that the diagonal arrow is also fully compatible with extensions.
\end{proof}

\sssec{}
Given a reductive group $\mathsf{H}$ and $Y\in\Schft$, let $\Rep(\mathsf{H})^{\QLisse}(Y)$ be the category of $\mathsf{H}$-representations, i.e., $\cO_{\mathsf{H}}$-comodules, on $\QLisse(Y)$. We have a canonical equivalence
\[\Rep(\mathsf{H})^{\QLisse}(Y)=\Rep(\mathsf{H})\otimes\QLisse(Y).\]
The collection of such categories form a sheaf $\mathbf{Rep}(H)^{\QLisse}$ used in Sect. \ref{sssec: list of sheaves}.
\medskip

Taking the categories of $H^I$-representations, i.e., $\cO_{H^I}$-comodules, on each categories of the diagram \eqref{eq: Abft QL LS diagram} we get a commutative diagram
\begin{equation}\label{eq: Rep H Abft QL LS}
    \begin{tikzcd}
	{\Rep(H^I)^{\Ab,\loc,\ft}(X'^I)} & \\
	{\left(\Rep(H)^{\QLisse}(X')\right)^{\otimes I}} & {\left(\Rep(H)^{\loc}(X')\right)^{\otimes I}.}
	\arrow[hook, from=1-1, to=2-1]
	\arrow[hook', from=1-1, to=2-2]
	\arrow[hook, from=2-1, to=2-2]
\end{tikzcd}
\end{equation}

\begin{lemma}\label{lem: Rep H Abft QL full compatibility}
    The vertical arrow of the diagram \eqref{eq: Rep H Abft QL LS} is fully compatible with extensions.
\end{lemma}
\begin{proof}
    We have a canonical equivalence
\[\left(\Rep(H)^{\loc}(X')\right)^{\otimes I} =\Rep(H^I)^{\loc}(X'^I),\]
and this category is equipped with a $t$-structure uniquely characterized by the requirement that the forgetful functor
\[\Rep(H^I)^{\loc}(X'^I)\rightarrow \LS(X'^I)\]
is $t$-exact with respect to the natural $t$-structure on the target. Its abelian heart is $\Rep(H^I)^{\Ab,\loc}(X'^I)$, which has $\Rep(H^I)^{\Ab,\loc,\ft}(X'^I)$ as a Serre subcategory. Therefore, by Example \ref{ex: t str full compatibility} the diagonal arrow of \eqref{eq: Rep H Abft QL LS} is fully compatible with extensions. We conclude using Lemma \ref{lem: full compatibility propagate}.
\end{proof}

\sssec{}
Suppose we are in a situation where we have a collection of DG categories $\cC_1,\cdots,\cC_r, \cD\in\DGcat$ and a functor
\[F: \cC_1\otimes\cdots\otimes \cC_r\rightarrow \cD \quad \te{in }\DGcat.\]
Suppose that we have a collection of $E$-linear abelian categories $\cA_1,\cdots,\cA_r,\cB$ equipped with fully faithful functors
\[\cA_1\hookrightarrow \cC_1,\cdots, \cA_r\hookrightarrow \cC_r,\quad\te{and}\quad \cB\hookrightarrow\cD,\]
all fully compatible with extensions. 

\begin{lemma}\label{lem: induced exactness}
    In the above situation, suppose that the composition
    \[\cA_1\times \cdots\times\cA_r\rightarrow \cC_1\otimes\cdots\otimes\cC_r\xrightarrow{F}\cD\]
    lands in $\cB$. Then the induced map $\cA_1\times\cdots\times\cA_r\rightarrow \cB$
    is exact on each variables.
\end{lemma}

\subsection{Characterization lemmas}\label{ssec: characterization lemmas}
 Our construction of $\mathbf{S}^{\QLisse}_{\{*\}}$ was quite involved, and hence so will be that of $\mathbf{S}^{\QLisse}_{I}$, which is going to be constructed using $\mathbf{S}^{\QLisse}_{\{*\}}$. So, it is very useful to have a characterization of these functors in terms of perverse ULA Satake equivalence. 
\medskip

A perverse ULA Satake functor $\mathbf{S}^{\Ab}_I$ is a $\mathbf{LS}^{\Ab,\ft}_I$-linear functor out of $\mathbf{Rep}(H^I)^{\Ab,\loc,\ft}_{\cG^{\boxtimes I}_Z}$, and this motivates us to establish a characterization of a $\mathbf{QLisse}^{\otimes I}$-linear functor out of $\left(\mathbf{Rep}(H)^{\QLisse}_{\cG_Z}\right)^{\otimes I}$ \footnote{$\mathbf{S}^{\QLisse}_I$ is an example of such.} by the restriction to $\mathbf{Rep}(H^I)^{\Ab,\loc,\ft}_{\cG^{\boxtimes I}_Z}$ as a $\mathbf{LS}^{\Ab,\ft}_I$-linear functor. This is done in Lemma \ref{lem: the characterization lemma} under a mild hypothesis, which will be harmless in our context. Finally, we record the following standard result:

\begin{lemma}
    Suppose $X\ne\pjl$ and $X'\in X_{\et}$. Then $X'\ne\pjl$.
\end{lemma}

\sssec{}\label{sssec: QLisse and ind ft under not P1}
Therefore, under the assumption $X\ne\pjl$, each $X'\in X_{\et}$ cannot have a connected component isomorphic to $\pjl$, and hence $\QLisse(X')=\LS^{\indft}(X')$.

    \sssec{The characterization lemma}
    Now we formulate the characterization lemma announced at the beginning of this subsection.
     Recall the inclusion
    \[\Rep(H^I)^{\Ab,\loc,\ft}(X^I)\hookrightarrow \left(\Rep(H)^{\QLisse}(X)\right)^{\otimes I}\]
    considered in the diagram \eqref{eq: Rep H Abft QL LS}.
    Both sides naturally sheafify to yield a subsheaf
    \[\mathbf{Rep}(H^I)^{\Ab,\loc,\ft}\rightarrow \left(\textbf{Rep}(H)^{\QLisse}\right)^{\otimes I}.\]
    Since the subsheaf relation is stable under twist, we have a subsheaf
    \begin{equation}\label{eq: Ab ft twisted in QLisse twisted}
        \mathbf{Rep}(H^I)^{\Ab,\loc,\ft}_{\cG^{\boxtimes I}_Z}\rightarrow \left(\textbf{Rep}(H)^{\QLisse}_{\cG_Z}\right)^{\otimes I}.
    \end{equation}

\begin{lemma}\label{lem: the characterization lemma}
    Assume either $X\ne \pjl$ or that $X=\pjl$ and $\cG_Z$ is trivializable. Suppose we have a $\mathbf{QLisse}^{\otimes I}$-module object $\mathbf{F}\in \Shv(X_{\et};\DGcat)$. A $\mathbf{QLisse}^{\otimes I}$-linear functor
    \begin{equation}\label{eq: functor out of Rep H QLisse}
        \left(\mathbf{Rep}(H)^{\QLisse}_{\cG_Z}\right)^{\otimes I}\rightarrow \mathbf{F}
    \end{equation}
        is uniquely recovered from its restriction 
        \begin{equation}\label{eq: functor out of Rep H I Ab ft}
            \mathbf{Rep}(H^I)^{\Ab,\loc,\ft}_{\cG^{\boxtimes I}_Z}\rightarrow \mathbf{F}
        \end{equation}
        along \eqref{eq: Ab ft twisted in QLisse twisted} as a $\mathbf{LS}^{\Ab,\ft}_I$-linear functor. 
    \end{lemma}

    \sssec{}\label{sssec: char paradigm} For the proof of Lemma \ref{lem: the characterization lemma} we consider the following categorical setup.
    Suppose that we have DG categories $\bA\in\DGCatMon$ and $\bM\in\DGcat$ with an action of $\bA$ on $\bM$. Assume furthermore that there are full subcategories $\bA^{\circ}\subseteq \bA$ and $\bM^{\circ}\subseteq \bM$ stable under the action of $\Vect^{\heartsuit,\on{f.d.}}$ satisfying the following properties:
    \begin{enumerate}
        \item $\bA^{\circ}$ (resp. $\bM^{\circ}$) compactly generates $\bA$ (resp. $\bM$);
        \item the monoidal structure of $\bA$ restricts to a monoidal structure on $\bA^{\circ}$ and the action of $\bA$ on $\bM$ restricts to an action of $\bA^{\circ}$ on $\bM^{\circ}$.
    \end{enumerate}
    Given another $\bA$-module object $\bM'\in \DGcat$ we have a restriction functor
    \begin{equation}\label{eq: restriction functor}
        \Fun_{\bA\lMod}(\bM,\bM')\rightarrow \Fun_{\bA^{\circ}\lMod}(\bM^{\circ}, \bM'),
    \end{equation}
    where $\bA\lMod$ on the source is internal to $\DGcat$ while $\bA^{\circ}\lMod$ on the target is internal to $\CatEc$. Then the functor \eqref{eq: restriction functor} is fully faithful.

    \sssec{Proof under the assumption $X\ne\pjl$}
        It suffices to show that on each $X'\in X_{\et}$ where $\cG_Z$ is trivializable, the restriction functor
        \[\Fun_{\QLisse(X')^{\otimes I}\lMod}\left(\left(\Rep(H)^{\QLisse}_{\cG_Z}(X')\right)^{\otimes I}, \mathbf{F}(X')\right)\rightarrow \Fun_{\LS^{\Ab,\ft}(X'^I)\lMod}\left(\Rep(H^I)^{\Ab,\loc,\ft}_{\cG^{\boxtimes I}_Z}(X'^I),\mathbf{F}(X')\right)\]
        is fully-faithful. This is an application of the paradigm of Sect. \ref{sssec: char paradigm} to
        \[\bA=\QLisse(X')^{\otimes I},\quad \bA^{\circ}=\LS^{\Ab,\ft}(X'^I),\quad \mathbf{M}=\left(\Rep(H)^{\QLisse}_{\cG_Z}(X')\right)^{\otimes I},\ \ \te{and}\ \ \mathbf{M}^{\circ}=\Rep(H^I)^{\Ab,\loc,\ft}_{\cG^{\boxtimes I}_Z}(X'^I).\]
        
        Let us check that these data satisfy the assumptions of Sect. \ref{sssec: char paradigm}. For this purpose, we may assume $\cG_Z$ is trivial. Throughout the proof, we will continue to use the fact that both inclusions
        $\bA^{\circ}\hookrightarrow \bA$ and $\bM^{\circ}\hookrightarrow\bM$
        are fully compatible with extensions. (see Corollary \ref{cor: Abft QL full compatibility} and Lemma \ref{lem: Rep H Abft QL full compatibility})
        \medskip
        
        First of all, we have $\bA^{\circ}$ as a monoidal subcategory of $\bA$ (Proposition \ref{prop: Ab ft I into QLisse tensor}). Second, we claim that $\bA$ is compactly generated by $\bA^{\circ}$. The containment $\bA^{\circ}\subseteq \bA^c$ follows from the fact that any object of $\bA^{\circ}$ can be obtained as a finite extension of external product of a form $\boxtimes_{i\in I}\cL_i$ for $\cL_i\in \LS^{\Ab,\ft}(X')$ where each $\cL_i$ is irreducible. (see the proof of Lemma \ref{lem: Ab ft on product into tensor of ind ft}) Since $\bA$ is compactly generated by the external products $\boxtimes_{i\in I}\cL_i$ where $\cL_i\in \LS^{\ft}(X')$ (see Sect. \ref{sssec: QLisse and ind ft under not P1}), it remains to show that every object in this form can be expressed as a finite extension of objects in $\bA^{\circ}$. This follows from the fact that each $\cL_i$ is expressed as a finite extension of local systems $\cL^j_i\in\LS^{\Ab,\ft}(X')$. Third, the action of $\bA$ on $\bM$ obviously restricts to an action of $\bA^{\circ}$ on $\bM^{\circ}$.
        \medskip

        Finally, we claim that $\bM$ is compactly generated by $\bM^{\circ}$. In addition to the assumption that $\cG_Z$ is trivial, we may assume $X'$ is connected.
        The category $\bM$ is compactly generated by the collection of objects
        \begin{equation}\label{eq: compact generators Rep QLisse}
            V\boxtimes (\boxtimes_{i\in I} \cL_i),\quad V\in\Rep(H^I)^{\heartsuit,\on{f.d.}},\quad \cL\in \LS^{\Ab,\ft}(X').
        \end{equation}
        For the containment $\bM^{\circ}\subseteq \bM^c$, we claim that each object of $\bM^{\circ}$ can be written as a finite extension of objects in a form \eqref{eq: compact generators Rep QLisse}. This will be proved in Lemma \ref{lem: Rep H Ab ft structure}. Since all the objects \eqref{eq: compact generators Rep QLisse} belong to $\bM^{\circ}$, we conclude that $\bM^{\circ}$ compactly generates $\bM$. \qed

        \begin{lemma}\label{lem: Rep H Ab ft structure}
            Suppose we have a reductive group $\mathsf{H}$ over $E$ and a path-connected topological space $Y$. Consider the abelian category $\Rep(\mathsf{H})^{\Ab,\loc,\ft}(Y)$ of $H$-representations on finite rank local systems on $Y$. Every object of $\Rep(\mathsf{H})^{\Ab,\loc,\ft}(Y)$ can be expressed as a finite direct sum of objects in a form
            \[V\boxtimes \cL,\quad V\in\Rep(\mathsf{H})^{\heartsuit,\on{f.d.}}, \ \ \cL\in\LS^{\Ab,\ft}(Y).\]
            If $\mathsf{H}=\prod_{u\in U}\mathsf{H}_u$ and $Y=\prod_{u\in U} Y_u$ for a finite set $U$ and reductive groups $\{\mathsf{H}_u\}$, every object of $\Rep(\mathsf{H})^{\Ab,\loc,\ft}(Y)$ can be expressed as a finite extension of objects in a form
            \[(\boxtimes_u V_u)\boxtimes (\boxtimes_{u}\cL_u)\in\Rep(\mathsf{H})^{\heartsuit,\on{f.d.}}, \quad V_u\in\Rep(\mathsf{H}_u)^{\Ab,\ft}, \ \ \cL_u\in\LS^{\Ab,\ft}(Y_u).\]
        \end{lemma}
        \begin{proof}
        The second assertion follows immediately from the first assertion, because a finite-rank irreducible local systems on $Y$ is a box tensor product of irreducible local systems on $Y_u$'s, and finite-dimensional irreducible representations of $\mathbf{H}$ is a box tensor product of irreducible representations of $\mathsf{H}_u$'s. Now we prove the first assertion.
        \medskip
        
        We identify
        \[\Rep(\mathsf{H})^{\Ab,\loc,\ft}(Y)=\Fun(Y, \Rep(\mathsf{H})^{\heartsuit,\on{f.d.}})=\Fun(\Pi_1(Y), \Rep(\mathsf{H})^{\heartsuit,\on{f.d.}}).\]
        In the last equality we used the fact that the target is an ordinary category. Suppose we are given an object $\mathcal{R}\in \Rep(\mathsf{H})^{\Ab,\loc,\ft}(Y)$. Choose a basepoint of $Y$ so that we can take the underlying object $V\in \Rep(\mathsf{H})^{\heartsuit,\on{f.d.}}$. Then the data of lifting $V$ to $\mathcal{R}$ is a pointed map from $\Pi_1(Y)$ to $B\on{Aut}(V)$, or equivalently a group homomorphism $\pi_1(Y)\rightarrow \on{Aut}(V)$. Using the semisimplicity of $\mathsf{H}$-representations, we can decompose
        \[V=V^{\oplus n_1}_1\oplus \cdots\oplus V^{\oplus n_r}_{r}\]
        where $\{V_s\}$ is the collection of distinct irreducible representations of $\mathsf{H}$. Then from the previous discussion, the data of lifting $V$ is a group homomorphism $\pi_1(Y)\rightarrow \prod_s \on{GL}_{n_s}(E)$. This is the collection 
        $\rho_s: \pi_1(Y)\rightarrow \on{GL}_{n_s}(E)$
        of $n_s$-dimensional representations, and under the monodromy equivalence this corresponds to the collection $\{\cL_s\}\in \LS^{\Ab,\ft}(Y)$ of local systems. Now we have
        \[\mathcal{R}=\bigoplus_{1\le s\le r} V_s\boxtimes \cL_s\quad \te{in } \Rep(\mathsf{H})^{\Ab,\loc,\ft}(Y),\]
        hence we conclude.
        \end{proof}
    
        \sssec{Proof under the assumption $X=\pjl$ and $\cG_Z$ is trivializable}\label{sssec: X=pjl and trivializable}
        In this case, the characterization lemma holds true for a \emph{different} reason. The previous reasoning does not work anymore because $\LS^{\Ab,\ft}(\pjl)$ does \emph{not} compactly generate $\QLisse(\pjl)$.
        \medskip

        However, we observe that a trivializable $Z_H(E)^{\on{tors}}$-gerbe on $\pjl$ is \emph{canonically} pulled back from a $Z_H(E)^{\on{tors}}$-gerbe on a point. More precisely, Let $\on{Ge}^{0}_{Z_H(E)^{\on{tors}}}(X)$ be the full subcategory of $\on{Ge}_{Z_H(E)^{\on{tors}}}(X)$ spanned by trivializable gerbes. The pullback functor of gerbes along $X\rightarrow \pt$
induces a map
\begin{equation}\label{eq: pullback to trivializable}
    \on{Ge}_{Z_H(E)^{\on{tors}}}(\pt)\rightarrow \on{Ge}^0_{Z_H(E)^{\on{tors}}}(X). \end{equation}
Comparing \'etale cohomology groups, we deduce that \eqref{eq: pullback to trivializable} is an equivalence. Indeed, we have
\[\pi_i(\on{Ge}_{Z_H(E)^{\on{tors}}}(X))=H^{2-i}_{\et}(X; Z_H(E)^{\on{tors}})=\begin{cases} Z_H(E)^{\on{tors}} & \te{if } i=0,2,\\
0 & \te{otherwise,}
\end{cases}\]
and therefore
\[\pi_i(\on{Ge}^0_{Z_H(E)^{\on{tors}}}(\pt))=\begin{cases} Z_H(E)^{\on{tors}} & \te{if } i=2,\\
0 & \te{otherwise.}
\end{cases}\]
From this description, we know that \eqref{eq: pullback to trivializable} induces an isomorphism on each $\pi_i$, hence an equivalence.
\medskip

Let $\cG^0\in \on{Ge}_{Z_H(E)^{\on{tors}}}(\pt)$ be the object corresponding to $\cG_Z$ under the equivalence \eqref{eq: pullback to trivializable}. Then $\mathbf{Rep}(H)^{\QLisse}_{\cG_Z}$ is a free $\mathbf{QLisse}$-module, i.e., we have
\[\mathbf{Rep}(H)^{\QLisse}_{\cG_Z}= \ul{\Rep(H)_{\cG^0}}\otimes \mathbf{QLisse}(X)\quad \te{in }\Shv(X_{\et};\DGcat),\]
where $\Rep(H)_{\cG^0}$ is the twist of $\Rep(H)$ by $\cG^0$, and $\ul{\Rep(H)_{\cG^0}}$ is the corresponding constant sheaf. Therefore, we have an equivalence between the space of $\mathbf{QLisse}^{\otimes I}$-linear functors \eqref{eq: functor out of Rep H QLisse} and plain functors
\begin{equation}\label{eq: Rep H 0 to F}
    \ul{\left(\left(\Rep(H)_{\cG^0}\right)^{\otimes I}\right)}^{\on{pre}}\rightarrow \mathbf{F}
\end{equation}
between presheaves, where the left-hand side is the constant presheaf. The category $\Rep(H)_{\cG^0}$ is compactly generated by the full subcategory $\Rep(H^I)^{\Ab,\ft}_{(\cG^0)^{\boxtimes I}}$ spanned by representations on finite-dimensional vector spaces, and hence \eqref{eq: Rep H 0 to F} is uniquely recovered by the further restriction
\begin{equation*}
\ul{\left(\Rep(H^I)^{\Ab,\ft}_{(\cG^0)^{\boxtimes I}}\right)}^{\on{pre}}\rightarrow \mathbf{F}.
\end{equation*}
However, the data of \eqref{eq: functor out of Rep H I Ab ft} determine the map above by restriction.  \qed

    \sssec{Multiple version}
    Now we consider a \emph{multiple} version of Lemma \ref{lem: the characterization lemma}. Suppose $I\in\fSet$ is equipped with a partition $I=\coprod_{u\in U}I_u$. We have a chain of $\Ectensbig_{u\in U}\mathbf{LS}^{\Ab,\ft}_{I_u}$-linear maps
\begin{equation}\label{eq: char multiple restriction map}
   \Ectensbig_{u\in U}\mathbf{Rep}(H^{I_u})^{\Ab,\loc,\ft}_{\cG^{\boxtimes I_u}_Z}\rightarrow \mathbf{Rep}(H^{I})^{\Ab,\loc,\ft}_{\cG^{\boxtimes I}_Z} \hookrightarrow \left(\mathbf{Rep}(H)^{\QLisse}_{\cG_Z}\right)^{\otimes I}  
\end{equation}
in $\Shv(X_{\et};\CatEc)$, where the first arrow is the external tensor product. Applying the proof of Lemma \ref{lem: the characterization lemma} verbatim, we obtain:

    \begin{lemma}\label{lem: the characterization lemma multiple}
    Assume either $X\ne \pjl$ or that $X=\pjl$ and $\cG_Z$ is trivializable.
        Suppose we have $I\in\fSet$ equipped with a partition $I=\coprod_{u\in U} I_u$.
        Suppose we have a $\mathbf{QLisse}^{\otimes I}$-module object $\mathbf{F}\in \Shv(X_{\et};\DGcat)$. A $\mathbf{QLisse}^{\otimes I}$-linear functor
        \[\left(\mathbf{Rep}(H)^{\QLisse}_{\cG_Z}\right)^{\otimes I}\rightarrow \mathbf{F}\]
        is uniquely recovered from its restriction 
        \[\Ectensbig_{u\in U}\mathbf{Rep}(H^{I_u})^{\Ab,\loc,\ft}_{\cG^{\boxtimes I_u}_Z}\rightarrow \mathbf{F}\]
        along \eqref{eq: char multiple restriction map} as a $\Ectensbig_{u\in U}\mathbf{LS}^{\Ab,\ft}_{I_u}$-linear functor in $\Shv(X_{\et}; \CatEc)$. 
    \end{lemma}

\subsection{A $\mathbf{LS}$-linear version of factorizable Satake: the general case}\label{ssec: LS linear general}
In this subsection, we assume either $X\ne \pjl$ or that $X=\pjl$ and $\cG_Z$ is trivializable.

\sssec{}
Let us forget the algebra structure of the functor $\bS^{\QLisse}_{\{*\}}$ (see \eqref{eq: Rep QLisse to Shv all}) and consider the functor as a $\mathbf{QLisse}$-linear functor. By construction, this is compatible with the perverse ULA equivalence, i.e., we have a commutative diagram
\[\begin{tikzcd}
	{\mathbf{Rep}(H)^{\QLisse}_{\cG_Z}} & {\mathbf{Shv}^q(\Hecke^{\loc}_{\{*\}})} \\
	{\mathbf{Rep}(H)^{\Ab,\loc,\ft}_{\cG_Z}} & {\textbf{Sat}^q_{\{*\}}}
	\arrow["{\textbf{S}_{\{*\}}}", from=1-1, to=1-2]
	\arrow[from=2-1, to=1-1]
	\arrow["{\textbf{S}^{\Ab}_{\{*\}}}"', from=2-1, to=2-2]
	\arrow[from=2-2, to=1-2]
\end{tikzcd}\]
as $\mathbf{LS}^{\Ab,\ft}$-modules. By Lemma \ref{lem: the characterization lemma}, this compatibility uniquely characterizes $\bS^{\QLisse}_{\{*\}}$ as a $\mathbf{QLisse}$-linear functor.

\sssec{}
Now we proceed to the construction of the functors $\bS^{\QLisse}_I$'s for $I\in\fSet$ in general. Once we choose a linear ordering $\sigma$ of $I$ thought of as a bijection $\sigma:\{1,2,\cdots,|I|\}\rightarrow I$, we get an ordered partition of $I$, and hence an external convolution product functor
\[\on{conv}_{\sigma}:\Shv^q(\Hecke^{\loc}_X)^{\otimes I}\rightarrow \Shv^q(\Hecke^{\loc}_{X^I}),\]
together with its natural $\LS(X^I)$-linearity. Hence, we get a $\QLisse(X)^{\otimes I}$-linear functor
\begin{equation*}
    \left(\Rep(H)^{\QLisse}_{\cG_Z}(X)\right)^{\otimes I}\xrightarrow{(\cS_{\{*\}})^{\otimes I}}\Shv^q(\Hecke^{\loc}_X)^{\otimes I}\xrightarrow{\on{conv}_{\sigma}} \Shv^q(\Hecke^{\loc}_{X^I}).
\end{equation*}
which we denote as $\cS^{\QLisse,\sigma}_{X,I}$. This sheafifies to yield a $\mathbf{QLisse}^{\otimes I}$-linear functor
\begin{equation}
    \textbf{S}^{\QLisse,\sigma}_I:\left(\mathbf{Rep}(H)^{\QLisse}_{\cG_Z}\right)^{\otimes I}\rightarrow \mathbf{Shv}^q(\Hecke^{\loc}_{X^I}).
\end{equation}
We claim the following compatibility with the perverse ULA version:
\begin{thm}\label{thm: characterization for I}
    For $I\in\fSet$, we have the following commutative diagram of $\mathbf{LS}^{\Ab,\ft}_I$-modules
    \begin{equation}\label{eq: characterization for I}
        \begin{tikzcd}
	{\left(\mathbf{Rep}(H)^{\QLisse}_{\cG_Z}\right)^{\otimes I}} & {\mathbf{Shv}^q(\Hecke^{\loc}_{I})} \\
	{\mathbf{Rep}(H^I)^{\Ab,\loc,\ft}_{\cG^{\boxtimes I}_Z}} & {\mathbf{Sat}^q_{X^I},}
	\arrow["{\mathbf{S}^{\QLisse,\sigma}_I}", from=1-1, to=1-2]
	\arrow["\iota^B_I", from=2-1, to=1-1]
	\arrow["{\mathbf{S}^{\Ab}_I}"', from=2-1, to=2-2]
	\arrow["\iota^A_I"', from=2-2, to=1-2]
\end{tikzcd}
    \end{equation}
    where the horizontal arrows $\iota^A_I$ and $\iota^B_I$ are inclusion maps.
\end{thm}

Thanks to Lemma \ref{lem: the characterization lemma}, the above assertion uniquely characterizes $\textbf{S}^{\QLisse,\sigma}_I$ as a $\mathbf{QLisse}^{\otimes I}$-linear functor. Since $\cS^{\Ab}_I$ does not depend on $\sigma$ we conclude:
\begin{cor}
    The $\mathbf{QLisse}^{\otimes I}$-linear functor $\mathbf{S}^{\QLisse,\sigma}_I$ does not depend on the choice of $\sigma$.
\end{cor}

\sssec{}\label{sssec: construction of S QL I}
Hence, we drop $\sigma$ from the notation and just denote the functor as 
\[\textbf{S}^{\QLisse}_I:\left(\mathbf{Rep}(H)^{\QLisse}_{\cG_Z}\right)^{\otimes I}\rightarrow \mathbf{Shv}^q(\Hecke^{\loc}_{I}).\]
Using the $\mathbf{LS}^{\otimes I}$-linearity on the target, we can tensor up the source along $\mathbf{QLisse}^{\otimes I}\rightarrow \mathbf{LS}^{\otimes I}$ to get a $\mathbf{LS}^{\otimes I}$-linear functor
\[\textbf{S}_I:\left(\mathbf{Rep}(H)^{\loc}_{\cG_Z}\right)^{\otimes I}\rightarrow \mathbf{Shv}^q(\Hecke^{\loc}_{I}).\]

\sssec{}
For the rest of this subsection, we prove Theorem \ref{thm: characterization for I}. Since its proof is a bit involved, we give an overview here. 
\medskip

From the compatibility of perverse ULA version of Satake with external convolution product, we have a diagram \footnote{Strictly speaking, in this diagram and the discussion to follow, the products of abelian categories should be replaced by $\Ectensbig$, but as in Sect. \ref{nota: introducing oneCatEc} we just take $\times$ remembering the obvious $\Vect^{\heartsuit,\on{f.d.}}$-multilinearity.}
\begin{equation}\label{eq: big diagram characterization}
    \begin{tikzcd}
	& {\left(\mathbf{Rep}(H)^{\QLisse}_{\cG_Z}\right)^{\otimes I}} & {\textbf{Shv}^{q}(\Hecke^{\loc}_I)} & \\
	{\mathbf{Rep}(H^I)^{\Ab,\ft}_{\cG^{\boxtimes I}_Z}} &&& {\mathbf{Sat}^q_{I}} \\
	& {\left(\mathbf{Rep}(H)^{\Ab,\ft}_{\cG_Z}\right)^{\times I}} & {\left(\mathbf{Sat}^q_{\{*\}}\right)^{\times I}.}
	\arrow["{\mathbf{S}^{\QLisse,\sigma}_I}", from=1-2, to=1-3]
	\arrow["{\iota^B_I}", from=2-1, to=1-2]
	\arrow["{\mathbf{S}^{\Ab}_I}", from=2-1, to=2-4]
	\arrow["{\iota^A_I}"', from=2-4, to=1-3]
	\arrow["{\iota^B}", from=3-2, to=1-2]
	\arrow["\boxtimes", from=3-2, to=2-1]
	\arrow["{\left(\mathbf{S}^{\Ab}_{\{*\}}\right)^{\times I}}"', from=3-2, to=3-3]
	\arrow["{\iota^A}"', from=3-3, to=1-3]
	\arrow["{\boxtimes^f}"', from=3-3, to=2-4]
\end{tikzcd}
\end{equation}
Let us summarize what we know about this diagram:
\begin{enumerate}
    \item both the left and right triangle commute;
    \item the lower square commutes from the compatibility of perverse ULA version of Satake with external convolution product;
    \item and the pentagon that only excludes the left-most vertex $\mathbf{Rep}(H^I)^{\Ab,\ft}_{\cG^{\boxtimes I}_Z}$ commutes from the construction of $\mathbf{S}^{\QLisse,\sigma}_I$;
\end{enumerate}
and all these commute as $(\mathbf{LS}^{\Ab,\ft})^{\times I}$-modules.
To prove Theorem \ref{thm: characterization for I}, we have to provide a commutativity of the upper square as $\mathbf{LS}^{\Ab,\ft}_I$-modules.
\medskip

We give a bit imprecise argument to convey the idea of the proof and comment on how to turn it into an actual proof. The idea of the proof is that the arrow $\boxtimes$ in the triangle on the left almost looks like the operation of tensoring up along the symmetric monoidal map
\begin{equation}\label{eq: boxtimes LSAbft}
    \boxtimes:(\mathbf{LS}^{\Ab,\ft})^{\times I}\rightarrow \mathbf{LS}^{\Ab,\ft}_I.
\end{equation}
From this point of view, the $\mathbf{LS}^{\Ab,\ft}_I$-linear commutativity of the upper square is equivalent to the $(\mathbf{LS}^{\Ab,\ft})^{\times I}$-linear commutativity of the pentagon in (3) above, because of the induction-forgetful adjunction
\begin{equation}\label{eq: induction forgetful adj}
(\mathbf{LS}^{\Ab,\ft})^{\times I}\lMod\leftrightarrows \mathbf{LS}^{\Ab,\ft}_I\lMod.
\end{equation}

However, this argument has several flaws. First of all, it is not true that the arrow $\boxtimes$ in the triangle on the left is the operation of tensoring up along \eqref{eq: boxtimes LSAbft}. Indeed, it is not true that $\Ectens_{i\in I}\Rep(H)^{\Ab,\ft}$ is equivalent to $\Rep(H^I)^{\Ab,\ft}$. For such a statement to work, we should work in the setting of presentable categories. Therefore, it is necessary to $\Ind$-complete the abelian categories in \eqref{eq: big diagram characterization} and replace the products, e.g., $(\mathbf
{Rep}(H)^{\Ab,\ft}_{\cG_Z})^{\times I}$, by appropriate Lurie tensor products of presentable categories.

\medskip
Second, even after $\Ind$-completing abelian categories in the diagram \eqref{eq: big diagram characterization}, the resulting arrows from abelian categories to derived categories, e.g., $\iota^B$ and $\iota^B_I$, do \emph{not} preserve colimits. Indeed, those arrows are similar to an inclusion of the abelian heart into a DG category equipped with a $t$-structure, and such an inclusion does not preserve colimits in general. The failure of colimit-preservation prevents us from using higher algebra in presentable categories; in particular, we can no longer use the induction–forgetful adjunction mentioned above. To remedy this, we prove that the composition $\mathbf{S}^{\QLisse,\sigma}_I\circ \iota^B_I$ factors through $\iota^A_I$, so that we can work only with abelian categories and exact functors between them. Note that an exact functor between abelian categories turns into a colimit-preserving functor after $\Ind$-completion. 

\sssec{}
Now we discuss the actual proof.
We remedy the second issue first:
\begin{lemma}\label{lem: landing in Sat}
    The composition $\mathbf{S}^{\QLisse,\sigma}_I\circ \iota^B_I$ lands in the subsheaf $\mathbf{Sat}^q_I$, and hence factors through $\iota^A_I$.
\end{lemma}
\begin{proof}
    It suffices to verify this assertion sectionwise for $X'\in X_{\et}$ where $\cG_Z$ gets trivializable.
    The subcategory $\Sat^q_{X'^I}\subseteq \Shv^q(\Hecke^{\loc}_{X'^I})$ is closed under extensions. (Corollary \ref{cor: Sat full compatibility}) Therefore, from Lemma \ref{lem: Rep H Ab ft structure} it suffices to show that objects in $\Rep(H^I)^{\Ab,\ft}_{\cG^{\boxtimes I}_Z}(X'^I)$ that lie in the essential image of $\boxtimes$ (the upper-left arrow in \eqref{eq: big diagram characterization}) are sent to the full subcategory $\Sat^q_{X'^I}$ under $\mathbf{S}^{\QLisse,\sigma}_I\circ \iota^B_I$. This is obvious, because the hexagon consisting of outer circuits of \eqref{eq: big diagram characterization} commutes. 
\end{proof}

\sssec{}
From Lemma \ref{lem: landing in Sat}, we get a canonical $\mathbf{LS}^{\Ab,\ft}_I$-linear functor 
\begin{equation}\label{eq: Psi sigma I}
    \Psi^{\sigma}_I: \mathbf{Rep}(H^I)^{\Ab,\ft}_{\cG^{\boxtimes I}_Z}\rightarrow \mathbf{Sat}^q_I.
\end{equation}
We have to prove that this canonically agrees with $\mathbf{S}^{\Ab}_I$ as $\mathbf{LS}^{\Ab,\ft}_I$-linear functors.  In the diagram
\begin{equation}\label{eq: smaller square characterization}
    \begin{tikzcd}
	{\mathbf{Rep}(H^I)^{\Ab,\ft}_{\cG^{\boxtimes I}_Z}} & {\mathbf{Sat}^q_{I}} \\
	{\left(\mathbf{Rep}(H)^{\Ab,\ft}_{\cG_Z}\right)^{\times I}} & {\left(\mathbf{Sat}^q_{\{*\}}\right)^{\times I}}
	\arrow["{\mathbf{S}^{\Ab}_I}"', from=1-1, to=1-2]
	\arrow["{\Psi^{\sigma}_I}", shift left=3, from=1-1, to=1-2]
	\arrow["\boxtimes", from=2-1, to=1-1]
	\arrow["{\left(\mathbf{S}^{\Ab}_{\{*\}}\right)^{\times I}}"', from=2-1, to=2-2]
	\arrow["{\boxtimes^f}"', from=2-2, to=1-2]
\end{tikzcd}
\end{equation}
we can make squares in two different ways - one with $\mathbf{S}^{\QLisse,\sigma}_I$ and the other with $\Psi^{\sigma}_I$ - what we know for now is that they both commute as $(\mathbf{LS}^{\Ab,\ft})^{\times I}$-modules.

\sssec{}\label{sssec: Ind completion Ec}
Now we address the first issue by passing to $\Ind$-completions. Consider the category $\Pr^L$ of presentable categories with colimit-preserving functors between them (\cite[Definition 5.5.3.1]{HTT}) and let $\Pr^L_{E^{\Ab}}:=\Vect^{\heartsuit}_E\lMod(\Pr^L)$. We have a functor
\begin{equation}\label{eq: Ind completion for oneCatEc}
    \Ind: \CatEc^{\omega}\rightarrow \PrL_{E^{\Ab}},\quad C\mapsto \Ind(C)
\end{equation}
where $\CatEc^{\omega}$ is the 1-full subcategory of $\CatEc$ consisting of $\Vect^{\heartsuit,\on{f.d.}}$-modules $C$ which has finite colimits and retracts, whose action $\Vect^{\heartsuit,\on{f.d.}}\times C\rightarrow C$ is right-exact on each variables, and only right-exact functors between module categories are allowed. Any $E$-linear abelian category $A$ lies in $\CatEc^{\omega}$, and $\Ind(A)$ is a $E$-linear Grothendieck abelian category. 
\medskip

Given a $\Vect^{\heartsuit,\on{f.d.}}$-multilinear functor $C_1\times \cdots \times C_r\rightarrow C$, we say that this functor is \textbf{multi right-exact} (resp. \textbf{multi exact}) if this is right-exact (resp. exact) on each variable. We record the following trivial observation: 

\begin{lemma}\label{lem: functoriality of Ind}
    A $\Vect^{\heartsuit,\on{f.d.}}$-multilinear functor $C_1\times \cdots \times C_r\rightarrow C$ that is multi right-exact induces a map
\[\Ind(C_1)\otimes\cdots \otimes \Ind(C_r)\rightarrow\Ind(C)\quad\te{in }\PrL_{E^{\Ab}}.\]
\end{lemma}

\sssec{Ind-finite abelian local systems}\label{sssec: ind-finite ab loc}
We list two main examples of Grothendieck abelian categories.
\begin{enumerate}
    \item For $Y\in\Schft$ consider the category $\LS^{\Ab,\indft}(Y):=\Ind(\LS^{\Ab,\ft}(Y))$.
    \item For an algebraic group $\mathsf{H}$ over $E$, we take the heart $\Rep(\mathsf{H})^{\heartsuit}$ of $\Rep(\mathsf{H})$.
\end{enumerate}
Both categories belong to $\Pr^L_{E^{\Ab}}$, and in fact they belong to $\CAlg(\PrL_{E^{\Ab}})$ thanks to Lemma \ref{lem: functoriality of Ind}. Indeed, they are $\Ind$-completions of $\LS^{\Ab,\ft}(Y)$ and $\Rep(\mathsf{H})^{\heartsuit,\on{f.d.}}$ respectively, and their monoidal structures are multi exact.
Hence we can take a tensor product
\begin{equation}\label{eq: H representation on indft}
    \Rep(\mathsf{H})^{\heartsuit}\otimes \LS^{\Ab,\indft}(Y) \quad \te{in } \CAlg(\on{Pr}^L_{E^{\Ab}}),
\end{equation}
which gets identified with the category of $\cO_{\mathsf{H}}$-comodules in $\LS^{\Ab,\indft}(Y)$. From the latter description, the category \eqref{eq: H representation on indft} is the $\Ind$-completion of $\Rep(H)^{\Ab,\loc,\ft}(Y)$.
\medskip

The assignment
\[Y\mapsto \LS^{\Ab,\indft}(Y),\quad (\Schft)^{\op}\rightarrow \on{Pr}^L_{E^{\Ab}}\]
satisfies descent, and therefore we have sheaves
\[\mathbf{LS}^{\Ab,\indft}: X'\mapsto \LS^{\Ab,\indft}(X') \quad \te{and}\quad \mathbf{LS}^{\Ab,\indft}_I: X'\mapsto \LS^{\Ab,\indft}(X'^I)\]
in $\CAlg(\Shv(X_{\et};\Pr^L_{E^{\Ab}}))$. From the dualizability of $\Rep(H)^{\heartsuit}$ and $\Rep(H^I)^{\heartsuit}$ we have sheaves
\begin{equation}\label{eq: Rep H indft}
    \textbf{Rep}(H)^{\Ab,\loc,\indft}: X'\mapsto \Rep(H)^{\heartsuit}\otimes \LS^{\Ab,\indft}(X')
\end{equation}
and
\begin{equation}\label{eq: Rep H indft I}
    \textbf{Rep}(H^I)^{\Ab,\loc,\indft}: X'\mapsto \Rep(H^I)^{\heartsuit}\otimes \LS^{\Ab,\indft}(X'^I)
\end{equation}
in $\CAlg(\Shv(X_{\et};\Pr^L_{E^{\Ab}}))$.
\medskip

Both $\LS^{\Ab,\indft}(X')$ and $\Rep(H)^{\heartsuit}$ are compactly generated and hence dualizable, so the assignment
\[X'\mapsto \Rep(H^I)^{\heartsuit}\otimes\LS^{\Ab,\indft}(X')^{\otimes I}\]
still satisfies \'etale descent. Using $(\Rep(H)^{\heartsuit})^{\otimes I}=\Rep(H^I)^{\heartsuit}$, this gets identified with the $I$-fold tensor product of $\mathbf{Rep}(H)^{\Ab,\loc,\indft}$ taken in $\CAlg(\Shv(X_{\et};\Pr^L_{E^{\Ab}}))$.
The external tensor product of local systems induces \footnote{We used Lemma \ref{lem: functoriality of Ind} here.} a map
\begin{equation}\label{eq: external product indft}
    (\mathbf{Rep}(H)^{\Ab,\loc,\indft})^{\otimes I}\xrightarrow{\boxtimes} \mathbf{Rep}(H^I)^{\Ab,\loc,\indft}\quad \te{in }\CAlg(\Shv(X_{\et};\on{Pr}^L_{E^{\Ab}})).
\end{equation}

The same argument goes through when $H$ is trivial, in which case we identify the assignment $X'\mapsto \left(\LS^{\Ab,\indft}(X')\right)^{\otimes I}$ with the $I$-fold tensor product of $\mathbf{LS}^{\Ab,\indft}$ and we have a map
\begin{equation}\label{eq: ext tensor prod LS indft}
    (\mathbf{LS}^{\Ab,\indft})^{\otimes I}\rightarrow \mathbf{LS}^{\Ab,\indft}_I \quad \te{in }\CAlg(\Shv(X_{\et};\on{\Pr}^L_{E^{\Ab}})).
\end{equation}

The sheaf \eqref{eq: Rep H indft} (resp. \eqref{eq: Rep H indft I}) admits a natural symmetric monoidal twist by a $Z_H(E)^{\on{tors}}$ (resp. $Z_{H^I}(E)^{\on{tors}}$) - gerbe, and hence the gerbe $\cG_Z$ defines the symmetric monoidal twists
\[\textbf{Rep}(H)^{\Ab,\loc,\indft}_{\cG_Z}\quad\te{and}\quad \textbf{Rep}(H^I)^{\Ab,\loc,\indft}_{\cG^{\boxtimes I}_Z}\quad \te{in }\CAlg(\Shv(X_{\et};\on{Pr}^L_{E^{\Ab}})).\]
Since the twist by $\cG_Z$ intertwines both sides of \eqref{eq: external product indft} we get a map
\begin{equation}\label{eq: external product indft twisted}
    (\mathbf{Rep}(H)^{\Ab,\loc,\indft}_{\cG_Z})^{\otimes I}\xrightarrow{\boxtimes} \mathbf{Rep}(H^I)^{\Ab,\loc,\indft}_{\cG^{\boxtimes I}_Z}\quad \te{in }\CAlg(\Shv(X_{\et};\on{Pr}^L_{E^{\Ab}})).
\end{equation}

\sssec{$\Ind$-completion of a sheaf}
 Given a sheaf $\cF$ valued in $\CatEc$ whose restriction functors are right-exact, we define a presheaf $\Ind^{\on{pre}}(\cF)$ valued in $\PrL_{E^{\Ab}}$ by taking $\Ind$-completion \eqref{eq: Ind completion for oneCatEc}, and we sheafify to get the sheaf $\Ind(\cF)$ valued in $\PrL_{E^{\Ab}}$. Since \'etale pullbacks are perverse $t$-exact, the sheaf $\mathbf{Sat}^q_I$ is valued in $\CatEc$ meets the condition above, and therefore we get $\Ind(\mathbf{Sat}^q_I)$ by the above procedure. Some sheaf of categories introduced previously are related to each other in terms of $\Ind$-completion:
\[\mathbf{LS}^{\Ab,\indft}=\Ind(\mathbf{LS}^{\Ab,\ft}),\quad \mathbf{LS}^{\Ab,\indft}_I=\Ind(\mathbf{LS}^{\Ab,\ft}_I),\quad \Ind\left(\mathbf{Rep}(H^I)^{\Ab,\ft}_{\cG^{\boxtimes I}_Z}\right) =\mathbf{Rep}(H^I)^{\Ab,\indft}_{\cG^{\boxtimes I}_Z}.\]

\sssec{}\label{sssec: sheafy version of multi right exact}
Such an $\Ind$-completion for sheaves enjoys an obvious analog of Lemma \ref{lem: functoriality of Ind}. Namely, if we have sheaves $\cF_i$'s and $\cF$ valued in $\CatEc$ whose restriction maps are right-exact, and a $\Vect^{\heartsuit,\on{f.d.}}$-multilinear map $\cF_1\times\cdots\times \cF_r\rightarrow \cF$ between sheaves that is section-wise multi right-exact, i.e., for any $X'\in X_{\et}$ the multilinear map
\[\cF_1(X')\times\cdots\times \cF_r(X')\rightarrow \cF(X')\]
is multi right-exact, then we have a canonical map
\[\Ind(\cF_1)\otimes\cdots\otimes\Ind(\cF_r)\rightarrow\Ind(\cF)\quad \te{in }\Shv(X_{\et};\PrL_{E^{\Ab}}).\]

\sssec{Passing to $\Ind$-completions}
We claim that using Lemma \ref{lem: functoriality of Ind} we can pass from \eqref{eq: smaller square characterization} to a $(\mathbf{LS}^{\Ab,\indft})^{\otimes I}$-linear diagram
\begin{equation}\label{eq: smaller diagram ind completion}
    \begin{tikzcd}[column sep=6em]
	{\mathbf{Rep}(H^I)^{\Ab,\indft}_{\cG^{\boxtimes I}_Z}} & {\Ind(\mathbf{Sat}^q_{I})} \\
	{\left(\mathbf{Rep}(H)^{\Ab,\indft}_{\cG_Z}\right)^{\otimes I}} & {\left(\Ind(\mathbf{Sat}^q_{\{*\}})\right)^{\otimes I}}
	\arrow["{\Ind(\mathbf{S}^{\Ab}_I)}"', from=1-1, to=1-2]
	\arrow["{\Ind(\Psi^{\sigma}_I)}", shift left=3, from=1-1, to=1-2]
	\arrow["\Ind(\boxtimes)", from=2-1, to=1-1]
	\arrow["{\left(\Ind(\mathbf{S}^{\Ab}_{\{*\}})\right)^{\otimes I}}"', from=2-1, to=2-2]
	\arrow["{\Ind(\boxtimes^f)}"', from=2-2, to=1-2]
\end{tikzcd}
\end{equation}
in $\Shv(X_{\et};\PrL_{E^{\Ab}})$ such that two squares commute. 

\sssec{Constructing linearity}\label{sssec: constructing linearity}
We explain the $(\mathbf{LS}^{\Ab,\indft})^{\otimes I}$-linearity on sheaves and maps in \eqref{eq: smaller diagram ind completion}. Let $X^{\triv}_{\et}$ be the full subsite of $X_{\et}$ spanned by $X'$'s on which $\cG_Z$ gets trivializable. By Sect. \ref{ssec: restriction to triv} it is enough to construct linearities over $X^{\triv}_{\et}$.
\medskip

We first construct linearity on sheaves. Using the techniques developed in Sect. \ref{ssec: full compatibility with extensions} we can easily show that the $(\mathbf{LS}^{\Ab,\ft})^{\times I}$-action on each sheaf appearing in the diagram \eqref{eq: smaller square characterization} is section-wise multi right-exact \footnote{By an action being multi right-exact, we mean that all the action maps involved are multi right-exact.} (in fact, they are section-wise multi exact). For example, to show that the action map
\[\LS^{\Ab,\ft}(X')\times \Rep(H)^{\Ab,\loc,\ft}_{\cG_Z}(X')\rightarrow \Rep(H)^{\Ab,\loc,\ft}_{\cG_Z}(X')\]
is exact on each variables, we apply Lemma \ref{lem: induced exactness} to the situation \footnote{That $\Rep(H)^{\Ab,\loc,\ft}_{\cG_Z}(X')\hookrightarrow \Rep(H)^{\loc}_{\cG_Z}(X')$ is fully compatible with extensions follows from Example \ref{ex: t str full compatibility}. The twist by $\cG_Z$ does not affect the full compatibility because $\cG_Z$ is trivializable on $X'$.}
\[\begin{tikzcd}
	{\LS(X')\otimes \Rep(H)^{\loc}_{\cG_Z}(X')} & {\Rep(H)^{\loc}_{\cG_Z}(X')} \\
	{\LS^{\Ab,\ft}(X')\times \Rep(H)^{\Ab,\loc,\ft}_{\cG_Z}(X')} & {\Rep(H)^{\Ab,\loc,\ft}_{\cG_Z}(X').}
	\arrow[from=1-1, to=1-2]
	\arrow[from=2-1, to=1-1]
	\arrow[from=2-1, to=2-2]
	\arrow[from=2-2, to=1-2]
\end{tikzcd}\]
Therefore, using Sect. \ref{sssec: sheafy version of multi right exact} we get natural $(\mathbf{LS}^{\Ab,\indft})^{\otimes I}$ linearities on each sheaves.
\medskip

To obtain linearity for functors, it suffices to show that all the functors in the diagram \eqref{eq: smaller square characterization} are section-wise multi right-exact (in fact, they are section-wise multi exact). For vertical arrows, it can be done similarly as before using Lemma \ref{lem: induced exactness} and Corollary \ref{cor: Sat full compatibility}. Since the perverse ULA Satake equivalence is exact, it only remains to show the exactness of $\Psi^{\sigma}_I$. Since $\cG_Z$ is trivializable on $X'$, Lemma \ref{lem: Rep H Abft QL full compatibility} tells us that the inclusion
\[\Rep(H^I)^{\Ab,\loc,\ft}_{\cG^{\boxtimes I}_Z}(X'^I)\hookrightarrow \left(\Rep(H)^{\QLisse}_{\cG_Z}(X')\right)^{\otimes I}\]
is fully compatible with extensions, and hence $\iota^B_I$ sends a short exact sequence in $\Rep(H^I)^{\Ab,\loc,\ft}_{\cG^{\boxtimes I}_Z}(X'^I)$ is sent to a cofiber sequence in $\left(\Rep(H)^{\QLisse}_{\cG_Z}(X')\right)^{\otimes I}$. Therefore, $\cS^{\QLisse,\sigma}_I\circ \iota^B_I$ sends a short exact sequence in $\Rep(H^I)^{\Ab,\loc,\ft}_{\cG^{\boxtimes I}_Z}(X'^I)$ to a cofiber sequence in $\Shv^q(\Hecke^{\loc}_{X'^I})$. Since this lands in $\Sat^q_{X'^I}$ (Lemma \ref{lem: landing in Sat}), using Corollary \ref{cor: Sat full compatibility} we conclude that this cofiber sequence is a short exact sequence in $\Sat^q_{X'^I}$. In conclusion, $\Psi^{\sigma}_I$ sends a short exact sequence to a short exact sequence, hence the proof of exactness.

\begin{rmk}
    One readily checks that the map $\Ind(\boxtimes)$ in the diagram \eqref{eq: smaller diagram ind completion} is canonically identified with \eqref{eq: external product indft twisted}.
\end{rmk}

 \sssec{}
 Since all sheaves and maps together with linearity arose from \eqref{eq: smaller square characterization} by means of $\Ind$-completion, the commutativity of \eqref{eq: smaller square characterization} as $(\mathbf{LS}^{\Ab,\ft})^{\times I}$-modules gives rise to a commutativity of (both squares of) \eqref{eq: smaller diagram ind completion} as $(\mathbf{LS}^{\Ab,\indft})^{\otimes I}$-modules.

\sssec{Induction-forgetful adjunction}
In the diagram \eqref{eq: smaller square characterization}, the arrows $\Psi^{\sigma}_I$ and  $\mathbf{S}^{\Ab}_I$ have natural $\mathbf{LS}^{\Ab,\ft}_I$-linearity lifting the $(\mathbf{LS}^{\Ab,\ft})^{\times I}$-linearity along the symmetric monoidal map
\[(\mathbf{LS}^{\Ab,\ft})^{\times I}\rightarrow \mathbf{LS}^{\Ab,\ft}_I\]
of box tensor product. Passing to $\Ind$-completions we get the symmetric monoidal map \eqref{eq: ext tensor prod LS indft}, and the action of $\mathbf{LS}^{\Ab,\ft}_I$ on both $\mathbf{Rep}(H^I)^{\Ab,\ft}_{\cG^{\boxtimes I}_Z}$ and $\mathbf{Sat}^q_I$ are sectionwise multi right-exact \footnote{This can be done as in Sect. \ref{sssec: constructing linearity}, using the techniques developed in Sect. \ref{ssec: full compatibility with extensions}.}. Combining these, we conclude that the arrows $\Psi^{\sigma}_I$ and $\mathbf{S}^{\Ab}_I$ in the $\Ind$-completed diagram \eqref{eq: smaller diagram ind completion} have natural $\mathbf{LS}^{\Ab,\indft}_I$-linearity lifting the $(\mathbf{LS}^{\Ab,\indft})^{\otimes I}$-linearity. We finally observe that we have an equivalence
\[\mathbf{Rep}(H^I)^{\Ab,\indft}_{\cG^{\boxtimes I}_Z}= \left(\mathbf{Rep}(H)^{\Ab,\indft}_{\cG_Z}\right)^{\otimes I}\bigotimes_{(\mathbf{LS}^{\Ab,\indft})^{\otimes I}} \mathbf{LS}^{\Ab,\indft}_I\]
and hence the box tensor product map
\[\left(\mathbf{Rep}(H)^{\Ab,\indft}_{\cG_Z}\right)^{\otimes I}\rightarrow \mathbf{Rep}(H^I)^{\Ab,\indft}_{\cG^{\boxtimes I}_Z}\]
is the unit map of the induction-forgetful adjunction\footnote{Both module categories are internal to $\PrL_{E^{\Ab}}$.}
\[(\mathbf{LS}^{\Ab,\indft})^{\otimes I}\lMod\leftrightarrows \mathbf{LS}^{\Ab,\indft}_I\lMod.\]
Therefore, both the $\mathbf{LS}^{\Ab,\indft}_I$-linear maps $\Ind(\Psi^{\sigma}_I)$ and $\Ind(\mathbf{S}^{\Ab}_I)$ in diagram \eqref{eq: smaller diagram ind completion} are uniquely determined as the adjoint to the $(\mathbf{LS}^{\Ab,\indft})^{\otimes I}$-linear map $\Ind(\boxtimes^f)\circ (\Ind(\mathbf{S}^{\Ab}_{*}))^{\otimes I}$, and hence we have a canonical identification between them. Restricting this equivalence to the subsheaf $\mathbf{Rep}(H^I)^{\Ab,\loc,\ft}_{\cG^{\boxtimes I}_Z}$ and $\mathbf{LS}^{\Ab,\ft}_I$-linearity, we prove Theorem \ref{thm: characterization for I}. \qed[Theorem \ref{thm: characterization for I}]

\section{Factorizable Satake as a natural transformation}\label{sec: Satake as a nat transf}
As mentioned in the introduction to Sect. \ref{sec: metaplectic Satake functor}, we continue to construct the `hard' part of the natural transformation. In Sect. \ref{sec: metaplectic Satake functor} we have obtained a family of functors $\{\cS_I\}_{I\in\fSet}$ where each $\cS_I$ is a $\LS(X^I)$-linear functor
\[\left(\Rep(H)^{\loc}_{\cG_Z}(X)\right)^{\otimes I}\rightarrow \Shv^q(\Hecke^{\loc}_{X^I}).\]
The goal of this section is to enhance this family of functors to a natural transformation between functors of the form
\begin{equation}\label{eq: fSet to fSet Alg}
    \fSet\rightarrow \fSet + \Alg(\LS^{\otimes}\lMod),
\end{equation}
where the target category consists of pairs $[I,\bA]$ where $I\in\fSet$ and $\bA\in\Alg(\LS(X)^{\otimes I}\lMod)$.

\sssec{}\label{sssec: Alg LS construction}
Suppose that we are given a diagram
\[\cC: \fSet\rightarrow \DGCatSymMon, \quad I\mapsto \cC_I.\]
We can consider a category $\fSet + \Alg(\cC\lMod)$ consisting of pairs $[I,\bA]$ where $\bA\in\Alg(\cC_I\lMod)$. Given a map $f:I\rightarrow J$ in $\fSet$ we have a symmetric monoidal $\cC_I\rightarrow \cC_J$ and hence a restriction map $\Alg(\cC_J\lMod)\rightarrow \Alg(\cC_I\lMod)$. In the category $\fSet + \Alg(\cC\lMod)$, a map between objects $[I,\bA_I]\rightarrow [J,\bA_J]$ is a data of a map $f:I\rightarrow J$ in $\fSet$ together with a map $\bA_I\rightarrow \bA_J$ in $\Alg(\cC_I\lMod)$. This construction can be done more rigorously using the language of coCartesian fibration, as we did in Sect. \ref{sssec: fSet LSMon}. 
\medskip

As in Sect. \ref{sssec: respects the projection}, we say that a functor $\fSet\rightarrow \fSet+\Alg(\cC\lMod)$ respects the projection if it forgets to the identity functor of $\fSet$. We use the notation
\[\Fun^{\circ}(\fSet,\fSet+\Alg(\cC\lMod)).\]
for the category of such functors. 
\medskip

Applying this construction to the diagram
\[\LS(X)^{\otimes}: \fSet\rightarrow \DGCatSymMon,\quad I\mapsto \LS(X)^{\otimes I},\]
 we get the category $\fSet+\Alg(\LS(X)^{\otimes}\lMod)$, which has appeared as the target of \eqref{eq: fSet to fSet Alg}. 
The main result of this section is:
\begin{thm}\label{thm: fact Satake nat transf}
    Assume either $X\ne \pjl$ or that $X=\pjl$ and $\cG_Z$ is trivializable. The collection of functors $\{\cS_I\}_{I\in\fSet}$ admits a natural enhancement into a natural transformation between two functors 
\[\fSet\rightarrow \fSet+ \Alg(\LS(X)^{\otimes}\lMod)\]
that respect the projection:
the first functor maps
\begin{equation}\label{eq: Satake functorial B}
    I\mapsto \left(\Rep(H)^{\loc}_{\cG_Z}(X)\right)^{\otimes I}
\end{equation}
and the second functor maps
\begin{equation}\label{eq: Satake functorial A}
    I\mapsto \Shv^q(\Hecke^{\loc}_{X^I}),
\end{equation}
and the natural transformation goes from \eqref{eq: Satake functorial B} to \eqref{eq: Satake functorial A}.
\end{thm}

\sssec{}
 The construction of each functor \eqref{eq: Satake functorial A} and \eqref{eq: Satake functorial B} is also a subject matter of this section. While the first functor is easy to construct, some thought is needed to construct the second functor in a homotopy-coherent manner. This issue of homotopy coherence gets even worse when we try to construct the natural transformation.
\medskip

To overcome this difficulty, we introduce an operad $\Conv$ that concisely encodes both the internal and external convolution structures of the collection $\{\Shv^q(\Hecke^{\loc}_{X^I})\}_{I\in\fSet}$. (Sect. \ref{ssec: an operad Conv}) In fact, we show that this operad admits a more familiar description: a $\Conv$-algebra object is the same as a functorial assignment of associative algebras to finite sets. $\Conv$ is an ordinary colored operad whose colors are finite sets, and we construct a $\Conv$-algebra object in $\DGcat$ whose evaluation at the color $I$ is $\Shv^q(\Hecke^{\loc}_{X^I})$. Moreover, we endow this $\Conv$-algebra with an action of a $\Conv$-algebra $\LS(X)^{\otimes}$ arising from the symmetric monoidal structure of $\LS(X)$. The analogous enhancement and module structure on the spectral side is constructed easily. These are done in Sect. \ref{ssec: Conv algebras from fact Satake}. After that, we provide a mechanism to pass from $\LS(X)^{\otimes}$-modules in $\Conv(\DGcat)$ to functors in a form \eqref{eq: fSet to fSet Alg}. (Sect. \ref{ssec: construction of a nat transf}) This mechanism lends to the following key idea: the desired natural transform is a shadow of a $\LS(X)^{\otimes}$-module map between $\Conv$-algebras enhancing the $\LS$-linear Satake functors $\cS_{X,I}$'s.
\medskip

In Sect. \ref{ssec: Perv ULA Satake in Conv language} we obtain the required $\LS(X)^{\otimes}$-module map, hence completing the proof of Theorem \ref{thm: fact Satake nat transf}. The idea is to use perverse ULA Satake equivalence. We reformulate the perverse ULA Satake equivalence in the language of $\Conv$-algebras and conclude from compatibilities between perverse ULA version and $\QLisse$-linear version of Satake functor (Sect. \ref{ssec: LS linear general}) together with characterization Lemmas (Sect. \ref{ssec: characterization lemmas}). The idea is that the homotopy coherence needed to enhance $\cS_{I}$ to a $\Conv$-algebra map is uniquely determined by the homotopy coherence built into the enhancement of $\cS^{\Ab}_I$ to a $\Conv$-algebra map.

\subsection{An operad $\Conv$}\label{ssec: an operad Conv}
We introduce an operad $\Conv$ with two goals in mind:
\begin{enumerate}
    \item to concisely encode convolution structures (both internal and external) and functoriality of the assignment $I\mapsto \Shv^q(\Hecke^{\loc}_{X^I})$ and compatibility between those at once;
    \item and to effectively deal with compatibilities needed for the natural transformation.
\end{enumerate}
This is an ordinary colored operad, and we use the notations of \cite[Definition 2.1.1.1]{HA} to describe it.
\begin{defn}\label{def: operad Conv}
    An operad $\Conv$ is an ordinary colored operad described as follows.
    \begin{enumerate}
        \item The collection of colors is the collection of finite sets.
        \item Suppose we are given a finite set $U$, an $U$-indexed collection of colors $\{I_u\}_{u\in U}$, and a color $J$, the set of morphisms $\Mul_{\Conv}(\{I_u\}_{u\in U}, J)$ is the set of pairs $(\sigma, \{f_{u}\}_{u\in U})$ where $\sigma$ is a total ordering of $U$ and each $f_u$ is a map $I_u\rightarrow J$.
        \item For each color $I$ we have a distinguished element
        $\id_I\in \Mul_{\Conv}(\{I\},I)$,
        whose data $(\sigma,f)$ is the unique total ordering of a singleton together with $\id_I:I\rightarrow I$.
        \item Suppose we have a map of finite sets $\phi:U\rightarrow V$ having fibers $\{U_v\}_{v\in V}$, a collection of colors $\{I_u\}_{u\in U}$, a collection of colors $\{J_v\}_{v\in V}$, and a color $K$. The composition map
        \begin{equation}\label{eq: composition Conv}
            \prod_{v\in V} \Mul_{\Conv}(\{I_u\}_{u\in U_v}, J_v)\times \Mul_{\Conv}(\{J_v\}_{v\in V}, K)\rightarrow \Mul_{\Conv}(\{I_u\}_{u\in U}, K)
        \end{equation}
        is defined as follows. Suppose for each $v\in V$ we have 
        \[(\sigma^v, \{f^v_u\}_{u\in U_v})\in \Mul_{\Conv}(\{I_u\}_{u\in U_v}, J_v)\]
        and
        $(\sigma^0, \{f_v\})\in \Mul_{\Conv}(\{J_v\}_{v\in V}, K)$.
        Then the resulting
        $(\sigma,\{f_u\})\in \Mul_{\Conv}(\{I_u\}_{u\in U}, K)$
        is defined so that each map $f_u: I_u\rightarrow K$ is the composition
        \[I_u\xrightarrow{f^{\phi(u)}_u}J_{\phi(u)}\xrightarrow{f_{\phi(u)}}K,\]
        and the total ordering $\sigma$ of $U$ is determined lexicographically: given $u_1,u_2\in U$, the order relation is defined as:
        \begin{equation}\label{eq: lexicographic}
        u_1\stackrel{\sigma}{\le} u_2\iff \phi(u_1)\stackrel{\sigma^0}{<}\phi(u_2)\te{ or } \left[\phi(u_1)=\phi(u_2)=v\te{ and } u_1\stackrel{\sigma^v}{\le} u_2\right].
        \end{equation}
    \end{enumerate}
\end{defn}
\begin{ex}
    We give an example of composition of operations. Consider the example $U=\{1,2,3\}$ and $V=\{1,2\}$ with a map $\phi:U\rightarrow V$ mapping $1,2\mapsto 1$ and $3\mapsto 2$. Consider the finite sets $I_1$, $I_2$, $I_3$, $J_1$, $J_2$, and $K$ considered as colors. Choose maps $I_1\xrightarrow{f^1_1} J_1$, $I_2\xrightarrow{f^1_2} J_1$, $I_3\xrightarrow{f^2_3} J_2$, $J_1\xrightarrow{f_1} K$, and $J_2\xrightarrow{f_2} K$. Let $\sigma$ be the ordering $2<1$ of $\{1,2\}$, and let $e$ be the unique ordering of a singleton. Then we have elements
    \[(\sigma,\{f^1_1,f^1_2\})\in \Mul_{\Conv}(\{I_1,I_2\},J_1),\ \ (e,\{f^2_3\})\in\Mul_{\Conv}(I_3,J_2),\ \ (\sigma, \{f_1,f_2\})\in \Mul_{\Conv}(\{J_1,J_2\},K).\]
    Then the composition map \eqref{eq: composition Conv} gives rise to  an element of $\Mul_{\Conv}(\{I_1,I_2,I_3\},K)$. The data of total ordering of $U=\{1,2,3\}$ is $3<2<1$, and the data of maps are 
    \[I_1\xrightarrow{f^1_1}J_1\xrightarrow{f_1}K,\quad I_2\xrightarrow{f^1_2}J_1\xrightarrow{f_1}K,\te{ and } I_3\xrightarrow{f^2_3}J_2\xrightarrow{f_2}K.\]
\end{ex}

\begin{rmk}
    For the operad $\Conv$, the colors are themselves finite sets. Thus, the role of a color (e.g. $I_u$'s, $J_v$'s and $K$) in the definition of the operad and the role of a finite set (e.g. $U$ and $V$) may become easy to confuse. For this reason, given a finite set, we will distinguish our terminology according to whether it is being regarded as a color of the operad or simply as a finite set. For example, for $I_u$ as in Definition \ref{def: operad Conv}, we write $I_u\in\Conv$, not $I_u\in\fSet$. Following the notation of \cite[Remark 2.1.1.12]{HA}, we use the notation $\Conv_{\langle 1\rangle}$ for the underlying category of $\Conv$. This is canonically equivalent to $\fSet$.
\end{rmk}

\sssec{}\label{sssec: Conv algebra object}
Until the end of this subsection, let $\bC$ be a symmetric monoidal category.
A data of a $\Conv$-algebra $\cA$ in $\bC$ consists of the following data: 
\begin{enumerate}
    \item a collection $\{\cA_I\}_{I\in\Conv_{\langle 1\rangle}}\subseteq \bC$ of objects, each of which is an evaluation of $\cA$ at each color;
    \item for any operation $(\sigma,\{f_u\}_{u\in U})\in \Mul_{\Conv}(\{I_u\}_{u\in U},J)$ as in (2) of Definition \ref{def: operad Conv} there should be an associated map $\bigotimes_{u\in U}\cA_{I_u}\xrightarrow{(\sigma,\{f_u\}_{u\in U})}\cA_J$, and the collection of such is supplied with higher coherence data.
\end{enumerate}

\sssec{}
Recall the definition of the operad $\E_1$ as an ordinary colored operad (\cite[Definition 4.1.1.1]{HA}). It has a single color $\mf{a}$, and given a finite set $U$, the set of operations $\Mul_{\E_1}(\{\mf{a}\}_{u\in U},\mf{a})$ is the set of linear ordering $\sigma$'s on $U$. Composition of operations, i.e., linear orderings, is determined in a lexicographical manner as in \eqref{eq: lexicographic}. We observe:

\begin{lemma}\label{lem: Conv and E1}
    For each $I\in\Conv$ there is a canonical map $\E_1\rightarrow \Conv$ between operads that sends a color $\mf{a}\in \E_1$ to the color $I\in\Conv$. Therefore, given a $\Conv$-algebra $\cA$ in $\bC$, each $\cA_I$ has a natural lift to $\Alg(\bC)$. This construction is functorial in $I$, and hence we get a functor $\Conv_{\langle 1\rangle}\rightarrow \Alg(\bC)$, which is $I\mapsto \cA_I$ on the level of the underlying objects.
\end{lemma}
\begin{proof}
    Given $I\in\Conv$, we construct a map $\E_1\rightarrow\Conv$ between colored operads by the following rule: on the level of colors it maps $\mf{a}\mapsto I$, and on the level of operations it maps
    \begin{equation}\label{eq: E1 to Conv for I}
        \Mul_{\E_1}(\{\mf{a}\}_{u\in U},\mf{a})\rightarrow \Mul_{\Conv}(\{I\}_{u\in U},I),\quad \sigma\mapsto (\sigma,\{\id_I\}_{u\in U}).
    \end{equation}
    As for the second claim pertaining to functoriality in $I$, we consider $\Conv_{\langle 1\rangle}$ as an operad via the construction of \cite[Remark 2.1.1.4]{HA}. A $\Conv_{\langle 1\rangle}$-algebra object is equivalent to a functor out of $\Conv_{\langle 1\rangle}$. Therefore, it suffice to construct a functor
    \[\Conv\lAlg(\bC)\rightarrow\Conv_{\langle 1\rangle}\lAlg(\Alg(\bC)),\]
    and hence suffices to construct a bifunctor of operads
    $\Conv_{\langle 1\rangle}\times \E_1\rightarrow \Conv$.
    On the level of colors this maps $(I,\mf{a})\mapsto I$, and given a finite set $U$ and $f\in\Mul_{\Conv\langle 1\rangle}(I, J)$ and $\sigma'\in\Mul_{\E_1}(\{\mf{a}\}_{u\in U},\mf{a})$, then the product $f\times \sigma'$ of operations gets mapped to the operation
    \[(\sigma', \{f\}_{u\in U})\in \Mul_{\Conv}(\{I\}_{u\in U}, J).\]
    Since the bifunctor restricts to \eqref{eq: E1 to Conv for I} on $\{I\}\times\E_1$, the resulting functor $\Conv_{\langle 1\rangle}\rightarrow \Alg(\bC)$ maps $I\rightarrow\cA_I$.
\end{proof}

\begin{prop}
    The functor 
    \begin{equation}\label{eq: Conv and E1}
        \Conv\lAlg(\bC)\rightarrow \Fun(\Conv_{\langle 1\rangle},\Alg(\bC))
    \end{equation}
    constructed in Lemma \ref{lem: Conv and E1} is an equivalence.
\end{prop}
\begin{proof}
    We exhibit an explicit inverse to the functor above. Suppose we have $(I\mapsto\cA_I)\in \Fun(\Conv_{\langle 1\rangle},\Alg(\bC))$. We describe what the corresponding $\Conv$-algebra $\cA$ is. In light of Sect. \ref{sssec: Conv algebra object} the evaluation of $\cA$ at each color $I$ is $\cA_I$. Moreover, the operation map $\bigotimes_{u\in U}\cA_{I_u}\xrightarrow{(\sigma,\{f_u\}_{u\in U})}\cA_J$ is the composition
    \[\bigotimes_{u\in U}\cA_{I_u}\simeq \bigotimes_{1\le s\le |U|}\cA_{I_{\sigma(s)}}\rightarrow \bigotimes_{1\le s\le |U|}\cA_{J}\xrightarrow{\on{mult}}\cA_J.\]
    Here, we view the linear ordering $\sigma$ as a bijection $\{1,2,\cdots,|U|\}\rightarrow U$, the first map is reordering, and the second map is the tensor product of maps $\cA_{I_{\sigma(s)}}\rightarrow \cA_J$, each coming from the functoriality of $I\mapsto \cA_I$.
\end{proof}

\sssec{}
Any ordinary category can be thought of as an operad, see \cite[Remark 2.1.1.4]{HA}. Hence, the category $\Conv_{\langle 1\rangle}$ gives rise to an operad, which we still denote by $\Conv_{\langle 1\rangle}$. We have $\Conv_{\langle 1\rangle}\lAlg(\bC)=\Fun(\fSet,\bC)$. Since the equivalence \eqref{eq: Conv and E1} holds for any symmetric monoidal category $\bC$, we conclude:

\begin{cor}
    We have a canonical equivalence $\Conv\simeq \Conv_{\langle 1\rangle}\otimes \E_1$ of operads.
\end{cor}

\subsection{$\Conv$-algebras arising from factorizable Satake}\label{ssec: Conv algebras from fact Satake}
In this subsection, we enhance both sides of the factorizable Satake into $\Conv$-algebras equipped with a module structure.

\sssec{}\label{sssec: local Hecke Conv algebra in Stk}
Consider the span category $\Corr(\PreStklft, \on{sch})$ of prestacks. Its objects are $\cY\in\PreStklft$ and a morphism from $\cY_1$ to $\cY_2$ is a diagram
\begin{equation*}\label{eq: Corr diagram nontwisted}
    \begin{tikzcd}
	{\cY_{1,2}} & {\cY_1} \\
	{\cY_2}
	\arrow["g", from=1-1, to=1-2]
	\arrow["f"', from=1-1, to=2-1]
\end{tikzcd}
\end{equation*}
where $f$ is required to be ind-schematic. We endow this category with a Cartesian symmetric monoidal structure. A rigorous definition of this category can be found in \cite[Part III, Chapter 7, Sect. 1]{GR}.

\sssec{}\label{sssec: Hecke loc as a Conv algebra}
 We observe that the collection $\{\Hecke^{\loc}_X\}_{I\in \Conv}$ of stacks can be naturally lifted to a $\Conv$-algebra object $\cHecke^{\loc}_X$ in $\Corr(\PreStklft, \on{sch})$. Given a color $I$, the evaluation $(\cHecke^{\loc}_X)_I$ is $\Hecke^{\loc}_{X^I}$. To unburden the notation, let us omit $X$ from the subscript and denote it as $\Hecke^{\loc}_I$. Suppose we are given an operation
    \[(\sigma,\{f_u\}_{u\in U})\in\Mul_{\Conv}(\{I_u\}_{u\in U},J).\]
    In light of Sect. \ref{sssec: Conv algebra object} we have to describe the corresponding map 
    \[\prod_{u\in U}\Hecke^{\loc}_{{I_u}}\rightarrow \Hecke^{\loc}_{J} \quad \te{in } \Corr(\PreStklft, \on{sch}).\]
    The correspondence is
    \begin{equation}\label{eq: operation for cHecke}
        \begin{tikzcd}
	{\wt{\Hecke}^{\loc}_{I_{\sigma(1)},\cdots,I_{\sigma(|U|)}\rightarrow J}} & {\prod_{1\le s\le |U|}\Hecke^{\loc}_{I_{\sigma(s)}}} \\
	{\Hecke^{\loc}_J.}
	\arrow[from=1-1, to=1-2]
	\arrow[from=1-1, to=2-1]
\end{tikzcd}
    \end{equation}
Here the total ordering $\sigma$ of $U$ is considered as a bijection $\sigma:\{1,\cdots,|U|\}\rightarrow U$, and the stack $\wt{\Hecke}^{\loc}_{I_{\sigma(1)},\cdots,I_{\sigma(|U|)}\rightarrow J}$ is a stack over $X^J$ whose points classify the tuples
\begin{equation}\label{eq: point of iterated Hecke loc}
(\ul{x},\cP^0,\cdots,\cP^{|U|},\alpha^1,\cdots,\alpha^{|U|}),
\end{equation}
where $\ul{x}\in X^J$, its restriction along $f_{\sigma(s)}$ is $\ul{x}_{s}\in X^{I_{\sigma(u)}}$, $\cP^s$'s are $G$-bundles on $D_{\ul{x}}$, and each $\alpha^s$ is an isomorphism between $\cP^{s-1}$ and $\cP^{s}$ outside of $\Gamma_{\ul{x}_s}$.
\medskip

The vertical arrow of the correspondence sends the point \eqref{eq: point of iterated Hecke loc} to $(\ul{x}, \cP^0,\cP^{|U|},\alpha^{\on{con}})\in \Hecke^{\loc}_{X^J}$, where $\alpha^{\on{con}}:\cP^0|_{\cD^{\circ}_{\ul{x}}}\xrightarrow{\sim} \cP^{|U|}|_{\cD^{\circ}_{\ul{x}}}$ is the composition of restrictions $\alpha^s|_{\cD^{\circ}_{\ul{x}}}: \cP^{u-1}|_{\cD^{\circ}_{\ul{x}}}\xrightarrow{\sim} \cP^{u}|_{\cD^{\circ}_{\ul{x}}}$. The horizontal arrow is the product of maps 
\[\wt{\Hecke}^{\loc}_{I_{\sigma(1),\cdots,I_{\sigma(|U|)}}\rightarrow J}\rightarrow \Hecke^{\loc}_{I_{\sigma(s)}}\]
for each $1\le s\le |U|$, which sends the point \eqref{eq: point of iterated Hecke loc} to the point
\[(\ul{x}_s,\cP^{s-1}|_{\cD_{\ul{x}_s}}, \cP^{s}|_{\cD_{\ul{x}_s}},\alpha^s|_{\cD^{\circ}_{\ul{x}_s}})\in \Hecke^{\loc}_{I_{\sigma(s)}}.\]

To make sense of this description of $\cHecke^{\loc}_X$ as a $\Conv$-algebra, we have to take care of the higher coherence data. This can be done by repeating the argument of \cite[Sect. 1.2.7]{ZhaoNew}, hence we conclude:
\begin{prop}
    There is a canonical $\Conv$-algebra object $\cHecke^{\loc}_X$ in the symmetric monoidal category $\Corr(\PreStklft, \on{sch})$ admitting the description of evaluation at colors and operation maps as in Sect. \ref{sssec: Hecke loc as a Conv algebra}.
\end{prop}

\begin{rmk}
    To be strict, $\cHecke^{\loc}_X$ does not make sense as a $\Conv$-algebra object in $\Corr(\PreStklft,\on{sch})$, because each $\Hecke^{\loc}_{X^I}$ does not lie in $\PreStklft$ and the vertical map of \eqref{eq: operation for cHecke} is ind-schematic, but not schematic. However, we use the method of approximating by finite pieces outlined in Sect. \ref{sssec: subtlety Hecke} and work as if $\cHecke^{\loc}_{X}$ makes sense as a $\Conv$-algebra in $\Corr(\PreStklft,\on{sch})$. The same applies to every variants of $\cHecke^{\loc}_X$ to be introduced later.
\end{rmk}

\sssec{}\label{sssec: local Hecke Conv algebra in Stk Ge}
Consider the category \[\PreStklft+\on{Ge}_A:= (\PreStklft)_{/\Bet^2 A}\] 
of stacks equipped with an $A$-gerbe, and the correspondence category $\Corr(\PreStklft+\on{Ge}_A, \on{sch})$. Its objects are $(\cY,\cG)\in\PreStklft+\on{Ge}_A$ where $\cY$ is a prestack and $\cG$ is an $A$-gerbe on it, and a morphism from $(\cY_1,\cG_1)$ to $(\cY_2,\cG_2)$ is a diagram
\begin{equation}\label{eq: morphism Corr PreStk Ge}
    \begin{tikzcd}
	{(\cY_{1,2},\cG_{1,2})} & {(\cY_1,\cG_1)} \\
	{(\cY_2,\cG_2),}
	\arrow["g", from=1-1, to=1-2]
	\arrow["f"', from=1-1, to=2-1]
\end{tikzcd}
\end{equation}
where the map $f$ on the level of prestacks is ind-schematic. The category is equipped with a symmetric monodal structure which is Cartesian on the level of prestacks and the external sum on the level of gerbes.

\sssec{}
From the compatibilities of gerbes on $\Hecke^{\loc}_{X^I}$ with respect to various versions of convolution diagrams established in \cite{ZhaoNew}, we can lift a $\Conv$-algebra $\cHecke^{\on{loc}}_X$ in $\Corr(\PreStklft, \on{sch})$ to a $\Conv$-algebra $(\cHecke^{\on{loc}}_X,\cG^{\on{loc}})$ in $\Corr(\PreStklft+\on{Ge}_A, \on{sch})$. As before, we describe some piece of its structure.
The evaluation $(\cHecke^{\on{loc}}_X,\cG^{\on{loc}})_I$ at the color $I$ is a tuple $(\Hecke^{\loc}_{X^I},\cG^I)$. In the correspondence \eqref{eq: operation for cHecke}, the pullback of the external sum of gerbes $\cG^{I_{\sigma(u)}}$'s on the stacks $\Hecke^{\loc}_{I_{\sigma(s)}}$'s to $\wt{\Hecke}^{\loc}_{I_{\sigma(1)},\cdots,I_{\sigma(|U|)}\rightarrow J}$ is canonically identified with the pullback of $\cG^{J}$ on $\Hecke^{\loc}_J$, and hence \eqref{eq: operation for cHecke} naturally lifts to a morphism
\[\prod_{u\in U}(\Hecke^{\loc}_{{I_u}},\cG^{I_u})\rightarrow (\Hecke^{\loc}_{J},\cG^J) \quad \te{in } \Corr(\Stklft+\on{Ge}_A).\]

We sum up this discussion as:

\begin{prop}
    We have a $\Conv$-algebra $(\cHecke^{\loc}_X,\cG^{\loc})$ in the symmetric monoidal category $\Corr(\PreStklft+\on{Ge}_A, \on{sch})$ admitting the description of evaluation at colors and operation maps as in Sect. \ref{sssec: local Hecke Conv algebra in Stk Ge}.
\end{prop}

\sssec{}
$(X,0)$ is a cocommutative coalgebra object in $\Stk+\on{Ge}_A$, and hence using the symmetric monoidal functor $(\Stk+\on{Ge}_A)^{\op}\rightarrow \Corr(\Stk+\on{Ge}_A)$ we get $(X,0)$ as a commutative algebra object in $\Corr(\Stk+\on{Ge}_A)$. For example, the multiplication map $(X,0)\times (X,0)\rightarrow (X,0)$ is given by the correspondence
\[\begin{tikzcd}
	(X,0) & {(X\times X,0)} \\
	(X,0).
	\arrow["\Delta", from=1-1, to=1-2]
	\arrow["{\id_X}"', from=1-1, to=2-1]
\end{tikzcd}\]

\sssec{}\label{sssec: Conv from E infty}
Suppose $\bC$ is a symmetric monoidal category. We have a canonical functor
\[\CAlg(\bC)\rightarrow \Conv\lAlg(\bC),\quad \cC\mapsto \cC^{\otimes}.\]
The evaluation of $\cC^{\otimes}$ at the color $I$ is the $I$-fold tensor product $\cC^{\otimes I}$, i.e., $(\cC^{\otimes})_I=\cC^{\otimes I}$. Using the notation of Sect. \ref{sssec: Conv algebra object}, the operation $\bigotimes_{u\in U}(\cC^{\otimes})_{I_u}\xrightarrow{(\sigma,f)}(\cC^{\otimes})_J$ is the multiplication map $\bigotimes_{u\in U}\cC^{\otimes I_u}\rightarrow\cC^{\otimes J}$ of $\cC$.
\medskip

The functor above has a natural symmetric monoidal structure. Therefore, using the equivalence $\CAlg(\CAlg(\bC))\simeq \CAlg(\bC)$ we get a functor
\begin{equation}\label{eq: CAlg to CAlg P Alg}
    \CAlg(\bC)\rightarrow \CAlg(\Conv\lAlg(\bC)).
\end{equation}

Applying this construction to $\bC=\Corr(\Stk+\on{Ge}_A)$, and $\cC=(X,0)$ we have $(X,0)^{\otimes}$ as an $\E_{\infty}$-object in $\Conv(\Corr(\Stk+\on{Ge}_A))$. We readily check:

\begin{prop}\label{prop: cHecke as an X module}
    $(\cHecke^{\loc}_X,\cG^{\loc})$ is a $(X,0)^{\otimes}$-module in the symmetric monoidal category $\Conv\lAlg(\Corr(\Stk+\on{Ge}_A))$.
\end{prop}

\sssec{}
The formation of twisted sheaves gets enhanced to the right-lax symmetric monoidal functor
\[\Shv^{\on{tw}}_{\on{Corr}}:\on{Corr}\left(\PreStklft+\on{Ge}_A, \on{sch}\right)\rightarrow \DGcat.\]
For a map \eqref{eq: morphism Corr PreStk Ge} this functor assigns the composition
\[\Shv(\cY_1)\xrightarrow{g^*} \Shv(\cY_{12})\xrightarrow{f_!}\Shv(\cY_2).\]
We refer readers to Sect. \ref{ssec: functoriality of twisted sheaves} for details.
\medskip

Therefore, we get a right-lax symmetric monoidal functor
\[\Shv^{\on{tw}}_{\on{Corr}}:\Conv\lAlg\left(\on{Corr}\left(\PreStklft+\on{Ge}_A, \on{sch}\right)\right)\rightarrow\Conv\lAlg(\DGcat).\]

\sssec{}
We define a $\Conv$-algebra $\Shv^q(\cHecke^{\loc}_X)$ in $\DGcat$ as the image of the $\Conv$-algebra $(\cHecke^{\loc}_X,\cG^{\loc})$ under the above functor. The structure of this $\Conv$-algebra efficiently packages the functoriality and various convolution operations on $\Shv^q(\Hecke^{\loc}_{X^I})$, together with compatibilities between them. Let us discuss this point in detail. 

\sssec{Monoidal structure}\label{sssec: ShvqHeckeloc monoidal str}
    The $\Conv$-algebra structure of $\Shv^q(\cHecke^{\loc}_X)$ encodes the monoidal structure of each category $\Shv^q(\Hecke^{\loc}_{X^I})$. Indeed, the correspondence \eqref{eq: operation for cHecke} for the operation \[(\triv, \{\id_I\}_{u\in \{1,2,\cdots,k\}})\in \Mul_{\Conv}(\{I\}_{u\in \{1,2,\cdots,k\}},\{I\})\]
    is the $k$-fold convolution diagram of $\Shv^q(\Hecke^{\loc}_{X^I})$.
    \medskip

    For each color $I$, there is a unique operation whose target is $I$ and $U=\varnothing$. In this case, the correspondence \eqref{eq: operation for cHecke} reads
    \[\begin{tikzcd}
	{(X^I,0)} & {(\pt,0)} \\
	{(\Hecke^{\loc}_{X^I},\cG^I),}
	\arrow[from=1-1, to=1-2]
	\arrow[from=1-1, to=2-1]
\end{tikzcd}\]
where the vertical map sends $\ul{x}$ to the tuple $(\ul{x},\cP_{\triv},\cP_{\triv},\id)$ where $\cP_{\triv}$ means the trivial $G$-bundle. The pull-push operation along this correspondence agrees with the unit map
    \[\Vect\rightarrow \Shv^q(\Hecke^{\loc}_{X^I})\]
    of the convolution monoidal structure.

\sssec{External convolution product}
    In the formation of the external convolution product, we started with an ordered partition $I=I_1\coprod I_2\cdots\coprod I_k$.
    The external convolution product functor 
    \begin{equation*}
    \Shv^q(\Hecke^{\loc}_{X^{I_1}})\otimes\cdots\otimes \Shv^q(\Hecke^{\loc}_{X^{I_k}})\rightarrow\Shv^q(\Hecke^{\loc}_{X^{I}})
\end{equation*}
is encoded in the operation
\[(\triv, \{\iota_u\}_{u\in \{1,2,\cdots,k\}})\in\Mul_{\Conv}(\{I_u\}_{u\in \{1,2,\cdots,k\}}, I).\]
Here we used the $U=\{1,2,\cdots,k\}$ and its trivial linear ordering, and we use the inclusion $\iota_u: I_u\hookrightarrow I$.

\sssec{Functoriality}
Finally, to encode the functoriality
\[\Shv^q(\Hecke^{\loc}_{X^I})\rightarrow \Shv^q(\Hecke^{\loc}_{X^I})\]
for $I\rightarrow J$, we let $U=\{*\}$ and consider the operation
\[(\triv, I\rightarrow J)\in\Mul_{\Conv}(\{I\},J).\]

\sssec{}
Recall the construction of $(X,0)^{\otimes}$ as an $\E_{\infty}$-object of $\Conv(\Corr(\Stklft+\on{Ge}_A))$. This maps to an $\E_{\infty}$-object of $\Conv\lAlg(\DGcat)$, which we denote temporarily as $\Shv^{\on{tw}}_{\Corr}((X,0)^{\otimes})$. On the other hand, we have an $\E_{\infty}$-object $\Shv(X)^{\otimes}$ of $\Conv(\DGcat)$ obtained by applying the paradigm of Sect. \ref{sssec: Conv from E infty} to the situation $\bC=\DGcat$, $\cP=\Conv$, and $\cC=\Shv(X)$. We have a canonical equivalence
\[\Shv^{\on{tw}}_{\Corr}((X,0)^{\otimes})\simeq \Shv(X)^{\otimes}\quad \te{in } \CAlg(\Conv\lAlg(\DGcat)).\]
In conclusion, applying $\Shv^{\on{tw}}_{\Corr}$ to the module structure of Proposition \ref{prop: cHecke as an X module} we conclude:

\begin{prop}
    $\Shv^q(\cHecke^{\loc}_X)$ has a natural $\Shv(X)^{\otimes}$-module structure in $\Conv\lAlg(\DGcat)$.
\end{prop}

\begin{cor}\label{cor: A side as a LS module in Conv}
    $\Shv^q(\cHecke^{\loc}_X)$ has a natural $\LS(X)^{\otimes}$-module structure in $\Conv\lAlg(\DGcat)$.
\end{cor}

\sssec{}\label{sssec: map from LS to Rep H loc}
 We have a symmetric monoidal map $\LS(X)\rightarrow \Rep(H)^{\loc}_{\cG_Z}(X)$, which serves as a unit map where $\Rep(H)^{\loc}_{\cG_Z}(X)$ is thought of as a commutative $\LS(X)$-linear algebra. Using the construction in Sect. \ref{sssec: Conv from E infty} we get a map
 \begin{equation}\label{eq: LS to twRepH}
\LS(X)^{\otimes}\rightarrow \left(\Rep(H)^{\loc}_{\cG_Z}(X)\right)^{\otimes}\quad \te{in } \CAlg(\Conv\lAlg(\DGcat)).
 \end{equation}
In particular, we conclude:

\begin{prop}\label{prop: B side as a LS module in Conv}
    $\left(\Rep(H)^{\loc}_{\cG_Z}(X)\right)^{\otimes}$ has a natural $\LS(X)^{\otimes}$-module structure in $\Conv\lAlg(\DGcat)$.
\end{prop}

\sssec{}
We can do the same construction using $\QLisse$-sheaves instead. Namely, from the $\E_{\infty}$-objects $\QLisse(X)$ and $\Rep(H)^{\QLisse}(X)$ of $\DGcat$ we produce
\[\QLisse(X)^{\otimes},\quad \left(\Rep(H)^{\QLisse}(X)\right)^{\otimes}\quad\te{in } \CAlg(\Conv\lAlg(\DGcat)),\]
and from the $\E_{\infty}$-map 
$\QLisse(X)\rightarrow \Rep(H)^{\QLisse}(X)$
we get:

\begin{prop}\label{prop: B QLisse as a QLisse module in Conv}
    $\left(\Rep(H)^{\QLisse}_{\cG_Z}(X)\right)^{\otimes}$ has a natural $\QLisse(X)^{\otimes}$-module structure in $\Conv\lAlg(\DGcat)$.
\end{prop}

\sssec{}\label{sssec: CAlg Conv algebras}
The constructions above admit an obvious sheafy analog. Namely, we can apply the paradigm of Sect. \ref{sssec: Conv from E infty} to the situation $\cP=\Conv$ and $\bC=\Shv(X_{\et};\DGcat)$ to get a functor
\[\CAlg(\Shv(X_{\et};\DGcat))\rightarrow \CAlg(\Conv\lAlg(\Shv(X_{\et};\DGcat))).\]
Applying this functor to the $\E_{\infty}$-sheaves
\[\textbf{LS}, \ \ , \textbf{QLisse}, \ \ \textbf{Rep}(H)^{\loc}_{\cG_Z},\ \te{and} \ \textbf{Rep}(H)^{\QLisse}_{\cG_Z},\]
we get
\[\textbf{LS}^{\otimes}, \ \ , \textbf{QLisse}^{\otimes}, \ \ \left(\textbf{Rep}(H)^{\loc}_{\cG_Z}\right)^{\otimes },\ \te{and} \ \left(\textbf{Rep}(H)^{\QLisse}_{\cG_Z}\right)^{\otimes}\]
in $\CAlg(\Conv\lAlg(\Shv(X_{\et};\DGcat)))$, and we get:

\begin{cor}
    $\left(\mathbf{Rep}(H)^{\loc}_{\cG_Z}\right)^{\otimes}$ has a natural $\mathbf{LS}^{\otimes}$-module structure in $\Conv\lAlg(\Shv(X_{\et};\DGcat))$. The same holds when we change the sheaf theory from $\LS$ to $\QLisse$.
\end{cor}

\sssec{}\label{sssec: rel tensor prod}
In the situation of Sect. \ref{sssec: Conv from E infty}, the functor \eqref{eq: CAlg to CAlg P Alg} preserves relative tensor products under the assumption that the tensor product of $\bC$ commutes with  colimits. Therefore, since we have
\[\textbf{Rep}(H)^{\loc}_{\cG_Z}=\textbf{Rep}(H)^{\QLisse}_{\cG_Z}\otimes_{\textbf{QLisse}}\textbf{LS}\]
in $\CAlg(\Shv(X_{\et}))$ we get
\[\left(\textbf{Rep}(H)^{\loc}_{\cG_Z}\right)^{\otimes}= \left(\textbf{Rep}(H)^{\QLisse}_{\cG_Z}\right)^{\otimes}\otimes_{\textbf{QLisse}^{\otimes}}\textbf{LS}^{\otimes}\]
in $\Conv\lAlg(\CAlg(\Shv(X_{\et})))= \CAlg(\Shv(X_{\et};\Conv\lAlg))$.

\subsection{Passage from $\Conv$-algebras to a natural transformation}\label{ssec: construction of a nat transf}\label{ssec: passage to nat transf}
We claim that we have the following mechanism to pass from $\LS(X)^{\otimes}$-module objects in $\Conv\lAlg(\DGcat)$ to functors in a form \eqref{eq: fSet to fSet Alg}:
\begin{prop}\label{prop: module in Conv to fSet Alg}
Suppose $\bC$ is a symmetric monoidal category, and $\cC$ is an $\E_{\infty}$-algebra object in $\Conv\lAlg(\bC)$.
There is a canonical functor \footnote{See Sect. \ref{sssec: Alg LS construction} for the $\Fun^{\circ}$ notation.}
\begin{equation}\label{eq: passage functor}
    \cC\lMod\left(\Conv\lAlg(\bC)\right)\rightarrow \Fun^{\circ}(\fSet, \fSet + \Alg(\cC\lMod)).
\end{equation}
\end{prop}

\sssec{}\label{sssec: LMcomm}
For the proof of Proposition \ref{prop: module in Conv to fSet Alg} we introduce another operad. Recall the colored operad $\LM$ from \cite[Sect. 4.2.1]{HA}. An $\LM$-algebra object is a pair $(\cA,\cM)$ where $\cA$ is an algebra object and $\cM$ is a left $\cA$-module. We consider its variant $\LMcomm$ whose algebra object is a pair $(\cC,\cM)$ where $\cC$ is a commutative algebra object and $\cM$ is a left $\cC$-module.

\begin{defn}\label{defn: LMcomm}
    An operad $\LMcomm$ is an ordinary colored operad described as follows.
    \begin{enumerate}
        \item The collection of colors is $\{\mf{c},\mf{m}\}$.
        \item Suppose we are given a finite set $U$, a collection of colors $\{C_u\}_{u\in U}$, and a color $C$. We have
        \[\Mul_{\LMcomm}(\{C_u\}_{u\in U}, C)=\begin{cases}
            \{*\} & \te{if } C=\mf{c} \te{ and }C_u=\mf{c}\ \forall u\in U,\\
            \varnothing & \te{if } C=\mf{c} \te{ and } \exists u\in U \te{ s.t. } C_u=\mf{m},\\
            \{\ast\} & \te{if } C=\mf{m} \te{ and } \exists! u\in  U \te{ s.t. } C_u=\mf{m},\\
            \varnothing & \te{otherwise.}
        \end{cases}\]
        \item Since the set of operations is at most a singleton, the distinguished element $\id_C$ for each color and the composition maps are defined in the unique way.
    \end{enumerate}
\end{defn}

\sssec{Proof of Proposition \ref{prop: module in Conv to fSet Alg}}
    We have a sequence of equivalences
    \begin{align}\label{eq: key operadic manipulation}
        \LMcomm\lAlg(\Conv\lAlg(\bC))&=\Conv\lAlg(\LMcomm\lAlg(\bC)) \\
        &\stackrel{\eqref{eq: Conv and E1}}{=} \Fun(\Conv_{\langle 1\rangle}, \Alg(\LMcomm\lAlg(\bC))),\nonumber
    \end{align}
    The functor \eqref{eq: key operadic manipulation} fits into a commutative diagram
    \begin{equation}\label{eq: operads comm diagram}
        \begin{tikzcd}
	{\LMcomm\lAlg(\Conv\lAlg(\bC))} & {\Fun(\Conv_{\langle 1\rangle}, \Alg(\LMcomm\lAlg(\bC)))} \\
	{\CAlg(\Conv\lAlg(\bC))} & {\Fun(\Conv_{\langle 1\rangle}, \CAlg(\bC))}
	\arrow["{\eqref{eq: key operadic manipulation}}", from=1-1, to=1-2]
	\arrow[from=1-1, to=2-1]
	\arrow[from=1-2, to=2-2]
	\arrow["U"',"\sim", from=2-1, to=2-2]
\end{tikzcd}
    \end{equation}
Here the vertical arrows are induced by forgetful functors $\LMcomm\rightarrow\CAlg$ and \[\Alg(\LMcomm\lAlg(\bC))\rightarrow\LMcomm\lAlg(\bC)\rightarrow\CAlg(\bC)\] respectively, and the lower horizontal arrow $U$ is the composition
\[\CAlg(\Conv\lAlg(\bC))\stackrel{\eqref{eq: Conv and E1}}{=}\CAlg(\Fun(\Conv_{\langle 1\rangle},\Alg(\bC)))=\Fun(\Conv_{\langle 1\rangle}, \CAlg(\bC)).\]

The category $\cC\lMod(\Conv\lAlg(\bC))$ fits into a Cartesian diagram
\[\begin{tikzcd}
	{\cC\lMod(\Conv\lAlg(\bC))} & {\LMcomm\lAlg(\Conv\lAlg(\bC))} \\
	\pt & {\CAlg(\Conv\lAlg(\bC)).} 
	\arrow[from=1-1, to=1-2]
	\arrow[from=1-1, to=2-1]
	\arrow[from=1-2, to=2-2]
	\arrow["\cC"', from=2-1, to=2-2]
\end{tikzcd}\]
Moreover, since the category $\Alg(\LMcomm\lAlg(\bC))$ consists of pairs $(\bc,\ba)$ where $\bc\in\CAlg(\bC)$ and $\ba\in\Alg(\bc\lMod(\bC))$,
the category $\Fun^{\circ}(\fSet,\fSet+\Alg(\cC\lMod))$ fits into a Cartesian diagram
\[\begin{tikzcd}
	{\Fun^{\circ}(\fSet,\fSet+\Alg(\cC\lMod))} & {\Fun(\Conv_{\langle 1\rangle}, \Alg(\LMcomm\lAlg(\bC))).} \\
	\pt & {\Fun(\Conv_{\langle 1\rangle}, \CAlg(\bC))}
	\arrow[from=1-1, to=1-2]
	\arrow[from=1-1, to=2-1]
	\arrow[from=1-2, to=2-2]
	\arrow["{(I\mapsto \cC_I)}"', from=2-1, to=2-2]
\end{tikzcd}\]
Since the object $\cC\in\CAlg(\Conv\lAlg(\bC))$ is mapped to $(I\mapsto \cC_I)\in \Fun(\Conv_{\langle 1\rangle},\CAlg(\bC))$ under $U$, the desired functor \eqref{eq: passage functor} is obtained as an induced map on pullbacks. \qed

\begin{rmk}\label{rmk: evaluation at I}
    From construction, for each $I$ we have a commutative diagram
    \[\begin{tikzcd}
	{\cC\lMod(\Conv\lAlg(\bC))} & {\Fun^{\circ}(\fSet,\fSet+\Alg(\cC\lMod))} \\
	{\cC_I\lMod(\Alg(\bC))} & {\Alg(\cC_I\lMod(\bC))}
	\arrow[from=1-1, to=1-2]
	\arrow[from=1-1, to=2-1]
	\arrow[from=1-2, to=2-2]
	\arrow["\sim",from=2-1, to=2-2]
\end{tikzcd}\]
where the left vertical arrow is the arrow induced by the operadic map $\E_1\rightarrow \Conv$ mapping $\mf{a}\mapsto I$ (see Lemma \ref{lem: Conv and E1}) and the right vertical arrow is evaluation at $I$.
\end{rmk}

\sssec{}
Now we discuss a proof strategy of Theorem \ref{thm: fact Satake nat transf}.
Applying Proposition \ref{prop: module in Conv to fSet Alg} to the situation
$\bC=\DGcat$ and  $\cC=\LS(X)^{\otimes}$,
we get a canonical functor
\begin{equation}\label{eq: passage functor LS module}
    \LS(X)^{\otimes}\lMod(\Conv\lAlg(\DGcat))\rightarrow \Fun^{\circ}(\fSet, \fSet+ \Alg(\LS(X)^{\otimes}\lMod)).
\end{equation}

In Subsection \ref{ssec: Conv algebras from fact Satake} we have constructed two $\LS(X)^{\otimes}$-module $\Conv$-algebras
\[\left(\Rep(H)^{\loc}_{\cG_Z}(X)\right)^{\otimes}\quad \te{and}\quad \Shv^q(\cHecke^{\loc}_X)\]
(see Corollary \ref{cor: A side as a LS module in Conv} and Proposition \ref{prop: B side as a LS module in Conv}).
Each of those gives rise to a functor $\fSet\rightarrow \fSet+\Alg(\LS(X)^{\otimes}\lMod)$ respecting the projection, and using Remark \ref{rmk: evaluation at I} their evaluation at $I$ is 
\[\left(\Rep(H)^{\loc}_{\cG_Z}(X)\right)^{\otimes I}\quad \te{and}\quad \Shv^q(\Hecke^{\loc}_{X^I}),\]
respectively. This accomplishes the first part of Theorem \ref{thm: fact Satake nat transf}, i.e., constructing functors of the form \eqref{eq: Satake functorial B} and \eqref{eq: Satake functorial A}. We are left with constructing a natural transformation between them. To do so, we again use the functor \eqref{eq: passage functor LS module} to reduce the problem to the following statement:

\begin{thm}\label{thm: fact Satake Conv at X}
    We have a canonical map
    \[\Rep(H)^{\loc}_{\cG_Z}(X)^{\otimes}\rightarrow \Shv^q(\cHecke^{\loc}_X)\]
    as $\LS(X)^{\otimes}$-modules in $\Conv\lAlg(\DGcat)$.
\end{thm}

\sssec{}
    For convenience, let us omit $\DGcat$ from the notations when it is clear from the context. Observe that both $\Conv$-algebras in the question sheafify, i.e., we have sheaves
    \begin{equation}\label{eq: Conv algebra sheafify}
        \left(\mathbf{Rep}(H)^{\loc}_{\cG_Z}\right)^{\otimes} \quad \te{and}\quad \mathbf{Shv}^q(\cHecke^{\loc})
    \end{equation}
    in $\Shv(X_{\et};\Conv\lAlg)$ whose section at $X'\in X_{\et}$ are
    \[\left(\Rep(H)^{\loc}_{\cG_Z}(X')\right)^{\otimes}\quad\te{and}\quad \Shv^q(\cHecke^{\loc}_{X'}),\]
    respectively. The commutative algebra object $\LS(X)^{\otimes}$ also sheafifies, i.e., we have an object $\mathbf{LS}^{\otimes}$ in $\CAlg(\Shv(X_{\et};\Conv\lAlg))$ whose section at $X'\in X_{\et}$ is $\LS(X')^{\otimes}\in\Conv\lAlg$. The $\LS(X)^{\otimes}$-module structures sheafify, i.e., we have a canonical $\mathbf{LS}^{\otimes}$-module structure on both sheaves \eqref{eq: Conv algebra sheafify}. Now it suffice to construct a canonical map
    \begin{equation}\label{eq: fact Satake ultimate}
        \left(\mathbf{Rep}(H)^{\loc}_{\cG_Z}\right)^{\otimes} \rightarrow \mathbf{Shv}^q(\cHecke^{\loc}_X)
    \end{equation}
    between $\mathbf{LS}^{\otimes}$-modules. Note that we have constructed this map `colorwise', i.e., for each $I$ the evaluation of \eqref{eq: fact Satake ultimate} at $I$ is the functor $\textbf{S}_I$ from Subsection \ref{ssec: LS linear fact Satake}. We want to lift this collection to a map between $\Conv$-algebras. 
   Using the relative tensor product formula in Sect. \ref{sssec: rel tensor prod}, it reduces to showing the following statement:
   
   \begin{thm}\label{thm: QLisse Satake ultimate}
      We have a canonical $\mathbf{QLisse}^{\otimes}$-linear map
\begin{equation}
    \left(\mathbf{Rep}(H)^{\QLisse}_{\cG_Z}\right)^{\otimes I}\rightarrow \mathbf{Shv}^q(\cHecke^{\loc})
\end{equation}
whose evaluation at the each color $I$ is $\mathbf{S}^{\QLisse}_{I}$.
   \end{thm}

   Its proof will be given in Sect. \ref{ssec: proof of QLisse ultimate}.
    
    \subsection{Perverse ULA Satake in the language of $\Conv$-algebras}\label{ssec: Perv ULA Satake in Conv language}
    In this subsection, we enhance the perverse ULA version of Satake equivalence (Theorem \ref{thm: factorizable Satake}) to an equivalence between $\Conv$-algebras. In fact, we introduce another operad $\Convcomm$ - which can be understood as a commutative version of $\Conv$ - equipped with a canonical operadic map $\Conv\rightarrow \Convcomm$ and enhance the equivalence to that between $\Convcomm$-algebras. This enhancement serves as an important step towards the proof of Theorem \ref{thm: QLisse Satake ultimate}. Indeed, the $\mathbf{QLisse}$-linear (and $\mathbf{LS}$-linear) Satake functors are uniquely determined by their compatibility with perverse ULA Satake equivalence. (see Lemma \ref{lem: the characterization lemma multiple} and Theorem \ref{thm: characterization for I})

    \begin{nota}
        Throughout this subsection, we will consider various algebras in $\Shv(X_{\et};\CatEc)$. To unburden the notation, we discard $\CatEc$ from the notation and write simply as $\Shv(X_{\et})$.
    \end{nota}

    \sssec{}
    Consider the object
    \[\mathbf{LS}^{\Ab,\ft}_?\in\CAlg(\Conv\lAlg(\Shv(X_{\et})))\]
    whose underlying functor is $I\mapsto \mathbf{LS}^{\Ab,\ft}_I$. Its structure as an $\E_{\infty}$-object internal to $\Conv$-algebras is realized as a subfunctor of $\mathbf{LS}^{\otimes}$ (see Sect. \ref{sssec: CAlg Conv algebras}). 
    \medskip
    
    For each $I$,
    \begin{equation}\label{eq: Perv ULA sheafify}
        \mathbf{Rep}(H^I)^{\Ab,\loc,\ft}_{\cG^{\boxtimes I}_Z}\quad \te{and}\quad \textbf{Sat}^q_I \quad\te{in }\Shv(X_{\et})
    \end{equation}
    are subsheaves of
    \[\left(\mathbf{Rep}(H)^{\QLisse}_{\cG_Z}\right)^{\otimes I} \quad \te{and}\quad \mathbf{Shv}^q(\cHecke^{\loc}_I) \quad\te{in }\Shv(X_{\et};\DGcat),\]
    respectively. The latter objects form the $\Conv$-algebras 
    \begin{equation}\label{eq: two Conv algebras QLisse}
        \left(\mathbf{Rep}(H)^{\QLisse}_{\cG_Z}\right)^{\otimes} \quad \te{and}\quad \mathbf{Shv}^q(\cHecke^{\loc}),
    \end{equation}
    and the $\Conv$-operations restrict to the collection of subsheaves \eqref{eq: Perv ULA sheafify}. Therefore, we have objects
    \begin{equation}\label{eq: two Conv algebras Perv ULA}
        \mathbf{Rep}(H^?)^{\Ab,\loc,\ft}_{\cG^{\boxtimes ?}_Z}\quad \te{and}\quad \textbf{Sat}^q \quad\te{in }\Conv\lAlg(\Shv(X_{\et}))
    \end{equation}
    whose evaluation at $I$ are \eqref{eq: Perv ULA sheafify}.
    \medskip
    
    Moreover, there is a map 
    \[\mathbf{LS}^{\Ab,\ft}_?\rightarrow \mathbf{QLisse}^{\otimes} \quad \te{in } \CAlg(\Conv\lAlg(\Shv(X_{\et})))\] induced from the inclusions $\LS^{\Ab,\ft}(X'^I)\hookrightarrow\QLisse(X')^{\otimes I}$, and hence we get a $\mathbf{LS}^{\Ab,\ft}_{?}$-action on \eqref{eq: two Conv algebras QLisse}. This action preserves the subsheaves \eqref{eq: Perv ULA sheafify}, and hence we get canonical $\mathbf{LS}^{\Ab,\ft}_{?}$-actions on both of \eqref{eq: two Conv algebras QLisse}.
    \medskip

    We claim that the perverse ULA Satake equivalence admits the following enhancement:
    \begin{thm}\label{thm: Perv ULA Satake ultimate}
        There is a canonical $\mathbf{LS}^{\Ab,\ft}_?$-linear equivalence between the $\Conv$-algebras \eqref{eq: two Conv algebras Perv ULA}, whose evaluation at each color $I$ is the  perverse ULA Satake equivalence $\mathbf{S}^{\Ab}_I$. \footnote{Recall from Sect. \ref{sssec: monoidal formulation of PULA Satake} that $\mathbf{S}^{\Ab}_I$ is an equivalence between sheaf of monoidal categories, not between symmetric monoidal categories.}
    \end{thm}

    \sssec{}
    Perverse ULA Satake is not just an equivalence between monoidal categories, but between symmetric monoidal categories. Moreover, the operations on the $\Conv$-algebras \eqref{eq: Perv ULA sheafify} can be alternatively described in terms of external fusion products and functoriality in $I$, which is independent of the choice of ordering $\sigma$. Therefore, it is natural to expect that the desired equivalence of Theorem \ref{thm: Perv ULA Satake ultimate} should be a consequence of an equivalence between algebras with more structures. This motivates us to define the following commutative version of $\Conv$:

    \begin{defn}\label{defn: operad Conv comm}
        An operad $\Convcomm$ is an ordinary colored operad described as follows.
    \begin{enumerate}
        \item The collection of colors is the collection of finite sets.
        \item Suppose we are given a finite set $U$, an $U$-indexed collection of colors $\{I_u\}_{u\in U}$, and a color $J$. The set of operations $\Mul_{\Conv}(\{I_u\}_{u\in U}, J)$ is the set of collections $\{f_{u}\}_{u\in U}$ where each $f_u$ is a map $I_u\rightarrow J$.
        \item For each color $I$ we have a distinguished element
        \[\{\id_I\}\in \Mul_{\Conv}(\{I\},I)\]
        and the composition law is defined to be the composition of functions.
    \end{enumerate}
    \end{defn}

    \sssec{}
    We have an obvious map of operads $\Conv\rightarrow\Convcomm$, and hence for a symmetric monoidal category $\bC$ we have a natural symmetric monoidal functor
    \begin{equation*}
        \Convcomm\lAlg(\bC)\rightarrow \Conv\lAlg(\bC).
    \end{equation*}

    \sssec{}
    We have the following analog of assertions in Sect. \ref{ssec: an operad Conv}:
    
    \begin{lemma}\label{lem: Convcomm and Einfty}
    For each $I\in\Convcomm$ there is a canonical map $\E_{\infty}\rightarrow \Convcomm$ between operads that sends the unique color of $\E_{\infty}$ to the color $I\in\Conv$. Therefore, given a $\Convcomm$-algebra $\cA$ in $\bC$, each $\cA_I$ has a natural lift to $\CAlg(\bC)$. This construction is functorial in $I$, and hence we get a functor $\Conv_{\langle 1\rangle}\rightarrow \CAlg(\bC)$, which is $I\mapsto \cA_I$ on the level of the underlying objects.
\end{lemma}

\begin{prop}
    The functor
    \begin{equation}\label{eq: Convcomm and Einfty}
        \Convcomm\lAlg(\bC)\rightarrow \Fun(\Conv_{\langle 1\rangle},\CAlg(\bC))
    \end{equation}
    constructed in Lemma \ref{lem: Convcomm and Einfty} is an equivalence.
\end{prop}

\begin{cor}
    We have a canonical equivalence $\Convcomm\simeq \Conv_{\langle 1\rangle}\otimes\E_{\infty}$ of operads.
\end{cor}

\begin{rmk}\label{rmk: lifting Conv to Convcomm}
    The data needed to lift a $\Conv$-algebra to a $\Convcomm$-algebra is the identification between $\Conv$-operations corresponding to $(\sigma, \{f_u\}_{u\in U})$ for different choices of linear ordering $\sigma$. Let $\cA$ be a $\Conv$-algebra and suppose it is equipped with a lift to a $\Convcomm$-algebra. Fix the choice of $\{f_u\}$ and denote the $\Conv$-operation on $\cA$ depending on $\sigma$ as
    \[F_{\sigma}: \bigotimes_{u\in U} \cA_{I_u}\rightarrow \cA_J.\]
    
    The space of linear orderings on $U$ is a torsor over the symmetric group $S_{|U|}$. Given $\tau\in S_{|U|}$ and a linear ordering $\sigma$, the data of lifting $\cA$ to a $\Convcomm$-algebra should give us a natural isomorphism 
    \begin{equation}\label{eq: symmetric group repn}
        V_{\sigma\rightarrow\tau\sigma}: F_{\sigma}\rightarrow F_{\tau\sigma}
    \end{equation}
    compatible with the torsor structure, e.g., a data of commutativity
    \[\begin{tikzcd}
	{F_{\sigma}} & {F_{\tau_1\sigma}} \\
	& {F_{\tau_2\tau_1\sigma}.}
	\arrow["{V_{\sigma\rightarrow\tau_1\sigma}}", from=1-1, to=1-2]
	\arrow["{V_{\sigma\rightarrow \tau_2\tau_1\sigma}}"', from=1-1, to=2-2]
	\arrow["{V_{\tau_1\sigma\rightarrow \tau_2\tau_1\sigma}}", from=1-2, to=2-2]
\end{tikzcd}\]
\end{rmk}

    \sssec{}
    We will lift the situation of $\textbf{LS}^{\Ab,\ft}_{?}$ acting on the $\Conv$-algebras \eqref{eq: two Conv algebras Perv ULA} to the same situation with $\Convcomm$-algebras.
    \medskip
    
    The formation of $\Conv$-algebra $\mathbf{Rep}(H^?)^{\Ab,\loc,\ft}_{\cG^{\boxtimes ?}_Z}$ and the commutative $\Conv$-algebra $\mathbf{LS}^{\Ab,\ft}_{?}$ easily gets lifted to $\Convcomm$-algebras. Indeed, the former is the restriction of $\mathbf{QLisse}^{\otimes}$-action on $\left(\mathbf{Rep}(H)^{\QLisse}_{\cG_Z}\right)^{\otimes}$ as $\Conv$-algebras, and this is induced from the analogous structure on $\E_{\infty}$-algebras. Thus, our assertion follows from the fact that $\Conv\rightarrow \E_{\infty}$ factors as $\Conv\rightarrow\Convcomm\rightarrow \E_{\infty}$.
    To summarize:

    \begin{prop}
        There is a canonical lift of $\mathbf{LS}^{\Ab,\ft}_{?}$ to an object of $\CAlg(\Convcomm\lAlg(\Shv(X_{\et})))$ and a canonical lift of $\mathbf{Rep}(H^?)^{\Ab,\loc,\ft}_{\cG^{\boxtimes ?}_Z}$ to an object of $\Convcomm\lAlg(\Shv(X_{\et}))$, and furthermore the lift of the action of $\mathbf{LS}^{\Ab,\ft}_{?}$ on $\mathbf{Rep}(H^?)^{\Ab,\loc,\ft}_{\cG^{\boxtimes ?}_Z}$ as $\Conv$-algebras to an action as $\Convcomm$-algebras.
    \end{prop}

    \sssec{} 
    From the fact that the external convolution product can be identified with the external fusion product on the perverse ULA sheaves, we can lift $\textbf{Sat}^q$ to an object of $\Convcomm$. Its underlying functor is still $I\mapsto \mathbf{Sat}^q_{I}$, but given an operation
    \[\{f_u\}_{u\in U}\in\Mul_{\Convcomm}(\{I_u\}_{u\in U}, J)\]
    the corresponding structure morphism
    \begin{equation}\label{eq: Convcomm operation for Sat}
        \Ectens\mathbf{Sat}^q_{I_u}\rightarrow \mathbf{Sat}^q_{J}
    \end{equation}
    is the composition
    \[\Ectens\mathbf{Sat}^q_{I_u}\rightarrow \mathbf{Sat}^q_{\coprod I_u}\rightarrow  \mathbf{Sat}^q_{J}\]
    where the first arrow is external fusion product and the second arrow is the functoriality of Satake category with respect to $(f_u): \coprod I_u\rightarrow J$ (see Sect. \ref{sssec: A side functoriality}). The compatibilities needed for its $\Convcomm$-algebra structure were already provided by the $\Conv$-algebra structure of $\mathbf{Sat}^q$. Moreover, the action of $\mathbf{LS}^{\Ab,\ft}_?$ on $\mathbf{Sat}^q$ as $\Conv$-algebras also naturally lifts to that as $\Convcomm$-algebras, because the required compatibilities were already encoded there on the level of $\Conv$-algebras. To summarize:

    \begin{prop}
        There is a canonical lift of $\mathbf{Sat}^q$ to an object of $\Convcomm\lAlg(\Shv(X_{\et}))$, and furthermore the lift of the action of $\mathbf{LS}^{\Ab,\ft}_{?}$ on $\mathbf{Sat}^q$ as $\Conv$-algebras to an action as $\Convcomm$-algebras.
    \end{prop}

    \sssec{Sign modification}
    To have the perverse ULA Satake equivalence as an equivalence between symmetric monoidal categories, we have to alter the commutativity constraint on both sides, and thus get a symmetric monoidal equivalence
    \[\plp{\mathbf{Sat}^q_{X^I}}\simeq \elp{\textbf{Rep}}(H^I)^{\Ab,\loc,\ft}_{\cG^{\boxtimes I}_Z}.\]
    Hence we have to construct two $\Convcomm$-algebras whose underlying functor (in the sense of Lemma \ref{lem: Convcomm and Einfty}) gets identified with
    \begin{equation}\label{eq: underlying functor after alter}
        \Conv_{\langle 1\rangle}\rightarrow \CAlg(\Shv(X_{\et})),\quad I\mapsto \plp{\mathbf{Sat}^q_{X^I}}
    \end{equation}
    and
    \begin{equation}\label{eq: underlying functor B after alter}
        \Conv_{\langle 1\rangle}\rightarrow \CAlg(\Shv(X_{\et})),\quad I\mapsto \elp{\textbf{Rep}}(H^I)^{\Ab,\loc,\ft}_{\cG^{\boxtimes I}_Z}
    \end{equation}
    respectively.
    \medskip
    
    On the automorphic side, we twist the $\mathbf{LS}^{\Ab,\ft}_?$-module $\Convcomm$-algebra $\mathbf{Sat}^q$ by a $\Z/2$-grading. Each $\mathbf{Sa}t^q_{I}$ has a natural $\mathbf{LS}^{\Ab,\ft}_I$-linear $\Z/2$-grading provided by the parity of connected component, and the $\Convcomm$-operation \eqref{eq: Convcomm operation for Sat} preserves the $\Z/2$-grading. Therefore, we get a canonical $\Z/2$-graded lift of $\mathbf{Sat}^q$ as a $\mathbf{LS}^{\Ab,\ft}_?$-module in $\Shv(X_{\et};\Convcomm\lAlg)$. Now we can apply the paradigm of Sect. \ref{sssec: twist by grading} to the situation
    \[\cX=X_{\et},\quad \cC=\mathbf{LS}^{\Ab,\ft}_?\quad\te{and}\quad \bO=\Convcomm\lAlg(\CatEc).\]
    Since $E^{\times}\simeq \Omega^2(\CatEc)^{\simeq}$, a $\Convcomm$-monoidal functor
    \begin{equation}\label{eq: twisting data comm altering}
        \cG^{\Z/2}:\Z/2\rightarrow\Gamma(X_{\et}; B^2 E^{\times})
    \end{equation}
    gives rise to a twisting data. Our choice of $\cG^{\Z/2}$ is induced along $\Z/2\hookrightarrow E^{\times}$ from a certain $\E_{\infty}$-map $\Z/2\rightarrow B^2(\Z/2)$. This $\E_{\infty}$-map in terms of a symmetric monoidal extension of $\Z/2$ by $B(\Z/2)$. Its underlying monoidal category is $\cM=\Z/2\times B(\Z/2)$ and the commutativity constraint is 
    \[c_{i,j}=(-1)^{ij},\quad i,j\in \Z/2=\pi_0(\cM).\]
    
    Applying the twisting construction to $\mathbf{Sat}^q$ using the aforementioned $\Z/2$-graded lift and twisting data \eqref{eq: twisting data comm altering} we get 
    \[\plp{\mathbf{Sat}^q} \in\mathbf{LS}^{\Ab,\ft}_?\lMod(\Shv(X_{\et};\Convcomm\lAlg))\] whose underlying functor is \eqref{eq: underlying functor after alter}. Moreover, $\plp{\mathbf{Sat}^q}$ together with its module structure still lifts the $\Conv$-algebra $\mathbf{Sat}^q$ together with its module structure.
\medskip

    The sign modification on the spectral side is done completely analogously. The sign component $\epsilon$ induces a $\Z/2$-grading on each $\textbf{Rep}(H^I)^{\Ab,\loc,\ft}_{\cG^{\boxtimes I}_Z}$ (see Sect. \ref{sssec: spectral sign}), and this gives rise to a $\Z/2$-graded lift of $\textbf{Rep}(H^?)^{\Ab,\loc,\ft}_{\cG^{\boxtimes ?}_Z}$ as a $\mathbf{LS}^{\Ab,\ft}_?$-module in $\Shv(X_{\et};\Convcomm\lAlg)$. Applying the same twisting construction with $\cG^{\Z/2}$ we get
    \[\plp{\textbf{Rep}}(H^?)^{\Ab,\loc,\ft}_{\cG^{\boxtimes ?}_Z} \in \mathbf{LS}^{\Ab,\ft}_?\lMod(\Shv(X_{\et};\Convcomm\lAlg))\]
    whose underlying functor is \eqref{eq: underlying functor B after alter}.

    \begin{rmk}\label{rmk: twist by sign}
        Let us discuss in more detail how does this $\Z/2$-graded twist affect the $\Convcomm$-algebra structure. If we started with a $\Z/2$-graded $\Convcomm$-algebra $\cA$, then the twisted $\Convcomm$-algebra $\plp{\cA}$ is another lift of the underlying $\Conv$-algebra of $\cA$. From Remark \ref{rmk: lifting Conv to Convcomm} the $\Convcomm$-algebra can be understood as its underlying $\Conv$-algebra together with the collection of natural isomorphisms \eqref{eq: symmetric group repn}. Therefore, we have to explain how to twist the collection \eqref{eq: symmetric group repn} to get another data
        \[\plp{V}_{\sigma\rightarrow\tau\sigma}:F_{\sigma}\rightarrow F_{\tau\sigma}.\]
        
        For each object
        $x\in \otimes_{u\in U}\cA_{I_u}$ the natural isomorphism \eqref{eq: symmetric group repn} evaluates to an isomorphism \[F_{\sigma}(x)\xrightarrow{V_{\sigma\rightarrow\tau\sigma}(x)} F_{\tau\sigma}(x).\]
        From the $\Z/2$-grading, we have a decomposition $\cA_I=\cA^0_I\oplus \cA^1_I$ for each color $I$. Suppose we have a partition $U=U^0\coprod U^1$ and that $x$ lies in the component \[(\otimes_{u\in U^0}\cA^0_{I_u})\bigotimes (\otimes_{u\in U^1}\cA^1_{I_u}).\] Then the evaluation of $\plp{V}_{\sigma\rightarrow\tau\sigma}$ is the composition
        \[F_{\sigma}(x)\xrightarrow{V_{\sigma\rightarrow\tau\sigma}(x)} F_{\tau\sigma}(x)\xrightarrow{\on{sgn}_{U_1}(\tau)} F_{\tau\sigma}(x),\]
        where $\on{sgn}_{U_1}(\tau)\in \{\pm 1\}$ records the parity of transpositions between elements of $U^1$.
    \end{rmk}
    
    \sssec{}
    The following $\Convcomm$-version of perverse ULA Satake equivalence implies the $\Conv$-version Theorem \ref{thm: Perv ULA Satake ultimate}:
    \begin{thm}\label{thm: Perv ULA Satake Convcomm}
        We have a canonical $\mathbf{LS}^{\Ab,\ft}_{?}$-linear equivalence
        \begin{equation}\label{eq: equiv as Convcomm algs}
            \plp{\mathbf{Sat}^q}\simeq \elp{\mathbf{Rep}}(H^?)^{\Ab,\loc,\ft}_{\cG^{?}_Z}\quad \te{in} \ \Convcomm\lAlg(\Shv(X_{\et};\CatEc))
        \end{equation}
         whose evaluation at each color $I$ is the perverse ULA Satake equivalence.
    \end{thm} 

    Using the equivalence \eqref{eq: Convcomm and Einfty}, this is exactly Theorem \ref{thm: factorizable Satake} together with the functoriality in $I$, see \eqref{eq: fSet compatibility diagram} and also the remark below that.

    \subsection{Proof of Theorem \ref{thm: QLisse Satake ultimate}}\label{ssec: proof of QLisse ultimate}
    Now we prove Theorem \ref{thm: QLisse Satake ultimate}, and thus complete the construction of the natural transformation as in Theorem \ref{thm: fact Satake nat transf}.

    \begin{proof}
     Assume either $X\ne \pjl$ or that $X=\pjl$ and $\cG_Z$ is trivializable, as in the theorem. We have constructed the desired map colorwise, i.e., for each color $I$ we had the canonical $\mathbf{QLisse}^{\otimes I}$-linear map
    \[\mathbf{S}^{\QLisse}_I: \left(\mathbf{Rep}(H)^{\QLisse}_{\cG_Z}\right)^{\otimes I}\rightarrow \mathbf{Shv}^q(\Hecke^{\loc}_{I}).\]
    (see Sect. \ref{sssec: construction of S QL I}) Now we have to lift this collection of functors to a $\mathbf{QLisse}^{\otimes}$-linear $\Conv$-algebra map.
    \medskip

    Note that we have achieved the same goal in the perverse ULA setting (Theorem \ref{thm: Perv ULA Satake ultimate}), i.e., we have upgraded the collection of $\mathbf{LS}^{\Ab,\loc,\ft}_I$-linear map
    \[\mathbf{S}^{\Ab}_I: \mathbf{Rep}(H^I)^{\Ab,\ft}_{\cG^{\boxtimes I}_Z}\rightarrow \mathbf{Sat}^q_I\]
    to a $\mathbf{LS}^{\Ab,\ft}_{?}$-linear $\Conv$-algebra map.
    \medskip

    All the data we have to supply for the aforementioned lift of $\{\mathbf{S}^{\QLisse}_I\}$ are (higher) homotopies between $\mathbf{QLisse}^{\otimes I}$-linear functors out of $\left(\mathbf{Rep}(H)^{\QLisse}_{\cG_Z}\right)^{\otimes I}$. We have obtained the following three results so far:
    \begin{enumerate}
        \item the characterization of such functor in terms of $\mathbf{LS}^{\Ab,\ft}_I$-linear (or $\Ectensbig_{u\in U}\mathbf{LS}^{\Ab,\ft}_{I_u}$-linear) functors out of $\mathbf{Rep}(H^I)^{\loc}_{\cG^{\boxtimes I}_Z}$ (Lemmas \ref{lem: the characterization lemma}, \ref{lem: the characterization lemma multiple});
        \item the $\mathbf{LS}^{\Ab,\ft}_I$-linear compatibility between $\mathbf{S}^{\Ab}_I$ and $\mathbf{S}^{\QLisse}_I$ (Theorem \ref{thm: characterization for I});
        \item a coherent collection of (higher) homotopies between $\Conv$-operations in the perverse ULA setup built in Theorem \ref{thm: Perv ULA Satake ultimate}.
    \end{enumerate}
    Combining all these, the desired (higher) homotopies for Theorem \ref{thm: QLisse Satake ultimate} are supplied in a coherent manner by those in the perverse ULA setup, hence the proof.   
    \end{proof}

    \sssec{}
    Let us illustrate the logic of the proof in a simple example. Suppose we are given an operation $(\sigma,\{f_u\}_{u\in U})\in\Mul_{\Conv}(\{I_u\}_{u\in U},J)$, which we may assume $U=\{1,2,\cdots,k\}$ and $\sigma=\triv$ without loss of generality. We need to supply a canonical commutativity data of the diagram
    \begin{equation}\label{eq: Conv compatibility}
    \begin{tikzcd}
	(\textbf{Rep}(H)^{\QLisse}_{\cG_Z})^{\otimes J} & \textbf{Shv}^q(\Hecke^{\loc}_{J}) \\
	\bigotimes_{1\le u\le k} (\textbf{Rep}(H)^{\QLisse}_{\cG_Z})^{\otimes I_u}& \bigotimes_{1\le u\le k} \textbf{Shv}^q(\Hecke^{\loc}_{I_u})
	\arrow["{\mathbf{S}^{\QLisse}_J}", from=1-1, to=1-2]
	\arrow[from=2-1, to=1-1]
	\arrow["{\bigotimes \mathbf{S}^{\QLisse}_{I_u}}"', from=2-1, to=2-2]
	\arrow[from=2-2, to=1-2]
\end{tikzcd}
    \end{equation}
as $\bigotimes_{1\le u\le k}\mathbf{QLisse}^{\otimes I_u}$-modules. By Lemma \ref{lem: the characterization lemma multiple} it suffices to identify the two circuits as $\Ectensbig_{1\le u\le k}\mathbf{LS}^{\Ab,\ft}_{I_u}$-modules after pulling back along
\begin{equation}\label{eq: compatibility check map}
    \Ectensbig_{1\le u\le k} \textbf{Rep}(H^{I_u})^{\Ab,\loc,\ft}_{\cG^{\boxtimes I_u}_Z}  \rightarrow \bigotimes_{1\le u\le k} (\textbf{Rep}(H)^{\QLisse}_{\cG_Z})^{\otimes I_u}.
\end{equation}

Consider the bigger diagram
\[
\begin{tikzcd}[column sep=1.1em, row sep=1.1em]
	(\textbf{Rep}(H)^{\QLisse}_{\cG_Z})^{\otimes J} && \mathbf{Shv}^q(\Hecke^{\loc}_J) & \\
	\bigotimes_{1\le u\le k} (\textbf{Rep}(H)^{\QLisse}_{\cG_Z})^{\otimes I_u} && \bigotimes_{1\le u\le k} \textbf{Shv}^q(\Hecke^{\loc}_{I_u}) \\
	& \textbf{Rep}(H^{J})^{\Ab,\loc,\ft}_{\cG^{\boxtimes J}_Z} && \textbf{Sat}^q_{J} \\
	& \Ectensbig_{1\le u\le k} \textbf{Rep}(H^{I_u})^{\Ab,\loc,\ft}_{\cG^{\boxtimes I_u}_Z} && \Ectensbig_{1\le u\le k} \textbf{Sat}^q_{I_u}
	\arrow["{\mathbf{S}^{\QLisse}_J}", from=1-1, to=1-3]
	\arrow[from=2-1, to=1-1]
	\arrow["{\bigotimes \mathbf{S}^{\QLisse}_{I_u}}"', from=2-1, to=2-3]
	\arrow[from=2-3, to=1-3]
	\arrow[from=3-2, to=1-1]
	\arrow["{\mathbf{S}^{\Ab}_{J}}", from=3-2, to=3-4]
	\arrow[from=3-4, to=1-3]
	\arrow[from=4-2, to=2-1]
	\arrow[from=4-2, to=3-2]
	\arrow["{\Ectensbig_{1\le u\le k}\mathbf{S}^{\Ab}_{I_u}}"', from=4-2, to=4-4]
	\arrow[from=4-4, to=2-3]
	\arrow[from=4-4, to=3-4]
\end{tikzcd}
\]
We have proved that all the faces but the one on the back commute:
\begin{enumerate}
    \item both of the side faces commute tautologically; 
    \item the commutativity of the upper face and the lower face is provided by Theorem \ref{thm: characterization for I};
    \item the commutativity of the front face is encoded in Theorem \ref{thm: Perv ULA Satake ultimate}.
\end{enumerate}

Moreover, all of those commute as $\Ectensbig_{u\in U}\mathbf{LS}^{\Ab,\ft}_{I_u}$-modules. Combining the commutativity of five faces as above, the pullback of two circuits of the diagram \eqref{eq: Conv compatibility} along \eqref{eq: compatibility check map} are equipped with a canonical $\Ectensbig_{u\in U}\mathbf{LS}^{\Ab,\ft}_{I_u}$-linear equivalence.

\sssec{}
Since Theorem \ref{thm: QLisse Satake ultimate} implies Theorem \ref{thm: fact Satake nat transf}, and hence we complete the `hard part' of the construction of the natural transformation needed as an input for Theorem \ref{thm: main categorical thm}, up to taking global sections. (see the introduction of Sect. \ref{sec: metaplectic Satake functor})
This was achieved by the following chain of implications:
\[\te{Thm } \ref{thm: Perv ULA Satake Convcomm}\implies \te{Thm } \ref{thm: Perv ULA Satake ultimate} \implies \te{Thm }\ref{thm: QLisse Satake ultimate}\implies \te{Thm }\ref{thm: fact Satake Conv at X}\implies \te{Thm }\ref{thm: fact Satake nat transf},\]
where in the second arrow the characterization lemmas (Lemmas \ref{lem: the characterization lemma}, \ref{lem: the characterization lemma multiple}) were needed, and in the last arrow the categorical mechanism of Proposition \ref{prop: module in Conv to fSet Alg} was needed.

\section{Factorizable action on $\cA^q_G$ arising from geometric Satake}\label{sec: factorizable action}
In this section, we finish the construction of the spectral action in the Betti metaplectic setting:
\begin{thm}\label{thm: spectral action}
    There is a canonical action of $\QC(\Loc^{\cG_Z}_H)$ on $\cA^q_G$.
\end{thm}

\sssec{}
In the introduction of Sect. \ref{sec: metaplectic Satake functor} we have argued that it suffices to construct a natural transformation $\eqref{eq: funct CX}\rightarrow \eqref{eq: funct AG}$ between functors
$\fSet\rightarrow \fSet+\LSMon$, which can be thought of as a factorizable action of the local system $\cL_{(X^{\topl,\sing},\Phi(\cG_Z),\Rep(H))}$ \footnote{See the notation from Sect. \ref{sssec: examples from gerbes}.} of symmetric monoidal DG categories on $\cA^q_G$.
\medskip

After taking global section, the functor  \eqref{eq: funct CX} turns into the functor 
\begin{equation}\label{eq: source global}
    \fSet\rightarrow\fSet+\Alg(\LS(X)^{\otimes}\lMod),\quad I\mapsto (\Rep(H)^{\loc}_{\cG_Z}(X))^{\otimes I}
\end{equation}
and the functor \eqref{eq: funct AG} turns into the functor
\begin{equation}\label{eq: target global}
    \fSet\rightarrow\fSet+\Alg(\LS(X)^{\otimes}\lMod),\quad I\mapsto \End(\cA^q_G)\otimes\LS(X)^{\otimes I}.
\end{equation}
Through Sect. \ref{sec: metaplectic Satake functor} and \ref{sec: Satake as a nat transf} we were proceeding towards the construction of a natural transformation $\eqref{eq: source global}\rightarrow \eqref{eq: target global}$ arising from the metaplectic geometric Satake, and the `hard part' was done in those sections. In this section, we do the other part to complete the construction.

\sssec{}
In general, we lose information in the process of taking global sections, and hence a natural transformation on the level of global sections does not give rise to a natural transformation on the level of local systems. However, in our situation of passing from $\eqref{eq: funct CX}\rightarrow\eqref{eq: funct AG}$ to $\eqref{eq: source global}\rightarrow\eqref{eq: target global}$, we can use some special features of the underlying topology of $X$ to prove that we do not lose any information. (see Proposition \ref{prop: equiv LS global}) Therefore, we prove Theorem \ref{thm: spectral action} once we construct $\eqref{eq: source global}\rightarrow\eqref{eq: target global}$.

\sssec{Extension of an action}
The desired natural transformation $\eqref{eq: source global}\rightarrow \eqref{eq: target global}$ when evaluated at the color $I$ takes a form
\[\left(\Rep(H)^{\loc}_{\cG_Z}(X)\right)^{\otimes I}\rightarrow \End(\cA^q_G)\otimes\LS(X)^{\otimes I}\quad\te{in }\Alg(\LS(X)^{\otimes I}\lMod),\]
and such data is equivalent to the extension of the free $\LS(X)^{\otimes I}$-action on $\cA^q_G\otimes\LS(X)^{\otimes I}$ to an action of $(\Rep(H)^{\loc}_{\cG_Z}(X))^{\otimes I}$, along the symmetric monoidal map $\LS(X)^{\otimes I}\rightarrow (\Rep(H)^{\loc}_{\cG_Z}(X))^{\otimes I}$. This is a variant of \cite[Sect. B.5.5]{GKRV}. 
\medskip

Let us briefly explain where the first action comes from. Combining the metaplectic Satake functor with the metaplectic Hecke action we get a canonical action of $(\Rep(H)^{\loc}_{\cG_Z}(X))^{\otimes I}$ on $\Shv^q(\Bun_G\times X^I)$. In Sect. \ref{ssec: nilpotent singular support} we prove that $\cA^q_G\otimes \LS(X)^{\otimes I}$ is a full subcategory of $\Shv^q(\Bun_G\times X^I)$ and moreover is preserved under the action of $(\Rep(H)^{\loc}_{\cG_Z}(X))^{\otimes I}$. These results are metaplectic analogs of \cite[Theorem 6.1.1]{NY} and \cite[Theorem A.3.9,  B.5.2]{GKRV} and the proof is almost identical. The $\LS(X)^{\otimes I}$-linearity of the constructions involved yields an obvious extension of the action as described above.

\sssec{}
To get the desired natural transformation, we need a `smarter version' of this argument which takes care of homotopy coherences as we vary $I$. We faced a similar kind of challenge when we upgraded the collection $\{\cS_I\}_{I\in\fSet}$ of Satake functor to a natural transformation, and we introduced the operads $\Conv$ and $\Convcomm$ to resolve this. Similarly, we introduce operads $\ConvMod$ and $\ConvcMod$ in Sect. \ref{ssec: an operad ConvMod} to achieve the aforementioned goal.
\medskip

We highlight the key features of the operad $\ConvMod$.
\begin{enumerate}
    \item We can think of a $\Conv$-algebra as a generalization of an associative algebra. In that aspect, a $\ConvMod$-algebra is encoding a $\Conv$-algebra together with its right module.
    \item As $\Conv$-operations encode the external convolution products of Hecke categories, $\ConvMod$-operations encode the Hecke action on $\Shv^q(\Bun_G)$.
\end{enumerate}

Because of the first feature, we write a $\ConvMod$-algebra as $(\cM;\cA)$ where $\cM$ is the `module part' and $\cA$ is the `algebra part'. $\cA$ itself is a $\Conv$-algebra. We refer readers to Sect. \ref{sssec: ConvMod algebra object} for details. In Sect. \ref{ssec: ConvMod algebra from Hecke}, we construct a canonical $\ConvMod$-algebra encoding the Hecke action from geometry, in the way similar to  Sect. \ref{ssec: Conv algebras from fact Satake}.
\medskip

 Through Sect. \ref{ssec: ConvMod algebras and Nilp} and \ref{ssec: passage ConvcMod nat transf} we realize the `smarter version' of the extension of action argument above. In Sect. \ref{ssec: ConvMod algebras and Nilp}, we construct a canonical functor 
 \[(\cA^q_G\otimes\LS(X)^{\otimes};\LS(X)^{\otimes})\rightarrow (\cA^q_G\otimes\LS(X)^{\otimes};(\Rep(H)^{\loc}_{\cG_Z}(X))^{\otimes})\]
 between $\ConvcMod$-algebras. Indeed, this can be understood as the extension of the free $\LS(X)^{\otimes}$-action on $\cA^q_G\otimes\LS(X)^{\otimes}$ to an action of $(\Rep(H)^{\loc}_{\cG_Z}(X))^{\otimes}$, along the $\Convcomm$-algebra map $\LS(X)^{\otimes}\rightarrow (\Rep(H)^{\loc}_{\cG_Z}(X))^{\otimes}$. In Sect. \ref{ssec: passage ConvcMod nat transf} we develop a categorical mechanism to pass from such a data to a natural transformation $\eqref{eq: source global}\rightarrow\eqref{eq: target global}$. This section is analogous to Sect. \ref{ssec: passage to nat transf}.

 \sssec{}
 In Sect. \ref{ssec: construction of the spectral action} we put all the ingredients together and prove Theorem \ref{thm: spectral action}.

\subsection{Nilpotent singular support}\label{ssec: nilpotent singular support}
In this subsection, we prove that the Hecke action provided by the factorizable metaplectic geometric Satake functors $\cS_{X,I}$
preserves the nilpotent singular support condition on $\Bun_G$.

\begin{note}
    Our convention is that the algebraic geometry on the automorphic side is always a classical complex algebraic geometry, see Sect. \ref{sssec: AG convention}. Throughout this subsection, $\cY$ is a classical complex algebraic stack.
\end{note}

\sssec{}\label{sssec: singsupp review}
We use the theory of singular support developed in \cite[Sect. A.3]{GKRV}. In \emph{loc. cit.}, singular support is considered as a purely set-theoretic notion. In particular, a singular support condition is generally determined by a conical Lagrangian \emph{subset} of $T^*\cY$.
 \medskip
 
 Given a conic Lagrangian subset $\cN\subset T^*\cY$ we have a subset $\cN(\C)\subset (T^*\cY)(\C)$ of $\C$-points defined as
 \[\cN(\C)=\{(y,\xi): y\in \cY(\C), \ \xi\in T^*_y\cY\cap \cN\},\]
 and the set $\cN(\C)$ completely determines $\cN$. Thus, we can restrict our concern to $\C$-points when addressing questions about singular support.

\sssec{Singular support for twisted sheaves}
We briefly indicate how singular support is defined for twisted sheaves.
Suppose we are given $\cY$, an abelian group $A$, a gerbe $\cG\in\on{Ge}_A(\cY)$, and a character $\chi:A\hookrightarrow E^{\times}$. As in Sect. \ref{sssec: twisted sheaves} we obtain the category $\Shv^{\cG}(\cY)$ of twisted shaves over $\cY$. Suppose $\cF\in \Shv^{\cG}(\cY)$. To specify the singular support $\cN\subset T^*\cY$ of $\cF$, it suffices to specify the subset $\cN(\C)\subset (T^*\cY)(\C)$. For a point $(y,\xi)\in (T^*\cY)(\C)$ we choose an \'etale neighborhood $(\mc{U},u)\rightarrow (\cY,y)$ of $y$ so that the pullback of $\cG$ to $\mc{U}$ is trivializable. Upon choosing the trivialization, we get a sheaf $\cF|_{\mc{U}}\in \Shv(\mc{U})$. Let $\cN_{\mc{U}}\subset T^*\mc{U}$ be the singular support of $\cF|_{\mc{U}}$. We let
\[(y,\xi)\in \cN(\C) \iff (u,\xi)\in \cN_{\mc{U}}(\C),\]
where on the right-hand side we used the isomorphism $T^*_y\cY\simeq T^*_u\mc{U}$. Whether a given point $(y,\xi)\in (T^*\cY)(\C)$ lies in $\cN(\C)$ or not does not depend on the choice of \'etale neighborhood $\mc{U}$ and the trivialization of $\cG$ on $\mc{U}$.

\sssec{}\label{sssec: sing supp via total space}
Here is another way of defining a singular support. Let $\on{Tot}(\cG)$ be the total space of the gerbe $\cG$. This space is equipped with a natural action of $\Bet A$. The category $\Shv(\Bet A)$ is equipped with the convolution symmetric monoidal structure, and acts on $\Shv(\on{Tot}(\cG))$ via convolution $\ast$. Furthermore, $\Vect_{-\chi}$ is the module category of $\Shv(\Bet A)$ whose underlying DG category is $\Vect$ and the action is given by the symmetric monoidal
\begin{equation}\label{eq: chi eigenfunctor}
    (\Shv(\Bet A),\ast) = (\Rep(A),\ast)=(\QC(\Spec(EA)),\otimes)\xrightarrow{-\chi}\Vect.
\end{equation}
In Appendix \ref{sec: twisted sheaves} we construct a canonical equivalence
\begin{equation}\label{eq: twisted sheaves tensor prod}
    \Shv^{\cG}(\cY)=\Shv(\on{Tot}(\cG))\otimes_{\Shv(\Bet(A))}\Vect_{-\chi}
\end{equation}
    in $\DGcat$. Here 
where the last arrow comes from the pullback along $-\chi$, interpreted as an $E$-point of $\Spec(EA)$. The equivalence \eqref{eq: twisted sheaves tensor prod} realizes $\Shv^{\cG}(\cY)$ as a full subcategory of $\Shv(\on{Tot}(\cG))$.
\medskip

Given a sheaf $\cF\in\Shv^{\cG}(\cY)$ we view it as a sheaf over $\on{Tot}(\cG)$ by the above procedure, and consider the singular support $\cN_{\Tot}\subset T^*(\Tot(\cG))$. Since $\on{Tot}(\cG)$ is a relative Deligne-Mumford stack over $\cY$, the natural projection induces a bijection 
\[(T^*\Tot(\cG))(\C)\xrightarrow{1:1} (T^*\cY)(\C).\]
Therefore we obtain a subset of $(T^*\cY)(\C)$, and this defines the singular support $\cN\subset T^*\cY$ of $\cF\in\Shv^{\cG}(\cY)$.

\sssec{}
Suppose $\cG$ is an $A$-gerbe over $\cY$. Given a conical Zariski-closed subset $\cN\subset T^*\cY$ we consider the full subcategory 
\[\Shv^{\cG}_{\cN}(\cY)\subset \Shv^{\cG}(\cY)\]
of twisted sheaves whose singular support is contained in $\cN$.

\sssec{The metaplectic automorphic category}\label{sssec: AqG}
As in \cite[Sect. 4.2]{BZN} we consider the conical subset $\Nilp \subset T^*\Bun_G$ called the \emph{global nilpotent cone}, defined as the zero fiber of the Hitchin fibration. Then the \textbf{metaplectic automorphic category} is defined as
\[\cA^q_G:=\Shv^q_{\Nilp}(\Bun_G),\]
i.e., the category of $\cG^{\on{glob}}$-twisted sheaves on $\Bun_G$ whose singular support is contained in $\Nilp$. We often use the term `nilpotent singular support' for the condition that singles out $\cA^q_G$ inside of $\Shv^q(\Bun_G)$.

\sssec{}
Thanks to the above two descriptions of singular support, the usual functorial properties of singular support carry over to twisted sheaves. More precisely, we have the following twisted analog of \cite[Sect. 5.4]{KS}:
\begin{prop}\label{prop: SingSupp functoriality}
    Suppose we have a map $f:(\cY_1,\cG_1)\rightarrow (\cY_2,\cG_2)$ in $\Stklft+\on{Ge}_A$ where both $\cY_i$'s are algebraic stacks. We have a correspondence
    \[T^*(\cY_1)\xleftarrow{df} T^*(\cY_2)\times_{\cY_2} \cY_1 \xrightarrow{f_{\nat}}T^*(\cY_2).\]
    \begin{enumerate}
        \item If $f$ is smooth and $\cF\in\Shv^{\cG_2}(\cY_2)$, we have 
        \[\on{SingSupp}(f^!\cF)=(df) (f_{\nat})^{-1} \on{SingSupp}(\cF).\]
        \item If $f$ is proper and $\cF\in\Shv^{\cG_1}(\cY_1)$, we have
        \[\on{SingSupp}(f_!\cF)\subseteq (f_{\nat})(df)^{-1}\on{SingSupp}(\cF).\]
    \end{enumerate}
\end{prop}

\sssec{Tensor product formula}
Using the description of Sect. \ref{sssec: sing supp via total space} we get the following twisted analog of \cite[Theorem A.3.9]{GKRV}: 

\begin{prop}
    Suppose $Y$ is a smooth scheme and $\cG$ is an $A$-gerbe on $\cY$. We get an $A$-gerbe on $\cY\times Y$ by pulling back from $\cY$, and hence the category $\Shv^{\cG}(\cY\times Y)$ of twisted sheaves. Suppose $\cN\subset T^*\cY$ is a conical Zariski closed subset of half-dimension, and let
    \[\cN\times Y\subset T^*(\cY\times Y)\]
    be the product of $\cN$ and the zero section $Y\subset T^*Y$.
    The external tensor product functor
    \begin{equation}\label{eq: ext tens prod zero section}
        \Shv^{\cG}_{\cN}(\cY)\otimes \LS(Y)\rightarrow \Shv^{\cG}_{\cN\times Y}(\cY\times Y)
    \end{equation}
    is an equivalence.
\end{prop}

\begin{rmk}
    The external tensor product functor
    \[\Shv^{\cG}(\cY)\otimes \Shv(Y)\rightarrow \Shv^{\cG}(\cY\times Y)\]
    restricts to \eqref{eq: ext tens prod zero section}. Indeed, this follows from the functoriality of singular support along smooth pullbacks (see Proposition \ref{prop: SingSupp functoriality}) together with the fact that tensoring with a non-trivial local system does not change the singular support.
\end{rmk}

\sssec{}
For a finite set $I$, let $\Shv^q(\Bun_G\times X^I)$ be the category of twisted sheaves on $\Bun_G\times X^I$, twisted by the pullback of $\cG^{\on{glob}}$ along the projection map.

\begin{cor}\label{cor: singling out AqG otimes LS}
    For a finite set $I$, the external tensor product functor
    \[\cA^q_G\otimes \LS(X^I)\rightarrow \Shv^q(\Bun_G\times X^I)\]
    is fully-faithful, and its essential image is $\Shv^q_{\Nilp\times X^I}(\Bun_G\times X^I)$.
\end{cor}

\sssec{Preservation of nilpotent singular support}
The key reason why one should impose nilpotent singular support condition for (non-metaplectic) Betti Langlands is that the Hecke action - which needs a choice of a point $x$ - depends locally constantly on $x\in X$. A sheaf $\cF\in\Shv(Y)$ being locally constant is equivalent to the sheaf having singular support equal to the zero section $Y\subset T^*Y$. The idea of \cite{NY} (and followed by \cite{GKRV, AGKRRV}) is to prove the local constancy by showing that the Hecke correspondence introduces only the zero section along the direction of the curve when restricted to sheaves with nilpotent singular support. The latter fact is called the \textbf{preservation of nilpotent singular support}.
\medskip

This theme of preservation of nilpotent singular support and the resulting local constancy of Hecke action carries over to the metaplectic case. We will cover those in Sect. \ref{sssec: preservation metaplectic}.

\sssec{Non-metaplectic case}\label{sssec: preservation non-metaplectic}
In the context of non-metaplectic geometric Langlands, we have the Hecke correspondence
\begin{equation}\label{eq: Hecke correspondence for I}
\begin{tikzcd}
	& {\Hecke^{\loc}_{X^I}} & \\
	{\Bun_G\times X^I} & {\Hecke_{X^I}} & {\Bun_G\times X^I,}
	\arrow["{\mf{r}}"', from=2-2, to=1-2]
	\arrow["{p_1}"', from=2-2, to=2-1]
	\arrow["{p_2}", from=2-2, to=2-3]
\end{tikzcd}
\end{equation}
such that given a point $(\ul{x},\cP^0,\cP^1,\alpha)$ as in the notation of Sect. \ref{sssec: stacks on the A side}, it is sent to $(\cP^0,\ul{x})$ under $p_1$ and sent to $(\cP^1,\ul{x})$ under $p_2$, and $\mf{r}$ is given by the restriction to $\cD_{\ul{x}}$. This gives rise to the functor
\[\Shv(\Bun_G\times X^I)\otimes \Shv(\Hecke^{\loc}_{X^I})\rightarrow \Shv(\Bun_G\times X^I),\]
which naturally gets enhanced to the action of $\Shv(\Hecke^{\loc}_{X^I})$ on $\Shv(\Bun_G\times X^I)$. Combined with the geometric Satake functor
\begin{equation}\label{eq: Satake non metaplectic}
    \Rep(\chG)^{\otimes I}\rightarrow \Shv(\Hecke^{\loc}_{X^I})
\end{equation}
we get an action of $\Rep(\chG)^{\otimes I}$ on $\Shv(\Bun_G\times X^I)$.

\begin{thm}[\cite{NY} Theorem 5.1,\cite{GKRV} Theorem B.5.2]
The action of $\Rep(\chG)^{\otimes I}$ preserves the full subcategory $\Shv_{\Nilp\times X^I}(\Bun_G\times X^I)$.
\end{thm}

\sssec{Metaplectic case}\label{sssec: preservation metaplectic}
The correspondence \eqref{eq: Hecke correspondence for I} has a natural compatibility with the gerbes on each stack. More precisely, on $\Bun_G\times X^I$ we have a pullback of $\cG^{\on{glob}}$ which we denote by $\pi^*\cG^{\on{glob}}$ and on $\Hecke^{\loc}_{X^I}$ we have $\cG^I$. Then we have a canonical isomorphism
\[p_1^*(\pi^*\cG^{\on{glob}})+ \mf{r}^*\cG^I \simeq p^*_2(\pi^*\cG^{\on{glob}}),\]
and hence a functor
\[\Shv^q(\Bun_G\times X^I)\otimes \Shv^q(\Hecke^{\loc}_{X^I})\rightarrow \Shv^q(\Bun_G\times X^I).\]
This naturally gets enhanced to the action of $\Shv^q(\Hecke^{\loc}_{X^I})$ on $\Shv^q(\Bun_G\times X^I)$. Combined with the $\LS(X^I)$-linear metaplectic Satake functor
\[\left(\Rep(H)^{\loc}_{\cG_Z}(X)\right)^{\otimes I}\rightarrow \Shv^q(\Hecke^{\loc}_{X^I})\]
(see Theorem \ref{thm: fact Satake nat transf}) we get a canonical action of $\left(\Rep(H)^{\loc}_{\cG_Z}(X)\right)^{\otimes I}$ on $\Shv^q(\Bun_G\times X^I)$. The following assertion follows verbatim from the proof of \cite[Theorem B.5.2]{GKRV}:

\begin{thm}\label{thm: preservance Nilp metaplectic}
    The action of $\left(\Rep(H)^{\loc}_{\cG_Z}(X)\right)^{\otimes I}$ on $\Shv^q(\Bun_G\times X^I)$ preserves the full subcategory $\Shv^q_{\Nilp\times X^I}(\Bun_G\times X^I)$.
\end{thm}

\sssec{}
We briefly comment on the proof of the above theorem. First, we reduce to $I=\{*\}$ using the associativity of the Hecke action. The proof of its non-metaplectic version given in \cite{GKRV} does not use something specific about the Satake kernels, i.e., the objects in the essential image of \eqref{eq: Satake non metaplectic}. The proof only used the fact that if $X=\afl$,\footnote{Here we temporarily break the projectivity assumption on $X$, but the Satake functor is defined for any smooth curve $X$.} the Satake kernels spread constantly along the curve direction of $\Hecke^{\loc}_{\afl}$. More precisely, consider the natural projection $\Hecke^{\loc}_{\afl}\rightarrow \afl$ and the fiber $\Hecke^{\loc}_0$ at $0\in \afl$. Translation along $\afl$ gives rise to a canonical identification $\Hecke^{\loc}_{\afl}= \Hecke^{\loc}_0\times \afl$, and every Satake kernel arises as a pullback from a sheaf on $\Hecke^{\loc}_0$ along the projection. For convenience, let us call this property the `constancy along $\afl$'.
    \medskip

Now we can see the reason why the same proof carries over to the metaplectic case. Since twisted sheaves enjoy the same functoriality as sheaves (see Proposition \ref{prop: SingSupp functoriality}) we use the same argument as in the non-metaplectic case to reduce Theorem \ref{thm: preservance Nilp metaplectic} to showing the constancy of metaplectic Satake kernels along $\afl$. If $X=\afl$, the metaplectic parameter is constant along the curve, i.e., a pullback of $q_{\pt}\in \Gamma_e(\Bet G; B^4 A(1))$. Suppose $\cG_{Z,\pt}\in \on{Ge}_{Z_H(E)}(\pt)$ be the metaplectic dual gerbe corresponding to $q_{\pt}$. From the functoriality of the construction of metaplectic dual gerbe, $\cG_Z$ is the pullback  of $\cG_{Z,\pt}$ along $X\rightarrow \pt$, possibly up to the factor of the theta-characteristic $\omega^{1/2}_X$. Upon choosing a square root of $\omega_X$ we have 
\[\Rep(H)^{\loc}_{\cG_Z}(\afl)\simeq \Rep(H)_{\cG_{Z,\pt}}\otimes \LS(\afl)= \Rep(H)_{\cG_{Z,\pt}}.\]
Combining the above, we conclude the constancy of metaplectic Satake kernels along $\afl$.

\subsection{An operad $\ConvMod$}\label{ssec: an operad ConvMod}
We introduce an operad $\ConvMod$ which encodes the notion of a `right module' for a $\Conv$-algebra. This is an ordinary colored operad, and we use the notations of \cite[Definition 2.1.1.1]{HA} to describe it. It should be clear from the definition below that this is a hybrid of the operads $\LM$ and $\Conv$.
\begin{defn}\label{def: operad ConvMod}
    An operad $\ConvMod$ is an ordinary colored operad described as follows.
    \begin{enumerate}
        \item The collection of colors is the disjoint union 
        \[ \{I^{\mf{m}}\}_{I\in\fSet}\coprod \{I^{\mf{a}}\}_{I\in\fSet}.\]
        The subscript $\mf{a}$ indicates that the color takes the role of an algebra while the subscript $\mf{m}$ indicates that the color takes the role of a module.\footnote{This notation is borrowed from the notation for colors in the definition of the $\LM$-operad, see \cite[Sect. 4.2.1]{HA}.} For this reason, a color $I^{\mf{a}}$ will be said to be an `algebra type' while a color $I^{\mf{m}}$ will said to be a `module type'. Note that we have a forgetful map
        \[\psi:\{\te{colors of }\ConvMod\}\rightarrow \{\te{colors of }\Conv\},\quad I^{\mf{m}}, I^{\mf{a}}\mapsto I.\]
        
        \item Suppose we are given a finite set $U$, an $U$-indexed collection of colors $\{\cI_u\}_{u\in U}$, and a color $\cJ$. The set of operations $\Mul_{\ConvMod}(\{\cI_u\}_{u\in U}, \cJ)$ is a subset of $\Mul_{\Conv}(\{\psi(\cI_u)\}_{u\in U},\psi(\cJ))$ determined as follows:
        \medskip
        
        - Suppose that all of $\cI_u$'s are of algebra type. If $\cJ$ is of a module type, then the set of operations is empty. If $\cJ$ is of an algebra type, we have
        \[\Mul_{\ConvMod}(\{\cI_u\}_{u\in U}, \cJ)=\Mul_{\Conv}(\{\psi(I_u)\}_{u\in U}, \psi(J)).\]
        Such operations are said to be of `multiplication type'.
    \medskip
    
        - Suppose that exactly one of $\cI_u$'s is of a module type, say  $\cI_{u_0}$ is the only module type. Then we have
        \[\Mul_{\ConvMod}(\{\cI_u\}_{u\in U}, \cJ)=\left\{(\sigma,\{f_u\}_{u\in U})\in \Mul_{\Conv}(\{\psi(\cI_u)\}_{u\in U}, \psi(\cJ)): u_0=\on{min}(U;\sigma)\right\},\]
        where $\on{min}(U;\sigma)$ means the minimum of $U$ considered as an ordered set subject to linear order $\sigma$. Such operations are said to be an `action type'.
        \medskip

        - The set of operations is empty otherwise.
        
        \item The composition law of $\Conv$-operations preserves the subset of $\ConvMod$-operations, hence defining the composition law of $\ConvMod$.
    \end{enumerate}
\end{defn}

\sssec{}\label{sssec: ConvMod algebra object}
Given a symmetric monoidal category $\bC$, a data of a $\ConvMod$-algebra in $\bC$ will be notated as $(\cM;\cA)$, and it consists of the following data: 
\begin{enumerate}
    \item a collection $\{\cA_I\}_{I\in\fSet}\subseteq \bC$ of objects, each of which is an evaluation of $(\cM;\cA)$ at each algebra type color;
    \item a collection $\{\cM_I\}_{I\in\fSet}\subseteq \bC$ of objects, each of which is an evaluation of $(\cM;\cA)$ at each module type color;
    \item for a multiplication type operation $(\sigma,f)\in \Mul_{\ConvMod}(\{\cI_u\}_{u\in U},\cJ)$ as in (2) of Definition \ref{def: operad ConvMod}, there should be an associated map $\bigotimes_{u\in U}\cA_{\psi(\cI_u)}\xrightarrow{(\sigma,f)}\cA_{\psi(\cJ)}$
    \item for an action type operation $(\sigma,f)\in \Mul_{\ConvMod}(\{\cI_u\}_{u\in U},\cJ)$ as in (2) of Definition \ref{def: operad ConvMod}, there should be an associated map $\cM_{\psi(\cI_{u_0})}\otimes\bigotimes_{u\in U\setminus u_0}\cA_{\psi(\cI_u)}\xrightarrow{(\sigma,f)}\cM_{\psi(\cJ)}$, and the collection of such structure maps is supplied with higher coherence data.
\end{enumerate}

Note that we have a natural operad map $\Conv\rightarrow\ConvMod$ and hence a forgetful map $\ConvMod\lAlg(\bC)\rightarrow \Conv\lAlg(\bC)$. In the above notation, the second argument $\cA$ in $(\cM;\cA)$ specifies what the underlying $\Conv$-algebra is.

\sssec{Building blocks}\label{sssec: building blocks}
For $I\rightarrow J$ in $\fSet$ the $\ConvMod$-operation gives rise to maps
\[\cA_I\rightarrow \cA_J,\quad \cM_I\rightarrow \cM_J,\]
which we call the \emph{functoriality maps}. We also have an action type operation \[\cM_I\otimes \cA_I\rightarrow \cA_I\] corresponding to the case $U=\{1,2\}$ and $\sigma$ is the trivial linear ordering, and $f_1=f_2=\id_I$.
    \medskip

    In fact, any action type $\ConvMod$-operation can be expressed as a combination of the above three operations. Indeed, using the composition  law of $\ConvMod$-operations, each action type map
\[\cM_{\psi(\cI_{u_0})}\otimes\bigotimes_{u\in U\setminus u_0}\cA_{\psi(\cI_u)}\xrightarrow{(\sigma,f)}\cM_{\psi(\cJ)}\]
can be expressed as the composition
\[\cM_{\psi(\cI_{u_0})}\otimes\bigotimes_{u\in U\setminus u_0}\cA_{\psi(\cI_u)}\rightarrow \cM_{\psi(\cJ)}\otimes\bigotimes_{u\in U\setminus u_0}\cA_{\psi(\cJ)}\rightarrow\cM_{\psi(\cJ)}\]
where the first map comes from the functoriality maps and the second map is the iteration of $\cM_{\psi(\cJ)}\otimes \cA_{\psi(\cJ)}\rightarrow \cA_{\psi(\cJ)}$.

\sssec{} 
We also introduce a commutative variant $\ConvcMod$, which is an ordinary colored operad encoding the notion of a module for a $\Convcomm$-algebra.

\begin{defn}
    An operad $\ConvcMod$ is an ordinary colored operad described as follows.
    \begin{enumerate}
        \item Its collection of colors is the same as that of $\ConvMod$.
        \item Suppose we are given a finite set $U$, an $U$-indexed collection of colors $\{\cI_u\}_{u\in U}$, and a color $\cJ$, we describe the set of operations $\Mul_{\ConvcMod}(\{\cI_u\}_{u\in U}, \cJ)$. We divide into three cases:
        
        - if all of $\cI_u$'s and $\cJ$ are algebra type;
    
        - if $\cJ$ is module type and exactly one of $\cI_u$'s is a module type;

        - all the other cases.
        
        The set of operations is empty in the third case. For both the first and second case we have
        \[\Mul_{\ConvcMod}(\{\cI_u\}_{u\in U},\cJ)=\Mul_{\Convcomm}(\{\psi(\cI_u)\}_{u\in U}, \psi(\cJ)).\]
        
        \item The composition law of $\Convcomm$-operations preserves the subset of $\ConvcMod$-operations, hence defining the composition law of $\ConvcMod$.
    \end{enumerate}
\end{defn}
\begin{rmk}
    The terminologies `module type' and `algebra type' for colors and `action type' and `multiplication type' for $\ConvMod$-operations have obvious analogs for $\ConvcMod$.
\end{rmk}

\sssec{}
    We have the following analog of assertions in Sect. \ref{ssec: an operad Conv}:
    
    \begin{lemma}\label{lem: ConvcMod and LMcomm}
    For each $I\in\Conv_{\langle 1\rangle}$ there is a canonical map $\LMcomm\rightarrow \ConvcMod$ between operads that sends a color $\mf{c}\in \LMcomm$ to the algebra type color $I^{\mf{a}}\in\ConvcMod$ and a color $\mf{m}\in\LMcomm$ to the module type color $I^{\mf{m}}\in \ConvcMod$. Therefore, given a $\ConvcMod$-algebra $(\cM;\cC)$ in $\bC$, each $(\cC_I,\cM_I)$ has a natural lift to $\LMcomm\lAlg(\bC)$. This construction is functorial in $I$, and hence we get a functor $\Conv_{\langle 1\rangle}\rightarrow \LMcomm\lAlg(\bC)$, which maps $I\mapsto (\cC_I,\cM_I)$.
\end{lemma}

\begin{prop}
    The functor
    \begin{equation}\label{eq: ConvcMod and LMcomm}
        \ConvcMod\lAlg(\bC)\rightarrow \Fun(\Conv_{\langle 1\rangle},\LMcomm\lAlg(\bC))
    \end{equation}
    constructed in Lemma \ref{lem: ConvcMod and LMcomm} is an equivalence.
\end{prop}

\begin{cor}
    We have a canonical equivalence $\ConvcMod\simeq \Conv_{\langle 1\rangle}\otimes\LMcomm$ of operads.
\end{cor}

\sssec{}\label{sssec: lifting ConvMod to ConvcMod}
We have a commutative diagram of maps between operads
\[\begin{tikzcd}
	\ConvcMod & \ConvMod \\
	\Convcomm & \Conv.
	\arrow[from=1-2, to=1-1]
	\arrow[from=2-1, to=1-1]
	\arrow[from=2-2, to=1-2]
	\arrow[from=2-2, to=2-1]
\end{tikzcd}\]
Given a symmetric monoidal category $\bC$ we have a resulting commutative diagram
\begin{equation}\label{eq: ConvcMod and ConvMod}
    \begin{tikzcd}
	{\ConvcMod(\bC)} & {\ConvMod(\bC)} \\
	{\Convcomm(\bC)} & {\Conv(\bC).}
	\arrow[from=1-1, to=1-2]
	\arrow[from=1-1, to=2-1]
	\arrow[from=1-2, to=2-2]
	\arrow[from=2-1, to=2-2]
\end{tikzcd}
\end{equation}
This diagram is a Cartesian diagram: this is an analog of the fact that the notion of a module of an algebra only depends on the $\E_1$-structure of the algebra.

\sssec{}
In what follows, we list some constructions of $\ConvMod$-algebras analogous to constructions of modules for associative algebras. Each of those three has a natural commutative analog for $\ConvcMod$-algebras.

\sssec{Pullback module}\label{sssec: pullback of ConvMod}
Suppose that we have a $\ConvMod$-algebra $(\cM;\cA)$ and a $\Conv$-algebra map $\cA'\rightarrow\cA$. This naturally induces a $\ConvMod$-algebra $(\cM;\cA')$, equipped with a $\ConvMod$-algebra map $(\cM;\cA')\rightarrow (\cM;\cA)$. This is an analog of pulling back the module structure.

\sssec{Free module}\label{sssec: free module construction}
Suppose $\bC$ is a symmetric monoidal category, $M\in\bC$, and $\cA\in\Conv\lAlg(\bC)$. Then we have a canonical object $(M\otimes \cA;\cA)\in\ConvMod\lAlg(\bC)$ whose evaluation at each color is
\[(M\otimes \cA;\cA)_{I^{\mf{a}}}= \cA_I,\quad (M\otimes \cA;\cA)_{I^{\mf{m}}}=M\otimes\cA_{I}\quad \te{for } I\in\fSet, \]
and $\ConvMod$-operations are defined in an obvious way using the $\Conv$-operations of $\cA$. This can be thought of as a `free module' construction for $\Conv$-algebras. This yields a functor
\[\bC\times \Conv\lAlg(\bC)\rightarrow \ConvMod\lAlg(\bC),\quad M,\cA\mapsto (M\otimes\cA;\cA).\]
As a special case, let $\cA=\cC^{\otimes}$ for an object $\cC\in\CAlg(\bC)$. We get a $\ConvMod$-algebra $(M\otimes\cC^{\otimes};\cC^{\otimes})$. 

\sssec{Submodule}\label{sssec: submodule construction}
Suppose we have a $\ConvMod$-algebra $(\cM;\cA)$ in $\DGcat$. Suppose that for each $I\in\fSet$ we have a full sub-DG category $\cM'_{I}\subseteq \cM_I$, such that the operations $\cM_I\rightarrow \cM_J$ and $\cM_I\otimes \cA_I\rightarrow \cM_I$ from Sect. \ref{sssec: building blocks} restrict to $\cM'_I\rightarrow\cM'_J$ and $\cM'_I\otimes \cA_I\rightarrow \cM'_I$. In this case, we have a $\ConvMod$-algebra $(\cM';\cA)$ whose evaluation at algebra type colors are $\cA_I$'s but evaluation at module type colors are $\cM'_I$'s. Since any action type $\ConvMod$-operation can be factorized into the combination of the operations
\[\cM_I\rightarrow \cM_J,\quad \cA_I\rightarrow \cA_J, \ \te{ and } \ \cM_I\otimes \cA_I\rightarrow \cM_I,\]
 the collection of subcategories $\{\cM'_I\}$ is preserved under any action type $\ConvMod$-action. This can be thought as: $\cM'$ is the `submodule' of $\cM$ under the $\cA$-action.

\sssec{}\label{sssec: ConvMod to Conv1}
Finally, we observe that there is a functor
\begin{equation*}
    \ConvMod\lAlg(\bC)\rightarrow \Fun(\Conv_{\langle 1\rangle},\bC)
\end{equation*}
which only remembers the `module part' of the $\ConvMod$-algebra. Given a $\ConvMod$-algebra $(\cM;\cA)$ we get an assignment $I\mapsto \cM_I$. Moreover, this is functorial in $I$: if we have $I\rightarrow J$ in $\fSet$, then the corresponding map $\cM_I\rightarrow \cM_J$ is the action type $\ConvMod$-operation for $U=\{*\}$, $\sigma$ is the unique linear ordering of $U$, and $\{f_u\}_{u\in U}=\{I\rightarrow J\}$.

\subsection{$\ConvMod$-algebras arising from the Hecke correspondence}\label{ssec: ConvMod algebra from Hecke} Recall from Sect. \ref{sssec: preservation non-metaplectic} and \ref{sssec: preservation metaplectic} the Hecke correspondence \eqref{eq: Hecke correspondence for I}
and the canonical compatibility data of gerbes appearing in the correspondence. Equivalently, we can consider the Hecke correspondence diagram 
\begin{equation}\label{eq: Hecke in a Corr form}
    \Bun_G\times X^I\times \Hecke_{X^I}\xleftarrow{(p_1,\id)} \Hecke_{X^I}\xrightarrow{p_2} \Bun_G\times X^I,
\end{equation}
and the pull-push operation along this correspondence. The stack $\Hecke_{X^I}$ has a canonical $A$-gerbe $\wt{\cG}^I$ given as a pullback of $\cG^I$ along the canonical map
\begin{equation}\label{eq: Hecke to Hecke loc}
\Hecke_{X^I}\rightarrow\Hecke^{\loc}_{X^I},
\end{equation}
and we have a canonical equivalence
\[\wt{\cG}^I+ p^*_1\cG^{\glob} \simeq p^*_2\cG^{\glob}.\]
This compatibility enables us to define a functor
\begin{equation}\label{eq: Shv Hecke action on Shv BunG}
    \Shv^q(\Bun_G\times X^I)\otimes \Shv^q(\Hecke_{X^I})\rightarrow \Shv^q(\Bun_G\times X^I),
\end{equation}
which in fact gets enhanced into the right action of $\Shv^q(\Hecke_{X^I})$ on $\Shv^q(\Bun_G\times X^I)$. In this subsection, we will encode this structure using the language of $\ConvMod$-algebras.

\sssec{}\label{sssec: Hecke as a Conv algebra}
Similar to Sect. \ref{sssec: Hecke loc as a Conv algebra}, we can enhance the collection $\{\Hecke_{X^I}\}_{I\in\fSet}$ to a $\Conv$-algebra $\cHecke_X$ in $\Corr(\PreStklft)_{\on{sch};\on{all}}$. The evaluation at the color $I$ is $(\cHecke_X)_I = \Hecke_{X^I}$, and given a $\Conv$-operation $(\sigma,\{f_u\})\in \Mul_{\Conv}(\{I_u\}, J)$, the corresponding map
 \[\prod_{u\in U}\Hecke_{{I_u}}\rightarrow \Hecke_{J} \quad \te{in } \Corr(\PreStklft, \on{sch})\]
    is the correspondence
    \begin{equation}\label{eq: operation for cHecke stack}
        \begin{tikzcd}
	{\wt{\Hecke}_{I_{\sigma(1)},\cdots,I_{\sigma(|U|)}\rightarrow J}} & {\prod_{1\le s\le |U|}\Hecke_{I_{\sigma(s)}}} \\
	{\Hecke_J.}
	\arrow[from=1-1, to=1-2]
	\arrow[from=1-1, to=2-1]
\end{tikzcd}
    \end{equation}
Here the total ordering $\sigma$ of $U$ is considered as a bijection $\sigma:\{1,\cdots,|U|\}\rightarrow U$, and the stack $\wt{\Hecke}_{I_{\sigma(1)},\cdots,I_{\sigma(|U|)}\rightarrow J}$ is a stack over $X^J$ whose points classify the tuples 
\begin{equation}\label{eq: point of iterated Hecke}
(\ul{x},\cP^0,\cdots,\cP^{|U|},\alpha^1,\cdots,\alpha^{|U|}),
\end{equation}
where $\ul{x}\in X^J$, whose restriction along $f_{\sigma(i)}$ is denoted $\ul{x}_{i}\in X^{I_{\sigma(i)}}$, $\cP^i$'s are $G$-bundles on $X$, and each $\alpha^i$ is an isomorphism between $\cP^{i-1}$ and $\cP^{i}$ outside of $\Gamma_{\ul{x}_i}$. 
\medskip

The vertical arrow of the correspondence sends the point \eqref{eq: point of iterated Hecke} to $(\ul{x}, \cP^0,\cP^{|U|},\alpha^{\on{con}})\in \Hecke_{X^J}$, where $\alpha^{\on{con}}:\cP^0|_{X\setminus \Gamma_{\ul{x}}}\xrightarrow{\sim} \cP^{|U|}|_{X\setminus \Gamma_{\ul{x}}}$ is the composition of restrictions $\alpha^s|_{X\setminus \Gamma_{\ul{x}}}: \cP^{u-1}|_{X\setminus \Gamma_{\ul{x}}}\xrightarrow{\sim} \cP^{u}|_{X\setminus \Gamma_{\ul{x}}}$. The horizontal arrow is the product of maps 
\begin{equation*}
    p_s:\wt{\Hecke}_{I_{\sigma(1),\cdots,I_{\sigma(|U|)}}\rightarrow J}\rightarrow \Hecke_{I_{\sigma(s)}}
\end{equation*}
for each $1\le s\le |U|$, which sends the point \eqref{eq: point of iterated Hecke} to the point
\[(\ul{x}_s,\cP^{s-1}, \cP^{s},\alpha^s|_{X\setminus\Gamma_{\ul{x}_s}})\in \Hecke_{I_{\sigma(s)}}.\]

\sssec{}\label{sssec: Shv q cHecke}
The collection of maps \eqref{eq: Hecke to Hecke loc} naturally gets enhanced to a map $\cHecke^{\loc}_X\rightarrow\cHecke_{X}$ between $\Conv$-algebras.\footnote{We are viewing the map \eqref{eq: Hecke to Hecke loc} as a map $\Hecke^{\loc}_{X^I}\rightarrow \Hecke_{X^I}$ in the category of correspondences.} This, together with the fact that the collection $\{(\Hecke^{\loc}_{X^I},\cG^I)\}_{I\in\fSet}$ gets enhanced to a $\Conv$-algebra in $\Corr(\PreStklft+\on{Ge}_A,\on{sch})$, we conclude that the collection
$\{(\Hecke_{X^I},\wt{\cG}^I)\}_{I\in\fSet}$ gets naturally enhanced to a $\Conv$-algebra in the same category of correspondences. We denote the resulting $\Conv$-algebra by $(\cHecke_X, \wt{\cG})$. By construction we have a $\Conv$-algebra map
\begin{equation}\label{eq: Conv map Hecke loc to Hecke}
    (\cHecke^{\loc}_X,\cG)\rightarrow (\cHecke_X,\wt{\cG}).
\end{equation}

Taking the right-lax symmetric monoidal functor $\Shv^{\on{tw}}_{\Corr}$ (see Theorem \ref{thm: twShvCorr functor}) to $(\cHecke_X, \wt{\cG})$ we get a $\Conv$-algebra in $\DGcat$, which we denote by $\Shv^q(\cHecke_X)$. By construction we have
\begin{equation}\label{eq: Shv q on Hecke loc to Shv q on Hecke}
    \Shv^q(\cHecke^{\loc}_X)\rightarrow \Shv^q(\cHecke_X)\quad \te{in }\Conv\lAlg(\DGcat).
\end{equation}
When evaluated at $I$, this is the pullback of sheaves along \eqref{eq: Hecke to Hecke loc}.
\medskip

Similar to Proposition \ref{prop: cHecke as an X module}, we have a natural $(X,0)^{\otimes}$-module structure on $(\cHecke_X,\wt{\cG})$, and the map \eqref{eq: Conv map Hecke loc to Hecke} has a natural $(X,0)^{\otimes}$-linearity. Therefore we conclude:
\begin{prop}
    The map \eqref{eq: Shv q on Hecke loc to Shv q on Hecke} has a natural $\Shv(X^{\otimes})$-linearity, and hence a $\LS(X)^{\otimes}$-linearity.
\end{prop}

\sssec{Monoidal structure}\label{sssec: monoidal structure ShvqHecke}
In the spirit of Sect. \ref{sssec: ShvqHeckeloc monoidal str}, the $\Conv$-algebra structure of $\Shv^q(\cHecke_{X})$ encodes the monoidal structure of the category $\Shv^q(\Hecke_{X^I})$. Indeed, the correspondence \eqref{eq: operation for cHecke stack} for the operation \[(\triv, \{\id_I\}_{u\in \{1,2,\cdots,k\}})\in \Mul_{\Conv}(\{I\}_{u\in \{1,2,\cdots,k\}},\{I\})\]
    is the $k$-fold convolution diagram of $\Shv^q(\Hecke_{X^I})$.
    \medskip

    For each color $I$, there is a unique operation whose target is $I$ and $U=\varnothing$. In this case, the correspondence  \eqref{eq: operation for cHecke stack} reads
    \begin{equation}\label{eq: pointing cHeckeX}
        \begin{tikzcd}
	{(\Bun_G\times X^I,0)} & {(\pt,0)} \\
	{(\Hecke_{X^I},\cG^I),}
	\arrow[from=1-1, to=1-2]
	\arrow[from=1-1, to=2-1]
\end{tikzcd}
    \end{equation}
where the vertical arrow sends a point $(\cP,\ul{x})\in\Bun_G\times X^I$ to the point $(\ul{x},\cP,\cP,\id)$.
The pull-push operation along this correspondence agrees with the unit map
    \[\Vect\rightarrow \Shv^q(\Hecke_{X^I})\]
    of the convolution monoidal structure.

\sssec{}\label{sssec: a ConvMod algebra Hecke on BunG}
We have a $\ConvMod$-algebra $(\Bun_G\times X^{\otimes}; \cHecke_{X})$ in $\Corr(\PreStklft)_{\on{sch};\on{all}}$ whose evaluation at the algebra type color $I^{\mf{a}}$ is $\Hecke_{X^I}$, and its evaluation at the module type color $I^{\mf{m}}$ is $\Bun_G\times X^I$. Under the forgetful functor from $\ConvMod$-algebras to $\Conv$-algebras, the $\ConvMod$-algebra $(\Bun_G\times X^{\otimes}; \cHecke_{X})$ forgets to the $\Conv$-algebra $\cHecke_X$. This explains the multiplicative type operations of $(\Bun_G\times X^{\otimes}; \cHecke_{X})$.
\medskip

Suppose we are given an action type $\ConvMod$-operation \[(\sigma,\{f_u\}_{u\in U})\in \Mul_{\ConvMod}(\{I^{\mf{m}}_{u_0}\}\sqcup\{I^{\mf{a}}_u\}_{u\in U\setminus u_0}, J^{\mf{m}}).\]
    Recall the notations from Sect. \ref{sssec: Hecke as a Conv algebra}. We have $\sigma(1)=u_0$ by definition. Then the corresponding map
    \[\Bun_G\times X^{I_{u_0}}\times \prod_{u\in U\setminus u_0}\Hecke_{{I_u}}\rightarrow \Bun_G\times X^J \quad \te{in } \Corr(\PreStklft, \on{sch})\]
    is the correspondence
    \begin{equation}\label{eq: Hecke ConvMod correspondence}
        \begin{tikzcd}
	{\wt{\Hecke}_{I_{\sigma(1)},\cdots,I_{\sigma(|U|)}\rightarrow J}} & {\Bun_G\times X^{I_{\sigma(1)}}\times \prod_{2\le s\le |U|}\Hecke_{I_{\sigma(s)}}} \\
	{\Bun_G\times X^J.}
	\arrow[from=1-1, to=1-2]
	\arrow[from=1-1, to=2-1]
\end{tikzcd}
    \end{equation}
     The vertical arrow of the correspondence sends the point \eqref{eq: point of iterated Hecke} to $(\cP^{|U|},\ul{x})\in \Bun_G\times X^J$. The horizontal arrow is the product of maps 
     $\{p_s\}_{2\le s \le |U|}$ from Sect. \ref{sssec: Hecke as a Conv algebra} and a map
     \[\wt{\Hecke}_{I_{\sigma(1),\cdots,I_{\sigma(|U|)}}\rightarrow J}\rightarrow \Bun_G\times X^{I_{\sigma(1)}}\]
    sending the point \eqref{eq: point of iterated Hecke} to $(\cP^{0},\ul{x}_{1})$.
\medskip

    The discussion above can be summarized into the following proposition, which can be proved using the argument of \cite[Sect. 1.2.7]{ZhaoNew}:

\begin{prop}
    There is a canonical $\ConvMod$-algebra object $(\Bun_G\times X^{\otimes}; \cHecke_I)$ in the symmetric monoidal category $\Corr(\PreStklft, \on{sch})$ admitting the description of evaluation at colors and operation maps as in Sect. \ref{sssec: a ConvMod algebra Hecke on BunG}.
\end{prop}

\sssec{}
As pointed out earlier, the Hecke correspondence \eqref{eq: Hecke in a Corr form} is compatible with the formation of gerbes $\wt{\cG}_I$ and $\cG^{\on{glob}}$. This compatibility admits a natural iterated version (see \cite[Lemma 7.1.11]{ZhaoNew}). For example, the correspondence \eqref{eq: Hecke ConvMod correspondence} gets equipped with a natural compatibility data pertaining to the gerbes $\cG^{\on{glob}}$ and $\{\wt{\cG}^{I_{\sigma(u)}}\}_{2\le u\le |U|}$.
\medskip

Hence, we have a $\ConvMod$-algebra object 
\begin{equation}\label{eq: ConvMod BunG Hecke gerbe}
    ((\Bun_G\times X^{\otimes},\pi^*\cG^{\on{glob}}); (\cHecke_{X},\wt{\cG}))
\end{equation}
in the category $\Corr(\PreStklft+\on{Ge}_A,\on{sch})$ whose evaluation at the algebra type color $I^{\mf{a}}$ is $(Hecke_{X^I},\wt{\cG}_I)$ and evaluation at the module type color $I^{\mf{m}}$ is $(\Bun_G\times X^I, \pi_I^*\cG^{\on{glob}})$ where $\pi_I: \Bun_G\times X^I\rightarrow \Bun_G$ is the projection. On the level of underlying spaces, the $\ConvMod$-operations are just the same as that of $\cHecke_X$. Its multiplication type $\ConvMod$-operation is determined by the $\Conv$-operation of $(\cHecke_{X},\wt{\cG})$, and its action type $\ConvMod$-operation is determined by the aforementioned compatibility data of Hecke correspondence. To summarize, we have:

\begin{prop}
    There is a canonical $\ConvMod$-algebra object \eqref{eq: ConvMod BunG Hecke gerbe} in the category $\Corr(\PreStklft+\on{Ge}_A,\on{sch})$ satisfying the following properties:
    \begin{enumerate}
        \item its evaluation at the algebra type color $I^{\mf{a}}$ is $(\Hecke_{X^I},\wt{\cG}^I)$ and its evaluation at the module type color $I^{\mf{m}}$ is $(\Bun_G\times X^I,\pi^*_I\cG^{\on{glob}})$;
        \item it forgets to the $\ConvMod$-algebra $(\Bun_G\times X^{\otimes}; \cHecke_{X})$ under the symmetric monoidal forgetful map \[\PreStklft+\on{Ge}_A\rightarrow \PreStklft;\]
        \item it forgets to the $\Conv$-algebra $(\cHecke_X, \wt{\cG})$ under the forgetful map
        \[\ConvMod\lAlg(\mhyphen)\rightarrow \Conv\lAlg(\mhyphen);\]
        \item its action type $\ConvMod$-operation on the level of gerbes is given by the data of compatibility of gerbes under the (iterated) Hecke correspondence.
    \end{enumerate}
\end{prop}

\sssec{}
Taking the right-lax symmetric monoidal functor $\Shv^{\on{tw}}_{\Corr}$ to the $\ConvMod$-algebra \eqref{eq: ConvMod BunG Hecke gerbe}, we get a $\ConvMod$-algebra in $\DGcat$, denoted by \[\left(\Shv^q(\Bun_G\times X^{\otimes}); \Shv^q(\cHecke_{X})\right).\] By construction, it forgets to the $\Conv$-algebra $\Shv^q(\cHecke_{X})$ introduced in Sect. \ref{sssec: Shv q cHecke}. Its evaluation at the module type color $I^{\mf{m}}$ is the category of twisted sheaves on $\Bun_G\times X^I$ twisted by the gerbe $\pi_I^*\cG^{\on{glob}}$, denoted by
$\Shv^q(\Bun_G\times X^I)$.

\begin{rmk}
    By construction, the functor \eqref{eq: Shv Hecke action on Shv BunG} is the structure map of the $\ConvMod$-algebra above corresponding to the action type operation \[(\triv, \{\id_I,\id_I\})\in \Mul_{\ConvMod}(\{I^{\mf{m}},I^{\mf{a}}\}, I^{\mf{m}}).\]
    Here $U=\{1,2\}$ and $\triv$ is the trivial linear ordering of it.
\end{rmk}

\sssec{}
Applying the construction of Sect. \ref{sssec: free module construction} to the situation
\[\bC=\Corr(\PreStklft)_{\on{sch};\on{all}},\quad \cC=X, \quad M=\Bun_G\]
we get a $\ConvMod$-algebra object $(\Bun_G\times X^{\otimes}; X^{\otimes})$ in the category of correspondences. 
\medskip

For each $I\in\fSet$ we have a correspondence
\[\begin{tikzcd}
	{\Bun_G\times X^I} & {X^I} \\
	{\Hecke_{X^I}}
	\arrow[from=1-1, to=1-2]
	\arrow[from=1-1, to=2-1]
\end{tikzcd}\]
where the horizontal map is the projection and the vertical map encodes the trivial modification. Moreover, the collection of such correspondences gets enhanced to a map $X^{\otimes}\rightarrow \cHecke_{X}$ of $\Conv$-algebras in the category of correspondences. This map induces a map
\[(\Bun_G\times X^{\otimes}; X^{\otimes})\rightarrow (\Bun_G\times X^{\otimes}; \cHecke_{X^I})\]
of $\ConvMod$-algebras in the category of correspondences. The source, which a priori obtained as a pullback, agrees with the free module construction applied to $\Bun_G$ and $X^{\otimes}$.

\sssec{}\label{sssec: unit canonical pointing gerbe}
The above admits a natural analog including gerbes. The pullback of $\wt{\cG}_I$ along the map $e_I$ is canonically trivialized, and hence $e_I$ gets enhanced to
$e_I: (X^I,0)\rightarrow (\Hecke_{X^I},\wt{\cG}_I)$.
This collection of maps gets enhanced to a $\Conv$-algebra map 
\begin{equation}\label{eq: can pointing cHecke gerbe}
    (X,0)^{\otimes}\rightarrow (\cHecke_{X},\wt{\cG}).
\end{equation}
Pulling back the $\ConvMod$-algebra \eqref{eq: ConvMod BunG Hecke gerbe} we get a $\ConvMod$-algebra denoted \begin{equation}\label{eq: BunG X by X ConvMod}
    ((\Bun_G\times X^{\otimes}, \pi^*\cG^{\on{glob}}); (X,0)^{\otimes}).
\end{equation}
One checks that this agrees with the free module
\[((\Bun_G,\cG^{\on{glob}})\times (X,0)^{\otimes}; (X,0)^{\otimes}).\]

From the pullback construction, we have a canonical map
\[((\Bun_G\times X^{\otimes},\pi^*\cG^{\on{glob}}); (X,0)^{\otimes})\rightarrow ((\Bun_G\times X^{\otimes},\pi^*\cG^{\on{glob}}); (\cHecke_X,\wt{\cG}))\]
of $\ConvMod$-algebras in the category of correspondences.  Taking the right-lax symmetric monoidal functor $\Shv^{\on{tw}}_{\Corr}$ we get a map
\begin{equation}\label{eq: pullback to the action of Shv X}
    (\Shv^q(\Bun_G\times X^{\otimes});\Shv(X^{\otimes}))\rightarrow (\Shv^q(\Bun_G\times X^{\otimes}); \Shv^q(\cHecke_X))
\end{equation}
in $\ConvMod\lAlg(\DGcat)$. Here, the source is defined to be the image of \eqref{eq: BunG X by X ConvMod} under $\Shv^{\on{tw}}_{\Corr}$. Its evaluation at the algebra type color $I^{\mf{a}}$ is $\Shv(X^I)$.

\sssec{}
Assume either $X\ne\pjl$ or $X=\pjl$ and $\cG_Z$ is trivializable.
We have a composition
\[(\Rep(H)^{\loc}_{\cG_Z}(X))^{\otimes}\rightarrow \Shv^q(\cHecke^{\loc}_X)\rightarrow \Shv^q(\cHecke_X) \quad \te{in }\Conv\lAlg(\DGcat)\]
where the first arrow is the $\Conv$-algebra form of $\LS$-linear factorizable Satake and the second arrow is \eqref{eq: Shv q on Hecke loc to Shv q on Hecke}. Sect. \ref{sssec: pullback of ConvMod} combined with the $\Conv$-algebra map above yields us an object
\begin{equation}\label{eq: Shv linear Satake kernel ConvMod}
    \left(\Shv^q(\Bun_G\times X^\otimes);(\Rep(H)^{\loc}_{\cG_Z}(X))^{\otimes}\right)\in \ConvMod\lAlg(\DGcat).
\end{equation}
From Sect. \ref{sssec: lifting ConvMod to ConvcMod}, this has a canonical lift to a $\ConvcMod$-algebra.

\subsection{$\ConvcMod$-algebras and nilpotent singular support}\label{ssec: ConvMod algebras and Nilp}
In this subsection, we construct the canonical functor
\begin{equation}\label{eq: ConvMod extension}
    (\cA^q\otimes \LS(X)^{\otimes};\LS(X)^{\otimes})\rightarrow \left(\cA^q_G\otimes \LS(X)^{\otimes}; (\Rep(H)^{\loc}_{\cG_Z}(X))^{\otimes}\right)
\end{equation}
in $\ConvcMod\lAlg(\DGcat)$ under a mild assumption on $X$ and $\cG_Z$. We first have to define both sides.

\sssec{}
The source of \eqref{eq: ConvMod extension} is defined via the free module construction (Sect. \ref{sssec: free module construction}). For the construction of the target we prove:

\begin{prop}\label{prop: constructing ConvMod AqG}
    Suppose either $X\ne\pjl$ or $X=\pjl$ and $\cG_Z$ is trivializable.
    There is a canonical object
    \begin{equation}\label{eq: ConvMod AqG twRepH}
    \left(\cA^q_G\otimes \LS(X)^{\otimes}; (\Rep(H)^{\loc}_{\cG_Z}(X))^{\otimes}\right)\in \ConvcMod\lAlg(\DGcat)
\end{equation}
whose evaluation at the algebra type color $I^{\mf{a}}$ is $(\Rep(H)^{\loc}_{\cG_Z}(X))^{\otimes I}$ and the evaluation at the module type color $I^{\mf{m}}$ is $\cA^q_G\otimes \LS(X)^{\otimes I}$, arising as a submodule of \eqref{eq: Shv linear Satake kernel ConvMod} in the sense of Sect \ref{sssec: submodule construction}.
\end{prop}
\begin{proof}
    We aim to apply the paradigm of Sect \ref{sssec: submodule construction} to the situation where $(\cM;\cA)$ is \eqref{eq: Shv linear Satake kernel ConvMod} and each $\cM'_I\subseteq \cM_I$ is the inclusion
    \[\Shv^q_{\Nilp\times X^I}(\Bun_G\times X^I)\subset \Shv^q(\Bun_G\times X^I).\]
    By Corollary \ref{cor: singling out AqG otimes LS}, we have $\cM'_I=\cA^q_G\otimes \LS(X)^{\otimes I}$. To use the submodule construction, we have to prove the following two assertions:
    \begin{enumerate}
        \item given $I\rightarrow J$ in $\fSet$, the action type operation 
        \begin{equation}\label{eq: Shv q BunGXI functoriality}
            \Shv^q(\Bun_G\times X^I)\rightarrow \Shv^q(\Bun_G\times X^J)
        \end{equation}
        restricts to $\Shv^q_{\Nilp\times X^I}(\Bun_G\times X^I)\rightarrow \Shv^q_{\Nilp\times X^J}(\Bun_G\times X^J)$;
        \item for each $I$, the action type operation
        \begin{equation}\label{eq: action of met Sat kernel}
            \Shv^q(\Bun_G\times X^I)\otimes (\Rep(H)^{\loc}_{\cG_Z}(X))^{\otimes I}\rightarrow \Shv^q(\Bun_G\times X^I)
        \end{equation}
        restricts to
        \[\Shv^q_{\Nilp\times X^I}(\Bun_G\times X^I)\otimes (\Rep(H)^{\loc}_{\cG_Z}(X))^{\otimes I}\rightarrow \Shv^q_{\Nilp\times X^I}(\Bun_G\times X^I).\]
    \end{enumerate}
    
    The correspondence \eqref{eq: Hecke ConvMod correspondence} corresponding to the operation \eqref{eq: Shv q BunGXI functoriality} is
    \[\begin{tikzcd}
	{\Bun_G\times X^J} & {\Bun_G\times X^I} \\
	{\Bun_G\times X^J}
	\arrow["{\id\times\Delta}", from=1-1, to=1-2]
	\arrow["\id"', from=1-1, to=2-1]
\end{tikzcd}\]
where $X^J\xrightarrow{\Delta} X^I$ is the diagonal map induced from $I\rightarrow J$, and therefore \eqref{eq: Shv q BunGXI functoriality} is just a pullback along $\id\times \Delta$. We have an obvious commutative diagram
\[\begin{tikzcd}
	{\Shv^q(\Bun_G\times X^I)} & {\Shv^q(\Bun_G\times X^J)} \\
	{\Shv^q(\Bun_G)\otimes \Shv(X^I)} & {\Shv^q(\Bun_G)\otimes \Shv(X^J).}
	\arrow["{(\id\times\Delta)^*}", from=1-1, to=1-2]
	\arrow[from=2-1, to=1-1]
	\arrow["{\id\otimes \Delta^*}"', from=2-1, to=2-2]
	\arrow[from=2-2, to=1-2]
\end{tikzcd}\]
This, combined with Corollary \ref{cor: singling out AqG otimes LS} verifies the first assertion.
\medskip

    The map \eqref{eq: action of met Sat kernel} is the action of metaplectic Satake kernels on $\Shv^q_{\Nilp\times X^I}(\Bun_G\times X^I)$, and hence the second assertion is exactly the content of Theorem \ref{thm: preservance Nilp metaplectic}.
\end{proof}

\sssec{}
In the remainder of this subsection, we will prove:
\begin{thm}\label{thm: extension of action}
    Suppose either $X\ne \pjl$ or $X=\pjl$ and $\cG_Z$ is trivializable. There is a canonical map \eqref{eq: ConvMod extension} 
    in $\ConvcMod\lAlg(\DGcat)$,
    where the source arises from a free module construction, and the target is \eqref{eq: ConvMod AqG twRepH}.
\end{thm}

\sssec{}\label{sssec: two constructions}
We can form a $\ConvcMod$-algebra \footnote{We added the superscript $\on{Sat}$ to indicate that this was constructed using the metaplectic Satake.}
\begin{equation}\label{eq: ConvcMod Sat}
    (\cA^q_G\otimes \LS(X)^{\otimes};\LS(X)^{\otimes})^{\on{Sat}}
\end{equation}
via the pullback of \eqref{eq: ConvMod AqG twRepH} along the $\Convcomm$-algebra map $\LS(X)^{\otimes}\rightarrow (\Rep(H)^{\loc}_{\cG_Z}(X))^{\otimes}$. This comes naturally with a $\ConvcMod$-algebra map to the target of \eqref{eq: ConvMod extension}, and moreover this has the same evaluation at colors and the underlying $\Convcomm$-algebra with the free module 
\begin{equation}\label{eq: ConvcMod Free}
    (\cA^q_G\otimes \LS(X)^{\otimes};\LS(X)^{\otimes})^{\on{Free}}.
\end{equation}
Therefore, it suffices to prove the following assertion:

\begin{prop}\label{prop: two constructions agree}
    We have a canonical equivalence between the $\ConvcMod$-algebras \eqref{eq: ConvcMod Sat} and \eqref{eq: ConvcMod Free}.
\end{prop}

\begin{rmk}\label{rmk: ConvMod suffices}
    From Sect. \ref{sssec: lifting ConvMod to ConvcMod}, it suffices to construct an equivalence as $\ConvMod$-algebras whose underlying $\Conv$-algebra map $\LS(X)^{\otimes}\rightarrow\LS(X)^{\otimes}$ is the identity.
\end{rmk}

\sssec{}
This assertion is not surprising: it is a consequence of $\LS(X)^{\otimes I}$-linearity of the metaplectic Satake functor and the metaplectic Hecke correspondence. However, some care is needed to provide a formal proof, because the construction of $(\cA^q_G\otimes \LS(X)^{\otimes};\LS(X)^{\otimes})^{\Sat}$ goes through the metaplectic Satake, whose construction was indirect. Our strategy is to give an alternative description of the above that does not involve geometric Satake, but only the geometry on the automorphic side. 

\sssec{Pointing of a $\Conv$-algebra}\label{sssec: pointing of a Conv alg}
Suppose $\bC$ is a symmetric monoidal category with the monoidal unit $1_{\bC}$, and $\cC\in\Conv\lAlg(\bC)$. For each color $I$, there is a unique $\Conv$-operation with target $I$ with $U=\varnothing$. This operation yields a canonical map $1_{\bC}\rightarrow \cC_I$. The collection of such maps gets enhanced to
\[1_{\bC}^{\otimes}\rightarrow \cC\quad\te{in }\Conv\lAlg(\bC).\]
This can be understood as a canonical pointing of a $\Conv$-algebra.
\medskip

Now suppose that $\cA\in\Alg(\Conv\lAlg(\bC))$ and that $\cC$ is an $\cA$-module internal to $\Conv\lAlg(\bC)$. Then the canonical pointing yields $\cA\rightarrow\cC$ in $\Conv\lAlg(\bC)$. This construction is organized into a functor
\begin{equation*}
    \cA\lMod(\Conv\lAlg(\bC))\rightarrow \Conv\lAlg(\bC)_{\cA/}.
\end{equation*}

\sssec{}\label{sssec: alternative description}
Applying the construction of Sect. \ref{sssec: pointing of a Conv alg} to the $\LS(X)^{\otimes}$-linearity  of the factorizable Satake functor (Theorem \ref{thm: fact Satake Conv at X}) and the $\Shv(X^{\otimes})$-linearity of the functor \eqref{eq: Shv q on Hecke loc to Shv q on Hecke}, we get the following commutative diagram in $\Conv\lAlg(\DGcat)$:
\[\begin{tikzcd}
	{\left(\Rep(H)^{\loc}_{\cG_Z}(X)\right)^{\otimes}} & {\Shv^q(\cHecke^{\loc}_X)} & {\Shv^q(\cHecke_X)} \\
	{\LS(X)^{\otimes}} & {\Shv(X^{\otimes}).}
	\arrow[from=1-1, to=1-2]
	\arrow[from=1-2, to=1-3]
	\arrow[from=2-1, to=1-1]
	\arrow[from=2-1, to=2-2]
	\arrow[from=2-2, to=1-2]
	\arrow[from=2-2, to=1-3]
\end{tikzcd}\]

From the commutativity of this diagram, we can understand the construction of \eqref{eq: ConvcMod Sat} as the following three-step procedure:
\begin{enumerate}
    \item start from $(\Shv^q(\Bun_G\times X^{\otimes}); \Shv^q(\cHecke_{X}))$ and pull it back along $\Shv(X^{\otimes})\rightarrow \Shv^q(\cHecke_{X})$ to get
    \begin{equation}\label{eq: ShvqBunG ShvX}
        (\Shv^q(\Bun_G\times X^{\otimes}); \Shv(X^{\otimes}));
    \end{equation}
    \item pull it back along $\LS(X)^{\otimes}\rightarrow \Shv(X^{\otimes})$ to get $(\Shv^q(\Bun_G\times X^{\otimes}); \LS(X)^{\otimes})$;
    \item use the submodule construction to get $(\cA^q_G\otimes\LS(X)^{\otimes};\LS(X)^{\otimes})^{\on{Sat}}$.
\end{enumerate}

\begin{lemma}\label{lem: pullback equals free}
    The $\ConvMod$-algebra \eqref{eq: ShvqBunG ShvX} is canonically identified with the image of the free module $\ConvMod$-algebra
    \begin{equation}\label{eq: free module BunG X gerbe}
        ((\Bun_G,\cG^{\on{glob}})\times (X,0)^{\otimes}; (X,0)^{\otimes})
    \end{equation}
    under the right-lax symmetric monoidal functor $\Shv^{\on{tw}}_{\Corr}$. Moreover, this identification induces the identity of $\Shv(X^{\otimes})$ on the level of $\Conv$-algebras.
    \end{lemma}
    \begin{proof}
    The construction given in Sect. \ref{sssec: pointing of a Conv alg} enjoys an obvious functoriality with respect to a right-lax symmetric monoidal functor $\bC_1\rightarrow\bC_2$. Therefore, it suffices to prove the following: consider a $\Conv$-algebra map $\eta^{\on{lin}}:(X,0)^{\otimes}\rightarrow (\cHecke_X,\wt{\cG})$ obtained from the construction of Sect. \ref{sssec: pointing of a Conv alg}, and then the pullback of the $\ConvMod$-algebra $((\Bun_G,\pi^*\cG^{\on{glob}});(\cHecke_X,\wt{\cG}))$ along $\eta^{\on{lin}}$ is canonically identified with the free module \eqref{eq: free module BunG X gerbe}, and moreover this identification induces the identity of $(X,0)^{\otimes}$ on the level of $\Conv$-algebras. 
    \medskip

        The canonical pointing of $(\cHecke_X,\wt{\cG})$ is described by the correspondences \eqref{eq: pointing cHeckeX}. Using this we see that the $\Conv$-algebra map $\eta^{\on{lin}}$ agrees with \eqref{eq: can pointing cHecke gerbe} from Sect. \ref{sssec: unit canonical pointing gerbe}. We conclude from the discussion in the same section.
        \end{proof}

\sssec{Proof of Proposition \ref{prop: two constructions agree}}
    Combining the above discussion, \eqref{eq: ConvcMod Sat} as a $\ConvMod$-algebra agrees with the pullback and submodule construction applied to the image of the free module $\Conv$-algebra \eqref{eq: free module BunG X gerbe} under $\Shv^{\on{tw}}_{\Corr}$.
   It is clear that the second and third steps of Sect. \ref{sssec: alternative description} applied to this $\Conv$-algebra gives us the underlying $\ConvMod$-algebra of \eqref{eq: ConvcMod Free}. Therefore, we identified \eqref{eq: ConvcMod Sat} and \eqref{eq: ConvcMod Free} as $\ConvMod$-algebras. Moreover, from the second part of Lemma \ref{lem: pullback equals free}, this identification induces the identity of $\LS(X)^{\otimes}$ on the underlying $\Conv$-algebras. Using Remark \ref{rmk: ConvMod suffices} we conclude. \qed

\subsection{Passage from $\ConvcMod$-algebras to a natural transformation}\label{ssec: passage ConvcMod nat transf}
In this subsection, we explain a categorical mechanism to pass from the data \eqref{eq: ConvMod extension} to a natural transformation $\eqref{eq: source global}\rightarrow \eqref{eq: target global}$.

\begin{note}
    Throughout this subsection, let $\bC$ be a closed symmetric monoidal category where the tensor product commutes with colimits on each variable. Let $\cHom_{\bC}$ be the internal $\Hom$ in $\bC$. In practice, we will need $\bC=\DGcat$. 
\end{note}

\sssec{}
Given a category $\bD$, let $\Arr(\bD)$ be the category of arrows of $\bD$.
We have canonical functors
\[\bD\xrightarrow{e} \Arr(\bD),\quad \Arr(\bD)\xrightarrow{s} \bD,\ \te{and }\Arr(\bD)\xrightarrow{t} \bD,\]
where the first functor sends $D\in \cD$ to $(D\xrightarrow{\id} D)$, the second (resp. the third) functors records the source (resp. the target) of the arrow. For convenience, we call those the \emph{neutral}, \emph{source}, \emph{target} maps in order.

\sssec{}
Recall from Definition \ref{defn: LMcomm} that we have an operad $\LMcomm$ such that a $\LMcomm$-algebra object is a pair $(C,M)$ where $C\in\CAlg(\bC)$ and $M\in C\lMod(\bC)$. Consider the fiber product
\begin{equation}\label{eq: arrow fib prod 1}
    \Arr(\LMcomm\lAlg(\bC))\times_{\Arr(\bC)}\bC
\end{equation}
of categories where the map $\Arr(\LMcomm\lAlg(\bC))\rightarrow\Arr(\bC)$ is induced from the map 
\begin{equation}\label{eq: LMcomm forget to module}
    \LMcomm\lAlg(\bC)\rightarrow \bC,\quad (C,M)\mapsto M,
\end{equation}
and the map $\bC\rightarrow\Arr(\bC)$ is the neutral map. The above fiber product is the category of the following data: an object $M$, commutative algebra objects $C_1$ and $C_2$, and an arrow $(C_1,M)\rightarrow (C_2,M)$ in $\LMcomm\lAlg(\bC)$. 
\medskip

Fixing an object $M_0\in \bC$ we get the `free module' functor 
\begin{equation}\label{eq: free module functor LMcomm}
    \on{Free}_{M_0}:\CAlg(\bC)\rightarrow \LMcomm\lAlg(\bC),\quad C\mapsto (C,M_0\otimes C).
\end{equation}
We can take a bigger fiber product
\begin{equation}\label{eq: arrow fib prod 2}
    \CAlg(\bC)\times_{\LMcomm\lAlg(\bC)} \Arr(\LMcomm\lAlg(\bC))\times_{\Arr(\bC)}\bC
\end{equation}
where the second fiber product is the same as in \eqref{eq: arrow fib prod 1}, and for the first fiber product we use the free module functor $\on{Free}_{M_0}$ and the source map for $\LMcomm\lAlg(\bC)$. This fiber product is the category of the following data: commutative algebra objects $C_1$ and $C_2$, a lift of $M_0\otimes C_1$ to $(C_2, M_0\otimes C_1)\in\LMcomm\lAlg(\bC)$, and an arrow $(C_1, M_0\otimes C_1)\rightarrow (C_2, M_0\otimes C_1)$ in $\LMcomm\lAlg(\bC)$. 

\sssec{}
The category $\Alg(\LMcomm\lAlg(\bC))$ consists of the following data: a commutative algebra object $C$ and an object $A\in \Alg(C\lMod(\bC))$. It has a forgetful map to $\CAlg(\bC)$ which only remembers the commutative algebra object $C$. We take a fiber product
\begin{equation}\label{eq: arrow fib prod 3}
    \Arr(\Alg(\LMcomm\lAlg(\bC)))\times_{\Arr(\CAlg(\bC))}\CAlg(\bC)
\end{equation}
where we use the neutral map for $\CAlg(\bC)$ and the forgetful map. This fiber product is the category of the following data: a commutative algebra object $C$ and a map $A_1\rightarrow A_2$ in $\Alg(C\lMod(\bC))$.

\sssec{}
Using the descriptions of the fiber products \eqref{eq: arrow fib prod 2} and \eqref{eq: arrow fib prod 3} above, we observe:

\begin{lemma}\label{lem: functor btw arrow fib prods}
    We have a canonical functor from \eqref{eq: arrow fib prod 2} to \eqref{eq: arrow fib prod 3} which maps an arrow $(C_1,M_0\otimes C_1)\rightarrow (C_2,M_0\otimes C_1)$ in $\LMcomm\lAlg(\bC)$ to the arrow $C_2\rightarrow \cHom_{\bC}(M_0, M_0\otimes C_1)$ in $\Alg(C_1\lMod(\bC))$.
\end{lemma}

\begin{rmk}
    It may seem unnatural or too restrictive to include the free module functor in the construction of \eqref{eq: arrow fib prod 2}. In fact, even from the object $(C_1,M)\rightarrow (C_2,M)$ of \eqref{eq: arrow fib prod 1} we can produce an object in \eqref{eq: arrow fib prod 3} by considering $C_2\rightarrow \mathcal{E}nd_{C_1}(M)$ in $\Alg(C_1\lMod)$, where $\mathcal{E}nd_{C_1}$ is the internal endomorphism in $C_1\lMod(\bC)$. However, this construction lacks functoriality in $M$, and hence cannot be upgraded to a functor as in Lemma \ref{lem: functor btw arrow fib prods}. We avoid this issue by requiring that $M$ is obtained by a free module construction with a fixed $M_0$.
\end{rmk}

\sssec{}
The main result of this subsection is a $\ConvcMod$-version of Lemma \ref{lem: functor btw arrow fib prods}. To state this, we consider the $\ConvcMod$-version of each of \eqref{eq: arrow fib prod 2} and \eqref{eq: arrow fib prod 3}.
\medskip

Given a category $\bD$ we use the notation
\[\Conv_{\langle 1\rangle}(\bD):=\Fun(\Conv_{\langle 1\rangle},\bD)=\Fun(\fSet,\bD).\]

\sssec{}\label{sssec: construction fib prod 1}
We consider the fiber product
\begin{equation}\label{eq: ConvMod fib prod 1}
    \Convcomm\lAlg(\bC)\times_{\ConvcMod\lAlg(\bC)} \Arr(\ConvcMod\lAlg(\bC))\times_{\Arr(\Conv_{\langle 1\rangle}(\bC)))} \Conv_{\langle 1\rangle}(\bC),
\end{equation}
which can be thought of as the $\ConvcMod$-version of \eqref{eq: arrow fib prod 2}. We explain how the fiber products are formed:
\begin{enumerate}
    \item for the first fiber product, we use the source map for $\ConvcMod\lAlg(\bC)$ and the free module functor
    \[\on{Free}_{M_0}: \Convcomm\lAlg(\bC)\rightarrow \ConvcMod\lAlg(\bC),\quad \cC\mapsto (M_0\otimes \cC;\cC);\]
    \item for the second fiber product, we use the neutral map for $\Conv_{\langle 1\rangle}(\bC)$ and the forgetful functor
    \begin{equation}\label{eq: forget to module}
        \ConvcMod\lAlg(\bC)\rightarrow \Conv_{\langle 1\rangle}(\bC),\quad (\cM;\cC)\mapsto (I\mapsto \cM_I)
    \end{equation}
    discussed in Sect. \ref{sssec: ConvMod to Conv1}.
\end{enumerate}

\sssec{}\label{sssec: ConvMod fib 1 data}
The fiber product \eqref{eq: ConvMod fib prod 1} is the category of the following data:
\begin{enumerate}
    \item $\Convcomm$-algebras $\cC_1$ and $\cC_2$;
    \item a lift of $M_0\otimes \cC_1\in \Conv_{\langle 1\rangle}(\bC)$ to a $\ConvcMod$-algebra $(M_0\otimes \cC_1;\cC_2)$ along \eqref{eq: forget to module};
    \item and a $\ConvcMod$-map $(M_0\otimes \cC_1;\cC_1)\rightarrow (M_0\otimes \cC_1;\cC_2)$.
\end{enumerate}

\sssec{}
We also consider the fiber product
\begin{equation}\label{eq: ConvMod fib prod 2}
    \Arr\left(\Conv_{\langle 1\rangle}\left(\Alg(\LMcomm\lAlg(\bC))\right)\right)\times_{\Arr\left(\Conv_{\langle 1\rangle}(\CAlg(\bC))\right)} \Conv_{\langle 1\rangle}(\CAlg(\bC)),
\end{equation}
which can be thought of as the $\ConvcMod$-version of \eqref{eq: arrow fib prod 3}. In the fiber product, the forgetful functor
\[\Alg(\LMcomm\lAlg(\bC))\rightarrow \Alg(\CAlg(\bC))=\CAlg(\bC)\]
and the neutral map for $\Conv_{\langle 1\rangle}(\CAlg(\bC))$ is used. 

\sssec{}\label{sssec: ConvMod fib 2 data}
Recall the notation
 \[\Fun^{\circ}(\fSet,\fSet+\Alg(\cC\lMod))\]
 introduced in Sect. \ref{sssec: Alg LS construction}, defined for $\cC\in\Conv_{\langle 1\rangle}(\DGCatSymMon)$. More generally, we can replace $\DGcat$ by $\bC$ to get a category
 \begin{equation}\label{eq: Alg LS construction for C}
     \Fun^{\circ}(\fSet,\fSet+\Alg(\cC\lMod(\bC))).
 \end{equation}
The fiber product \eqref{eq: ConvMod fib prod 2} is the category of the following data:
\begin{enumerate}
    \item a map $\cC:\Conv_{\langle 1\rangle}\rightarrow \CAlg(\bC)$;
    \item an arrow in the category \eqref{eq: Alg LS construction for C}.
\end{enumerate}

\sssec{}
In the remainder of this subsection we will prove:

\begin{prop}\label{prop: extension mechanism}
    There is a canonical functor from \eqref{eq: ConvMod fib prod 1} to \eqref{eq: ConvMod fib prod 2}.
\end{prop}

\begin{proof}
    Using Lemma \ref{lem: Convcomm and Einfty} and Lemma \ref{lem: ConvcMod and LMcomm}, \eqref{eq: ConvMod fib prod 1} is equivalent to the following category:
    \begin{equation}\label{eq: mechanism proof intermediate}
        \Conv_{\langle
    1\rangle}(\CAlg(\bC))\times_{\Conv_{\langle 1\rangle}(\LMcomm\lAlg(\bC))}\Arr(\Conv_{\langle 1\rangle}(\LMcomm\lAlg(\bC)))\times_{\Arr(\Conv_{\langle 1\rangle}(\bC))}\Conv_{\langle 1\rangle}(\bC),
    \end{equation}
    where the fiber products are formed as:
    \begin{enumerate}
        \item for the first fiber product, we use the free module functor \eqref{eq: free module functor LMcomm} and the source map for $\Conv_{\langle 1\rangle}(\LMcomm\lAlg(\bC))$;
        \item for the second fiber product, we use the neutral map for $\Conv_{\langle 1\rangle}(\bC)$ and the forgetful map \eqref{eq: LMcomm forget to module}.
    \end{enumerate}

    Using the fact that $\Conv_{\langle 1\rangle}(\mhyphen)$ commutes with $\Arr(\mhyphen)$ and that $\Conv_{\langle 1\rangle}(\mhyphen)$ commutes with fiber products, we can rewrite \eqref{eq: mechanism proof intermediate} as \[\Conv_{\langle 1\rangle}\left( \CAlg(\bC)\times_{\LMcomm\lAlg(\bC)} \Arr(\LMcomm\lAlg(\bC))\times_{\Arr(\bC)}\bC\right),\]
    i.e., the $\Conv_{\langle 1\rangle}(\mhyphen)$ applied to \eqref{eq: arrow fib prod 2}. From Lemma \ref{lem: functor btw arrow fib prods} this admits a canonical functor to
    \[\Conv_{\langle 1\rangle}\left(\Arr(\Alg(\LMcomm\lAlg(\bC)))\times_{\Arr(\CAlg(\bC))}\CAlg(\bC)\right).\]
    Using again the commutativity of  $\Conv_{\langle 1\rangle}$ with $\Arr$ and fiber products, the above category is identified with \eqref{eq: ConvMod fib prod 2}, and hence the proof.
\end{proof}

\sssec{}\label{sssec: description image}
We give a more concrete description of the functor from \eqref{eq: ConvMod fib prod 1} to \eqref{eq: ConvMod fib prod 2}. Suppose we are given an object of the former, thought in terms of the list of data in Sect. \ref{sssec: ConvMod fib 1 data}. To specify what its image in \eqref{eq: ConvMod fib prod 2} is, we use the description in Sect. \ref{sssec: ConvMod fib 2 data}. Using the notation therein, the map $\cC:\Conv_{\langle 1\rangle}\rightarrow \CAlg(\bC)$ is the image of $\cC_1$ under the functor
\[\Convcomm\lAlg(\bC)\xrightarrow{\te{Lemma }\ref{lem: Convcomm and Einfty}} \Conv_{\langle 1\rangle}(\CAlg(\bC)),\]
which we denote by the same notation $\cC_1$. Another data is an arrow in the category
\begin{equation}\label{eq: Alg cC_1 Mod}
    \Fun^{\circ}(\fSet,\fSet+\Alg(\cC_1\lMod(\bC))).
\end{equation}
Its source is the functor
$I\mapsto [I, (\cC_2)_I]$
and the target is the functor
$I\mapsto [I, \cHom_{\bC}(M_0, M_0\otimes (\cC_1)_I)]$. We can describe the arrow colorwise: when evaluated at the color $I$ it should be a map
\begin{equation}\label{eq: arrow at level I}
    (\cC_2)_I\rightarrow \cHom_{\bC}(M_0, M_0\otimes (\cC_1)_I)\quad \te{in } \Alg((\cC_1)_I\lMod(\bC)).
\end{equation}
To describe this, note that the data
\[(M_0\otimes\cC_1;\cC_1)\rightarrow (M_0\otimes\cC_1;\cC_2)\quad\te{in }\ConvcMod\lAlg(\bC)\]
from Sect. \ref{sssec: ConvMod fib 1 data} is evaluated to
\begin{equation}\label{eq: arrow at level I 2}
    ((\cC_1)_I, M_0\otimes(\cC_1)_I)\rightarrow ((\cC_2)_I,M_0\otimes (\cC_1)_I) \quad \te{in }\LMcomm\lAlg(\bC)
\end{equation}
at $I$. Now, the arrow \eqref{eq: arrow at level I} is obtained from the arrow \eqref{eq: arrow at level I 2} via Lemma \ref{lem: functor btw arrow fib prods}.

\subsection{Construction of the spectral action}\label{ssec: construction of the spectral action}
Combining the results developed earlier, we prove Theorem \ref{thm: spectral action} in this subsection.
We observe:

\begin{lemma}
    If $X=\pjl$ and the metaplectic dual gerbe $\cG_Z$ is not trivializable, then the stack $\Loc^{\cG_Z}_H(X)$ of twisted $H$-local systems is empty.
\end{lemma}
\begin{proof}
    A $Z_H(E)^{\on{tors}}$-gerbe on $\pjl$ is classified by $a\in Z_H(E)^{\on{tors}}$, thought of as a second cohomology class. We may consider the following model of $\cG_Z$: it is trivialized on $\pjl\setminus 0$ and $\pjl\setminus \infty$, and the transition data between those charts is given by the $Z_H(E)^{\on{tors}}$-local system on $\pjl\setminus\{0,\infty\}$ with monodromy $a$. Then the Cartesian diagram
    \[\begin{tikzcd}
	{\Loc^{\cG_Z}_H(\pjl)} & {\Loc^{\cG_Z}_H(\pjl\setminus \infty)} \\
	{\Loc^{\cG_Z}_H(\pjl\setminus 0)} & {\Loc^{\cG_Z}_H(\pjl\setminus\{0,\infty\})}
	\arrow[from=1-1, to=1-2]
	\arrow[from=1-1, to=2-1]
	\arrow[from=1-2, to=2-2]
	\arrow[from=2-1, to=2-2]
\end{tikzcd}\]
gets identified with
\[
\begin{tikzcd}
	{\Loc^{\cG_Z}_H(\pjl)} & {\pt/ H} \\
	{\Bet H} & {H/H,}
	\arrow[from=1-1, to=1-2]
	\arrow[from=1-1, to=2-1]
	\arrow["a", from=1-2, to=2-2]
	\arrow["0"', from=2-1, to=2-2]
\end{tikzcd}
\]
where the arrow $0$ is induced from the unit map $\pt\rightarrow H$ while the arrow $a$ is induced from $\pt\xrightarrow{a} H$. This fiber product is empty unless $a=0$, or equivalently $\cG_Z$ is trivializable.
\end{proof}

\sssec{}\label{sssec: assumption}
Therefore, for the proof of Theorem \ref{thm: spectral action} we may assume either $X\ne\pjl$ or $X=\pjl$ and $\cG_Z$ is trivializable. We retain this assumption for the remainder of this subsection. Note that this is exactly the assumption for Theorem \ref{thm: extension of action}.

\sssec{}
    We note that Proposition \ref{prop: extension mechanism} serves as a categorical mechanism to pass from the data obtained in Theorem \ref{thm: extension of action} to a natural transformation between functors
    \begin{equation}\label{eq: fSet Alg LS}
        \fSet\rightarrow \fSet+\Alg(\LS(X)^{\otimes}\lMod)
    \end{equation}
    respecting the projection.
\medskip

Consider the situation
\[\bC= \DGcat,\quad \cC_1=\LS(X)^{\otimes},\quad \cC_2=\left(\Rep(H)^{\loc}_{\cG_Z}(X)\right)^{\otimes},\ \te{ and }\ M_0=\cA^q_G\]
and the resulting categories \eqref{eq: ConvMod fib prod 1} and \eqref{eq: ConvMod fib prod 2}. Theorem \ref{thm: extension of action} tells us that the arrow \eqref{eq: ConvMod extension} defines an object of the category \eqref{eq: ConvMod fib prod 1}.
\medskip

Therefore, applying Proposition \ref{prop: extension mechanism} we get an arrow in \eqref{eq: Alg LS construction for C}, which is a natural transformation between functors \eqref{eq: fSet Alg LS} in this case. Using Sect. \ref{sssec: description image}, the source functor is 
\begin{equation}\label{eq: cat input global source}
    \fSet\rightarrow\fSet+\Alg(\LS(X)^{\otimes}\lMod),\quad I\mapsto [I, (\Rep(H)^{\loc}_{\cG_Z}(X))^{\otimes I}],
\end{equation}
which appeared earlier in \eqref{eq: Satake functorial B}, and the target functor is
\[\fSet\rightarrow\fSet+\Alg(\LS(X)^{\otimes}\lMod),\quad I\mapsto [I, \Fun(\cA^q_G,\cA^q_G\otimes\LS(X)^{\otimes I})].\]

\begin{lemma}\label{lem: LS lin endo}
    Let $M\in\DGcat$ and $Y\in\Spc$. There is a canonical equivalence
    \[\End(M)\otimes\LS(Y)= \Fun(M,M\otimes\LS(Y)).\]
\end{lemma}
\begin{proof}
    For any $B\in\DGcat$ we have
    \[B\otimes \LS(Y)=B\otimes (\colim_Y \Vect) = \colim_Y B =\lim_Y B,\]
    where we passed to adjoints in the last equivalence. Hence, we have a sequence of equivalences
    \[\Fun(M,M\otimes\LS(Y))=\lim_Y \Fun(M, M) = \End(M)\otimes\LS(Y)\]
    as desired.    
\end{proof}

\sssec{}\label{sssec: cat input global}
Using the lemma above, we rewrite the target functor as
\begin{equation}\label{eq: cat input global target}
    \fSet\rightarrow\fSet+\Alg(\LS(X)^{\otimes}\lMod),\quad I\mapsto [I, \End(\cA^q_G)\otimes\LS(X)^{\otimes I}].
\end{equation}
Theorem \ref{thm: extension of action} yields a natural transformation from \eqref{eq: cat input global source} to \eqref{eq: cat input global target} respecting the projection.

\sssec{Twisted factorization homology}\label{sssec: applying twfact}
Recall the categorical formalism of the action of twisted factorization homology, developed in Sect. \ref{sec: action of twisted fact hlgy}.  We apply Theorem \ref{thm: main categorical thm} to the situation
\[Y=X^{\topl,\sing},\quad A=Z_H(E)^{\on{tors}}, \quad \cG=\cG_Z, \ \te{ and } \ \cC=\Rep(H).\]

As in Sect. \ref{sssec: main cat thm prep} we have functors
\begin{equation}\label{eq: cat thm input source}
    \fSet\rightarrow\fSet+\LSMon,\quad I\mapsto \left[I,\left(\int^{\plain}_{\cG_Z}\Rep(H)\right)^{\otimes I}\right]
\end{equation}
and 
\begin{equation}\label{eq: cat thm input target}
    \fSet\rightarrow\fSet+\LSMon,\quad I\mapsto [I, \End(\cA^q_G)\otimes\LS(Y^I)].
\end{equation}

A natural transformation from \eqref{eq: cat thm input source} to \eqref{eq: cat thm input target} respecting the projection is equivalent to a map
\[\int^{\CAlg}_{\cG_Z}\Rep(H)\rightarrow \End(\cA^q_G)\quad\te{in }\DGCatMon,\]
which is equivalent to the action of the twisted factorization homology $\int^{\CAlg}_{\cG_Z}\Rep(H)$ on $\cA^q_G$. In Theorem \ref{thm: QC over twlocsys as a tw fact hlgy} we identified $\int^{\CAlg}_{\cG_Z}\Rep(H)$ with $\QC(\Loc^{\cG_Z}_H)$ as symmetric monoidal DG categories. Therefore, we finish the proof of Theorem \ref{thm: spectral action} once we construct a canonical natural transformation from \eqref{eq: cat thm input source} to \eqref{eq: cat thm input target}. In fact, we prove:

\begin{prop}\label{prop: equiv LS global}
   There is a canonical equivalence between the space of natural transformation from \eqref{eq: cat input global source} to \eqref{eq: cat input global target} and the space of natural transformation from \eqref{eq: cat thm input source} to \eqref{eq: cat thm input target}, where all natural transformations considered here preserve the projection.
\end{prop}
\begin{proof}
    Taking global sections of local systems of monoidal categories over $(X^{\sing,\topl})^I$'s yields a canonical functor
    \[\Gamma: \Fun^{\circ}(\fSet,\fSet+\LSMon)\rightarrow \Fun^{\circ}(\fSet,\fSet+\Alg(\LS(X)^{\otimes}\lMod)).\]
    (see Sect. \ref{ssec: LSCat}) We observe that this functor sends \eqref{eq: cat thm input source} to \eqref{eq: cat input global source} and \eqref{eq: cat thm input target} to \eqref{eq: cat input global target}. Indeed, this is obvious for the latter case, and for the former case we use the chain of equivalences
    \begin{align*}
    \Gamma((X^{\topl,\sing})_B, \cL_{(X^{\topl,\sing}),\Phi(\cG_Z)_B,\Bet H})&=\int^{\plain}_{(X^{\topl,\sing},\Phi(\cG_Z))}\Rep(H)\\
    &=\QC\left(\Tot(\Phi(\cG_Z)_B)\times^{\Bet A}\Bet H\right)\\
    &=\Rep(H)^{\loc}_{\cG_Z}(X).
    \end{align*}
    For the first two equivalences we used Sect. \ref{sssec: examples from gerbes} and for the last equivalence we used Theorem \ref{thm: twisted representations}.
    \medskip

    If $X\ne\pjl$, all the $(X^{\topl,\sing})^I$'s are 1-truncated, and hence the functor $\Gamma$ above is an equivalence (Theorem \ref{thm: 1-affineness}), and hence the proof.
    \medskip

    If $X=\pjl$, there is a canonical gerbe $\cG^0\in \on{Ge}_{Z_H(E)^{\on{tors}}}(\pt)$ and a canonical identification of $\cG_Z$ with the pullback of $\cG^0$ along $X\rightarrow\pt$. (see Sect. \ref{sssec: X=pjl and trivializable}) We claim that both the space of natural transformations $\eqref{eq: cat thm input source}\rightarrow \eqref{eq: cat thm input target}$ and $\eqref{eq: cat input global source}\rightarrow\eqref{eq: cat input global target}$ are equivalent to the space of natural transformations between the functors $\fSet\rightarrow\DGCatMon$, where the source is $I\mapsto (\Rep(H)_{\cG^0})^{\otimes I}$ and the target is $I\mapsto \End(\cA^q_G)\otimes\LS(X)^{\otimes I}$. 
\medskip
   
    The functor \eqref{eq: cat thm input source} can be rewritten as
    \[\fSet\rightarrow \fSet+\LSMon,\quad I\mapsto \left[I,\ul{(\Rep(H)_{\cG^0})^{\otimes I}}_{X^I}\right].\]
    Note that the local system of monoidal categories are all pulled back from a point. Using the adjunction
    \[f^*:\DGCatMon=\LSMon(\pt)\leftrightarrows
    \LSMon((X^{\topl,\sing})^I):f_*\]
    from Sect. \ref{sssec: LSCat pull push}, we confirm the above claim for $\eqref{eq: cat thm input source}\rightarrow \eqref{eq: cat thm input target}$. 
    \medskip

    On the other hand, the functor \eqref{eq: cat input global source} can be rewritten as
    \[I\mapsto [I, (\Rep(H)_{\cG^0})^{\otimes I}\otimes\LS(X)^{\otimes I}].\]
    Using the induction-forgetful adjunction
    \[\DGCatMon\leftrightarrows\Alg(\LS(X)^{\otimes I}\lMod),\]
    we confirm the above claim for $\eqref{eq: cat input global source}\rightarrow\eqref{eq: cat input global target}$ as well, completing the proof.
\end{proof}

\sssec{}
From the above proposition, the canonical natural transformation $\eqref{eq: cat input global source} \rightarrow \eqref{eq: cat input global target}$ in Sect. \ref{sssec: cat input global} gives rise to a canonical natural transformation $\eqref{eq: cat thm input source} \rightarrow \eqref{eq: cat thm input target}$. From the discussion of Sect. \ref{sssec: applying twfact} we finish the proof of Theorem \ref{thm: spectral action}, and hence the construction of spectral action in Betti metaplectic Langlands.

\appendix
\section{Twisting constructions}\label{sec: twisting construction}
The procedure of twisting a sheaf of categories was given in \cite[Sect. 1.7]{GL}, and is written in more detail in \cite[Appendix A]{ZhaoNew}. In this section, we review the construction and rewrite it in a form that will be useful later. The method of twisting discussed in this section is largely borrowed from the references mentioned above. However, in this section we give a more systematic treatment of how to twist sheaves of algebras and modules. The main reason is that in the factorizable Satake input 
\[\left(\Rep(H)^{\loc}_{\cG_Z}(X)\right)^{\otimes I}\rightarrow \Shv^q(\Hecke^{\loc}_{X^I})\]
needed to obtain the spectral action, it is crucial that the functor is a morphism in $\Alg(\LS(X^I)\text{-}\Mod)$. In the construction of the functor above, we need to make use of twisting constructions, and hence we have to make sure that the twistings we apply interact nicely with the algebra and module structures.
\medskip

We develop the methodology of twisting algebras and modules throughout Sect. \ref{ssec: Yifei construction} and Sect. \ref{ssec: twisting modules}. We do it in a couple of stages, with increasing structure on sheaf of categories. In each setting, we discuss the data needed to twist the sheaf while preserving the structure, and provide a more practical way of generating the data needed using linearity. In fact, such twisting using linearity suffices for most of the twisting construction we need for this paper. In developing the twisting formalism, we introduce two key players in this paper: twisted sheaves and twisted representations, which also serve as our running examples.
\medskip

The most subtle part in the twisting constructions needed for this paper is the $\LS(X)$-linear monoidal twist of $\Shv^q(\Hecke^{\loc}_X)$. Although this fits into the framework in Sect. \ref{ssec: twisting modules}, this is more naturally understood using the framework of twisting by grading, which will be discussed in Sect. \ref{ssec: twisting from grading}. Under a certain $\Z$-linearity assumption, twisting by grading agrees with the previous method. Although this will be the case for most of the applications in our paper, when studying the perverse ULA version of Satake, e.g. in Sect. \ref{ssec: Perv ULA Satake in Conv language}, we have to alter the commutativity constraint, and in that case it is necessary to use the method of Sect. \ref{ssec: twisting from grading}.
\medskip

In Sect. \ref{ssec: restriction to triv} we give a more explicit description of the twisting procedure by restricting to a subsite where a given gerbe is trivializable. This serves as an important computational tool and yields us a very concrete understanding of twisted sheaves (Theorem \ref{thm: twisted sheaves tensor product formula}) and twisted representations (Theorem \ref{thm: twisted representations}). 

\subsection{Sheaf of categories of sheaves}\label{ssec: Yifei construction}
Let $\bO$ be a category with limits and $\cX$ be a site, where $(\cX,\bO)$ satisfies some conditions stated below.
Suppose that we have a sheaf $\cF\in \Shv(\cX;\bO)$ and a group-like $\E_2$-algebra object $\cA$ in $\Shv(\cX;\Cat)$. Given a data $\xi$ we introduce later, we get a twisted sheaf $\cF_{\cG}\in \Shv(\cX;\bO)$ for an $\cA$-gerbe $\cG\in \Gamma(\cX,B^2 \cA)$.

\sssec{Conditions on $(\cX,\bO)$}\label{sssec: conditions on cX bO}
The condition we impose on the site $\cX$ is that the topology is generated by a collection of morphisms $\kappa$ in $\cX$ satisfying
\begin{enumerate}
    \item $\kappa$ is stable under pullback;
    \item if we have $x_1\rightarrow x_2\rightarrow x$ in $\cX$ such that both $x_2\rightarrow x$ and $x_1\rightarrow x$ lie in $\kappa$, then $x_1\rightarrow x_2$ lies in $\kappa$.
\end{enumerate}
More precisely, given $x\in\cX$, a sieve $\cS\subseteq \cX_{/x}$ is a covering sieve if and only if it contains a morphism $x'\rightarrow x$ in $\kappa$. Then a presheaf on $\cX$ is a sheaf with respect to this topology if and only if it satisfies descent along the \v{C}ech nerve of any morphism in $\kappa$.
\medskip

Moreover, we assume the existence of the sheafification functor
\begin{equation}\label{eq: sheafification functor}
    (\mhyphen)^{\on{sh}}:\PSh(\cX;\bO)\rightarrow \Shv(\cX;\bO),
\end{equation}
i.e., the left adjoint to the inclusion functor
\[\Shv(\cX;\bO)\hookrightarrow \PSh(\cX;\bO).\]
We have a functor
\begin{equation}\label{eq: constant sheaf functor}
    \ul{(\mhyphen)}:\bO\rightarrow \Shv(\cX;\bO),\quad \cO\mapsto \ul{\cO}
\end{equation}
 where $\ul{\cO}$ is the sheafification of the constant presheaf with value $\cB$. We simply say that $\ul{\cO}$ is the constant sheaf with value $\cO$.
 \medskip
 
Our favorite example of such $(\cX,\bO)$ would be
\[\cX=(\Affft)_{/\cY} \te{ for } \cY\in\PreStklft \quad\te{and}\quad \bO=\DGcat,\] where $\kappa$ is the collection of \'etale surjections so that $\cX$ is equipped with the \'etale topology.
\begin{prop}\label{prop: sheaf of sheaf categories}
    The functor
\begin{equation}
    \ul{\Shv}(\cX;\bO): \cX^{\op}\rightarrow \Cat,\quad y\mapsto \Shv(\cX_{/y};\bO)
    \end{equation}
    satisfies descent, hence is an object in $\Shv(\cX;\Cat)$.
\end{prop}

\sssec{Twisting construction}\label{sssec: twisting construction}
The sheaf $\cF$ defines a pointing $\ul{\cF}$ of $\ul{\Shv}(\bO)$ as a morphism 
\begin{equation}\label{eq: pointing by cF}
    \ast\rightarrow \ul{\Shv}(\cX;\bO),\quad y\mapsto \cF|_{\cX_{/y}}
\end{equation}
in $\Shv(\cX;\Cat)$ where $\ast$ is the constant sheaf valued in $\pt$.
The data $\xi$ we need for the twisting construction is an
$\E_0$-map
\begin{equation}\label{eq: twisting data}
\xi:B^2 \cA\rightarrow \ul{\Shv}(\cX;\bO)
\end{equation}
in $\Shv(\cX;\Cat)$,
where the $\E_0$-structure on the target is $\ul{\cF}$.
Given such $\xi$, taking global sections we have a map
\[\Gamma(\cX;B^2 \cA)\rightarrow \Shv(\cX;\bO)\]
in $\E_0(\Cat)$, where the pointing on the target is $\cF$. We define the \textbf{twisted sheaf} $\cF_{\cG}\in \Shv(\cX;\bO)$ to be the image of $\cG\in \Gamma(\cX;B^2\cA)$ under this map.
\medskip

We call the data of $\xi$ the $\cA$-\textbf{linearity} on the sheaf $\cF$, because it is equivalent to the action of the group-like sheaf $B\cA$ on the sheaf $\cF$. In other words, for each $y\in \cX$ we have the identity map $\cF(y)\xrightarrow{\id}\cF(y)$, and there should be an $\E_2$-map
$\cA(y)\rightarrow \on{Aut}(\id_{\cF(y)})$
compatible with restriction of sections. These are justified by the following: through looping-delooping adjunction, the $\E_0$-map $\xi$ is equivalent to an $\E_1$-map
\[B\cA\rightarrow\Omega\ul{\Shv}(\cX;\bO)\]
and in turn an $\E_2$-map
\[\cA\rightarrow \Omega^2\ul{\Shv}(\cX;\bO)\]
where on the target we used \eqref{eq: pointing by cF} as the pointing of $\ul{\Shv}(\bO)$.

\sssec{Symmetric monoidal structure on the sheaf category}\label{sssec: sym mon str on Shv}
Let $(\cX,\bO)$ as in Sect. \ref{sssec: conditions on cX bO}.
Suppose $\bO$ is a symmetric monoidal category. This induces a symmetric monoidal structure on $\PSh(\cX;\bO)$. We further impose the condition:
\[\bO\te{ is cocomplete, and the tensor product of } \bO \te{ preserves all colimits in each variables.}\]
Then $\Shv(\cX;\bO)$ has a unique symmetric monoidal structure which makes the sheafification functor \eqref{eq: sheafification functor} symmetric monoidal. This implies that the constant sheaf functor \eqref{eq: constant sheaf functor} is symmetric monoidal.
\medskip

Given an operad $\cP$ we have equivalences
\begin{equation}\label{eq: PSh Shv and algebras}
    \cP\lAlg(\PSh(\cX;\bO))=\PSh(\cX;\cP\lAlg(\bO)),\quad \cP\lAlg(\Shv(\cX;\bO))=\Shv(\cX;\cP\lAlg(\bO))
\end{equation}
 compatible with forgetful maps to $\PSh(\cX;\bO)$ and $\Shv(\cX;\bO)$, respectively. Using the symmetric monoidality of \eqref{eq: sheafification functor} and identifications \eqref{eq: PSh Shv and algebras}, we get the sheafification functor
\[(\mhyphen)^{\on{sh}}:\PSh(\cX;\cP\lAlg(\bO))\rightarrow \Shv(\cX;\cP\lAlg(\bO)).\]
 In conclusion, if the pair $(\cX,\bO)$ satisfies the condition of Sect. \ref{sssec: conditions on cX bO}, then the pair $(\cX,\cP\lAlg(\bO))$ also satisfies the condition.

\sssec{Twisting a sheaf of algebras}\label{sssec: twisting algebra construction}
In the situation of Sect. \ref{sssec: sym mon str on Shv}, we may consider a sheaf of algebras $\cF\in \Shv(\cX;\cP\lAlg(\bO))$ and the twisting procedure as a sheaf valued in $\cP\lAlg(\bO)$. We proceed the same as in Sect. \ref{sssec: twisting construction} but replace $\bO$ by $\cP\lAlg(\bO)$. Hence the data we need for the twisting is an $\E_0$-map
\begin{equation}\label{eq: twisting data algebra}
    \xi:B^2\cA\rightarrow \ul{\Shv}(\cX;\cP\lAlg(\bO))
\end{equation}
in $\Shv(\cX;\Cat)$. Given $\cG\in\Gamma(\cX;B^2\cA)$ we get the twisted sheaf of algebras $\cF_{\cG}\in\Shv(\cX;\cP\lAlg(\bO))$.

\sssec{}\label{sssec: twisting and forgetting}
Twisting with monoidality is compatible with that without monoidality in an obvious way. To wit, consider the sheaf
$\cF^{\on{oblv}}\in\Shv(\cX;\bO)$ obtained by forgetting the algebra structure of $\cF$. Then
 \eqref{eq: twisting data algebra} composed with the forgetful map
\[\ul{\Shv}(\cX;\cP\lAlg(\bO))\rightarrow \ul{\Shv}(\cX;\bO)\]
yields an $\E_0$-map
\[\xi^{\on{oblv}}: B^2\cA\rightarrow \ul{\Shv}(\cX;\bO)\]
which serves as a twisting data of $\cF^{\on{oblv}}$. Hence given $\cG\in\Gamma(\cX;B^2\cA)$ we get the twisted sheaf $(\cF^{\on{oblv}})_{\cG}\in\Shv(\cX;\bO)$. The compatibility is that we get $(\cF^{\on{oblv}})_{\cG}$ from $\cF_{\cG}$ by forgetting the algebra structure.

\sssec{Twisting by linearity}\label{sssec: twisting by linearity}
In practice, the data of twisting arise from a simpler setup. Let $(\cX,\bO)$ as in Sect. \ref{sssec: sym mon str on Shv}. Suppose that $\cA$ arises as a constant sheaf $\ul{\sfA}$ for an abelian group $\sfA$ and that the sheaf $\cF\in\Shv(\cX;\bO)$ has a form $\ul{\cB}\otimes \cF_0$ for $\cB\in\bO$ and $\cF_0\in \Shv(\cX;\bO)$. In this case, an $\sfA$-linearity structure on the object $\cB$ induces a $\cA$-linearity structure on the sheaf $\cF$ by the functoriality of tensor product.
\medskip

In detail, assume that we have an $\E_0$-map $\xi^{\on{lin}}:B^2\sfA\rightarrow \bO$ where the pointing on the target is $\cB\in\bO$. The functor
\[\ul{(\mhyphen)}\otimes \cF_0: \bO\rightarrow \Shv(\cX_{/y};\bO)\]
induces a compatible family of $\E_0$-maps $\bO\rightarrow \Shv(\cX_{/y};\bO)$ in $\Cat$
for each $y\in\cX$ where the pointing on the left-hand side is $\cB$ and the pointing on the right-hand side is $\ul{\cB}\otimes \cF|_{\cX_{/y}}$. The collection of compositions
\[\xi^{\on{pre}}_y:B^2 \sfA\xrightarrow{\xi^{\on{lin}}} \bO\xrightarrow{\ul{(\mhyphen)}\otimes\cF_0} \Shv(\cX_{/ y};\bO)\]
for $y\in \cX$ defines an $\E_0$-map
\[\xi^{\on{pre}}:\ul{B^2\sfA}^{\on{pre}}\rightarrow \ul{\Shv}(\cX;\bO)\]
in $\PSh(\cX;\Cat)$, where the source is the constant presheaf valued in $B^2 \sfA$. We get the desired $\cA$-linearity from sheafifying $\xi^{\on{pre}}$.
\medskip

 Suppose we have $\cB\in\cP\lAlg(\bO)$ and assume that we have $\E_0$-map $\xi^{\on{lin}}:B^2\sfA\rightarrow \cP\lAlg(\bO)$ where the pointing on the target is $\cB$. Then, by the same logic, this determines an $\cA$-linearity on the $\cF$ as a sheaf of $\cP$-algebras.
\medskip

We introduce two examples of this construction, each of which plays a major role in our paper. Twisted sheaves will be studied further in Section \ref{sec: twisted sheaves}, and twisted representations will be studied further in Section \ref{sec: twisted representations}.

\sssec{Twisted sheaves}\label{sssec: twisted sheaves}
    Let $\bO=\DGcat$ and $\cB=\Vect$. Then the functor $\cB\otimes \mhyphen$ is an identity, so any sheaf $\cF$ is realized as $\ul{\cB}\otimes \cF_0$ where $\cB=\Vect$ and $\cF_0=\cF$. We take $\sfA$ to be a finite abelian group $A$ together with a character $\chi:A\rightarrow E^{\times}$. Then $\chi$ determines an $\E_0$-map $B^2 A\rightarrow \DGcat$ where the pointing on the target is $\Vect$, because any element of $E^{\times}$ determines an automorphism $\Vect\xrightarrow{\id}\Vect$ by scaling. Therefore, any sheaf $\cF$ valued in $\DGcat$ can be twisted by a point in $\Gamma(\cX,B^2 A)$.
\medskip

As an example, let $\cY\in \PreStklft$ and $\cX=(\Affft)_{/\cY}$ with \'etale topology. Then we take the sheaf $\cF=\mathbf{Shv}$ as
\begin{equation}\label{eq: bold Shv}
    \mathbf{Shv}:\left((\Affft)_{/\cY}\right)^{\op}\rightarrow \DGcat,\quad S\mapsto \Shv(S).
\end{equation}
We have $\Gamma(\cX;B^2\cA)=\on{Ge}_A(\cY)$.
Given a gerbe $\cG\in \on{Ge}_A(\cY)$ we get a twist $\mathbf{Shv}_{\cG}$, and we let
\[\Shv^{\cG}(\cY):=\Gamma((\Affft)_{/\cY}, \mathbf{Shv}_{\cG}).\]
This is the category of \textbf{twisted sheaves} on $\cY$.
\medskip

As another example, we take $\cF=\mathbf{LS}$ as
\begin{equation}\label{eq: bold LS}
    \mathbf{LS}: \left((\Affft)_{/\cY}\right)^{\op}\rightarrow \DGcat,\quad S\mapsto \LS(S).
\end{equation}
Given a gerbe $\cG\in\on{Ge}_A(\cY)$ we get a twist $\mathbf{LS}_{\cG}$, and we let
\[\LS^{\cG}(\cY):=\Gamma((\Affft)_{/\cY}, \mathbf{LS}_{\cG}).\]
This is the category of \emph{twisted local systems} on $\cY$.

\sssec{Twisted representations}\label{sssec: twisted representations}
Let $\bO=\DGcat$, $\cP=\E_{\infty}$, and $\cX=Y_{\et}$ be the small \'etale site of a scheme $Y\in\Schft$, so that $\cP\lAlg(\bO)=\DGCatSymMon$. We have a sheaf $\cF_0=\mathbf{LS}$ valued in $\DGCatSymMon$. Take $\cB=\Rep(H)\in\DGCatSymMon$ for an affine algebraic group $\mathsf{H}$ over $E$. 
Then $\cF$ is
\begin{equation}\label{eq: bold Rep H loc}
    \mathbf{Rep}(\mathsf{H})^{\loc}: Y^{\op}_{\et}\rightarrow \DGCatSymMon,\quad Y'\mapsto \Rep(\mathsf{H})\otimes\LS(Y').
\end{equation}
This is a sheaf because of the dualizability of $\Rep(\mathsf{H})$.
\medskip

We take $\sfA$ to be a torsion abelian group with a group homomorphism $\phi:\sfA\rightarrow Z_\mathsf{H}(E)$, then we have an $\E_0$-map $B^2 \sfA\rightarrow \DGCatSymMon$ where the pointing on the target is $\Rep(\mathsf{H})$. Indeed, from the natural $B\sfA$-action of $\Bet \mathsf{H}$  we have an $\E_0$-map 
\[B^2\sfA\rightarrow(\PreStklft)^{\op},\quad \ast\mapsto \Bet \mathsf{H},\]
and we compose this with 
\[\on{QC}:(\PreStklft)^{\op}\rightarrow\DGCatSymMon.\]

Given an \'etale gerbe $\cG\in \on{Ge}_{\sfA}(Y)$ we get a twist $\mathbf{Rep}(\mathsf{H})^{\loc}_{\cG}$, and we let
\[\Rep(\mathsf{H})^{\loc}_{\cG}(Y):=\Gamma(Y,\mathbf{Rep}(\mathsf{H})^{\loc}_{\cG}) \in\DGCatSymMon.\]
This is the category of \textbf{twisted} $\mathsf{H}$-\textbf{representations} on local systems on $Y$.

\subsection{Twisting and module structure}\label{ssec: twisting modules}
Now we establish the versions of twisting constructions mentioned above in the situation where the sheaf of categories has a module structure. For this purpose, we recall some standard categorical lemmas.

\sssec{}
Suppose we have a diagram
    \[F:I\rightarrow\CAlg(\Cat),\quad i\mapsto F(i)\]
    and $\wt{F}$ be the limit of this diagram. Consider $\wt{A}\in\Alg(\wt{F})$ and let $A(i)\in F(i)$ be the image of $\wt{A}$ under $\on{ev}_i:\wt{F}\rightarrow F(i)$. Since $\on{ev}_i$ is symmetric monoidal, $A(i)\in\Alg(F(i))$, and this functor induces
    \[\wt{A}\lMod(\wt{F})\rightarrow A(i)\lMod(F(i)).\]
    Hence we get
    \begin{equation}\label{eq: module cat comparison}
        \wt{A}\lMod(\wt{F})\rightarrow \lim_{i\in I}A(i)\lMod(F(i))
    \end{equation}
    in $\Cat$. Then we have:
    
\begin{lemma}\label{lem: module category limit}
    The functor \eqref{eq: module cat comparison} is an equivalence.
\end{lemma}

\begin{lemma}\label{lem: limit and algebra objects}
Suppose $\cP$ is an operad.
    The canonical functor
    \[\cP\lAlg(\wt{F})\rightarrow\lim_{I}\cP\lAlg(F(i))\]
    in $\CAlg(\Cat)$ is an equivalence.
\end{lemma}

\sssec{Sheaf of modules}\label{sssec: setup sheaf of modules}
Suppose $(\cX,\bO)$ as in Sect \ref{sssec: sym mon str on Shv}. Consider $\cC\in \Shv(\cX;\CAlg(\bO))$. Then $\cC$ is an $\E_{\infty}$-algebra object in the symmetric monoidal category $\Shv(\cX;\bO)$, hence we have a symmetric monoidal category
$\cC\lMod(\Shv(\cX;\bO))$.

\begin{prop}\label{prop: sheaf of module categories}
    The functor
    \begin{equation}\label{eq: sheaf of module sheaves}
        \ul{\Shv}(\cX;\cC\lMod(\bO)):\cX^{\on{op}}\rightarrow \Cat,\quad y\mapsto \cC|_{\cX_{/y}}\lMod\left(\Shv(\cX_{/y};\bO)\right)
    \end{equation}
    is a sheaf, and lifts canonically to an object in $\Shv(\cX;\CAlg(\Cat))$.
\end{prop}
\begin{proof}
    We first observe that \eqref{eq: sheaf of module sheaves}
    lifts to a presheaf valued in $\CAlg(\Cat)$, because restriction functors are symmetric monoidal. Since the forgetful functor $\CAlg(\Cat)\rightarrow \Cat$ is conservative, it suffices to show that \eqref{eq: sheaf of module sheaves}
    as a presheaf valued in $\Cat$ is a sheaf.
    This follows from Proposition \ref{prop: sheaf of sheaf categories} and Lemma \ref{lem: module category limit}.
\end{proof}

\begin{cor}
    Suppose $\cP$ is an operad. Then the functor
    \[\ul{\Shv}(\cX;\cP\lAlg(\cC\lMod(\bO))):\cX^{\op}\rightarrow\CAlg(\Cat), \quad y\mapsto \cP\lAlg\left(\cC|_{y}\lMod(\Shv(\cX_{/y};\bO))\right)\]
    is a sheaf.
\end{cor}
\begin{proof}
    This follows immediately from Proposition \ref{prop: sheaf of module categories} and Lemma \ref{lem: limit and algebra objects}.
\end{proof}

\begin{prop}
    We have a canonical functor
    \begin{equation}\label{eq: P Alg A Mod Shv to PSh}
        \cP\lAlg(\cC\lMod(\Shv(\cX;\bO)))\rightarrow \cP\lAlg(\cC\lMod(\PSh(\cX;\bO)))
    \end{equation}
    which forgets to
    \[\iota: \Shv(\cX;\bO)\hookrightarrow \PSh(\cX;\bO),\]
    i.e., making the diagram
\[\begin{tikzcd}
	{\cP\lAlg(\cC\lMod(\Shv(\cX;\bO)))} & {\cP\lAlg(\cC\lMod(\PSh(\cX;\bO)))} \\
	{\Shv(\cX;\bO)} & {\PSh(\cX;\bO)}
	\arrow[from=1-1, to=1-2]
	\arrow[from=1-1, to=2-1]
	\arrow[from=1-2, to=2-2]
	\arrow["\iota"', from=2-1, to=2-2]
\end{tikzcd}\]
    commute, where the vertical arrows are forgetful functors.
\end{prop}
\begin{proof}
    We have an adjunction \begin{equation}\label{eq: A Mod PSh Shv adjunction}
        \cC\lMod(\PSh(\cX;\bO))\leftrightarrows \cC\lMod(\Shv(\cX;\bO))
    \end{equation}
    that forgets to the adjunction
    \begin{equation}\label{eq: PSh Shv adjunction}
        (\mhyphen)^{\on{sh}}:\PSh(\cX;\bO)\leftrightarrows \Shv(\cX;\bO):\iota.
    \end{equation}

    In the adjunction \eqref{eq: A Mod PSh Shv adjunction}, the functor
    \[\cC\lMod(\PSh(\cX;\bO))\rightarrow \cC\lMod(\Shv(\cX;\bO))\]
    is symmetric monoidal, because the sheafification \eqref{eq: sheafification functor}
    is symmetric monoidal and preserves colimits, and hence preserves relative tensor products. Hence the adjunction \eqref{eq: A Mod PSh Shv adjunction} induces an adjunction
    \[\cP\lAlg(\cC\lMod(\PSh(\cX;\bO)))\leftrightarrows\cP\lAlg(\cC\lMod(\Shv(\cX;\bO))).\]
    By construction, this adjunction forgets to the adjunction \eqref{eq: PSh Shv adjunction}.
\end{proof}

\begin{cor}\label{cor: algebra structure on section}
    Let
\[\cF\in\cP\lAlg(\cC\lMod(\Shv(\cX;\bO))).\]
For any $y\in\cX$, the object $\cF(y)\in\bO$ canonically lifts to an object of $\cP\lAlg(\cC(y)\lMod(\bO))$.
\end{cor}
\begin{proof}
    We can consider the image of $\cF$ under \eqref{eq: P Alg A Mod Shv to PSh}. Since the functor does not change the underlying presheaf, it means that $\cF^{\on{oblv}}\in\Shv(\cX;\bO)$ has a canonical lift to $\cP\lAlg(\cC\lMod(\PSh(\cX;\bO)))$. Now the rest follows from the fact that the evaluation map
\[\on{ev}_y: \cC\lMod(\PSh(\cX;\bO))\rightarrow \cC(y)\lMod(\bO)\]
is symmetric monoidal.
\end{proof}

\sssec{Twisting a sheaf of modules}\label{sssec: twisting construction modules}
We remain in the context of Sect. \ref{sssec: setup sheaf of modules}. The data needed to twist $\cF\in\cC\lMod(\Shv(\cX;\bO))$ as a $\cC$-module is an $\E_0$-map
\begin{equation}\label{eq: twisting data modules}
    \xi:B^2\cA\rightarrow \ul{\Shv}(\cX;\cC\lMod(\bO))\quad \te{in} \ \Shv(\cX;\Cat).
\end{equation}
Given a gerbe $\cG\in\Gamma(\cX;B^2\cA)$ we get a twist $\cF_{\cG}\in \cC\lMod(\Shv(\cX;\bO))$ as in Sect. \ref{sssec: twisting construction}. For each $y\in \cX$ we have $\cF_{\cG}(y)\in\cC(y)\lMod(\bO)$.
\medskip

Given an operad $\cP$, the data needed to twist $\cF\in \cP\lAlg(\cC\lMod(\Shv(\cX;\bO)))$ as a sheaf of $\cC$-linear $\cP$-algebras is an $\E_0$-map
\begin{equation}\label{eq: twisting data algebras in modules}
    \xi: B^2\cA\rightarrow \ul{\Shv}(\cX;\cP\lAlg(\cC\lMod(\bO)))\quad \te{in} \ \Shv(\cX;\Cat).
\end{equation}
Given a gerbe $\cG\in\Gamma(\cX;B^2\cA)$ we get a twist $\cF_{\cG}\in \cP\lAlg(\cC\lMod(\Shv(\cX;\bO)))$ as in Sect. \ref{sssec: twisting construction}. From Corollary \ref{cor: algebra structure on section} we have
$\cF_{\cG}(y)\in\cP\lAlg(\cC(y)\lMod(\bO))$ for each $y\in\cX$.
\medskip

Twisting $\cC$-modules and $\cP$-algebra objects in $\cC$-module (Sect. \ref{sssec: twisting construction modules}) are compatible with the twisting without any additional structure (Sect. \ref{sssec: twisting construction}) in the sense of Sect. \ref{sssec: twisting and forgetting}.

\sssec{Twisting modules by linearity}\label{sssec: twisting modules by linearity}
Twisting by linearity (Sect. \ref{sssec: twisting by linearity}) has its version in both of the contexts discussed in Sect. \ref{sssec: twisting construction modules}.
Suppose that $\cA$ arises as a constant sheaf $\ul{\sfA}$ for an abelian group $\sfA$.
\medskip

First, assume that $\cF\in\cC\lMod(\Shv(\cX;\bO))$ has a form $\ul{\cB}\otimes\cF_0$ where $\cB\in\bO$ and $\cF_0\in \cC\lMod(\Shv(\cX;\bO))$.  Then the $\E_0$-map $\xi^{\on{lin}}:B^2\sfA\rightarrow \bO$, where the pointing on the target is $\cB$, determines the twisting data \eqref{eq: twisting data modules}, because we have a compatible family of functors
\[\bO\xrightarrow{\ul{(\mhyphen)}\otimes\cF_0|_{\cX_{/y}}} \cA\lMod(\Shv(\cX_{/y};\bO))\te{ for } y\in \cX.\]

Second, assume that $\cF\in\cP\lAlg(\cC\lMod(\Shv(\cX;\bO)))$ has a form $\ul{\cB}\otimes\cF_0$ where $\cF_0\in\cP\lAlg(\cC\lMod(\Shv(\cX;\bO)))$ and $\cB\in \cP\lAlg(\bO)$. Then the same data $\xi^{\on{lin}}$ determines the data \eqref{eq: twisting data algebras in modules}, because we have a compatible family of functors
\[\cP\lAlg(\bO)\xrightarrow{\ul{(\mhyphen)}\otimes\cF_0|_{\cX_{/y}}} \cP\lAlg(\cC\lMod(\Shv(\cX_{/y};\bO))),\te{ for } y\in\cX.\]
This functor is obtained by inserting $\cF_0$ to the second argument of the functor
\[\cP\lAlg(\bO)\times \cP\lAlg(\cC\lMod(\Shv(\cX_{/y};\bO)))\rightarrow \cP\lAlg(\cC\lMod(\Shv(\cX_{/y};\bO)))\]
induced from the symmetric monoidal functor
\begin{equation}\label{eq: tensoring with a constant sheaf}
    \bO\times\cC\lMod(\Shv(\cX_{/y};\bO))\rightarrow \cC\lMod(\Shv(\cX_{/y};\bO)),\quad (\mathbf{o},\cH)\mapsto \ul{\mathbf{o}}\otimes\cH.
\end{equation}

\sssec{}\label{sssec: twisting modules examples}
In Sect. \ref{sssec: twisted sheaves} and Sect. \ref{sssec: twisted representations} we introduced two key players in our paper: twisted sheaves and twisted representations. We explain how those two examples fit into our framework. Recall the definitions of \[\mathbf{Shv},\quad \mathbf{LS}, \quad  \te{and} \quad \mathbf{Rep}(\mathsf{H})^{\on{loc}}\]
from those two sections. The first one is a sheaf over $(\Affft)_{/\cY}$ and the others are sheaves over $Y_{\et}$.
\medskip

Let $\bO=\DGcat$, $\cB=\Vect$, and $\sfA$ a finite abelian group $A$ together with a character $\chi:A\rightarrow E^{\times}$. Arguing as in Sect. \ref{sssec: twisted sheaves}, for any $\cC\in\CAlg(\Shv(\cX;\DGcat))$ and $\cF\in\cC\lMod(\Shv(\cX;\DGcat))$, we can twist $\cF$ as a $\cC$-module given a gerbe $\cG\in\Gamma(\cX;B^2 A)$. As an example, we take $\cX=(\Affft)_{/\cY}$ with the \'etale topology. Then $\mathbf{Shv}$ is an $\E_{\infty}$-algebra object in $\Shv((\Affft)_{/\cY};\DGcat)$. Hence, take $(\cC,\cF)=(\mathbf{Shv},\mathbf{Shv})$ and conclude that $\mathbf{Shv}_{\cG}$ is a module over $\mathbf{Shv}$. In particular, $\Shv^{\cG}(\cY)$ is a module over $\Shv(\cY)$.
\medskip

Let $\bO=\DGcat$ and $\cB=\Rep(\mathsf{H})$ as an object of $\CAlg(\DGcat)$. Then 
\[\mathbf{Rep}(\mathsf{H})^{\loc}=\ul{\Rep(\mathsf{H})}\otimes \mathbf{LS},\]
which realizes it as an object of 
\begin{equation}\label{eq: CAlg LS Mod}
\CAlg\left(\mathbf{LS}\lMod\left(\Shv(Y_{\et};\DGcat)\right)\right).
\end{equation}
Using Sect. \ref{sssec: twisting modules by linearity} we conclude that the twist
$\mathbf{Rep}(\mathsf{H})^{\on{loc}}_{\cG}$
has a natural lift to \eqref{eq: CAlg LS Mod}.

\subsection{Twisting from grading}\label{ssec: twisting from grading}
We remain in the context of Sect. \ref{sssec: sym mon str on Shv}.
Suppose that a sheaf $\cF\in \Shv(\cX;\bO)$ has a graded $\cP$-monoidal structure. Then we will discuss a method of twisting $\cF$ as a graded $\cP$-monoidal sheaf.

\sssec{}\label{sssec: LKE sym mon}
Let $\bC$ be a symmetric monoidal category which is cocomplete and the tensor product commutes with all colimits. $\Upsilon$ be an abelian group considered as a discrete symmetric monoidal category. The left-Kan extension
\[\on{LKE}_{\Upsilon\rightarrow\pt}:\Fun(\Upsilon,\bC)\rightarrow \bC\]
is symmetric monoidal, where the symmetric monoidal structure of $\Fun(\Upsilon,\bC)$ is given by the Day convolution. Hence, given an operad $\cP$ we get an induced functor
\begin{equation}\label{eq: LKE graded monoidal}
    \on{LKE}_{\Upsilon\rightarrow\pt}:\cP\lAlg(\Fun(\Upsilon,\bC))\rightarrow\cP\lAlg(\bC).
\end{equation}
Note that the source category is identified with the category of right-lax $\cP$-monoidal functors $\Upsilon\rightarrow \bC$.

\sssec{Graded monoidal structure}\label{sssec: graded monoidal structure}
Suppose we have the data $(\cX,\bO)$, $\cC$, and $\cP$ as in Sect. \ref{sssec: twisting construction modules}, and $\Upsilon$ as before. Furthermore, suppose we have a lax $\cP$-monoidal functor
\[F:\Upsilon\rightarrow \cC\lMod(\Shv(\cX;\bO)).\]
Applying the construction of Sect. \ref{sssec: LKE sym mon} to $\bC=\cC\lMod(\Shv(\cX;\bO))$, we get \[\on{LKE}_{\Upsilon\rightarrow\pt}(F)\in\cP\lAlg\left(\cC\lMod\left(\Shv(\cX;\bO)\right)\right).\] 

\sssec{}\label{sssec: twist by grading outline}
Suppose $\cF\in\cP\lAlg(\bC\lMod(\Shv(\cX;\bO)))$ is realized as $\on{LKE}_{\Upsilon\rightarrow\pt}(F)$ for some $\Upsilon$ and $F$ as above. Then with a certain $\Upsilon$-family of gerbes (say $\cG^\Upsilon$) which we elaborate on later, we can twist the data $F$ to a right-lax $\cP$-monoidal
\begin{equation}\label{eq: twist of grading data}
    F_{\cG^{\Upsilon}}:\Upsilon\rightarrow\cC\lMod(\Shv(\cX;\bO)),
\end{equation}
and thus obtain
\[\cF_{\cG^{\Upsilon}}:=\on{LKE}_{\Upsilon\rightarrow\pt}(F_{\cG^{\Upsilon}})\in\cP\lAlg\left(\cC\lMod\left(\Shv(\cX;\bO)\right)\right).\]
We can think of this as the twist of $\cF$ by $\cG^{\Upsilon}$.

\sssec{}
Since $\bO$ is symmetric monoidal, $\bO^{\simeq}$ is an $\E_{\infty}$-object in $\Spc$. Applying loop functors at the monoidal identity, we get $\Omega^2\bO^{\simeq}$ as an $\E_{\infty}$-object in $\Spc$. We have $B^2\ul{\Omega^2\bO^{\simeq}}\in \CAlg(\Shv(\cX;\Spc))$ as the sheafification of the constant presheaf with value $\Omega^2\bO^{\simeq}$.
\medskip

By looping-delooping adjunction we have a canonical symmetric monoidal functor
\begin{equation}\label{eq: Upsilon family of gerbes}
    B^2\Omega^2\bO^{\simeq}\rightarrow\bO.
\end{equation}
Now we are in the situation of Sect. \ref{sssec: twisting modules by linearity}, where $\sfA=\Omega^2\bO^{\simeq}$, $\cB=1_{\bO}$, and $\xi^{\on{lin}}:B^2\sfA\rightarrow\bO$ is the underlying $\E_0$-map of \eqref{eq: Upsilon family of gerbes}, and therefore given a $\Omega^2\bO^{\simeq}$-gerbe, we can twist $\cC$-module sheaves. This construction is organized into a functor
\begin{equation}\label{eq: twisting module action functor}
    \Gamma(\cX;B^2\ul{\Omega^2\bO^{\simeq}})\times \cC\lMod(\Shv(\cX;\bO))\rightarrow\cC\lMod(\Shv(\cX;\bO)).
\end{equation}
In the special case where $\cC$ is the monoidal identity of $\Shv(\cX;\bO)$, we have
\begin{equation}\label{eq: twisting action functor}
\Gamma(\cX;B^2\ul{\Omega^2\bO^{\simeq}})\times \Shv(\cX;\bO)\rightarrow\Shv(\cX;\bO).
\end{equation}

\begin{lemma}
    The functor \eqref{eq: twisting module action functor} is symmetric monoidal.
\end{lemma}
\begin{proof}
    By \cite[Lemma A.2.2]{ZhaoNew} the functor \eqref{eq: twisting action functor} is symmetric monoidal. Inserting the monoidal identity on the second argument, we get a symmetric monoidal functor
    \[\Gamma(\cX;B^2\ul{\Omega^2\bO^{\simeq}})\rightarrow\Shv(\cX;\bO).\]
    The tensor product map
    \[\Shv(\cX;\bO)\times\cC\lMod(\Shv(\cX;\bO))\rightarrow \cC\lMod(\Shv(\cX;\bO))\]
    is symmetric monoidal, and hence we get a symmetric monoidal functor
    \[\Gamma(\cX;B^2\ul{\Omega^2\bO^{\simeq}})\times\cC\lMod(\Shv(\cX;\bO))\rightarrow \cC\lMod(\Shv(\cX;\bO)).\]
    From construction, this functor agrees with \eqref{eq: twisting module action functor} as a plain functor, hence the proof.
\end{proof}

\sssec{Twisting by grading}\label{sssec: twist by grading}
We remain in the context of Sect. \ref{sssec: twist by grading outline}. Now we claim that given a $\cP$-monoidal functor
\[\cG^{\Upsilon}:\Upsilon\rightarrow \Gamma(\cX;B^2\ul{\Omega^2\bO^{\simeq}}),\]
we can carry out the construction outlined in Sect. \ref{sssec: twist by grading outline}, i.e., we can construct a right-lax $\cP$-monoidal functor \eqref{eq: twist of grading data}. Indeed, the desired functor is obtained as the composition
\[\Upsilon\xrightarrow{(\cG^{\Upsilon}, F)}\Gamma(\cX;B^2\ul{\Omega^2\bO^{\simeq}})\times \cC\lMod(\Shv(\cX;\bO))\xrightarrow{\eqref{eq: twisting module action functor}}\cC\lMod(\Shv(\cX;\bO)).\]
Such a composition has a natural right-lax $\cP$-monoidal structure, because $\cG^{\Upsilon}$ and $F$ comes with a right-lax $\cP$-monoidal structure, and \eqref{eq: twisting action functor} is symmetric monoidal.

\sssec{Comparison to the previous method}\label{sssec: comparison of twisting methods}
If the twisting data $\cG^{\Upsilon}$ has an additional structure, the construction above fits into the framework of our earlier construction Sect. \ref{sssec: twisting construction modules}. Suppose we have an abelian group $\sfA$ equipped with a $\E_{\infty}$-map $\sfA\rightarrow \Omega^2\bO^{\simeq}$. Both $\Upsilon$ and $\Gamma(\cX; B^2\ul{\sfA})$ have a natural structure as a $\Z$-linear connective spectra, and hence we can consider the space
\[\Hom_{\Z}(\Upsilon,\Gamma(\cX;B^2\ul{\sfA}))\]
of $\Z$-linear functors between them. We have a canonical functor
\begin{equation}\label{eq: Hom Z to Fun rlax P}
    \Hom_{\Z}(\Upsilon,\Gamma(\cX;B^2\ul{\sfA}))\rightarrow \Fun_{\cP}(\Upsilon, \Gamma(\cX;B^2\ul{\Omega^2\bO^{\simeq}})).
\end{equation}
 Indeed,
this is the composition of
\[\Hom_{\Z}(\Upsilon,\Gamma(\cX;B^2\ul{\sfA}))\rightarrow \Fun_{\CAlg(\Cat)}(\Upsilon, \Gamma(\cX;B^2\ul{\Omega^2\bO^{\simeq}}))\rightarrow \Fun_{\cP}(\Upsilon, \Gamma(\cX;B^2\ul{\Omega^2\bO^{\simeq}})),\]
where the second arrow is the obvious one and the first arrow is the composition with the $\E_{\infty}$-map
\[\Gamma(\cX;B^2\ul{\sfA})\rightarrow\Gamma(\cX;B^2\ul{\Omega^2\bO^{\simeq}}).\]

The construction Sect. \ref{sssec: twist by grading} can be organized into the functor
\[\Fun^{\on{lax}}_{\cP}(\Upsilon, \Gamma(\cX;B^2\ul{\Omega^2\bO^{\simeq}}))\rightarrow \cP\lAlg\left(\cC\lMod\left(\Shv(\cX;\bO)\right)\right),\quad \cG^{\Upsilon}\mapsto \cF_{\cG^{\Upsilon}},\]
and this has a canonical $\E_0$-structure where the pointing on the target is $\cF$. Composing with \eqref{eq: Hom Z to Fun rlax P} we get an $\E_0$-map
\[\Hom_{\Z}(\Upsilon,\Gamma(\cX;B^2\ul{\sfA}))\rightarrow \cP\lAlg\left(\cC\lMod\left(\Shv(\cX;\bO)\right)\right).\]
Using the canonical equivalence
\[\Hom_{\Z}(\Upsilon,\Gamma(\cX; B^2\ul{\sfA}))=\Gamma(\cX;B^2\ul{\Hom_{\on{Grp}}(\Upsilon,\sfA)})\]
we get an $\E_0$-map
\begin{equation}\label{eq: twisting map Z linear grading}
    \Gamma(\cX;B^2\ul{\Hom_{\on{Grp}}(\Upsilon,\sfA)})\rightarrow \cP\lAlg\left(\cC\lMod\left(\Shv(\cX;\bO)\right)\right).
\end{equation}
This provides us an $\E_0$-map
\begin{equation*}
    \xi: B^2\ul{\Hom_{\on{Grp}}(\Upsilon,\sfA)}\rightarrow\ul{\Shv}\left(\cX;\cP\lAlg\left(\cC\lMod(\bO)\right)\right),
\end{equation*}
which is the twisting data \eqref{eq: twisting data algebras in modules} for $\cF$ where $\cA=\ul{\Hom_{\on{Grp}}(\Upsilon,\sfA)}$. By construction, \eqref{eq: twisting map Z linear grading} is the functor of twisting $\cF$ with the data $\xi$ through the procedure of Sect. \ref{sssec: twisting construction modules}.
\medskip

To sum up the discussion, suppose that we are in the situation of Sect. \ref{sssec: graded monoidal structure} and further suppose that we have an abelian group $\sfA$ equipped with an $\E_{\infty}$-map $\sfA\rightarrow\Omega^2\bO^{\simeq}$. Given a $\Z$-linear functor
\begin{equation}\label{eq: Z linear Upsilon family}
    \cG^{\Upsilon}_{\Z}:\Upsilon\rightarrow \Gamma(\cX;B^2\ul{A}),
\end{equation}
this gives rise to the data needed for the two twisting constructions in Sect. \ref{sssec: twist by grading} and Sect. \ref{sssec: twisting construction modules}. On one hand, $\cG^{\Upsilon}_{\Z}$ gives rise to the data $\cG^{\Upsilon}$ needed for the twisting construction of Sect. \ref{sssec: twist by grading}. On the other hand,  $\cG^{\Upsilon}_{\Z}$ gives rise to the data 
\begin{equation}\label{eq: formation of cG}
    \cG\in\Gamma(\cX;B^2\ul{\Hom_{\on{Grp}}(\Upsilon,\sfA)})
\end{equation}
and $\xi$ needed for the twisting construction of Sect. \ref{sssec: twisting construction modules}. Those two constructions result in the same twist of $\cF$, i.e., we have a canonical equivalence
$\cF_{\cG}=\cF_{\cG^{\Upsilon}}$.

\sssec{Twisting by linearity}\label{sssec: twisting graded algebras by linearity}
Suppose that we have an object
\[F_0\in\cP\lAlg\left(\cC\lMod\left(\Shv(\cX;\bO)\right)\right),\]
and that we have $\cB\in\cP\lAlg(\bO)$ equipped with its lift along the functor \eqref{eq: LKE graded monoidal}, i.e., a right-lax $\cP$-monoidal functor $\cB^{\Upsilon}:\Upsilon\rightarrow \bO$. Then we can form a right-lax $\cP$-monoidal functor
\[\ul{\cB^{\Upsilon}}\otimes F_0:\Upsilon\rightarrow \cC\lMod\left(\Shv(\cX;\bO)\right)\]
obtained as the composition
\[\Upsilon=\Upsilon\times\pt\xrightarrow{\cB^{\Upsilon} \times  F_0} \bO\times \cC\lMod\left(\Shv(\cX;\bO)\right)\xrightarrow{\eqref{eq: tensoring with a constant sheaf}}\cC\lMod\left(\Shv(\cX;\bO)\right).\]
This is indeed right-lax $\cP$-monoidal, because the first arrow is right-lax $\cP$-monoidal and the second arrow is symmetric monoidal. Hence, given a data of $\cG^{\Upsilon}$ as in Sect. \ref{sssec: twist by grading} we can twist $\ul{\cB^{\Upsilon}}\otimes F_0$.

\sssec{}\label{sssec: comparison twisting by linearity}
The construction of Sect. \ref{sssec: twisting graded algebras by linearity} agrees with the construction of Sect. \ref{sssec: twisting modules by linearity} in the sense of Sect. \ref{sssec: comparison of twisting methods}. We remain in the context of Sect. \ref{sssec: graded monoidal structure}, and let
\[\cF= \on{LKE}_{\Upsilon\rightarrow\pt}(F)\in\cP\lAlg\left(\cC\lMod\left(\Shv(\cX;\bO)\right)\right),\]
then we have
\begin{equation}\label{eq: tensor product and LKE}
    \ul{\cB}\otimes\cF= \on{LKE}_{\Upsilon\rightarrow\pt}(\ul{\cB^{\Upsilon}}\otimes F_0).
\end{equation}
Suppose that we have an abelian group $\sfA$ equipped with an $\E_{\infty}$-map $\sfA\rightarrow\Omega^2\bO^{\simeq}$, and a $\Z$-linear $\cG^{\Upsilon}_{\Z}$ as in \eqref{eq: Z linear Upsilon family}. Then we get $\cG$ as in \eqref{eq: formation of cG}. On the one hand, we can twist \eqref{eq: tensor product and LKE} via the procedure of Sect. \ref{sssec: twisting graded algebras by linearity}, and on the other hand, we can twist \eqref{eq: tensor product and LKE} via the procedure of Sect. \ref{sssec: twisting modules by linearity}. Those two agree.

\sssec{}\label{sssec: comparison twist of Rep H}
In view of the comparison in Sect. \ref{sssec: comparison twisting by linearity}, the construction of twisted representations
\[\mathbf{Rep}(\mathsf{H})^{\on{loc}}_{\cG}\in \CAlg\left(\mathbf{LS}\lMod\left(\Shv(Y_{\et};\DGcat)\right)\right)\]
can be understood via the framework of Sect. \ref{sssec: twisting graded algebras by linearity}. We take $\Upsilon$ to be the group $\wh{Z}_{\mathsf{H}}$ of characters of $Z_{\mathsf{H}}$. We have a functor
\begin{equation}\label{eq: graded monoidality on Rep H}
    \wh{Z}_H\rightarrow \DGcat,\quad \upsilon\mapsto \Rep(\mathsf{H})_{\upsilon}
\end{equation}
where $\Rep(\mathsf{H})_{\upsilon}$ is the full subcategory of $\Rep(\mathsf{H})$ where the center $Z_{\mathsf{H}}$ acts by the character $\upsilon$. Then \eqref{eq: graded monoidality on Rep H} admits a natural right-lax symmetric monoidal structure. Moreover, since $E^{\times}\simeq \Omega^2\DGcat$ we can take $\sfA=E^{\times}$ equipped with a canonical $\E_{\infty}$-map $\sfA\rightarrow \Omega^2\DGcat$. Since $E$ is algebraically closed, we have 
\[\Hom_{\Z}(\wh{Z}_{\mathsf{H}},E^{\times})= Z_{\mathsf{H}}(E).\]

\subsection{Restriction to trivializing covers}\label{ssec: restriction to triv}
In this subsection, we give a more hands-on description of the twisting procedure, which is closer to the description of \cite[Sect. 1.7]{GL}.

\sssec{}
Suppose we are given $\cX$, $\bO$, $\cA$ and $\xi$ as in Sect. \ref{sssec: conditions on cX bO}, and a gerbe $\cG\in \Gamma(\cX, B^2\cA)$.
Consider the full subcategory $\iota:\cX_{\triv}\hookrightarrow \cX$
spanned by objects $y\in \cX_{\triv}$ such that $\cG_y\in \Gamma(y,B^2\cA)$ is trivializable. $\cX_{\triv}$ inherits a topology from $\cX$. We consider the two sheaf categories
\[\Shv(\cX;\bO) \quad \te{and} \quad \Shv(\cX_{\triv};\bO).\]
There is a restriction functor
\begin{equation}\label{eq: Res triv general}
    \on{Res}_{\triv}:\Shv(\cX;\bO)\rightarrow \Shv(\cX_{\triv};\bO)
\end{equation}
and a functor
\begin{equation}\label{eq: Ext triv general}
    \on{Ext}_{\triv}:\Shv(\cX_{\triv};\bO)\rightarrow \Shv(\cX;\bO)
\end{equation}
defined as a composition
\begin{equation*}
    \Shv(\cX_{\triv};\bO)\rightarrow\PSh(\cX;\bO)\rightarrow \Shv(\cX;\bO),
\end{equation*}
where the first arrow is the right-Kan extension along the opposite of $\iota$, and the second arrow is the sheafification functor. We readily check the following. 

\begin{prop}\label{prop: sheaf on whole vs triv}
    The functors \eqref{eq: Res triv general} and \eqref{eq: Ext triv general} define an equivalence between
    \[\Shv(\cX;\bO) \te{ and } \Shv(\cX_{\triv};\bO).\]
\end{prop}

\sssec{}
Suppose that $\cX$, $\bO$, $\cA$, and $\xi$ are as in Sect. \ref{sssec: setup sheaf of modules}. Then both $\Shv(\cX;\bO)$ and $\Shv(\cX_{\triv};\bO)$ are symmetric monoidal, and the restriction functor \eqref{eq: Res triv general} is symmetric monoidal. Since $\CAlg(\Cat)\rightarrow\Cat$ is conservative, the equivalence of Proposition \ref{prop: sheaf on whole vs triv} is symmetric monoidal. 
\begin{cor}
Suppose $\cC\in\CAlg(\Shv(\cX;\bO))$ and $\cP$ be an operad.
    We have symmetric monoidal equivalences
    \[\cC\lMod\left(\Shv(\cX;\bO)\right)\simeq \cC|_{\cX_{\triv}}\lMod\left(\Shv(\cX_{\triv};\bO)\right)\]
    and
    \[\cP\lAlg\left(\cC\lMod\left(\Shv(\cX;\bO)\right)\right)\simeq \cP\lAlg\left(\cC|_{\cX_{\triv}}\lMod\left(\Shv(\cX_{\triv};\bO)\right)\right),\]
    both forget to the equivalence of Proposition \ref{prop: sheaf on whole vs triv}.
\end{cor}

\sssec{}\label{sssec: sheaf from a gerbe}
Hence, to describe the twist $\cF_{\cG}$ it suffices to describe a sheaf of categories over $\cX_{\triv}$. We define $\on{Tot}(y)\in \PSh(\cX;\Spc)$ by the pullback
\[\begin{tikzcd}
	{\on{Tot}(\cG_y)} & \ast \\
	y & {B^2\cA}
	\arrow[from=1-1, to=1-2]
	\arrow[from=1-1, to=2-1]
	\arrow[from=1-2, to=2-2]
	\arrow["\cG"', from=2-1, to=2-2]
\end{tikzcd}\]
and let $\Gamma(y,\cG)$ be the space of sections to $\on{Tot}(\cG_y)\rightarrow y$. Since $\Omega  B^2\cA = B\cA$, we have a canonical action $\Gamma(y,B\cA)$ on $\Gamma(y,\cG)$. The assignment 
\[\cX^{\op}\rightarrow \Spc,\quad y\mapsto \Gamma(y,\cG)\] is a sheaf on $\cX$, so in this way we can view the gerbe $\cG$ as an object of $\Shv(\cX,\Spc)$.
\medskip

If $y\in \cX_{\triv}$, then $\cG_y$ is trivializable hence this action makes $\Gamma(y,\cG)$ as a torsor over $\Gamma(y,B\cA)$. The data of $\xi$ defines an action of $\Gamma(y,B\cA)$ on $\cF(y)$, and the value of the twist $\cF_{\cG}$ at $y$ is the twist of $\cF(y)$ by the torsor $\Gamma(y,B\cA)$.

\sssec{}
Let us specialize to the case $\bO=\DGcat$ and let $y\in\cX_{\triv}$. From the previous discussion $\cF_{\cG}(y)$ is a twist of $\cF(y)$ by a torsor. If we forget the DG category structure $\cF_{\cG}(y)$ and view it as an $\infty$-category, then we have
\begin{equation}\label{eq: twist as a geometric realization}
    \cF_{\cG}(y)= \colim_{[i]\in \Delta^{\op}} \left(\cF(y)\times \Gamma(y,B\cA)^{\times i}\times \Gamma(y,\cG)\right)
\end{equation}
in $\Cat$. In the simplicial diagram, the transition maps are defined using the monoidal structure of $\Gamma(y,B\cA)$, and the action of $\Gamma(y,B\cA)$ on $\cF(y)$ and $\Gamma(y,\cG)$. If $\cG_y$ is trivial, then $\Gamma(y,\cG)=\Gamma(y,B\cA)$, hence the geometric realization is canonically $\cF(y)$.
The geometric realization on the right-hand side is a priori an object of $\Cat$, but it has a canonical DG category structure, such that for any trivialization of $\cG_y$ the induced equivalence $\cF_{\cG}(y)\simeq \cF(y)$ is an equivalence as DG categories. 
\medskip

The identification \eqref{eq: twist as a geometric realization} is compatible with restriction of sections. As a result, \eqref{eq: twist as a geometric realization} gets upgraded to an equivalence between sheaves on $\cX_{\triv}$; we have
\begin{equation}\label{eq: twist as a geom real sheaf ver}
    \cF_{\cG}= \colim_{[i]\in \Delta^{\op}} (\cF\times (B\cA)^{\times i}|_{\cX_{\triv}}\times \cG|_{\cX_{\triv}})
\end{equation}
in $\Shv(\cX_{\triv};\Cat)$.\footnote{A colimit in the category of sheaves is the sheafification of a colimit as presheaves. If $y\in\cX_{\triv}$, the geometric realization on the right-hand side as presheaves is already a sheaf when restricted to $\cX_{/y}$, hence the section at $y$ does not change under the sheafification. This accounts for why we get \eqref{eq: twist as a geometric realization} from \eqref{eq: twist as a geom real sheaf ver}.} Here we viewed $\cG$ as an object of $\Shv(\cX;\Spc)$ following the construction of Sect. \ref{sssec: sheaf from a gerbe}. By an argument similar to the previous paragraph, the geometric realization on the right-hand side has a natural enhancement to $\Shv(\cX_{\triv};\DGcat)$, and \eqref{eq: twist as a geom real sheaf ver} is an equivalence in $\Shv(\cX_{\triv},\DGcat)$.

\section{Twisted sheaves}\label{sec: twisted sheaves}
Let $A$ be a finite abelian group and let $\cY\in \PreStklft$. Let $\mc{G}$ be an \'etale $A$
-gerbe over $\cY$. Fix a group homomorphism $\chi: A\rightarrow E^{\times}$. In Sect. \ref{sssec: twisted sheaves} we defined the category $\Shv^{\cG}(\cY)$ of twisted sheaves over $\cY$, which is the key player in this paper. In this section, we study this notion in depth. In Sect. \ref{ssec: twiste sheaves monodromy} we give another interpretation of twisted sheaves as sheaves over the total space $\on{Tot}(\cG)$ of the gerbe with a fixed monodromy along the fibers of $\on{Tot}(\cG)\rightarrow \cY$, see Theorem \ref{thm: twisted sheaves tensor product formula}. Based on this result, we establish the six-functor formalism for twisted sheaves in Sect. \ref{ssec: functoriality of twisted sheaves}.

\sssec{Sheaf theory}\label{sssec: six functor extension} Throughout this section, except for Sect. \ref{ssec: t str twisting}, we work with Betti sheaves $\Shv$. However, all the results contained here work equally well with $\Ind$-constructible sheaves $\Shv^{\indconstr}$.

\sssec{Six-functor formalism}\label{sssec: six functor stacky}
Recall the notation $\Corr(C,E)$ from \cite{Sch} and the hypothesis on $C$ and $E$ in Sect. 2 of \emph{op. cit.}. $E$ is a collection of morphisms in the category $C$ along which we can pushforward. We often write a map in $\Corr(C,E)$ as a diagram
\[\begin{tikzcd}
	{\cC_{12}} & {\cC_1} \\
	{\cC_2}
	\arrow[from=1-1, to=1-2]
	\arrow[from=1-1, to=2-1]
\end{tikzcd}\]
where all the objects lie in $C$ and the vertical arrow belongs to $E$. We have a six-functor formalism \footnote{In the sense of \cite{Sch}}
\[\Shv: \Corr(\Affft,\on{all})\rightarrow\DGcat \]
which sends the correspondence $S_1\xleftarrow{f} S\xrightarrow{g} S_2$ to the functor $g_!\circ f^*:\Shv(S_1)\rightarrow\Shv(S_2)$. Using \cite[Proposition A.5.16]{Man} we have a six-functor formalism
\[\Shv:\Corr(\NearStklft,\on{affsch})\rightarrow \DGcat,\]
where $\on{affsch}$ is the collections of morphisms in $\NearStklft$ representable by affine schemes. Applying \cite[Theorem 5.19]{Sch} to this situation, there is a collection $\wt{E}$ of morphisms in $\NearStklft$ that contain $\on{affsch}$ such that is stable under disjoint unions, local on the target,
local on the source, and tame. From the closure property of $\wt{E}$, it contains schematic maps and also maps of the form $\Bet \sfA\rightarrow \pt$ where $\sfA$ is a finite abelian group. Using the fact that $\Shv(\mhyphen)$ is insensitive to 
\'etale sheafification, we have a six-functor formalism
\[\Shv:\Corr(\PreStklft, \wt{E}')\rightarrow \DGcat,\]
where $\wt{E}'$ is the collection of morphisms in $\PreStklft$ that map to morphisms in $\wt{E}$ under the \'etale sheafification functor $L_{\ft}:\PreStklft\rightarrow \NearStklft$. This is the functoriality of $\Shv$ we will use throughout this section. The collection $\wt{E}'$ still contains schematic maps and maps of the form $\Bet \sfA\rightarrow \pt$ as above. By the hypothesis on $\wt{E}'$ made in Sect. 2 of \cite{Sch}, it also contains projection maps $\cY\times\Bet \sfA\rightarrow\cY$ for $\cY\in\PreStklft$. Suppose the situation where we have an \'etale $\sfA$-gerbe $\cG$ on $\cY$. Then the total space $\on{Tot}(\cG)$ has a natural action of $\Bet \sfA$. The map $a:\Bet\sfA\times\on{Tot}(\cG)\rightarrow \on{Tot}(\cG)$ from the action is equivalent to the projection $\pr: \Bet\sfA\times \on{Tot}(\cG)\rightarrow \on{Tot}(\cG)$ by an automorphism of $\Bet\sfA\times \on{Tot}(\cG)$, and therefore lies in $\wt{E}'$. The action of $\Bet \sfA$ on $\on{Tot}(\cG)$ will play a key role in Sect. \ref{ssec: twiste sheaves monodromy}, and the fact that we can $!$-pushforward sheaves along such $a$ and $\pr$ will be used throughout.
\begin{rmk}
    Restricting the class $\wt{E}'$ to schematic maps we get the six-functor formalism
\[\Shv:\Corr(\PreStklft,\on{sch})\rightarrow\DGcat.\]
\end{rmk}

\subsection{The category $\Shv(\BA)$ and its action}
Consider the category $\Rep(A)$ of $E$-linear representations of the group $A$.
\begin{lemma}
    There is a canonical equivalence 
    \begin{equation}\label{eq: Shv BA and Rep}
        (\Shv(\Bet A),\otimes)\cong (\Rep(A),\otimes_E)\quad\te{in } \DGCatSymMon.
    \end{equation}
\end{lemma}

\begin{lemma}
    For $Y\in\PreStklft$ we have canonical equivalences
\[\Shv(A)\otimes\Shv(Y)\cong \Shv(A\times Y)\quad\te{and}\quad \Shv(\BA)\otimes\Shv(Y)\cong \Shv(\BA\times Y)\]
     in $\DGCatSymMon$, both given by the box tensor product.
\end{lemma}

\sssec{Convolution}
Since $A$ is abelian we have 
\begin{equation}\label{eq: multiplication on BA}
    m:\Bet A\times \Bet A\rightarrow \Bet A
\end{equation}
induced from the multiplication of $A$. We define the \emph{convolution} product on $\Shv(\BA)$, using the $!$-pushforward along \eqref{eq: multiplication on BA}. This turns $\Shv(BA)$ into a symmetric monoidal category $(\Shv(BA),\ast)$.
\medskip

We define \emph{convolution} monoidal structure $\ast$ on $\Rep(A)$ using a sequence of equivalences
\begin{equation}\label{eq: equivalences Rep A}
    \Rep(A)\cong EA\lMod\cong \QC(\Spec(EA)),
\end{equation}
where $EA$ means the group algebra.
The rightmost term $\QC(\Spec(EA))$ has a natural symmetric monoidal structure $\otimes_{EA}$, so we transfer this structure to $\Rep(A)$ and to get $(\Rep(A),\ast)$ in $\DGCatSymMon$.

\begin{lemma}
    We have a canonical equivalence
    \begin{equation}\label{eq: Shv BA and Rep convolution}
        (\Shv(\Bet A),\ast)\cong (\Rep(A),\ast) \cong (\QC(\Spec(EA)),\otimes_{EA})
    \end{equation}
    in $\DGCatSymMon$, and the first equivalence has the same underlying plain functor with \eqref{eq: Shv BA and Rep}.
\end{lemma}

\sssec{}\label{sssec: L reg}
If $V\in \Rep(A)$, then there is a decomposition
\begin{equation}\label{eq: Rep A decomposition}
    V=\oplus_{\chi'\in A^{\vee}} V_{\chi'}
\end{equation}
where $A^{\vee}=\Hom_{\on{Grp}}(A,E^{\times})$.
The $!$-pushforward functor $\Shv(\Bet A)\rightarrow\Vect$ along $\Bet A\rightarrow\pt$ is identified with $V\mapsto V_0$ under the decomposition \eqref{eq: Rep A decomposition}. 
\medskip

Hence, given $V,W\in\Rep(A)$, the convolution product $V\ast W\in \Rep(A)$ is defined in terms of the decomposition \eqref{eq: Rep A decomposition} as
$(V\ast W)_{\chi'}= V_{\chi'}\otimes W_{\chi'}$. The monoidal unit of $(\Shv(\Bet A),\ast)$ is
\[\cL_{\reg}:=(\pt\rightarrow\BA)_! \ul{E}_{\pt} \in \Shv(\BA),\]
and this corresponds to the regular representation of $A$ under the correspondence \eqref{eq: Shv BA and Rep convolution}.

\sssec{}\label{sssec: Vect chi}
A character $\chi'$ defines a point $\chi':\pt\rightarrow \Spec(EA)$ and hence we have a corresponding $(\QC(\Spec(EA)),\otimes_{EA})$-action on $\Vect$.
Under the equivalence \eqref{eq: Shv BA and Rep convolution} we have a corresponding $(\Shv(\Bet A),\ast)$-action on $\Vect$. We denote this module by $\Vect_{\chi'}$ and its action map by
$\on{act}_{\chi'}: \Shv(\Bet A)\otimes \Vect_{\chi'}\rightarrow \Vect_{\chi'}$. Under the identification $\Vect_{\chi'}\cong\Vect$ and using the decomposition \eqref{eq: Rep A decomposition}, this functor is $V\mapsto V_{\chi'}$. 
\medskip

We consider the skyscraper sheaf at $\chi': \pt\rightarrow \Spec(EA)$, which we denote $\delta_{\chi'}\in \QC(\Spec(EA))$. We have a corresponding object $\cL_{\chi'} \in \Shv(\Bet A)$. Then $\cL_{\chi'}$ is a rank-one local system with monodromy representation $\chi': A\rightarrow E^{\times}$. $\cL_{\chi'}$ has a canonical multiplicative structure; for instance we have a canonical isomorphism $m^!\cL_{\chi'}\cong\cL_{\chi'}\boxtimes\cL_{\chi'}$.  Under \eqref{eq: Rep A decomposition} $\cL_{\chi'}$ corresponds to $V\in\Rep(A)^{\heartsuit}$ with $V_{\chi'}=E$ and all other components are zero. The case $\chi'=-\chi$ will be particularly important for us.

\subsection{Twisted sheaves as sheaves with monodromy}\label{ssec: twiste sheaves monodromy}
From the data of $\cY$ and $\cG$ we have the total space $\on{Tot}(\cG)\in \PreStklft$ that fits into the pullback diagram
\[\begin{tikzcd}
	{\on{Tot}(\cG)} & \pt \\
	\cY & {\Bet^2 A}
	\arrow[from=1-1, to=1-2]
	\arrow[from=1-1, to=2-1]
	\arrow[from=1-2, to=2-2]
	\arrow["\cG"', from=2-1, to=2-2]
\end{tikzcd}\]
The canonical map $\on{Tot}(\cG)\rightarrow \cY$ is a $\Bet A$-bundle. The slogan of this section is:
\begin{displayquote}
$\Shv^{\cG}(\cY)$ is equivalent to the category of sheaves over the total space $\on{Tot}(\cG)$ with monodromy $-\chi$ along the fiber.
\end{displayquote}
The latter category consists of $\cF\in \Shv(\on{Tot}(\cG))$ together with twisted equivariant structure against the action of $\Bet A$, where the twist is given by the multiplicative local system $\cL_{-\chi}\in \Shv(\Bet A)$. While its precise formulation will be given in Sect. \ref{sssec: chi monodromic sheaves}, a piece of the data is an isomorphism
\begin{equation}\label{eq: twisted equivariance}
    a^*\cF\cong \cL_{-\chi}\boxtimes \cF
\end{equation}
in $\Shv(\Bet A\times \on{Tot}(\cG))$ where $a$ is the action map $\Bet A\times\on{Tot}(\cG)\rightarrow\on{Tot}(\cG)$.

\sssec{$(-\chi)$-monodromic sheaves}\label{sssec: chi monodromic sheaves}
The action of $\Bet A$ on $\on{Tot}(\cG)$ induces an action of $(\Shv(\Bet A),\ast)$ on $\Shv(\on{Tot}(\cG))$. For instance, the action map
\[\Shv(\Bet A)\otimes \Shv(\on{Tot}(\cG))\rightarrow \Shv(\on{Tot}(\cG))\]
is the composition
\[\Shv(\Bet A)\otimes \Shv(\on{Tot}(\cG))\xrightarrow{\simeq} \Shv(\Bet A\times \on{Tot}(\cG))\xrightarrow{a_!} \Shv(\on{Tot}(\cG))\]
where $a: \Bet A\times \on{Tot}(\cG)\rightarrow \on{Tot}(\cG)$ is the action map. We consider the DG category
\begin{equation}\label{eq: functor category definition}
    \Fun_{(\Shv(\Bet A),\ast)\lMod}(\Vect_{-\chi},\Shv(\on{Tot}(\cG)))
\end{equation}
of functors between module categories. ($\Vect_{-\chi}$ was defined in Sect. \ref{sssec: Vect chi}) We call \eqref{eq: functor category definition} the category of $\chi$-\emph{monodromic sheaves}.
\medskip

The category \eqref{eq: functor category definition} indeed attempts to capture the idea of $\cL_{-\chi}$-twisted equivariant sheaves on $\on{Tot}(\cG)$.
If we have a sheaf $\cF$ on the total space $\on{Tot}(\cG)$ there is a (plain) functor $\Vect_{-\chi}\cong \Vect\rightarrow \Shv(\on{Tot}(\cG))$ defined by this object $\cF$, and we check that the data \eqref{eq: twisted equivariance} yields a commutativity data of the square
\[\begin{tikzcd}
	{\Shv(\Bet A)\times\Vect_{-\chi}} & {\Shv(\Bet A)\times\Shv(\on{Tot}(\cG))} \\
	{\Vect_{-\chi}} & {\Shv(\on{Tot}(\cG))}
	\arrow["\cF", from=1-1, to=1-2]
	\arrow[from=1-1, to=2-1]
	\arrow[from=1-2, to=2-2]
	\arrow["\cF", from=2-1, to=2-2]
\end{tikzcd}\]
where the vertical arrows are action maps.
\medskip

The map $-\chi:\pt\rightarrow \Spec(EA)$ is perfect, so
$\Vect_{-\chi}$ is self-dual as a $(\Shv(\Bet A),\ast)$-module category (\cite[Corollary 4.8]{BZFN}). Hence, \eqref{eq: functor category definition} is equivalent to
\begin{equation}\label{eq: tensor product definition}
    \Shv(\on{Tot}(\cG))\otimes_{(\Shv(\Bet A),\ast)}\Vect_{-\chi}.
\end{equation}
We later prove in Subsection \ref{ssec: proof of tensor product formula}:
\begin{thm}\label{thm: twisted sheaves tensor product formula}
   There is a canonical equivalence of DG categories
   \[\Shv(\on{Tot}(\cG))\otimes_{(\Shv(\Bet A),\ast)}\Vect_{-\chi} \cong \Shv^{\cG}(\cY).\]
\end{thm}

\sssec{Twisted descent}\label{sssec: twisted descent} 
Another approach to $(-\chi)$-monodromic sheaves can be made through descent.
Consider the simplicial diagram
\begin{equation}\label{eq: simplicial stack}
    \Delta^{\op}\rightarrow \PreStklft,\quad [n]\mapsto (\Bet A)^{\times n}\times \on{Tot}(\cG)
\end{equation}
arising from the $\Bet A$-action on $\on{Tot}(\cG)$, and for each map $\alpha:[i]\rightarrow [j]$ in $\Delta$ we denote the transition map
\begin{equation}\label{eq: T alpha}
    T_{\alpha}: (\Bet A)^{\times j}\times \on{Tot}(\cG)\rightarrow (\Bet A)^{\times i}\times \on{Tot}(\cG).
\end{equation}
The colimit of \eqref{eq: simplicial stack}
is canonically equivalent to $\cY$ up to sheafification. Hence we get
\begin{equation}\label{eq: sheaves on Y simplicial}
    \Shv(\cY)=\lim_{[n]\in\Delta} \Shv((\Bet A)^{\times n}\times \on{Tot}(\cG))
\end{equation}
in $\DGcat$ where transition maps are given by $*$-pullbacks. An object of the limit on the right-hand side is the collection $\{\cF_n\}_{[n]\in\Delta}$ where $\cF_n\in \Shv((\Bet A)^{\times n}\times \on{Tot}(\cG))$ together with isomorphisms
$\cF_j\xrightarrow{\sim} T^*_{\alpha}\cF_i$ for each $\alpha:[i]\rightarrow [j]$ in $\Delta$ and higher coherence data. 
\medskip

The limit \eqref{eq: sheaves on Y simplicial} is how we define equivariant sheaves (see \cite[Appendix B]{BL}). For example, for
\[\alpha_0,\alpha_1:[0]\rightarrow [1],\quad \alpha_0(0)=1 \te{ and }\alpha_1(0)=0,\]
we have that $T_{\alpha_0}$ is the projection and $T_{\alpha_1}$ is the action map, hence we have isomorphisms $\cF_1\xrightarrow{\sim}\ul{E}_{\Bet A}\boxtimes\cF_0$ and $\cF_1\xrightarrow{\sim} a^*\cF_0$. Composing these, we get $a^*\cF_0\cong \ul{E}_{\Bet A}\boxtimes \cF_0$, which is a part of equivariance data of $\cF_0$. It is natural to expect that we can form a similar cosimplicial diagram
\begin{equation}\label{eq: TD chi}
    \on{TD}_{-\chi}:\Delta\rightarrow \DGcat,\quad [n]\mapsto \Shv((\Bet A)^{\times n}\times\on{Tot}(\cG))
\end{equation}
(with different transition maps!)
whose totalization would be equivalent to the category of twisted equivariant sheaves. In fact, it is not hard to guess what the appropriate cosimplicial diagram should be. For example, this limit should encode an isomorphism \eqref{eq: twisted equivariance}, so we can define the transition map for $\alpha_1$ as before but define the new transition map for $\alpha_0$ as
\[\Shv(\on{Tot}(\cG))\rightarrow \Shv(\Bet A\times \on{Tot}(\cG))\cong \Shv(\Bet A)\otimes \Shv(\on{Tot}(\cG)),\quad \cF\mapsto \cL_{-\chi}\boxtimes\cF.\]
Then the corresponding limit would encode isomorphisms $\cF_1\xrightarrow{\sim}\cL_{-\chi}\boxtimes \cF_0$ and  $\cF_1\xrightarrow{\sim}a^*\cF_0$, hence we get a data \eqref{eq: twisted equivariance} for twisted equivariance of $\cF_0$.
\medskip

We claim that the expression \eqref{eq: tensor product definition} captures this idea of twisted descent. We will construct an equivalence between \eqref{eq: tensor product definition} and a certain totalization that encodes twisted equivariance.

\begin{nota}
    In the category $\Delta$ we let $\pa_i:[n-1]\rightarrow [n]$ the unique map whose image is $[n]\setminus \{i\}$, and $\sigma_i:[n+1]\rightarrow [n]$ the unique map which is surjective and has the same value on $i$ and $i+1$. Hence for a simplicial set $X$ we have a face map $\pa_i:X_n\rightarrow X_{n-1}$ and a degeneracy map $\sigma_i:X_n\rightarrow X_{n+1}$. We denote those maps as $\pa^{[n]}_i$ and $\sigma^{[n]}_i$ when we want to make clear which $n$ we are considering.
\end{nota}

\sssec{}
By definition, the tensor product \eqref{eq: tensor product definition} is the geometric realization
\[\colim_{[n]\in\Delta^{\op}}\Vect_{-\chi}\otimes \Shv(\Bet A)^{\otimes n}\otimes \Shv(\on{Tot}(\cG))\]
in $\DGcat$. The transition maps are as follows:
\[\pa_0: \Vect_{-\chi}\otimes \Shv((\Bet A))^{\otimes n}\otimes\Shv(\on{Tot}(\cG)))\rightarrow \Vect_{-\chi}\otimes \Shv((\Bet A))^{\otimes n-1}\otimes \Shv(\on{Tot}(\cG))\]
is the tensor product $\on{act}_{-\chi}\otimes \Shv(\Bet A)^{\otimes n-1}\otimes\Shv(\on{Tot}(\cG))$, and 
\[\pa^{[n]}_i: \Vect_{-\chi}\otimes \Shv((\Bet A))^{\times n}\times\on{Tot}(\cG))\rightarrow \Vect_{-\chi}\otimes \Shv((\Bet A)^{\otimes n-1}\times \on{Tot}(\cG)) \]
for $i> 0$ is defined to be the tensor product $\Vect_{-\chi}\otimes (T_{\pa^{[n]}_i})_{!}$. For any $i$,
\[\sigma_i:  \Vect_{-\chi}\otimes \Shv((\Bet A))^{\times n}\times\on{Tot}(\cG))\rightarrow \Vect_{-\chi}\otimes \Shv((\Bet A)^{\otimes n+1}\times \on{Tot}(\cG))\]
is defined as the tensor product $\Vect_{-\chi}\otimes (T_{\sigma^{[n]}_i})_{!}$. (see Sect. \ref{sssec: Vect chi} and \eqref{eq: T alpha} for notations)
\medskip

Note that each of the transition maps has a right adjoint. This follows from \textbf{(1)} each $T_{\alpha}$ is \'etale, hence $(T_{\alpha})^*$ is right adjoint to $(T_{\alpha})_!$; \textbf{(2)} $\on{act}_{-\chi}$ has a right adjoint $\Vect_{-\chi}\xrightarrow{\id\otimes \cL_{-\chi}}\Vect_{-\chi}\otimes \Shv(\Bet A)$. Hence, we can pass to right adjoints (\cite[Proposition 2.5.7]{GR}) and conclude that the tensor product \eqref{eq: tensor product definition} is equivalent to
\[\lim_{[n]\in\Delta} \Vect_{-\chi}\otimes \Shv(\Bet A)^{\otimes n} \otimes \Shv(\on{Tot}(\cG))\]
in $\DGcat$. Using $\Vect_{-\chi}\cong \Vect$ this gets identified with
\begin{equation}\label{eq: monodromic sheaves simplicial}
    \lim_{[n]\in\Delta} \Shv(\Bet A)^{\otimes n} \otimes \Shv(\on{Tot}(\cG)) =\lim_{[n]\in\Delta}\Shv((\Bet A)^{\times n}\times \on{Tot}(\cG)).
\end{equation}

In Sect. \ref{sssec: twisted descent} we wanted to construct a cosimplicial diagram $\on{TD}_{-\chi}$ (see \eqref{eq: TD chi}) that yields the category of $(-\chi)$-monodromic sheaves. The right-hand side of \eqref{eq: monodromic sheaves simplicial} - which we proved to be equivalent to the category of $(-\chi)$-monodromic sheaves - is the desired $\on{TD}_{-\chi}$. 
We have made expectations for transition maps of $\alpha_0=\sigma_0$ and $\alpha_1=\sigma_1$ in $\on{TD}_{-\chi}$. (see Sect. \ref{sssec: twisted descent}) We claim that this agrees with transition maps of \eqref{eq: monodromic sheaves simplicial}. From our construction, the transition map of $\alpha_0$ in \eqref{eq: monodromic sheaves simplicial} is 
\[\Vect_{-\chi}\otimes\Shv(\on{Tot}(\cG))\xrightarrow{\id\otimes \cL_{-\chi}\otimes\id} \Vect_{-\chi}\otimes\Shv(\Bet A)\otimes \Shv(\on{Tot}(\cG)),\]
which gets identified with
\[\Shv(\on{Tot}(\cG))\rightarrow \Shv(\Bet A)\otimes \Shv(\on{Tot}(\cG)),\quad \cF\mapsto \cL_{-\chi}\boxtimes \cF\]
via $\Vect_{-\chi}\cong\Vect$, hence it agrees with the expectation we made. A similar argument shows the agreement for $\alpha_1$.

\subsection{Proof of Theorem \ref{thm: twisted sheaves tensor product formula}}\label{ssec: proof of tensor product formula}
In this subsection, we prove Theorem \ref{thm: twisted sheaves tensor product formula}, which asserts a canonical equivalence
\begin{equation}\label{eq: twisted and tensor product}
    \Shv(\on{Tot}(\cG))\otimes_{(\Shv(\Bet A),\ast)}\Vect_{-\chi} \cong \Shv^{\cG}(\cY).
\end{equation}
Our approach is to prove a sheafified version of this assertion.

\sssec{}
We observe that construction \eqref{eq: functor category definition} is upgraded to a sheaf.
The functor 
\begin{equation}\label{eq: functor cat defn sheaf}
\left((\Affft)_{/\cY}\right)^{\op}\rightarrow \DGcat,\quad S\mapsto \Fun_{(\Shv(\Bet A),\ast)\lMod}(\Vect_{-\chi},\Shv(\on{Tot}(\cG|_{S})))
\end{equation}
is a sheaf. Indeed, this is a combination of \textbf{(1)} the assignment $S\mapsto \Shv(\on{Tot}(\cG|_S))$ is a sheaf; \textbf{(2)} the forgetful functor $\Shv(\Bet A,\ast)\lMod\rightarrow \DGcat$ creates limits; and \textbf{(3)} $\Fun_{(\Shv(\Bet A),\ast)}(\Vect_{-\chi},\mhyphen)$ preserves limits. Hence the functor
\begin{equation}\label{eq: tensor product defn sheaf}
\left((\Affft)_{/\cY}\right)^{\op}\rightarrow \DGcat,\quad S\mapsto \Vect_{-\chi}\otimes_{(\Shv(\Bet A),\ast)}\Shv(\on{Tot}(\cG|_{S}))
\end{equation}
is a sheaf. Moreover, by a similar argument the global section of \eqref{eq: functor cat defn sheaf} is
\[\Fun_{(\Shv(\Bet A),\ast)\lMod}(\Vect_{-\chi},\Shv(\on{Tot}(\cG)),\]
hence the global section of \eqref{eq: tensor product defn sheaf} is the left-hand side of \eqref{eq: twisted and tensor product}.
\medskip

On the other hand, the functor
\begin{equation}\label{eq: twisted sheaf sheaf}
    \left((\Affft)_{/\cY}\right)^{\op}\rightarrow \DGcat,\quad S\mapsto \Shv^{\cG}(S)
\end{equation}
is a sheaf, and its global section is $\Shv^{\cG}(\cY)$ from construction. Hence, once we construct a canonical equivalence between sheaves \eqref{eq: tensor product defn sheaf} and \eqref{eq: twisted sheaf sheaf}, we prove Theorem \ref{thm: twisted sheaves tensor product formula}.

\sssec{}
By Sect. \ref{ssec: restriction to triv} it suffices to construct a canonical equivalence between their restriction to $(\Affft)_{/\cY}^{\triv}$, the full subcategory of $(\Affft)_{/\cY}$ spanned by objects $S\rightarrow \cY$ along which the pullback of the gerbe $\cG$ is trivializable. From \eqref{eq: twist as a geometric realization} the restriction of \eqref{eq: twisted sheaf sheaf} is
\begin{equation}\label{eq: twisted sheaf sheaf on triv}
    \left((\Affft)^{\triv}_{/\cY}\right)^{\op}\rightarrow \DGcat,\quad S\mapsto \colim_{[n]\in\Delta^{\op}}(\Shv(S)\times \Bet A(S)^{\times n}\times \Gamma(S,\cG)).
\end{equation}
On the other hand, using the bar construction of relative tensor product, the restriction of \eqref{eq: tensor product defn sheaf} is
\begin{align}\label{eq: tensor product defn sheaf triv}
    \left((\Affft)^{\triv}_{/\cY}\right)^{\op}\rightarrow \DGcat,& \\
    S\mapsto & \colim_{[n]\in\Delta^{\op}}\left(\Vect_{-\chi}\otimes \Shv(\Bet A)^{\otimes n}\otimes \Shv(\on{Tot}(\cG_S))\right).\nonumber
\end{align}
For the rest of this subsection, we prove:

\begin{prop}\label{prop: equiv btw sheaves triv}
    There is a canonical equivalence between sheaves from \eqref{eq: twisted sheaf sheaf on triv} to \eqref{eq: tensor product defn sheaf triv}.
\end{prop}

\sssec{}\label{sssec: GnS TnS}
We have to construct a compatible system of canonical maps
\begin{equation}\label{eq: map between bar constructions}
    \colim_{[n]\in\Delta^{\op}}(\Shv(S)\times \Bet A(S)^{\times n}\times \cG(S))\rightarrow \colim_{[n]\in\Delta^{\op}} \left(\Vect_{-\chi}\otimes \Shv(\Bet A)^{\otimes n}\otimes \Shv(\on{Tot}(\cG_S))\right)
\end{equation}
for $S\in (\Affft)^{\triv}_{/\cY}$. For notational convenience, we denote the left-hand side by $G_S$ and the right-hand side by $T_S$. This will be accomplished by constructing a compatible system of canonical maps
\begin{equation}\label{eq: trivialization box tensor map}
    \Shv(S)\times \Bet A(S)^{\times n}\times \cG(S)\rightarrow \Vect_{-\chi}\otimes \Shv(\Bet A)^{\otimes n}\otimes \Shv(\on{Tot}(\cG_S)).
\end{equation}
For notational convenience, we denote the left-hand side by $G^n_S$ and the right-hand side by $T^n_S$.

\begin{rmk}
    Each $\cG_S$ is trivializable, and if $\cG_S$ is trivial, we can concretely identify both the $G_S$ and $T_S$ to $\Shv(S)$. For this reason, the proof of Proposition \ref{prop: equiv btw sheaves triv} is more or less immediate after constructing a canonical map between those sheaves of categories.
\end{rmk}

\begin{rmk}
    Note that the colimit on the left-hand side is taken in $\Cat$ while the colimit on the right-hand side is taken in $\DGcat$. However, the colimit on the left-hand side is canonically an object of $\DGcat$, and the map \eqref{eq: map between bar constructions} should be a map between DG categories. We construct a system of maps \eqref{eq: trivialization box tensor map} in $\Cat$ so that the resulting map \eqref{eq: map between bar constructions} is a priori in $\Cat$, but we argue that this has a natural lift to $\DGcat$. (see Sect. \ref{sssec: showing equiv btw bar constructions})
\end{rmk}

\sssec{Construction of \eqref{eq: trivialization box tensor map}} 
We explain the construction of canonical functor $G^n_S\rightarrow T^n_S$. 
\medskip

We construct an $A^{n+1}$-gerbe $\cG^{(n)}_S$ over $S$ as the composition
\[S\xrightarrow{\cG_S} \Bet^2 A\rightarrow \Bet^2 A^{n+1},\]
where the second map is induced from the inclusion $A\hookrightarrow A^{n+1}$ to the first factor. Then we have $\on{Tot}(\cG^{(n)}_S)= \on{Tot}(\cG_S)\times (\Bet A)^{\times n}$, hence
\[G^n_S=\Shv(S)\times \Gamma(S,\cG^n_S) \quad \te{and}\quad T^n_S= \Shv(\on{Tot}(\cG^n_S))\]
using $\Vect_{-\chi}\cong\Vect$ for $T^n_S$.
We observe that $\Gamma(S,\cG^{(n)}_S)$ is the category of trivializations (say $\beta$)  of the $A^{n+1}$-gerbe $\cG^{(n)}_S$, hence $G^n_S$ is the category of pairs $(\cF,\beta)$ where $\cF\in\Shv(S)$ and $\beta$ is a trivialization of $\cG^{(n)}_S$. In more detail, we have an equivalence
\[\Bet(A)^{\times n}\times\cG(S)\rightarrow \{\te{trivializations of } \cG^{(n)}_S\}\]
which maps $(a_1,\cdots,a_n,p)$ to the trivialization
\[(\Bet A)^{\times n+1}\times S\rightarrow (\Bet A)^{\times n}\times\on{Tot}(\cG_S),\quad (a'_1,\cdots,a'_{n+1},s)\mapsto (a_1 a'_1,\cdots, a_n a'_n, a'_{n+1} p(s)),\]
where we view $p\in\cG(S)$ as a map $S\rightarrow \on{Tot}(\cG)$.
\medskip

Recall from Sect. \ref{sssec: L reg} we had the monoidal unit $\cL_{\on{reg}}\in (\Shv(\Bet A),\ast)$ and $\cL_{-\chi}\in \Shv(\Bet A)$, which corresponds to regular representation of $A$ and the character $-\chi$ of $A$, respectively. A trivialization $\beta$ of $\cG^{(n)}_S$ induces a map $\beta^{\on{triv}}: (\Bet A)^{\times n+1}\times S\rightarrow \cG^{(n)}_S$. We construct a functor $\Phi^n_S:G^n_S\rightarrow T^n_S$ as
\begin{equation*}
    \Phi^n_S(\cF,\beta)=\beta^{\on{triv}}_!(\cL_{\on{reg}}\boxtimes\cdots\boxtimes\cL_{\on{reg}}\boxtimes \cF)
\end{equation*}
where $\cL_{\on{reg}}$ appears $n+1$-times.

\sssec{}
To show that the collection $\{\Phi^n_S\}$ defines a map \eqref{eq: map between bar constructions} we have to show that the maps $\Phi^n_S$ intertwine transition maps in the formation of colimits. We have to check the compatibility for $\pa_i$ and $\sigma_i$'s. The cases of $\pa_i$'s for $i\ge 1$ and $\sigma_i$'s are treated uniformly, but we have to make a separate argument for $\pa_0$.

\sssec{}
Let $A_1$ and $A_2$ be finite abelian groups with a group homomorphism $c:A_1\rightarrow A_2$. For $\cP\in\on{Ge}_{A_1}(S)$ we have $\cP'\in \on{Ge}_{A_2}(S)$ by the composition $S\xrightarrow{\cP} \Bet^2 A_1\xrightarrow{c} \Bet^2 A_2$. For $\cK\in \Shv(\Bet A_1)$ we construct a functor
\[\Phi_{\cK}:\Shv(S)\times\cP(S)\rightarrow \Shv(\on{Tot}(\cP_S)),\quad (\cF,\beta)\mapsto \beta^{\on{triv}}_{!}(\cK\boxtimes\cF).\]
Using $c_!\cK \in \Shv(\Bet A_2)$ we can similarly construct
\[\Phi_{c_!\cK}:\Shv(S)\times\cP'(S)\rightarrow \Shv(\on{Tot}(\cP'_S)),\quad (\cF,\beta)\mapsto \beta^{\on{triv}}_{!}(c_!\cK\boxtimes\cF).\]
Note that $c$ induces canonical maps $c:\cP(S)\rightarrow \cP'(S)$ and $\on{Tot}(\cP)\rightarrow \on{Tot}(\cP')$, which are denoted by $c$. We readily check the following:

\begin{lemma}\label{lem: Phi map and change of group}
    In the above situation, the following diagram commutes:
    \[\begin{tikzcd}
	{\Shv(S)\times\cP(S)} & {\Shv(\on{Tot}(\cP))} \\
	{\Shv(S)\times \cP'(S)} & {\Shv(\on{Tot}(\cP'))}
	\arrow["{\Phi_{\cK}}", from=1-1, to=1-2]
	\arrow["{1\times c}"', from=1-1, to=2-1]
	\arrow["{c_!}", from=1-2, to=2-2]
	\arrow["{\Phi_{c_*\cK}}"', from=2-1, to=2-2]
\end{tikzcd}\]
\end{lemma}

\sssec{Compatibility except for $\pa_0$} Say that $\alpha$ is a map in $\Delta$ of the form either $\pa_i:[n-1]\rightarrow [n]$, $i>0$ or $\sigma_i:[n+1]\rightarrow [n]$. We need commutativity of either
\[\begin{tikzcd}
	{\Shv(S)\times\cG^{(n)}(S)} & {\Shv(\on{Tot}(\cG^{(n)}_S))} \\
	{\Shv(S)\times\cG^{(n-1)}(S)} & {\Shv(\on{Tot}(\cG^{(n-1)}_S))}
	\arrow["{\Phi^n_S}", from=1-1, to=1-2]
	\arrow["{\pa_i}"', from=1-1, to=2-1]
	\arrow["{\pa_i}", from=1-2, to=2-2]
	\arrow["{\Phi^{n-1}_S}"', from=2-1, to=2-2]
\end{tikzcd}\]
or
\[\begin{tikzcd}
	{\Shv(S)\times\cG^{(n+1)}(S)} & {\Shv(\on{Tot}(\cG^{(n+1)}_S))} \\
	{\Shv(S)\times\cG^{(n)}(S)} & {\Shv(\on{Tot}(\cG^{(n)}_S))}
	\arrow["{\Phi^{n+1}_S}", from=1-1, to=1-2]
	\arrow["{\sigma_i}", from=2-1, to=1-1]
	\arrow["{\Phi^{n-1}_S}"', from=2-1, to=2-2]
	\arrow["{\sigma_i}"', from=2-2, to=1-2]
\end{tikzcd}\]
Both are special cases of Lemma \ref{lem: Phi map and change of group}. For $\pa_i$, $i>0$ we let 
\[c:A_1= A^{\times n+1}\rightarrow A_2=A^{\times n},\quad (a_1,\cdots,a_{n+1})\mapsto (a_1,\cdots, a_{i-1},a_i a_{i+1},a_{i+2},\cdots,a_{n+1})\]  and let $\cP=\cG^{(n)}_S$ and $\cK=\cL_{\on{reg}}^{\boxtimes n+1}$. Then we have $c_!\cK=\cL^{\boxtimes n}_{\on{reg}}$, hence $\Phi^n_S=\Phi_{\cK}$ and $\Phi^{n-1}_S=\Phi_{c_!\cK}$ in the notation of Lemma \ref{lem: Phi map and change of group}.
\medskip

For $\sigma_i$ we make a similar argument with
\[A_1=A^{\times n+1}\rightarrow A_2=A^{\times n+2},\quad (a_1,\cdots, a_{n+1})\mapsto (a_1,\cdots,a_i,1,a_{i+1},\cdots,a_{n+1})\]
and $\cP=\cG^{(n)}_S$ and $\cK=\cL_{\on{reg}}^{\boxtimes n+1}$.

\sssec{Compatibility for $\pa_0$}
We aim to show the commutativity of
\[\begin{tikzcd}
	{\Shv(S)\times\cG^{(n)}(S)} & {\Shv(\on{Tot}(\cG^{(n)}_S))} \\
	{\Shv(S)\times\cG^{(n-1)}(S)} & {\Shv(\on{Tot}(\cG^{(n-1)}_S))}
	\arrow["{\Phi^n_S}", from=1-1, to=1-2]
	\arrow["{\pa_0}"', from=1-1, to=2-1]
	\arrow["{\pa_0}", from=1-2, to=2-2]
	\arrow["{\Phi^{n-1}_S}"', from=2-1, to=2-2]
\end{tikzcd}\]
We start with a reduction step.
We observe that each category in the diagram has a $(\Bet A)(S)^{\times n}$-action and each arrow is equivariant under this action.
An object $\cT\in (\Bet A)(S)^{\times n}$ is an $A^{\times n}$-torsor over $S$, and we can form an $A^{\times n+1}$-torsor $\cT'$ over $S$ using 
\[A^{\times n}\rightarrow A^{\times n+1},\quad (a_1,\cdots,a_n)\mapsto (1,a_1,\cdots,a_n).\]

Since $\cG^{(n)}$ is an $A^{\times n+1}$-gerbe, the category $\Shv(S)\times \cG^{(n)}(S)$ has a $\Bet A(S)^{\times n+1}$-action. Then the action of $\cT$ on $\Shv(S)\times \cG^{(n)}(S)$ is given by the action of $\cT'$.
The category $\Shv(S)\times\cG^{(n-1)}(S)$ has a $(\Bet A)(S)^{\times n}$-action because $\cG^{(n-1)}$ is an $A^{\times n}$-gerbe. 
\medskip

$(\Bet A)(S)^{\times n}$ also acts on $\Shv(\on{Tot}(\cG^{(n)}_S))$: an object $\cT$ determines an automorphism 
\[\on{act}_{\cT'}:\on{Tot}(\cG^{(n)})\rightarrow \on{Tot}(\cG^{(n)})\] which is the action of $\cT'$. We define the action of $\cT$ on $\Shv(\on{Tot}(\cG^{(n)}))$ to be the pushforward along $\on{act}_{\cT'}$. Likewise, $\cT$ acts on $\Shv(\on{Tot}(\cG^{(n-1)}))$ by the pushforward along 
\[\on{act}_{\cT}:\on{Tot}(\cG^{(n-1)})\rightarrow \on{Tot}(\cG^{(n-1)}).\]

Each arrow is equivariant under this $(\Bet A)(S)^{\times n}$-action. The equivariance for $\Phi^n_S$ and $\Phi^{n-1}_S$ comes from construction, and the equivariance of transition maps $\pa_0$ come from the fact that $\pa_0$ only affects the first component in the $A^{n+1}$-gerbe structure on $\cG^{(n)}$, while the $(\Bet A)(S)^{\times n}$-action is trivial on the first factor. Using this equivariance, we can just check the commutativity for objects of the form
\begin{equation}\label{eq: test obj}
    \cO_{\on{test}}:=(\cF,\cT,\on{triv},\cdots,\on{triv},p)\in \Shv(S)\times \Bet A(S)^{\times n}\times \cG(S).
\end{equation}

\sssec{}
The action of $\Bet A(S)$ on $\Shv(S)$ is as follows: an $A$-torsor $\cT: S\rightarrow \Bet A$ acts on $\cF\in\Shv(S)$ as
$\cF\cdot \cT= \cF\otimes (\cT^*\cL_{\chi})$. Hence 
\[\pa_0(\cO_{test})=(\cF\otimes (\cT^*\cL_{\chi}),\on{triv},\cdots,\on{triv},p)\in \Shv(S)\times \Bet A(S)^{n-1}\times\cG(S).\]
Consider the map
\[p_{\on{triv}}: (\Bet A)^{\times n}\times S\rightarrow (\Bet A)^{\times  n-1}\times\on{Tot}(\cG_S),\quad (a_1,\cdots, a_n,s)\mapsto (a_1,a_2,\cdots, a_{n-1}, a_{n}p(s))\]
and 
\[p'_{\on{triv}}: (\Bet A)^{\times n+1}\times S\rightarrow (\Bet A)^{\times  n}\times\on{Tot}(\cG_S),\quad (a_1,\cdots, a_{n+1},s)\mapsto (\cT(s)a_1,a_2,\cdots, a_{n}, a_{n+1}p(s)).\]
Then 
\begin{equation}\label{eq: upper circuit}
    \Phi^n_S(\cO_{\on{test}})=(p'_{\on{triv}})_!(\cL^{\boxtimes n+1}_{\on{reg}}\boxtimes \cF)
\end{equation}
and
\begin{equation}\label{eq: lower circuit}
    \Phi^{n-1}_S(\pa_0(\cO_{\on{test}}))=(p_{\triv})_!(\cL^{\boxtimes n}_{\on{reg}}\boxtimes (\cF\otimes \cT^*\cL_{\chi})).
\end{equation}
It remains to show:
\begin{claim}\label{claim: O test claim}
    The image of \eqref{eq: upper circuit} under $\pa_0$ is canonically equivalent to \eqref{eq: lower circuit} in $\Shv(\on{Tot}(\cG^{(n-1)}_S))$.
\end{claim}
\begin{proof}
    Let us work in the case $n=1$ just for notational convenience.
Consider the $*$-pullback functor $\tau:\Shv(\on{Tot}(\cG_S))\rightarrow\Shv(\Bet A\times S)$ along the equivalence
\[\Bet A\times S\rightarrow \on{Tot}(\cG_S),\quad (a_1,s)\mapsto  a_1 p(s).\]
Because $\tau$ is determined by the data of $\cO_{\on{test}}$, it suffices to check Claim \ref{claim: O test claim} after applying $\tau$ to both sides: it suffices to construct an equivalence
\begin{equation}\label{eq: assertion after tau}
    \tau \circ\pa_0\circ (p'_{\triv})_!(\cL_{\on{reg}}\boxtimes\cL_{\on{reg}}\boxtimes \cF)\cong \tau\circ (p_{\triv})_!(\cL_{\on{reg}}\boxtimes (\cF\otimes \cT^*\cL_{\chi})).
\end{equation}
We observe that this reduces the proof to the case $\cG_S$ is trivial. To elaborate, we consider the map
\[p'^0_{\triv}: \Bet A\times \Bet A \times S\rightarrow \Bet A\times \Bet A\times S,\quad (a_1, a_2, s)\mapsto (\cT(s)a_1, a_2, s).\]
The left-hand side of \eqref{eq: assertion after tau} is equivalent to $\pa_0\circ (p'^0_{\triv})_!(\cL_{\on{reg}}\boxtimes \cL_{\on{reg}}\boxtimes \cF)$ where we consider
\[\pa_0: \Shv(\Bet A)\otimes \Shv(\Bet A\times S)\xrightarrow{\on{act}_{\chi}\otimes \Shv(\Bet A\times S)}\Shv(\Bet A\times S)\]
using $\Vect_{-\chi}\cong \Vect$, and the right-hand side of \eqref{eq: assertion after tau} is equivalent to $\cL_{\on{reg}}\boxtimes (\cF\otimes\cT^*\cL_{\chi})$. Therefore, we need an equivalence
\begin{equation}\label{eq: equiv after reduction to triv}
    \pa_0\circ (p'^0_{\triv})_!(\cL_{\on{reg}}\boxtimes \cL_{\on{reg}}\boxtimes \cF)\cong \cL_{\on{reg}}\boxtimes (\cF\otimes\cT^*\cL_{\chi}).
\end{equation}
Note that this statement is just the special case of \eqref{eq: assertion after tau} when $\cG_S$ is trivial. 
\medskip

We have $(p'^0_{\triv})_!=((p'^0_{\triv})^{-1})^*$ and 
\[(p'^0_{\triv})^{-1}: \Bet A\times\Bet A\times S\rightarrow \Bet A\times\Bet A\times S,\quad (a_1,a_2,s)\mapsto (\cT(s)^{-1}a_1,a_2,s)\]
gets decomposed as
$(p'^0_{\triv})^{-1}=(m\times \Bet A\times S)\circ \Gamma$
where $m$ is the multiplication map of $\Bet A$ and
\[\Gamma: \Bet A\times\Bet A\times S\rightarrow \Bet A\times\Bet A\times\Bet A\times S,\quad (a_1,a_2,s)\mapsto (a_1, \cT(s)^{-1},a_2,s).\]
Since $m^*\cL_{\on{reg}}=\bigoplus_{\chi'} \cL_{\chi'}\boxtimes \cL_{\chi'}$ we have 
\begin{align*}
    (p'^0_{\triv})_!(\cL_{\on{reg}}\boxtimes \cL_{\on{reg}}\boxtimes \cF) &= ((p'^0_{\triv})^{-1})^* (\cL_{\on{reg}}\boxtimes \cL_{\on{reg}}\boxtimes \cF) \\
    &  
    =\bigoplus_{\chi'} \Gamma^*(\cL_{\chi'}\boxtimes\cL_{\chi'}\boxtimes \cL_{\on{reg}}\boxtimes\cF)= \bigoplus_{\chi'}(\cL_{\chi'}\boxtimes \cL_{\on{reg}}\boxtimes (\cF\otimes\cT^{*}\cL_{-\chi'})).
\end{align*}
For the last equality, we used the fact that the involution on $\Shv(\Bet A)$ induced by $A\rightarrow A$, $a\mapsto -a$ interchanges $\cL_{\chi'}$ and $\cL_{-\chi'}$.
We have
\[\pa_0(\cL_{\chi'}\boxtimes \cL_{\on{reg}}\boxtimes (\cF\otimes\cT^{*}\cL_{-\chi'}))=\begin{cases}
    \cL_{\on{reg}}\boxtimes (\cF\otimes\cT^{*}\cL_{\chi}) &\te{if } \chi'=-\chi,\\
    0 &\te{otherwise,}
\end{cases}\]
hence the equivalence \eqref{eq: equiv after reduction to triv}.
\end{proof}

\sssec{}\label{sssec: showing equiv btw bar constructions}
The above arguments confirm that the maps $\Phi^n_S$ intertwine transition maps of $\{G^n_S\}_{[n]\in\Delta^{\op}}$ and $\{T^n_S\}_{[n]\in\Delta^{\op}}$, thus defining a canonical functor \eqref{eq: map between bar constructions} in $\Cat$. Denote this functor by $\Phi_S:G_S\rightarrow T_S$. Note that $\cG_S$ is trivializable by assumption. Once we choose a trivialization $\cG_S\cong \triv_S$, this induces an equivalence $G_S\cong \Shv(S)$ in $\Cat$ and an equivalence $T_S\cong \Vect_{-\chi}\otimes\Shv(S)$ in $\DGcat$. (We elaborate on these equivalences in Sect. \ref{sssec: proof of Claim trivial case}) We later prove:
\begin{claim}\label{claim: Phi n in the trivial case}
    If $\cG_S$ is trivial, then $\Phi_S:G_S\rightarrow T_S$ is identified with the equivalence $\Shv(S)\rightarrow \Vect_{-\chi}\otimes\Shv(S)$ obtained from $\Vect\cong \Vect_{-\chi}$.
\end{claim}

\sssec{}
Assume the validity of Claim \ref{claim: Phi n in the trivial case} for now.
This shows that $\Phi_S$ is an equivalence of categories. 
The category
 $G_S$ has a canonical lift to $\DGcat$ where the structure as a DG category is transported from $\Shv(S)$ under the equivalence $G_S\cong \Shv(S)$ considered above. From the discussion above, we conclude that $\Phi_S$ canonically lifts to an equivalence of DG categories. It is clear that the collection of equivalences $\Phi_S$ is compatible with maps in $\left((\Affft)_{/\cY}\right)^{\op}$, and thus becomes upgraded to an equivalence $\eqref{eq: twisted sheaf sheaf on triv} \rightarrow \eqref{eq: tensor product defn sheaf triv}$ between sheaves. This concludes the proof of Proposition \ref{prop: equiv btw sheaves triv}, and thus Theorem \ref{thm: twisted sheaves tensor product formula}. \qed

 \sssec{Proof of Claim \ref{claim: Phi n in the trivial case}}\label{sssec: proof of Claim trivial case}
 If $\cG_S=\triv_S$, then
 $G^n_S=\Shv(S)\times (\Bet A(S))^{n+1}$ and $T^n_S=\Vect_{-\chi}\otimes \Shv(\Bet A)^{\otimes n+1}\otimes \Shv(S)$. Then the functor $\Phi^n_S$ is
 \[\Phi^n_S:G^n_S\rightarrow T^n_S,\quad 
 (\cF,(a_1,\cdots,a_{n+1}))\mapsto \beta^{\triv}_!(\cL_{\on{reg}}^{\boxtimes n+1}\boxtimes \cF)\]
 where $\beta=(a_1,\cdots,a_{n+1})\in \Bet A(S)^{n+1}$ and $\beta^{\triv}$ is the automorphism
 \[\beta^{\triv}:(\Bet A)^{\times n+1}\times S\rightarrow (\Bet A)^{\times n+1}\times S,\quad (a'_1,\cdots,a'_{n+1},s)\mapsto (a_1(s) a'_1,\cdots, a_{n+1}(s)a'_{n+1},s).\]

 We first observe that the geometric realization $G_S$ is canonically equivalent to $\Shv(S)$; the insertion maps are
 \[G^n_S=\Shv(S)\times (\Bet A(S))^{n+1}\rightarrow \Shv(S)\]
 given by the ($(n+1)$-fold) action of $\Bet A(S)$ on $\Shv(S)$. Similarly, the geometric realization $T_S$ is canonically equivalent to $\Vect_{-\chi}\otimes\Shv(S)$; the insertion maps are
 \[T^n_S=\Vect_{-\chi}\otimes \Shv(\Bet A)^{\otimes n+1}\otimes \Shv(S)\rightarrow \Vect_{-\chi}\otimes \Shv(S)\]
where the map is the ($(n+1)$-fold) action of $\Shv(\Bet A)$ on $\Vect_{-\chi}$ tensored with the identity endomorphism of $\Shv(S)$. We can further use $\Vect_{-\chi}\cong \Vect$ to canonically identify $T_S\cong\Shv(S)$.
 \medskip

 From the construction of $\Phi_S$ from the collection $\{\Phi^n_S\}$, for each $n$ we have a canonical commutative square
\[\begin{tikzcd}
	{G^n_S} & {T^n_S} \\
	{G_S} & {T_S}
	\arrow["{\Phi^n_S}", from=1-1, to=1-2]
	\arrow["{\on{ins}^G_n}"', from=1-1, to=2-1]
	\arrow["{\on{ins}^T_n}", from=1-2, to=2-2]
	\arrow["{\Phi_S}"', from=2-1, to=2-2]
\end{tikzcd}\]
 Since insertion maps arise from the action of $\Shv(\Bet A)$, we identify the insertion map for $\on{ins}^T_n$ as
\[\Shv(\Bet A)^{\otimes n+1}\otimes \Shv(S)\xrightarrow{-\chi^{\otimes n+1}\otimes \Shv(S)}\Shv(S),\]
where $\Shv(\Bet A)\xrightarrow{-\chi}\Vect$ is the functor that gets identified with the functor $\QC(\on{Spec}(EA))\rightarrow \QC(\pt)$ taking fiber at $-\chi \in \on{Spec}(EA)$. 
Since the insertion map $\on{ins}^G_n$ comes from the action of $\Bet A(S)$ on $\Shv(S)$, the composition
\[\Shv(S)\xrightarrow{(\id,\triv_S,\cdots,\triv_S)} \Shv(S)\times \Bet A(S)^{\times n+1} \xrightarrow{\on{ins}^G_m}\Shv(S)\]
is the identity. Combining these, we recover $\Phi_S$ as the composition
\[\Shv(S)\xrightarrow{(\id,\triv_S,\cdots,\triv_S)} \Shv(S)\times \Bet A(S)^{\times n+1}\xrightarrow{\Phi^n_S}\Shv(\Bet A)^{\otimes n+1}\otimes \Shv(S)\xrightarrow{-\chi^{\otimes n+1}\otimes \Shv(S)}\Shv(S).\]
This sends $\cF\in\Shv(S)$ to
\[(-\chi^{\otimes n+1}\otimes \Shv(S))(\cL^{\boxtimes n+1}_{\on{reg}}\otimes \cF),\]
which is canonically equivalent to $\cF$ because $-\chi(\cL_{\on{reg}})=E$. Therefore, we conclude that $\Phi_S$ is the identity functor of $\Shv(S)$. \qed

\begin{rmk}
    In fact, for the purpose of constructing a map $\eqref{eq: twisted sheaf sheaf on triv}\rightarrow \eqref{eq: tensor product defn sheaf triv}$, we could have used   
    any kernel of the form $\cL^{\boxtimes n}_{\on{reg}}\boxtimes \cL$ for $\cL\in\Shv(\Bet A)$ for the definition of $\Phi^n_S$, i.e., define
    \[\Phi^n_S=\beta^{\triv}_*(\cL_{\on{reg}}^{\boxtimes n}\boxtimes \cL\boxtimes \cF).\]
    The collection of these maps would still be compatible with transition maps, hence define $\Phi_S:G_S\rightarrow T_S$.  
    However, for the resulting $\Phi_S$ to be the identity functor $\Shv(S)\rightarrow \Shv(S)$, we need an equivalence $-\chi(\cL)\cong E$. If we chose $\cL=\cL_{-\chi'}$ for $\chi'\ne \chi$, then the resulting $\Phi_S$ would be zero.
\end{rmk}

\subsection{Functoriality of twisted sheaves}\label{ssec: functoriality of twisted sheaves}
Having proved Theorem \ref{thm: twisted sheaves tensor product formula}, we can study the functoriality of twisted sheaves without difficulty. 

\sssec{Category of correspondences}
Recall the symmetric monoidal category
\[\PreStklft+\on{Ge}_A := (\PreStklft)_{/\Bet^2 A}\]
from Sect. \ref{sssec: functoriality twshv brief}. We consider the category
\begin{equation}\label{eq: Corr category}
    \on{Corr}\left(\PreStklft+\on{Ge}_A, \on{sch}\right)
\end{equation}
of correspondences. Its objects are $(\cY,\cG)\in\PreStklft+\on{Ge}_A$, and morphisms are diagrams
\begin{equation}\label{eq: Corr diagram}
    \begin{tikzcd}
	{(\cY_{1,2},\cG_{1,2})} & {(\cY_1,\cG_1)} \\
	{(\cY_2,\cG_2)}
	\arrow["g", from=1-1, to=1-2]
	\arrow["f"', from=1-1, to=2-1]
\end{tikzcd}
\end{equation}
where the map $\cY_{1,2}\rightarrow \cY_2$ underlying $f$ should be schematic, while there is no restriction to $g$. The composition of \[\begin{tikzcd}
	{(\cY_{2,3},\cG_{2,3})} & {(\cY_2,\cG_2)} \\
	{(\cY_3,\cG_3)}
	\arrow[ from=1-1, to=1-2]
	\arrow[ from=1-1, to=2-1]
\end{tikzcd}\]
with \eqref{eq: Corr diagram} is given by the diagram
\[\begin{tikzcd}
	{(\cY_{1,3},\cG_{1,3})} & {(\cY_1,\cG_1)} \\
	{(\cY_3,\cG_3)}
	\arrow[from=1-1, to=1-2]
	\arrow[from=1-1, to=2-1]
\end{tikzcd}\]
where 
\[(\cY_{1,3},\cG_{1,3})=(\cY_{1,2},\cG_{1,2})\times_{(\cY_2,\cG_2)}(\cY_{2,3},\cG_{2,3}) \quad \te{in} \ \PreStklft+\on{Ge}_A.\]
Note that
$\cY_{1,3}=\cY_{1,2}\times_{\cY_2}\cY_{2,3}$ in $\PreStklft$. The category \eqref{eq: Corr category} is equipped with a natural symmetric monoidal structure where
\[(\cY_1,\cG_1)\otimes (\cY_2,\cG_2) = (\cY_1\times\cY_2, \cG_1\boxplus \cG_2).\]

\sssec{}
We are interested in the analogue of the six-functor formalism
\begin{equation}\label{eq: shv Corr sch}
    \Shv:\Corr(\PreStklft,\on{sch})\rightarrow\DGcat
\end{equation}
in the context of twisted sheaves. (see Sect. \ref{sssec: six functor stacky}) We will explain that Theorem \ref{thm: twisted sheaves tensor product formula} enables us to pass from the six-functor formalism
\begin{equation}\label{eq: six functor stacky}
    \Shv: \Corr(\PreStklft,\wt{E}')\rightarrow\DGcat
\end{equation}
introduced in Sect. \ref{sssec: six functor stacky} to the following results about twisted sheaves:

\begin{thm}\label{thm: twShvCorr functor}
    We have a right-lax symmetric monoidal functor
    \begin{equation}\label{eq: twisted shv Corr case 1}
       \Shv^{\on{tw}}_{\on{Corr}}:\on{Corr}\left(\PreStklft+\on{Ge}_A, \on{sch}\right)\rightarrow \DGcat 
    \end{equation}
    which sends an object $(\cY,\cG)$ to the DG category $\Shv^{\cG}(\cY)$.
\end{thm}

\sssec{Construction of $\Shv^{\on{tw}}_{\Corr}$}
We have a functor
\begin{equation}\label{eq: functor Tot}
    \on{Tot}(\mhyphen):\PreStklft+\on{Ge}_A\rightarrow \Bet A\lMod (\PreStklft)
\end{equation}
of taking total space of gerbes. Indeed, we have a gerbe $\cY\xrightarrow{\cG}\Bet^2 A$, there is a natural action of
\[\Bet A=\pt\times_{\Bet^2 A} \pt \quad \te{on} \quad \on{Tot}(\cG)= \cY\times_{\Bet^2 A}\pt.\]
The functor \eqref{eq: functor Tot} preserves fiber products, hence we get an induced functor
\[\on{Corr}(\PreStklft+\on{Ge}_A,\on{sch})\rightarrow \on{Corr}(\Bet A\lMod(\PreStklft),\on{sch}).\]
Given a prestack $\wt{\cY}$ with a $\Bet A$-action, then $\Shv(\wt{\cY})$ is a $(\Shv(\Bet A),\ast)$-module in $\DGcat$. Moreover, if we have $\wt{\cY}_1\rightarrow \wt{\cY}_2$ as $\Bet A$-modules, then both
\[f_*: \Shv(\wt{\cY}_1)\rightarrow \Shv(\wt{\cY}_2)\quad \te{and}\quad f^!:\Shv(\wt{Y}_2)\rightarrow \Shv(\wt{Y}_1)\]
are $(\Shv(\Bet A),\ast)$-module maps, and thus we have
\begin{equation}\label{eq: Corr BA mod to Shv BA mod}
    \on{Corr}(\Bet A\lMod(\PreStklft),\on{sch})\rightarrow (\Shv(\Bet A),\ast)\lMod(\DGcat).
\end{equation}

Composing these, we get
\begin{equation}\label{eq: Corr prestk gerbe to module cat}
    \on{Corr}(\PreStklft+\on{Ge}_A, \on{sch})\rightarrow (\Shv(\Bet A),\ast)\lMod(\DGcat).
\end{equation}
The desired functor $\Shv^{\on{tw}}_{\Corr}$ is obtained as the composition
\[ \on{Corr}(\PreStklft+\on{Ge}_A, \on{sch})\xrightarrow{\eqref{eq: Corr prestk gerbe to module cat}} (\Shv(\Bet A),\ast)\lMod(\DGcat)\xrightarrow{(\mhyphen)\otimes_{\Shv(\Bet A)}\Vect_{-\chi}} \DGcat.\]
By construction, this sends a pair $(Y,\cG)$ to
\[\Shv(\on{Tot}(\cG))\otimes_{\Shv(\Bet A)}\Vect_{-\chi},\]
which is canonically equivalent to $\Shv^{\cG}(\cY)$ by Theorem \ref{thm: twisted sheaves tensor product formula}.

\sssec{}
We can provide a more formal construction of \eqref{eq: Corr BA mod to Shv BA mod} using the six-functor formalism \eqref{eq: six functor stacky}. The stack $\Bet A$ is naturally a commutative algebra object in 
$\on{Corr}(\PreStklft,\on{sch})$, and $(\Shv(\Bet A),*)\in\DGCatSymMon$ is the image of this under \eqref{eq: six functor stacky}. We have a natural functor
\begin{equation}\label{eq: Corr BA to BA Corr}
    \on{Corr}(\Bet A\lMod(\PreStklft),\on{sch})\rightarrow \Bet A\lMod (\on{Corr}(\PreStklft)_{\on{sch},\on{all}}).
\end{equation}
It sends $\wt{\cY}\in\Bet A\lMod(\PreStklft)$ to $\wt{\cY}\in \on{Corr}(\PreStklft)$ equipped with an action
\[\Bet A\times \wt{\cY}\rightarrow \wt{\cY} \quad \te{in} \ \on{Corr}(\PreStklft)_{\on{DM},\on{all}}\]
given by the diagram
\[\begin{tikzcd}
	{\Bet A\times \wt{\cY}} & {\Bet A\times \wt{\cY}} \\
	{\wt{\cY}}
	\arrow["\id", from=1-1, to=1-2]
	\arrow["{\on{act}}"', from=1-1, to=2-1]
\end{tikzcd}\]
Since the functor \eqref{eq: six functor stacky} is right-lax symmetric monoidal we get
\begin{equation}\label{eq: BA Corr to module cat}
    \Bet A\lMod (\on{Corr}(\PreStklft,\on{sch}))\rightarrow (\Shv(\Bet A),\ast)\lMod(\DGcat)
\end{equation}
from it. We get the desired functor \eqref{eq: Corr BA mod to Shv BA mod} composing \eqref{eq: Corr BA to BA Corr} and \eqref{eq: BA Corr to module cat}.

\sssec{}
We briefly sketch that $\Shv^{\on{tw}}_{\on{Corr}}$ has a canonical right-lax symmetric monoidal structure. To see this, given $(\cY_1,\cG_1)$ and $(\cY_2,\cG_2)$ in $\PreStklft+\on{Ge}_A$ we have to construct a canonical functor
\begin{equation}\label{eq: right lax tw Corr}
    \Shv^{\cG_1}(\cY_1)\otimes \Shv^{\cG_2}(\cY_2)\rightarrow \Shv^{\cG_1\boxtimes\cG_2}(\cY_1\times\cY_2).
\end{equation}
From Theorem \ref{thm: twisted sheaves tensor product formula} the left-hand side is identified with
\begin{equation}\label{eq: tensor prod A2}
    \left(\Shv(\on{Tot}(\cG_1))\otimes \Shv(\on{Tot}(\cG_2))\right)\otimes_{\Shv(\Bet A^2)}\Vect_{(-\chi,-\chi)}.
\end{equation}
From the right-lax symmetric monoidality of \eqref{eq: shv Corr sch} we have a functor
\begin{equation}\label{eq: right lax applied}
    \eqref{eq: tensor prod A2}\rightarrow \Shv(\on{Tot}(\cG_1)\times\on{Tot}(\cG_2))\otimes_{\Shv(\Bet A^2)}\Vect_{(-\chi,-\chi)}.
\end{equation}
Consider the summation map $A^2\rightarrow A$. We have 
\[\on{Tot}(\cG_1\boxplus\cG_2)=(\on{Tot}(\cG_1)\times\on{Tot}(\cG_2))\times^{\Bet A^2}\Bet A,\]
and the $(-\chi,-\chi)$ is the pullback of $-\chi$ along the summation map. Therefore, \eqref{eq: right lax applied} is equivalent to
\[
\Shv(\on{Tot}(\cG_1\boxplus\cG_2))\otimes_{\Shv(\Bet A)}\Vect_{-\chi},
\]
which is in turn equivalent to 
\[\Shv^{\cG_1\boxplus\cG_2}(\cY_1\times\cY_2)\]
by Theorem \ref{thm: twisted sheaves tensor product formula}. Combining all these, we get the desired functor \eqref{eq: right lax tw Corr}.

\subsection{Perverse $t$-structure on twisted sheaves}\label{ssec: t str twisting}
In this subsection, we study the perverse $t$-structure on the category of twisted ind-constructible sheaves on algebraic stacks. Throughout this subsection, let $\cY$ be an algebraic stack with an affine diagonal.

\sssec{}
Let $Y$ be a quasi-compact scheme. We may consider the DG category
\[\Shv^{\indconstr}(Y):=\Ind(\Shv^{\constr}(Y)).\]

As in \cite[Sect. 1.1.5]{AGKRRV} there is a unique $t$-structure on $\Shv^{\indconstr}(Y)$ compatible with filtered colimits that extends the perverse $t$-structure on $\Shv^{\constr}(Y)$. We call this the \textbf{perverse} $t$-\textbf{structure} on $\Shv^{\indconstr}(Y)$. Its heart is identified with $\Ind(\Perv(Y))$.
\medskip

It is proved in \cite[Sect. E.1]{AGKRRV} that for such $Y$, the perverse $t$-structure on $\Shv^{\indconstr}(Y)$ is left-complete. Finally, from \cite[Proposition 2.13]{AGH} the perverse $t$-structure is accessible, i.e., $\Shv^{\indconstr}(Y)^{\le 0}$ is presentable.
To summarize, the perverse $t$-structure on $\Shv^{\indconstr}(Y)$ has the following properties:
\begin{enumerate}
    \item it is left-complete;
    \item it is compatible with filtered colimits;
    \item $\Shv^{\indconstr}(Y)^{\le 0}$ has all colimits; 
    \item its heart is $\Ind$-complete.
\end{enumerate}
In fact, the last property is implied by the second and the third:

\begin{lemma}\label{lem: heart Ind complete}
    Say that $\bC\in\DGcat$ is equipped with a $t$-structure which is compatible with filtered colimits and $\bC^{\le 0}$ has filtered colimits. Then its heart $\bC^{\heartsuit}$ is $\Ind$-complete and the inclusion $\bC^{\heartsuit}\hookrightarrow \bC$ preserves filtered colimits.
\end{lemma}

\sssec{}
We will generalize the aforementioned properties of perverse $t$-structure on $\Shv^{\indconstr}(Y)$ to twisted ind-constructible sheaves on algebraic stacks. The following categorical lemma will be used repeatedly.

\begin{lemma}\label{lem: limit t str}
    Consider a diagram
    \[I\rightarrow \DGcat,\quad i\mapsto \bC_i.\] Assume that each $\bC_i$ is endowed with a $t$-structure, and all of the transition functors $\bC_i\rightarrow \bC_j$ are $t$-exact. Let $\bC:=\lim_{i\in I} \bC_i$. Then we have:
    \begin{enumerate}
        \item the category $\bC$ acquires a unique $t$-structure such that the evaluation functors $\bC\rightarrow \bC_i$ are $t$-exact;
        \item if the $t$-structure on each $\bC_i$ is compatible with filtered colimits, then so is the one on $\bC$;
        \item if the $t$-structure on each $\bC_i$ is left-complete, then so is the one on $\bC$;
        \item if each $\bC^{\le 0}_i$ has filtered colimits, then so is the one on $\bC$.
    \end{enumerate}
\end{lemma}

\sssec{}
Let $\cY$ be an algebraic stack with an affine diagonal. From \cite[F.1.2]{AGKRRV}, we have 
\[\Shv^{\indconstr}(\cY)= \lim_{S\rightarrow \cY}\Shv^{\indconstr}(S)\]
where the index category is the 1-full subcategory of $(\Affft)_{/\cY}$, where we allow only smooth $S\rightarrow \cY$'s and smooth maps between those. In this limit expression, the transition functors are smooth pullbacks, hence perverse $t$-exact up to shifts. Therefore, we can apply Lemma \ref{lem: limit t str} to conclude the following.

\begin{prop}\label{prop: t str indconstr alg stk}
For any algebraic stack $\cY$, the DG category $\Shv^{\indconstr}(\cY)$ has a canonical perverse $t$-structure determined by the requirement that for any smooth map $S\rightarrow\cY$ from an affine scheme, the $*$-pullback $\Shv^{\indconstr}(\cY)\rightarrow \Shv^{\indconstr}(S)$ is perverse $t$-exact up to a dimension shift. This construction satisfies the following propreties.
    \begin{enumerate}
        \item If we have a smooth map $\cY'\rightarrow\cY$ between algebraic stacks, then the $*$-pullback
    \[\Shv^{\indconstr}(\cY)\rightarrow \Shv^{\indconstr}(\cY')\]
    is perverse $t$-exact up to a dimension shift.
    \item The perverse $t$-structure on $\Shv^{\indconstr}(\cY)$ is left-complete, compatible with filtered colimits, and its non-positive part has filtered colimits.
    \end{enumerate}
\end{prop}

\sssec{}\label{sssec: t str on twisted sheaves triv case}
Say that we have an \'etale $A$-gerbe $\cG$ over $\cY$ satisfying the assumption as before.  We aim to explain the perverse $t$-structure on this DG category $\Shv^{\indconstr,\cG}(\cY)$ of twisted sheaves.
\medskip

First, assume that $\cG$ is trivializable on $\cY$. There is an equivalence between $\Shv^{\indconstr}(\cY)$ and $\Shv^{\indconstr,\cG}(\cY)$ depending on the choice of trivialization of $\cG$. There is a canonical perverse $t$-structure on $\Shv^{\indconstr,\cG}(\cY)$ such that for any such choice of trivialization, the perverse $t$-structure is identified with the perverse $t$-structure of $\Shv^{\indconstr}(\cY)$.
 If we have $(\cY,\cG)$ where $\cG$ is trivializable and a smooth map $\cY'\rightarrow\cY$, then the $*$-pullback map
\[\Shv^{\indconstr,\cG}(\cY)\rightarrow \Shv^{\indconstr,\cG|_{\cY'}}(\cY')\]
is perverse $t$-exact up to a dimension shift.
\medskip

Second, we construct a perverse $t$-structure on any such $(\cY,\cG)$, removing the hypothesis that $\cG$ is trivializable. Given a pair $(Y,\cG)$, from the \'etale cover $\on{Tot}(\cG)\rightarrow\cY$ we get a descent
\begin{equation}\label{eq: descent tw indconstr}
    \Shv^{\indconstr,\cG}(\cY) = \lim_{[n]\in\Delta}\Shv^{\indconstr, \cG}(\on{Tot}(\cG)^n/\cY).
\end{equation}
On each $\on{Tot}(\cG)^n/\cY$ the gerbe $\cG$ is trivializable, hence each term in this limit diagram is equipped with a canonical perverse $t$-structure and all the transition maps are perverse $t$-exact up to a dimension shift. Therefore, we can apply Lemma \ref{lem: limit t str} to get a canonical perverse $t$-structure on $\Shv^{\indconstr,\cG}(\cY)$. Since every map $\on{Tot}(\cG)^n/\cY\rightarrow \cY$ factors through $\on{Tot}(\cG)\rightarrow\cY$, we have
\begin{equation}\label{eq: criterion pullback to tot}
\cF\in \Shv^{\indconstr,\cG}(\cY)^{?} \iff \cF|_{\on{Tot}(\cG)}\in \Shv^{\indconstr,\cG}(\on{Tot}(\cG))^{?},\quad ?= \le 0, \ge 0.
\end{equation}

We have the following analogue of Proposition \ref{prop: t str indconstr alg stk} for twisted sheaves:
\begin{prop}\label{prop: t str twisted indconstr alg stk}
    For any $(\cY,\cG)\in\PreStklft+\on{Ge}_A$ where $\cY$ is an algebraic stack, the DG category $\Shv^{\indconstr,\cG}(\cY)$ of twisted sheaves has a canonical perverse $t$-structure determined by the requirement \eqref{eq: criterion pullback to tot} together with the construction in Sect. \ref{sssec: t str on twisted sheaves triv case} for the trivializable case. This construction satisfies the following properties:
    \begin{enumerate}
        \item If we have a smooth map $\cY'\rightarrow \cY$ between algebraic stacks, then the $*$-pullback
        \[\Shv^{\indconstr, \cG}(\cY)\rightarrow \Shv^{\indconstr, \cG|_{\cY'}}(\cY')\]
        is perverse $t$-exact up to a dimension shift.
        \item The perverse $t$-structure on $\Shv^{\indconstr,\cG}(\cY)$ is left-complete, compatible with filtered colimits, and its non-positive part has filtered colimits.
    \end{enumerate}
\end{prop}

We have the following alternative characterization of the perverse $t$-structure on $\Shv^{\indconstr,\cG}(\cY)$:
\begin{prop}
    The perverse $t$-structure on $\Shv^{\indconstr,\cG}(\cY)$ is uniquely characterized by the requirement that for any smooth $(S \rightarrow \cY)\in (\Affft)_{/\cY}$ where $\cG$ is trivializable over $S$, the $*$-pullback
    \[\Shv^{\indconstr,\cG}(\cY)\rightarrow \Shv^{\indconstr}(S)\]
    is perverse $t$-exact up to a dimension shift, where we used any trivialization of $\cG|_S$ to form the functor.
\end{prop}

\sssec{Supplement on Sect. \ref{sssec: ind nature precise}}
In Sect. \ref{sssec: ind nature precise} we have made the following claim:

\begin{prop}\label{prop: t str on untwisted indconstr Hecke loc}
    For each $\lambda\in\Lambda^+$, the DG category $\left(\Shv^{\on{ind-constr},q}(\Hecke^{\loc,\le\lambda}_{X'})\right)_{-\cG_Z}$ has a canonical perverse $t$-structure, which is left-complete, compatible with filtered colimits, and its non-positive part has filtered colimits. Its abelian heart is $\Ind$-complete.
\end{prop}

\sssec{Twisting and $t$-structure}\label{sssec: twisting and t str}
Suppose that a sheaf $\cF\in\Shv(X_{\et};\DGcat)$ is equipped with the following structure: for each $X'\in X_{\et}$ the DG category $\cF(X')$ is equipped with a $t$-structure, and all the restriction maps $\cF(X')\rightarrow \cF(X'')$ are $t$-exact. Suppose $\sfA$ is a torsion abelian group and $\cF$ has an $\sfA$-symmetry. Then we have a twist $\cF_{\cG}$. We claim that each section of $\cF_{\cG}$ has a $t$-structure and all the restriction maps are $t$-exact. This is done similar to the construction of perverse $t$-structure on twisted $\Ind$-constructible sheaves.
\medskip

On $X'\in X^{\on{triv}}_{\et}$ we have a canonical $t$-structure on $\cF_{\cG}(X')$ such that for every trivialization $\cG|_{X'}\simeq \triv_{X'}$ the induced $\cF_{\cG}(X')\simeq \cF(X')$ is $t$-exact. For $X'\in X_{\et}$ in general, we choose $X''\rightarrow X'$ such that $X''\in X^{\triv}_{\et}$. Since $(X'')^n/X'$ in the \v{C}ech nerve always lie in $X^{\triv}_{
\et}$, we can use the limit expression
\[\cF_{\cG}(X')=\lim_{[n]\in\Delta}\cF_{\cG}((X'')^n/X') \quad \te{in }\DGcat\]
and Lemma \ref{lem: limit t str} to obtain a canonical $t$-structure on $\cF_{\cG}(X')$. This $t$-structure does not depend on the choice of $X''$. 
\medskip

Suppose that the $t$-structure on $\cF$ has the following additional property: for each $X'\in X_{\et}$ the $t$-structure on $\cF(X')$ is left-complete, compatible with filtered colimits, its non-positive part has all colimits. Then by Lemma \ref{lem: limit t str} the sheaf $\cF_{\cG}$ has the same property.

\sssec{Proof of Proposition \ref{prop: t str on untwisted indconstr Hecke loc}}
Recall that 
\[\Shv^q(\Hecke^{\loc,\le\lambda}_{X}):=\Shv^q(L^{n}G\backslash \Gr^{\le\lambda}_X)\]
for sufficiently large $n$. From Proposition \ref{prop: t str twisted indconstr alg stk} this has a perverse $t$-structure which is left-complete, compatible with filtered colimits, and its non-positive part has all colimits. We can replace $X$ by any $X'\in X_{\et}$ to form a sheaf of DG categories with a $t$-structure. Since an \'etale $*$-pullback is perverse $t$-exact, we are in the situation of the last paragraph of Sect. \ref{sssec: twisting and t str}. Therefore, the global section of the twist by $-\cG_Z$, i.e., $\left(\Shv^{\on{ind-constr},q}(\Hecke^{\loc,\le\lambda}_{X'})\right)_{-\cG_Z}$, has a natural perverse $t$-structure with the same properties. Now we conclude from Lemma \ref{lem: heart Ind complete}. \qed

\section{Twisted representations}\label{sec: twisted representations}
We recall the construction of Sect. \ref{sssec: twisted representations}. We consider an algebraic group $\mathsf{H}$ over $E$, a scheme $Y\in\Schft$, a torsion abelian group $\sfA$ with a group homomorphism $\sfA\rightarrow Z_\mathsf{H}(E)$, and an \'etale $\sfA$-gerbe $\cG$ over $Y$. Then we could define the category $\Rep(\mathsf{H})^{\loc}_{\cG}(Y)$ of twisted $\mathsf{H}$-representations. This is the main result of this section:
\begin{thm}\label{thm: twisted representations}
    Let $\mathsf{H}$, $Y$, $\sfA$, and $\cG$ be as above. We have a canonical equivalence
    \begin{equation}\label{eq: twrep QC}
        \QC\left(\on{Tot}(\Phi_B(\cG))\times^{\Bet{\sfA}}\Bet{\mathsf{H}}\right)=\Rep(\mathsf{H})^{\loc}_{\cG}(Y).
    \end{equation}
\end{thm}

\sssec{}
We explain the notation on the left-hand side.
Recall from Sect. \ref{ssec: Betti realization} that $\Phi_B(\cG)$ is an $\sfA$-gerbe over $(Y^{\topl,\sing})_B \in \Stk^E$. We can take its total space $\Tot(\Phi_B(\cG))$, which gets equipped with an action of $\Bet \sfA$, and 
\begin{equation}\label{eq: ass bundle BH}
    \on{Tot}(\Phi_B(\cG))\times^{\Bet \sfA}\Bet \mathsf{H}\in\Stk^E
\end{equation}
is obtained from the associated bundle construction. (see Sect. \ref{sssec: Associated bundle construction})

\begin{ex}\label{ex: equiv trivial gerbe}
    If $\cG$ is the trivial gerbe, then so is $\Phi_B(\cG)$. The right-hand side of \eqref{eq: twrep QC} is $\Rep(\mathsf{H})\otimes \LS(Y)$, and \eqref{eq: ass bundle BH} is $(Y^{\topl,\sing})_B\times \Bet \mathsf{H}$. Because $\QC((Y^{\topl,\sing})_B)=\LS(Y)$ is dualizable, we have
    \[\QC((Y^{\topl,\sing})_B\times \Bet \mathsf{H})=\LS(Y)\otimes \Rep(\mathsf{H}).\]
\end{ex}

\subsection{Proof of Theorem \ref{thm: twisted representations}}
Our proof of Theorem \ref{thm: twisted representations} is very similar to the proof of Theorem \ref{thm: twisted sheaves tensor product formula}, so we omit some details.

\sssec{}
First, we reduce to the case where $\sfA$ is finite. If we have a group homomorphism $\sfA'\rightarrow\sfA$ between abelian groups and the $\sfA$-gerbe $\cG$ is induced from an $\sfA'$-gerbe $\cG'$, then we can define both sides of \eqref{eq: twrep QC} using $\sfA'$ together with the composition $\sfA'\rightarrow \sfA\rightarrow Z_{\mathsf{H}}(E)$ and the gerbe $\cG'$. We get the same result using either choice, i.e., we have canonical equivalences
\[\on{Tot}(\Phi_B(\cG))\times^{\Bet \sfA}\Bet \mathsf{H}= \on{Tot}(\Phi_B(\cG'))\times^{\Bet \sfA'}\Bet \mathsf{H}\]
and
\[\Rep(\mathsf{H})^{\loc}_{\cG}(Y)=\Rep(\mathsf{H})^{\loc}_{\cG'}(Y).\]
Therefore, Theorem \ref{thm: twisted representations} for $\sfA'$ and $\cG'$ implies the theorem for $\sfA$ and $\cG$.
\medskip

An $\sfA$-gerbe on $Y$ is classified by a class in $H^2_{\et}(Y; \sfA)$, and this is the filtered colimit of $H^2_{\et}(Y;\sfA')$'s where $\sfA'$ ranges in the collection of finite subgroups of $\sfA$. Therefore, there should be some finite $\sfA'\subseteq \sfA$ such that $\cG$ is induced from an $\sfA'$-gerbe $\cG'$. From the discussion above, we reduce to the case $\sfA$ is finite.

\sssec{}
Next, we observe that both sides of \eqref{eq: twrep QC} naturally form sheaves over $Y_{\et}$, and we will show those two sheaves are equivalent to each other. From Sect. \ref{ssec: restriction to triv}, it suffices to construct the equivalence over $Y^{\triv}_{\et}$, which is the subsite of $Y_{\et}$ spanned by $Y'\rightarrow Y$'s such that $\cG|_{Y'}$ is trivializable.
\medskip

The sheafy version of the right-hand side is $\mathbf{Rep}(\mathsf{H})^{\loc}_{\cG}$, and its restriction to $Y^{\triv}_{\et}$ is
\begin{align}\label{eq: twRep sheaf on triv}
    (Y^{\triv}_{\et})^{\op}\rightarrow\DGcat,& \\
    Y'&\mapsto \colim_{[n]\in\Delta^{\op}}(\LS(Y')\otimes\Rep(\mathsf{H}))\times \Bet A(Y')^{\times n}\times \Gamma(Y',\cG). \nonumber
\end{align}
Following the notation of Sect. \ref{sssec: GnS TnS} we let
\[G^n_{Y'}:= (\LS(Y')\otimes\Rep(\mathsf{H}))\times \Bet A(Y')^{\times n}\times \Gamma(Y',\cG).\]

The sheafy version of the left-hand side is 
\begin{equation}\label{eq: LHS sheafy}
    Y'\mapsto \QC(\Tot(\Phi_B(\cG|_{Y'}))\times^{\Bet \sfA} \Bet \mathsf{H}).
\end{equation}
Since we have
\begin{equation}\label{eq: ass bdl geom real}
    \Tot(\Phi_B(\cG|_{Y'}))\times^{\Bet \sfA} \Bet \mathsf{H} = \colim_{[n]\in\Delta^{\op}}\Tot(\Phi_B(\cG|_{Y'}))\times \Bet \sfA^{\times n}\times \Bet \mathsf{H} \quad \te{in }\Stk^E,
\end{equation}
we have
\begin{equation}\label{eq: QC BH fib totalization}
    \QC(\Tot(\Phi_B(\cG|_{Y'}))\times^{\Bet \sfA} \Bet \mathsf{H}) 
    = \lim_{[n]\in\Delta} \QC(\Tot(\Phi_B(\cG|_{Y'}))\times \Bet \sfA^{\times n}\times \Bet \mathsf{H})
\end{equation}
where the pullback functors are used in the limit. 
\medskip

In fact, all the transition maps in the geometric realization \eqref{eq: ass bdl geom real} have the form 
\[\cY\times \Bet A_1\rightarrow \cY\times\Bet A_2\]
for finite abelian groups $A_1$, $A_2$, and a homomorhpism $A_1\rightarrow A_2$. For such a morphism the pullback-pushforward adjunction for quasi-coherent sheaves is a biadjunction.\footnote{This is the consequence of the \emph{wrong-way} Frobenius reciprocity for finite groups.} Therefore, we can pass \eqref{eq: QC BH fib totalization} to left adjoints to get
\begin{equation}\label{eq: QC BH fib geom real}
    \QC(\Tot(\Phi_B(\cG|_{Y'}))\times^{\Bet \sfA} \Bet \mathsf{H}) 
    = \colim_{[n]\in\Delta^{\op}} \QC(\Tot(\Phi_B(\cG|_{Y'}))\times \Bet \sfA^{\times n}\times \Bet \mathsf{H}),
\end{equation}
using pushforwards. Following the notation of Sect. \ref{sssec: GnS TnS} we let
\[T^n_{Y'}:= \QC(\Tot(\Phi_B(\cG|_{Y'}))\times \Bet \sfA^{\times n}\times \Bet \mathsf{H}).\]

\sssec{}
As in the proof of Theorem \ref{thm: twisted sheaves tensor product formula} we construct a family of canonical maps $G^n_{Y'}\rightarrow T^n_{Y'}$. A section $\tau\in \Bet \sfA(Y')^{\times n}\times \Gamma(Y',\cG)$ can be thought of as a trivialization of the $\sfA^{n+1}$-gerbe
\[(\cG|_{Y'})^{+n}: Y'\rightarrow Y\xrightarrow{\cG} \Btet \sfA\rightarrow \Btet\sfA^{n+1},\]
where the last arrow is induced from the inclusion $\sfA\rightarrow\sfA^{n+1}$ into the first component.  Therefore, $\tau$ induces an equivalence
\begin{align}\label{eq: composition pullback}
\Tot(\Phi_B(\cG|_{Y'}))\times \Bet \sfA^{\times n}\times \Bet \mathsf{H} &= \Tot(\Phi_B((\cG|_{Y'})^{+n}))\times \Bet \mathsf{H}\\
& \xrightarrow{\tau} ((Y')^{\topl,\sing})_B\times \Bet \sfA^{\times n+1}\times \Bet \mathsf{H}.\nonumber
\end{align}
Hence, given such $\tau$ and $\cF\in\LS(Y')\otimes\Rep(\mathsf{H})$, i.e., an object of $G^n_{Y'}$ we can consider the sheaf 
\begin{equation}\label{eq: twRep kernel}
    \cF\boxtimes \cL_{\reg}^{\boxtimes n+1}\in \QC(((Y')^{\topl,\sing})_B\times \Bet \mathsf{H})\otimes \QC(\Bet \sfA^{n+1})
\end{equation}
where $\cL_{\reg}\in\QC(\Bet \sfA)$ corresponds to the regular representation of $\sfA$, and get an object in $T^n_{Y'}$ by pulling back \eqref{eq: twRep kernel} along \eqref{eq: composition pullback}.

\sssec{}
As in the proof of Theorem \ref{thm: twisted sheaves tensor product formula}, we confirm that the collection of functors $G^n_{Y'}\rightarrow T^n_{Y'}$ is compatible with transition functors of the geometric realization and hence defines a canonical functor
\begin{equation}\label{eq: can ftr G to T}
    \mathbf{Rep}(\mathsf{H})^{\loc}_{\cG}(Y')\rightarrow \QC(\Tot(\Phi_B(\cG|_{Y'}))\times^{\Bet \sfA} \Bet \mathsf{H}).
\end{equation}
This is functorial in $Y'\in Y^{\triv}_{\et}$, hence defines a functor from the sheaf \eqref{eq: twRep sheaf on triv} to \eqref{eq: LHS sheafy}. This is an equivalence, because upon choosing a trivialization of $\cG|_{Y'}$, the functor \eqref{eq: can ftr G to T} is exactly the equivalence in Example \ref{ex: equiv trivial gerbe}. This completes the proof of Theorem \ref{thm: twisted representations}.


\begin{thebibliography}{99}
\bibitem[Ach21]{Ach}
Pramod N. Achar, \emph{Perverse Sheaves and Applications to Representation Theory}, 
  Mathematical Surveys and Monographs, vol. 258, American Mathematical Society, Providence, RI, 2021. MR 4310578

\bibitem[AGH19]{AGH}
Benjamin Antieau, David Gepner, and Jeremiah Heller, \emph{K-theoretic obstructions to bounded $t$-structures}, 
  Invent. Math. 216 (2019), no. 1, 241--300. MR 3925528

\bibitem[AGK+20]{AGKRRV}
Dima Arinkin, Dennis Gaitsgory, David Kazhdan, Sam Raskin, Nick Rozenblyum, and Yasha
Varshavsky, \emph{The stack of local systems with restricted variation and geometric langlands theory
with nilpotent singular support}, \url{https://arxiv.org/abs/2010.01906}, 2020, preprint.

\bibitem[BD01]{BD}
Jean-Luc Brylinski and Pierre Deligne, \emph{Central extensions of reductive groups by $\mathbf{K}_2$}, Publ. Math.
Inst. Hautes Études Sci. (2001), no. 94, 5–85. MR 1896177

\bibitem[BL94]{BL}
Joseph Bernstein and Valery Lunts, \emph{Equivariant sheaves and functors}, 
  Lecture Notes in Mathematics, vol. 1578, Springer-Verlag, Berlin, 1994. MR 1299527

\bibitem[BZFN10]{BZFN}
David Ben-Zvi, John Francis, and David Nadler, \emph{Integral transforms and Drinfeld centers in derived algebraic geometry}, 
  J. Amer. Math. Soc. 23 (2010), no. 4, 909--966. MR 2669705

\bibitem[BZN18]{BZN}
David Ben-Zvi and David Nadler, \emph{Betti geometric Langlands}, Algebraic geometry: Salt Lake
City 2015, Proc. Sympos. Pure Math., vol. 97, Amer. Math. Soc., Providence, RI, 2018, pp. 3–-41. MR 3821166

\bibitem[DI04]{DI04}
Daniel Dugger, Daniel Isaksen, \emph{Topological hypercovers and $\afl$-realizations}, Math. Z. 246 (2004), no. 4, 667–689. MR 2045835

\bibitem[FS24]{FS}
Laurent Fargues and Peter Scholze, \emph{Geometrization of the local Langlands correspondence}, \url{https://arxiv.org/abs/2102.13459}, 2024, preprint.

\bibitem[Gai15a]{Ga}
Dennis Gaitsgory, \emph{Sheaves of categories and the notion of 1-affineness}, 
  Stacks and categories in geometry, topology, and algebra, 127–-225.
Contemp. Math., 643,
American Mathematical Society, Providence, RI, 2015. MR 3381473

\bibitem[Gai15b]{Gai1}
Dennis Gaitsgory. Outline of the proof of the geometric Langlands conjecture for $\on{GL}_2$, Astérisque No. 370 (2015), 1–112. MR 3364744

\bibitem[GG18]{GG}
Wee Teck Gan and Fan Gao, \emph{The Langlands-Weissman program for Brylinski-Deligne extensions},
Ast\'erisque, no. 398, Soci\'et\'e Math\'ematique de France,
Paris, 2018, L-groups and the Langlands program for covering groups, pp. 187--275. MR 3802419

\bibitem[GGW18]{GGW}
Wee Teck Gan, Fan Gao, and Martin H. Weissman, \emph{L-groups and the Langlands program for covering groups: a historical introduction}, Ast\'erisque, no. 398, Soci\'et\'e Math\'ematique de France, Paris, 2018, L-groups and the Langlands program for
covering groups, pp. 1--31. MR 3802417

\bibitem[GHN17]{GHN}
David Gepner, Rune Haugseng, and Thomas Nikolaus, \emph{Lax colimits and free fibrations in $\infty$-categories}, 
  Doc. Math. 22 (2017), 1225--1266. MR 3690268

\bibitem[GKRV22]{GKRV}
Dennis Gaitsgory, David Kazhdan, Nick Rozenblyum, and Yakov Varshavsky, \emph{A toy model for the
Drinfeld-Lafforgue shtuka construction}, Indag. Math. (N.S.) 33 (2022), no. 1, 39--189. MR 4368480

\bibitem[GL18]{GL}
Dennis Gaitsgory and Sergey Lysenko, \emph{Parameters and duality for the metaplectic geometric Lang-
lands theory}, Selecta Math. (N.S.) 24 (2018), no. 1, 227--301, References are to the corrected version,
available at: \url{https://lysenko.perso.math.cnrs.fr}. MR 3769731

\bibitem[GR17a]{GR}
Dennis Gaitsgory and Nick Rozenblyum, \emph{A Study in Derived Algebraic Geometry. Vol. I. Correspondences and Duality}, 
  Mathematical Surveys and Monographs, vol. 221, American Mathematical Society, Providence, RI, 2017. MR 3701352

\bibitem[GR17b]{GR2}
Dennis Gaitsgory and Nick Rozenblyum, \emph{A Study in Derived Algebraic Geometry. Vol. II. Deformations, Lie Theory and Formal Geometry}, 
  Mathematical Surveys and Monographs, vol. 221, American Mathematical Society, Providence, RI, 2017. MR 3701353


\bibitem[KS90]{KS}
Masaki Kashiwara and Pierre Schapira, \emph{Sheaves on Manifolds}, 
  Grundlehren der mathematischen Wissenschaften, vol. 292, Springer-Verlag, Berlin, 1990. MR 1074006

\bibitem[Laf18]{Laf}
Vincent Lafforgue, \emph{Chtoucas pour les groupes r\'eductifs et param\'etrisation de Langlands globale},
J. Amer. Math. Soc. 31 (2018), no. 3, 719--891. MR 3787407



\bibitem[Lur09]{HTT}
Jacob Lurie, \emph{Higher Topos Theory}, 
  Annals of Mathematics Studies, vol. 170, Princeton University Press, Princeton, NJ, 2009. MR~2522659

  \bibitem[Lur17]{HA}
Jacob Lurie, \emph{Higher Algebra}, 
  available at \url{https://www.math.ias.edu/~lurie/}, 2017.

\bibitem[MORezk]{MORezk}
Charles Rezk (\url{https://mathoverflow.net/users/437/charles-rezk}), When do Kan extensions preserve limits/colimits?, URL (version: 2023-07-05): \url{https://mathoverflow.net/q/450214}

\bibitem[Man22]{Man}
Lucas Mann, \emph{A $p$-adic 6-functor formalism in rigid-analytic geometry},\url{arXiv:2206.02022}, 2022, preprint

\bibitem[NY19]{NY}
David Nadler and Zhiwei Yun, \emph{Spectral action in Betti geometric Langlands}, Israel J. Math. 232 (2019), no. 1, 299--349. MR 3990944



\bibitem[PPS25]{PPS}
James Pascaleff, Emanuele Pavia, and Nicol{\`o} Sibilla, \emph{Higher Koszul duality and $n$-affineness}, \url{arXiv:2504.16935}, 2025, preprint

\bibitem[Rei12]{Rei}
Ryan Cohen Reich, \emph{Twisted geometric Satake equivalence via gerbes on the factorizable Grass-
mannian}, Represent. Theory 16 (2012), 345--449. MR 2956088

\bibitem[Sch25]{Sch}
Peter Scholze, \emph{Six-functor formalisms}, \url{https://arxiv.org/abs/2510.26269}, 2025, preprint

\bibitem[Wei18]{Wei}
Martin Weissman, \emph{L-groups and parameters for covering groups}, Ast\'erisque, no. 398, Soci\'et\'e Math\'ematique
de France, 2018, L-groups and the Langlands program for covering groups, 33--186. MR 3802418

\bibitem[Zha22]{Zhao1}
Yifei Zhao, \emph{\'Etale metaplectic covers of reductive group schemes}, \url{https://arxiv.org/abs/2204.00610}, 2022, preprint.

\bibitem[Zha24]{ZhaoNew}
Yifei Zhao, \emph{Spectral decomposition of genuine cusp forms over global function fields}, Compos. Math. 160 (2024), no. 6, 1194--1260. MR 4741316

\bibitem[Zhu17]{Zhu}
Xinwen Zhu, \emph{An introduction to affine Grassmannians and the geometric Satake equivalence}, Geometry
of moduli spaces and representation theory, IAS/Park City Math. Ser., vol. 24, Amer. Math. Soc.,
Providence, RI, 2017, pp. 59–154. MR 3752460
\end{thebibliography}
\end{document}